\documentclass{article}
\usepackage{amssymb,amsmath,theorem,euscript}

\usepackage[cp1251]{inputenc}
\usepackage[english,russian]{babel}

\input{epsf}

\newcounter{sec}

\def\sm{\smallskip}

\newcounter{punct}[sec]

\def\punct{\refstepcounter{punct}{\arabic{sec}.\arabic{punct}.  }}

\def\COUNTERS{\addtocounter{sec}{1}
              \setcounter{punct}{0}
          \setcounter{equation}{0}
          \setcounter{theorem}{0}
                  }
\newtheorem{theorem}{Теорема}[sec]
\newtheorem{proposition}[theorem]{Предложение}
\newtheorem{lemma}[theorem]{Лемма}
\newtheorem{corollary}[theorem]{Следствие}

\newtheorem{definition}[theorem]{Определение}

\begin{document}

 \def\ov{\overline}
\def\wt{\widetilde}
\def\wh{\widehat}
 \newcommand{\rk}{\mathop {\mathrm {rk}}\nolimits}
\newcommand{\Aut}{\mathop {\mathrm {Aut}}\nolimits}
\newcommand{\Out}{\mathop {\mathrm {Out}}\nolimits}
\newcommand{\Abs}{\mathop {\mathrm {Abs}}\nolimits}
\renewcommand{\Re}{\mathop {\mathrm {Re}}\nolimits}
\renewcommand{\Im}{\mathop {\mathrm {Im}}\nolimits}
 \newcommand{\tr}{\mathop {\mathrm {tr}}\nolimits}
  \newcommand{\Hom}{\mathop {\mathrm {Hom}}\nolimits}
   \newcommand{\diag}{\mathop {\mathrm {diag}}\nolimits}
   \newcommand{\supp}{\mathop {\mathrm {supp}}\nolimits}
 \newcommand{\im}{\mathop {\mathrm {im}}\nolimits}
 \newcommand{\grad}{\mathop {\mathrm {grad}}\nolimits}
  \newcommand{\sgrad}{\mathop {\mathrm {sgrad}}\nolimits}
 \newcommand{\rot}{\mathop {\mathrm {rot}}\nolimits}
  \renewcommand{\div}{\mathop {\mathrm {div}}\nolimits}

\def\Br{\mathrm {Br}}
\def\Vir{\mathrm {Vir}}

 \def\Ham{\mathrm {Ham}}
\def\SL{\mathrm {SL}}
\def\Pol{\mathrm {Pol}}
\def\SU{\mathrm {SU}}
\def\GL{\mathrm {GL}}
\def\U{\mathrm U}
\def\OO{\mathrm O}
 \def\Sp{\mathrm {Sp}}
  \def\Ad{\mathrm {Ad}}
 \def\SO{\mathrm {SO}}
\def\SOS{\mathrm {SO}^*}
 \def\Diff{\mathrm{Diff}}
 \def\Vect{\mathfrak{Vect}}
\def\PGL{\mathrm {PGL}}
\def\PU{\mathrm {PU}}
\def\PSL{\mathrm {PSL}}
\def\Symp{\mathrm{Symp}}
\def\Cont{\mathrm{Cont}}
\def\End{\mathrm{End}}
\def\Mor{\mathrm{Mor}}
\def\Aut{\mathrm{Aut}}
 \def\PB{\mathrm{PB}}
\def\Fl{\mathrm {Fl}}
\def\Symm{\mathrm {Symm}} 
 \def\Herm{\mathrm {Herm}} 
  \def\SDiff{\mathrm {SDiff}} 
 
 \def\cA{\mathcal A}
\def\cB{\mathcal B}
\def\cC{\mathcal C}
\def\cD{\mathcal D}
\def\cE{\mathcal E}
\def\cF{\mathcal F}
\def\cG{\mathcal G}
\def\cH{\mathcal H}
\def\cJ{\mathcal J}
\def\cI{\mathcal I}
\def\cK{\mathcal K}
 \def\cL{\mathcal L}
\def\cM{\mathcal M}
\def\cN{\mathcal N}
 \def\cO{\mathcal O}
\def\cP{\mathcal P}
\def\cQ{\mathcal Q}
\def\cR{\mathcal R}
\def\cS{\mathcal S}
\def\cT{\mathcal T}
\def\cU{\mathcal U}
\def\cV{\mathcal V}
 \def\cW{\mathcal W}
\def\cX{\mathcal X}
 \def\cY{\mathcal Y}
 \def\cZ{\mathcal Z}
\def\0{{\ov 0}}
 \def\1{{\ov 1}}
 \def\frA{\mathfrak A}
 \def\frB{\mathfrak B}
\def\frC{\mathfrak C}
\def\frD{\mathfrak D}
\def\frE{\mathfrak E}
\def\frF{\mathfrak F}
\def\frG{\mathfrak G}
\def\frH{\mathfrak H}
\def\frI{\mathfrak I}
 \def\frJ{\mathfrak J}
 \def\frK{\mathfrak K}
 \def\frL{\mathfrak L}
\def\frM{\mathfrak M}
 \def\frN{\mathfrak N} \def\frO{\mathfrak O} \def\frP{\mathfrak P} \def\frQ{\mathfrak Q} \def\frR{\mathfrak R}
 \def\frS{\mathfrak S} \def\frT{\mathfrak T} \def\frU{\mathfrak U} \def\frV{\mathfrak V} \def\frW{\mathfrak W}
 \def\frX{\mathfrak X} \def\frY{\mathfrak Y} \def\frZ{\mathfrak Z} \def\fra{\mathfrak a} \def\frb{\mathfrak b}
 \def\frc{\mathfrak c} \def\frd{\mathfrak d} \def\fre{\mathfrak e} \def\frf{\mathfrak f} \def\frg{\mathfrak g}
 \def\frh{\mathfrak h} \def\fri{\mathfrak i} \def\frj{\mathfrak j} \def\frk{\mathfrak k} \def\frl{\mathfrak l}
 \def\frm{\mathfrak m} \def\frn{\mathfrak n} \def\fro{\mathfrak o} \def\frp{\mathfrak p} \def\frq{\mathfrak q}
 \def\frr{\mathfrak r} \def\frs{\mathfrak s} \def\frt{\mathfrak t} \def\fru{\mathfrak u} \def\frv{\mathfrak v}
 \def\frw{\mathfrak w} \def\frx{\mathfrak x} \def\fry{\mathfrak y} \def\frz{\mathfrak z} \def\frsp{\mathfrak{sp}}
 \def\bfa{\mathbf a} \def\bfb{\mathbf b} \def\bfc{\mathbf c} \def\bfd{\mathbf d} \def\bfe{\mathbf e} \def\bff{\mathbf f}
 \def\bfg{\mathbf g} \def\bfh{\mathbf h} \def\bfi{\mathbf i} \def\bfj{\mathbf j} \def\bfk{\mathbf k} \def\bfl{\mathbf l}
 \def\bfm{\mathbf m} \def\bfn{\mathbf n} \def\bfo{\mathbf o} \def\bfp{\mathbf p} \def\bfq{\mathbf q} \def\bfr{\mathbf r}
 \def\bfs{\mathbf s} \def\bft{\mathbf t} \def\bfu{\mathbf u} \def\bfv{\mathbf v} \def\bfw{\mathbf w} \def\bfx{\mathbf x}
 \def\bfy{\mathbf y} \def\bfz{\mathbf z} \def\bfA{\mathbf A} \def\bfB{\mathbf B} \def\bfC{\mathbf C} \def\bfD{\mathbf D}
 \def\bfE{\mathbf E} \def\bfF{\mathbf F} \def\bfG{\mathbf G} \def\bfH{\mathbf H} \def\bfI{\mathbf I} \def\bfJ{\mathbf J}
 \def\bfK{\mathbf K} \def\bfL{\mathbf L} \def\bfM{\mathbf M} \def\bfN{\mathbf N} \def\bfO{\mathbf O} \def\bfP{\mathbf P}
 \def\bfQ{\mathbf Q} \def\bfR{\mathbf R} \def\bfS{\mathbf S} \def\bfT{\mathbf T} \def\bfU{\mathbf U} \def\bfV{\mathbf V}
 \def\bfW{\mathbf W} \def\bfX{\mathbf X} \def\bfY{\mathbf Y} \def\bfZ{\mathbf Z} \def\bfw{\mathbf w}
 \def\R {{\mathbb R }} \def\C {{\mathbb C }} \def\Z{{\mathbb Z}} \def\H{{\mathbb H}} \def\K{{\mathbb K}}
 \def\N{{\mathbb N}} \def\Q{{\mathbb Q}} \def\A{{\mathbb A}} \def\T{\mathbb T} \def\P{\mathbb P} \def\G{\mathbb G}
 \def\bbA{\mathbb A} \def\bbB{\mathbb B} \def\bbD{\mathbb D} \def\bbE{\mathbb E} \def\bbF{\mathbb F} \def\bbG{\mathbb G}
 \def\bbI{\mathbb I} \def\bbJ{\mathbb J} \def\bbL{\mathbb L} \def\bbM{\mathbb M} \def\bbN{\mathbb N} \def\bbO{\mathbb O}
 \def\bbP{\mathbb P} \def\bbQ{\mathbb Q} \def\bbS{\mathbb S} \def\bbT{\mathbb T} \def\bbU{\mathbb U} \def\bbV{\mathbb V}
 \def\bbW{\mathbb W} \def\bbX{\mathbb X} \def\bbY{\mathbb Y} \def\kappa{\varkappa} \def\epsilon{\varepsilon}
 \def\phi{\varphi} \def\le{\leqslant} \def\ge{\geqslant}

\def\UU{\bbU}
\def\Mat{\mathrm{Mat}}
\def\tto{\rightrightarrows}

\def\F{\mathbf{F}}

\def\Gms{\mathrm {Gms}}
\def\Ams{\mathrm {Ams}}
\def\Isom{\mathrm {Isom}}

\def\Gr{\mathrm{Gr}}

\def\graph{\mathrm{graph}}

\def\O{\mathrm{O}}

\def\la{\langle}
\def\ra{\rangle}


 \def\ov{\overline}
\def\wt{\widetilde}

\renewcommand{\Re}{\mathop {\mathrm {Re}}\nolimits}
\def\Br{\mathrm {Br}}

 \def\Isom{\mathrm {Isom}}
 \def\Hier{\mathrm {Hier}}
\def\SL{\mathrm {SL}}
\def\SU{\mathrm {SU}}
\def\GL{\mathrm {GL}}
\def\U{\mathrm U}
\def\OO{\mathrm O}
 \def\Sp{\mathrm {Sp}}
  \def\GLO{\mathrm {GLO}}
 \def\SO{\mathrm {SO}}
\def\SOS{\mathrm {SO}^*}
 \def\Diff{\mathrm{Diff}}
 \def\Vect{\mathfrak{Vect}}
\def\PGL{\mathrm {PGL}}
\def\PU{\mathrm {PU}}
\def\PSL{\mathrm {PSL}}
\def\Symp{\mathrm{Symp}}
\def\ASymm{\mathrm{Asymm}}
\def\Asymm{\mathrm{Asymm}}
\def\Gal{\mathrm{Gal}}
\def\End{\mathrm{End}}
\def\Mor{\mathrm{Mor}}
\def\Aut{\mathrm{Aut}}
 \def\PB{\mathrm{PB}}
 \def\cA{\mathcal A}
\def\cB{\mathcal B}
\def\cC{\mathcal C}
\def\cD{\mathcal D}
\def\cE{\mathcal E}
\def\cF{\mathcal F}
\def\cG{\mathcal G}
\def\cH{\mathcal H}
\def\cJ{\mathcal J}
\def\cI{\mathcal I}
\def\cK{\mathcal K}
 \def\cL{\mathcal L}
\def\cM{\mathcal M}
\def\cN{\mathcal N}
 \def\cO{\mathcal O}
\def\cP{\mathcal P}
\def\cQ{\mathcal Q}
\def\cR{\mathcal R}
\def\cS{\mathcal S}
\def\cT{\mathcal T}
\def\cU{\mathcal U}
\def\cV{\mathcal V}
 \def\cW{\mathcal W}
\def\cX{\mathcal X}
 \def\cY{\mathcal Y}
 \def\cZ{\mathcal Z}
\def\0{{\ov 0}}
 \def\1{{\ov 1}}
 
 \def\frA{\mathfrak A}
 \def\frB{\mathfrak B}
\def\frC{\mathfrak C}
\def\frD{\mathfrak D}
\def\frE{\mathfrak E}
\def\frF{\mathfrak F}
\def\frG{\mathfrak G}
\def\frH{\mathfrak H}
\def\frI{\mathfrak I}
 \def\frJ{\mathfrak J}
 \def\frK{\mathfrak K}
 \def\frL{\mathfrak L}
\def\frM{\mathfrak M}
 \def\frN{\mathfrak N} \def\frO{\mathfrak O} \def\frP{\mathfrak P} \def\frQ{\mathfrak Q} \def\frR{\mathfrak R}
 \def\frS{\mathfrak S} \def\frT{\mathfrak T} \def\frU{\mathfrak U} \def\frV{\mathfrak V} \def\frW{\mathfrak W}
 \def\frX{\mathfrak X} \def\frY{\mathfrak Y} \def\frZ{\mathfrak Z} \def\fra{\mathfrak a} \def\frb{\mathfrak b}
 \def\frc{\mathfrak c} \def\frd{\mathfrak d} \def\fre{\mathfrak e} \def\frf{\mathfrak f} \def\frg{\mathfrak g}
 \def\frh{\mathfrak h} \def\fri{\mathfrak i} \def\frj{\mathfrak j} \def\frk{\mathfrak k} \def\frl{\mathfrak l}
 \def\frm{\mathfrak m} \def\frn{\mathfrak n} \def\fro{\mathfrak o} \def\frp{\mathfrak p} \def\frq{\mathfrak q}
 \def\frr{\mathfrak r} \def\frs{\mathfrak s} \def\frt{\mathfrak t} \def\fru{\mathfrak u} \def\frv{\mathfrak v}
 \def\frw{\mathfrak w} \def\frx{\mathfrak x} \def\fry{\mathfrak y} \def\frz{\mathfrak z} \def\frsp{\mathfrak{sp}}
 \def\bfa{\mathbf a} \def\bfb{\mathbf b} \def\bfc{\mathbf c} \def\bfd{\mathbf d} \def\bfe{\mathbf e} \def\bff{\mathbf f}
 \def\bfg{\mathbf g} \def\bfh{\mathbf h} \def\bfi{\mathbf i} \def\bfj{\mathbf j} \def\bfk{\mathbf k} \def\bfl{\mathbf l}
 \def\bfm{\mathbf m} \def\bfn{\mathbf n} \def\bfo{\mathbf o} \def\bfp{\mathbf p} \def\bfq{\mathbf q} \def\bfr{\mathbf r}
 \def\bfs{\mathbf s} \def\bft{\mathbf t} \def\bfu{\mathbf u} \def\bfv{\mathbf v} \def\bfw{\mathbf w} \def\bfx{\mathbf x}
 \def\bfy{\mathbf y} \def\bfz{\mathbf z} \def\bfA{\mathbf A} \def\bfB{\mathbf B} \def\bfC{\mathbf C} \def\bfD{\mathbf D}
 \def\bfE{\mathbf E} \def\bfF{\mathbf F} \def\bfG{\mathbf G} \def\bfH{\mathbf H} \def\bfI{\mathbf I} \def\bfJ{\mathbf J}
 \def\bfK{\mathbf K} \def\bfL{\mathbf L} \def\bfM{\mathbf M} \def\bfN{\mathbf N} \def\bfO{\mathbf O} \def\bfP{\mathbf P}
 \def\bfQ{\mathbf Q} \def\bfR{\mathbf R} \def\bfS{\mathbf S} \def\bfT{\mathbf T} \def\bfU{\mathbf U} \def\bfV{\mathbf V}
 \def\bfW{\mathbf W} \def\bfX{\mathbf X} \def\bfY{\mathbf Y} \def\bfZ{\mathbf Z} \def\bfw{\mathbf w}

 \def\R {{\mathbb R }} \def\C {{\mathbb C }} \def\Z{{\mathbb Z}} \def\H{{\mathbb H}} \def\K{{\mathbb K}}
 \def\N{{\mathbb N}} \def\Q{{\mathbb Q}} \def\A{{\mathbb A}} \def\T{\mathbb T} \def\P{\mathbb P} \def\G{\mathbb G}
 \def\bbA{\mathbb A} \def\bbB{\mathbb B} \def\bbD{\mathbb D} \def\bbE{\mathbb E} \def\bbF{\mathbb F} \def\bbG{\mathbb G}
 \def\bbI{\mathbb I} \def\bbJ{\mathbb J} \def\bbL{\mathbb L} \def\bbM{\mathbb M} \def\bbN{\mathbb N} \def\bbO{\mathbb O}
 \def\bbP{\mathbb P} \def\bbQ{\mathbb Q} \def\bbS{\mathbb S} \def\bbT{\mathbb T} \def\bbU{\mathbb U} \def\bbV{\mathbb V}
 \def\bbW{\mathbb W} \def\bbX{\mathbb X} \def\bbY{\mathbb Y} \def\kappa{\varkappa} \def\epsilon{\varepsilon}
 \def\phi{\varphi} \def\le{\leqslant} \def\ge{\geqslant}

\def\UU{\bbU}
\def\Mat{\mathrm{Mat}}
\def\tto{\rightrightarrows}

\def\Gr{\mathrm{Gr}}

\def\B{\bfB} 

\def\graph{\mathrm{graph}}

\def\O{\mathrm{O}}

\def\la{\langle}
\def\ra{\rangle}

\begin{center}
{\bf\Large
Topological groups and invariant measures}

\bigskip


\large
{\sc Yury Neretin}

\bigskip

Lecture notes in Russian%
\footnote{The first part is based on my lectures
	given in Moscow State University in November--Dececember 2015.
	The second part is based on my lectures in Independent
	University of Moscow in different years and in Bialowieza, Poland, in July 2014.
	Partially supported by the grant FWF, Project P25142.}. Draft
\end{center}

{\sc Chapter A. The Haar measure on locally compact groups}

\bigskip

\noindent
{\sc 1. Topological groups. Haar measure}
\newline
{\sc 2. Various proofs of existence of Haar measure.}
\newline
{\sc 3. Modular character. Invariant measures on homogeneous spaces.}
\newline
{\sc 4. Invariant measures on Grassmannians and orthogonal groups.}
\newline
{\sc Bibliography to \S\S1-4.}

\bigskip

\noindent
{\sc Chapter B. Miscelania} 

\bigskip 

\noindent
{\sc 5. Zoo of topological groups.}
\newline
{\sc 6. Fourier transform on compact groups.}
\newline
{\sc 7. Distribution of eigenvalues.}
\newline
{\sc 8. Symmetries of Gaussian measures.}
\newline
{\sc 9. Poisson measures.}
\newline
{\sc 10. Inverse limits. Virtual permutations} 
\newline
{\sc 11. Inverse limits. Unitary groups.}
\newline
{\sc 12. Closures of actions}

\sm

\noindent{\sc Addendum. Lebesgue spaces}

\newpage

\begin{center}
{\bf\Large
Топологические группы и инвариантные меры}

\bigskip

\large
{\sc Ю.А.Неретин}

\bigskip

Записки лекций%
\footnote{Первая часть основана на моих лекциях 
	в МГУ, ноябрь-декабрь 2014.
	Вторая часть - на лекциях, читавшихся в Московском
	независимом университете в разные годы и в Беловеже, Польша,
	в июле 2014. Частично поддержано грантом FWF, P25142.}

\end{center}

{\sc Глава A. Мера Хаара на локально компактных группах}

\bigskip

\noindent
{\sc 1. Топологические группы. Мера Хаара}
\newline
{\sc 2. Разные доказательства существования меры Хаара}
\newline
{\sc 3. Модулярный характер. Меры на однородных пространствах.}
\newline
{\sc 4. Инвариантные меры на грассманиане и ортогональной группе.}
\newline
{\sc Литература%
\footnote{Литература из этого списка цитируется аббревиатурами типа 
\cite{Zhe1}, после каждого параграфа 5-11 идет свой список литературы,
которая нумеруется цифрами, скажем, \cite{Shir}.} \S\S1-4}

\bigskip

\noindent
{\sc Глава B. Разное} 

\bigskip 

\noindent
{\sc 5. Зоопарк топологических групп}
\newline
{\sc 6. Преобразование Фурье на компактных группах}
\newline
{\sc 7. Распределение собственных чисел.}
\newline
{\sc 8. Симметрии гауссовых мер.}
\newline
{\sc 9. Меры Пуассона.}
\newline
{\sc 10. Обратные пределы. Виртуальные перестановки} 
\newline
{\sc 11. Обратные пределы. Унитарные группы}
\newline
{\sc 12. Замыкание действий}

\sm

\noindent
{\sc Добавление. Лебеговские пространства}

\newpage

\begin{center}
	\huge
Глава A. Мера Хаара на локально компактных группах
\end{center}

\section[Топологические группы. Мера Хаара]{Топологические группы. Мера Хаара}

\COUNTERS

{\bf\punct Топологические группы.}

\begin{definition}
Топологическая группа  $G$ -- это хаусдорфово топологическое пространство
с групповыми операциями, причем умножение $G\times G\to G$ непрерывно по совокупности
переменных, а взятие обратного элемента $G\to G$ непрерывно.
\end{definition}

Ниже в формулировках и доказательствах основных утверждений мы, как правило, будем
предполагать, что топология на $G$ имеет  счетную базу. 

\sm

{\sc Замечание.}
Эти условия можно ослабить и потребовать лишь:

1. раздельную непрерывность умножения $(x,y)\to xy$;

2. непрерывность умножения в точке $(e,e)$; 3. непрерывность отображения $x\mapsto x^{-1}$
в точке $e$.

 Это иногда бывает полезно при  проверке того, что данная группа $G$
является топологической,
 см, например \cite[2.6]{Kir}. 
\hfill $\lozenge$

\sm

В \S\S 1-4 для нас основной интерес представляют {\it локально компактные группы}. Напомним, что
топологическое пространство называется локально компактным, если любая его точка имеет окрестность с компактным  замыканием.

\sm

См. \cite[\S2]{Zhe1}, \cite[гл.~4]{Zhe2},   \cite[2.6]{Kir}, \cite[гл.3]{Pon}, \cite[гл.2]{HR1},
\cite[гл.3]{Bou1}.

\sm

{\bf\punct Начальный зоопарк топологических групп.%
\label{ss:zoo-nach}}

1) Любая группа, снабженная дискретной топологией, является топологической группой.

\sm 
 
2) {\it Классические группы.} 

\sm 

--- $\GL(n,\R)$ -- группа всех обратимых вещественных матриц порядка $n$.

\sm 

--- $\SL(n,\R)$ -- группа всех  вещественных матриц порядка $n$ c определителем 1.

\sm 

--- $\GL(n,\C)$, $\SL(n,\R)$ -- определяемые аналогично группы комплексных матриц.

\sm 

--- $\O(n)$ -- группа ортогональных ($gg^t=1$) вещественных матриц порядка $n$; $\SO(n)$ --
ее подгруппа, состоящая из матриц с определителем единица; 

\sm 

--- $\U(n)$ -- группа унитарных ($gg^*=1$) комплексных матриц порядка $n$; $\SU(n)$ --
ее подгруппа, состоящая из матриц с определителем единица; 

\sm 

--- $\O(p,q)$ -- {\it псевдоортогональная группа}, т.е., группа вещественных матриц порядка $(p+q)$, cохраняющих знаконеопределенное
скалярное произведение; иными словами $g\in \O(p,q)$, если%
\footnote{Через $1_p$ мы обозначаем единичную матрицу размера $p$; через $g^t$ -- транспонированную
	матрицу.}
$$
g\begin{pmatrix}1_p&0\\ 0&-1_q  \end{pmatrix} g^t= \begin{pmatrix}1_p&0\\ 0&-1_q  \end{pmatrix}
.$$

\sm

{\sc Замечание.} Обозначения происходят из аббревиатур немецко-английских слов general linear, special linear, orthogronal, unitary.
\hfill $\lozenge$

\sm

3) Некоторые матричные группы.

\sm

-- Группа матриц порядка 2 вида $\begin{pmatrix} a&b\\0&1 \end{pmatrix}$, где $a>0$.

\sm

--- {\it Группа верхнетреугольных матриц} порядка $n$, т.е. группа матриц вида
$$
\begin{pmatrix}t_{11}& t_{11}& t_{13}&\dots\\
               0&t_{22}&t_{23}&\dots\\
               0&0&t_{33}&\dots\\
               \vdots&\vdots&\vdots&\ddots\end{pmatrix}
,$$
где диагональные элементы $t_{jj}$  обратимы;

\sm

--- {\it Группа строго верхнетреугольных матриц} состоит их матриц того же вида с $t_{jj}=1$.

\sm

4) Группа диффеоморфизмов компактного многообразия, снабженная топологией равномерной 
сходимости всех частных производных.

\sm

5) Полная группа $\mathbf{U}(\infty)$ унитарных операторов гильбертова пространства.
На ней есть две естественных топологии, равномерная и слабая.
Подробнее эта группа обсуждается в (п.\ref{ss:unitary-infty}.

\sm

6) Любая подгруппа топологической группы является топологической группой.

\sm

7) Замкнутая подгруппа локально компактной группы является локально компактной группой
(а замкнутая подгруппа компактной группы компактна).

\sm

8) {\it Фактор-группа.} Пусть $G$ -- топологическая группа, а $H$ -- ее замкнутая нормальная
подгруппа. Тогда фактор-группа $G/H$, снабженная фактор-топологией, является топологической
группой.  Если $G$ - локально компактна, то $G/H$ локально компактна

\sm

9) Конечное или счетное произведение компактных групп - компактная группа. Конечное
произведение локально компактных групп -- локально компактная группа.

\sm

Для читателя, неудовлетворенного краткостью и очевидностью приведенного списка примеров,
мы  продолжим его в \S5.

Стоит отметить, что все вышеперечисленные группы, кроме унитарной группы, снабженной равномерной
операторной топологией (а также несчетных дискретных групп), имеют счетную базу топологии. 

\sm

{\bf\punct Теорема Хаара.}
{\it Левым сдвигом} $L_h:G\to G$, где $h\in G$, называется преобразование  
$L_h g=hg$. Аналогично, {\it правый сдвиг $R_h$} -- это преобразование
$R_hg=gh$.

\begin{theorem}
 На любой локально компактной топологической группе со счетной базой существует единственная
с точностью до постоянного множителя $\sigma$-конечная мера, инвариантная относительно
левых сдвигов.
\end{theorem}

Мы  предполагаем, что мера определена на $\sigma$-алгебре, порожденной открытыми множествами и что мера любого компактного множества конечна.

Эта мера называется {\it лево-инвариантной мерой Хаара.} 
Разумеется, право-ин\-ва\-ри\-ант\-ная мера
тоже существует, она получается из лево-ин\-ва\-ри\-ант\-ной преобразованием
$g\mapsto g^{-1}$.

\begin{theorem}
На компактной группе мера Хаара является конечной, она инвариантна также относительно
 преобразования $g\mapsto g^{-1}$.
\end{theorem}

{\sc Замечание.} Стоит отметить, что согласно теореме А.Вейля \cite{Wei},
 мера Хаара существует и единственна также на локально-компактных группах, не имеющих счетной базы, но она уже
не должна быть $\sigma$-конечной, и $\sigma$-алгебру тоже надо определять аккуратней%
\footnote{Множество называется бэровским, если оно компактно и представимо как счетное
пересечение открытых множеств.  Мера определяется на $\sigma$-алгебре, порожденной 
бэровскими множествами. Мера  бэровского множества полагается конечной. См. \cite{RS1},
IV.4.}.
\hfill $\lozenge$

\sm

Известные доказательства общей теоремы (см. \S2) довольно хитроумны, однако для конкретных
групп  меры Хаара несложно построить явно. Начнем с тривиальных примеров.

\sm

{\sc Пример.} Мера Хаара на $\R^n$ -- это мера Лебега. Мера Хаара на торе
$(\R/\Z)^n$ -- тоже мера Лебега. \hfill $\lozenge$

\sm

{\sc Пример.} Для дискретной группы мы полагаем, что мера каждой точки равна 1.
\hfill $\lozenge$

\sm

{\sc Пример.} Мера Хаара на мультипликативной группе положительных чисел равна
$x^{-1}dx$, где $dx$ -- мера Лебега (а чему равна мера Хаара на мультипликативной группе
комплексных чисел?).
\hfill $\lozenge$

\sm

{\sc Пример.} Рассмотрим счетное произведение $\Z_2^\infty$ групп $\Z_2$. Мера Хаара совпадает
с естественной мерой на счетном произведении двоеточий.
\hfill $\lozenge$

\sm

См. \cite[гл.7]{Bou2}, \cite[\S7-8]{Wei}, \cite[\S26]{Zhe2}, \cite{Nach},
\cite[\S29]{Pon}, \cite[\S15]{HR1},
\cite[гл.11]{Hal}.

\sm

{\bf\punct Мера Хаара на $\GL(n,\R)$.} Рассмотрим пространство $\Mat(n)$ вещественных матриц
порядка $n$. Группа $\GL(n,\C)\subset\Mat(n)$ является дополнением до гиперповерхности
$\det(X)=0$. Пусть $dX$ -- мера Лебега на $\Mat(n)$. Мы будем искать меру Хаара
в виде $\phi(X)\,dX$, где $\phi(X)$ -- функция плотности.

Сначала заметим, что для $A\in \GL(n,\R)$ выполнено
$$
d(AX)=\det(A)^n\cdot dX
$$
В самом деле, матрица $X$ составлена из своих столбцов $X^1$, \dots, $X^n$,
т.е. $\Mat(n)\simeq(\R^n)^n$. 
Преобразование $X\mapsto AX$ устроено как 
$$
(X^1,\dots,X^n)\mapsto (AX^1,\dots,AX^n).
$$
На каждой копии $\R^n$ мера Лебега умножается на $\det(A)$, а на $(\R^n)^n$ -- на 
$\det(A)^n$. 

Теперь мы можем написать множитель $\psi(X)$.

\begin{proposition}
Мера Хаара на $\GL(n,\R)$ равна $\det(X)^{-n}\,dX$.
\end{proposition}

{\sc Задача.} Чему равна мера Хаара на $\GL(n,\C)$?

\sm

Ниже нам понадобится еще одно упражнение в том же духе. Рассмотрим пространства
$\Symm(n)$ и $\ASymm(n)$, состоящие из вещественных симметричных и кососимметричных матриц.
Рассмотрим линейное преобразование
$$
\frA:X\mapsto AXA^t
$$
пространства всех матриц. Оно распадается в прямую сумму таких же преобразований
на $\Symm(n)$ и $\ASymm(n)$, обозначим их через $\frA_{\Symm(n)}$ и $\frA_{\ASymm(n)}$.

\begin{lemma}
\label{l:AXA}
$$
\det(\frA_{\Symm(n)})=\det(A)^{n+1},\qquad \det(\frA_{\ASymm(n)})=\det(A)^{n-1}.
$$
\end{lemma}

{\sc Доказательство.} Проверим, например, первое равенство.
Обозначим искомый определитель через $\cD(A)$. Очевидно, что
$\cD(AB)=\cD(A)\,\cD(B)$.  Поэтому
$\cD(CAC^{-1})=\cD(C)$. Следовательно,
 утверждение достаточно проверить для
диагональных матриц.

Рассмотрим диагональную матрицу $\Lambda$ с элементами
$\lambda_j$. Тогда матрица $\Lambda X\Lambda^t$ имеет матричные элементы
$x_{ij}\lambda_i\lambda_j$. Обозначим через $E_{ij}$ матричные единицы, т.е. матрицы, у которых
$ij$-ый элемент равен 1, а остальные равны нулю. Тогда базис $E_{ij}+E_{ji}$, $E_{ii}$
оказывается собственным для линейного преобразования $X\mapsto\Lambda X\Lambda^t$, 
а потому $D(\Lambda)$ равно произведению собственных чисел,
$$\qquad \qquad\qquad\qquad\quad
D(\Lambda)=\prod_{i<j} \lambda_i\lambda_j\prod_i\lambda_i^2=\det(\Lambda)^{n+1}
. \qquad \qquad\qquad\qquad\quad\lozenge
$$ 

\sm

{\bf\punct Дальнейшие примеры.} 
1) Рассмотрим группу строго верхнетреугольных матриц $T$. Мы утверждаем, что мера Хаара
совпадает с мерой Лебега $dT=\prod_{i<j} dt_{ij}$. Пусть, для определенности, размер матриц равен  4. Вместо меры мы будем писать дифференциальную форму старшей степени
$$
dt_{12}\wedge dt_{13}\wedge dt_{14}\wedge dt_{23}\wedge dt_{24}\wedge dt_{34}.
$$
Перемножая матрицы
\begin{multline*}
\begin{pmatrix}
1&a_{12}&a_{13}&a_{14}\\
0&1&a_{23}&a_{24}\\
0&0&1&a_{34}\\
0&0&0&1
\end{pmatrix}
\begin{pmatrix}
1&t_{12}&t_{13}&t_{14}\\
0&1&t_{23}&t_{24}\\
0&0&1&t_{34}\\
0&0&0&1
\end{pmatrix}=
\\=
\begin{pmatrix}
 1&t_{12}+a_{12} &t_{13}+a_{12}t_{23}+a_{13} &t_{14}+a_{12}t_{24}+a_{13}t_{34}+a_{14}\\
 0&1&t_{23}+a_{23} & t_{24}+a_{23}t_{34}+a_{24}\\
 0&0&1&t_{34}+a_{34}\\
 0&0&0&1
\end{pmatrix}.
\end{multline*}
мы получаем формулу для левого сдвига. Соответственно, наша дифференциальная форма переходит
в
\begin{multline*}
d(t_{12}+a_{12})\wedge d(t_{13}+a_{12}t_{23}+a_{13})\wedge
 d(t_{14}+a_{12}t_{24}+a_{13}t_{34}+a_{14})\wedge
\\ \wedge d(t_{23}+a_{23})\wedge d(t_{24}+a_{23}t_{34}+a_{24})\wedge d(t_{34}+a_{34}).
\end{multline*}
Так как $a_{ij}$ - постоянные, $da_{ij}=0$. Заметим, что в нашем произведении присутствуют 
сомножители 
$$
dt_{12}=d(t_{12}+a_{12}), \quad dt_{23}=d(t_{23}+a_{23}),\quad dt_{34}=d(t_{34}+a_{34}).
$$
Поэтому под другими знаками дифференциалов
переменные $t_{12}$, $t_{23}$, $t_{34}$ могут быть убраны, и мы приходим к
$$
dt_{12}\wedge dt_{13}\wedge
 d(t_{14}+a_{12}t_{24})\wedge
dt_{23}\wedge dt_{24}\wedge dt_{34}.
$$
Теперь в качестве множителей появились $dt_{13}$ и $dt_{24}$, и мы повторяем тот же довод.
В случае произвольного $n$ надо было бы сказать <<и так далее>>, а в нашем случае вычисление на этом кончилось.

Разумеется, дифференциальных форм можно было бы не писать, а заметить, что матрица Якоби
левого сдвига как преобразования пространства $\R^{n(n-1)/2}$ является строго верхнетреугольной.

Легко видеть также, что наша мера правоинвариантна.


\sm

2) Рассмотрим группу вещественных матриц вида 
$\begin{pmatrix}u&v\\0&1 \end{pmatrix}$, где $u>0$. Тогда при левом сдвиге 
$$
\begin{pmatrix}u&v\\0&1 \end{pmatrix}\mapsto
\begin{pmatrix}a&b\\0&1 \end{pmatrix}\begin{pmatrix}u&v\\0&1 \end{pmatrix}
=\begin{pmatrix}au&av+b\\0&1 \end{pmatrix}
$$
дифференциальная форма $du\wedge dv$ переходит в 
$$
d(au)\wedge d(av+b)= a^2 du\wedge dv
$$
Соответственно, левоинвариантная форма объема задается формулой 
$$u^{-2}du\wedge dv.$$

\sm

{\sc Задача.} Вычислить правоинвариантную меру Хаара на той же группе и убедиться, 
что она отлична от левоинвариантной. \hfill $\lozenge$

\sm 

{\sc Задача.} Вычислить левоинвариантную меру Хаара на группе верхнетреугольных вещественных
(комплексных) матриц.

\sm

3) Рассмотрим группу $\SL(n,\R)$. Для элемента $X\in \SL(n,\R)$ общего положения
матричный элемент $x_{nn}$ однозначно определяется остальными матричными элементами
$x_{ij}$ из уравнения $\det (X)=1$. Фактически, мы получаем координаты на открытом плотном подмножестве группы.

\sm

 {\sc Задача.} Обозначим через $[X]_{n-1}$ левый верхний угол матрицы $X$ размера $(n-1)\times (n-1)$.
   Двусторонне инвариантная мера Хаара на $\SL(n,\R)$ задается формулой
   $$
  \bigl( \det[X]_{n-1}\bigr)^{-1} \prod_{i,j:\, (i,j)\ne (n,n)} dx_{ij}
  .
   $$
   
{\bf\punct Мера Хаара на ортогональной группе.}
Во всех ранее разобранных случаях на группах легко вводились координаты. Координатное
выражение для меры Хаара на $\SO(n)$ мы напишем в \S4, и отсюда, в частности, будет
следовать ее существование. Но здесь также работает простой общий прием.

\begin{theorem}
На группе $\SO(n)$ существует единственная, с точностью до пропорциональности,
левоинвариантная дифференциальная форма максимальной степени.
\end{theorem}

Это высказывание очевидно. {\it Единственность}. В силу инвариантности, эта форма
определяется своим значением в точке 1. В силу максимальности степени, такое значение
единственно с точностью до пропорциональности.

{\it Существование}. Мы фиксируем значение формы в точке 1, а потом посредством
левых сдвигов переносим ее во все точки группы. \hfill $\square$

\sm 

Легко видеть, что в этом доказательстве мы почти не использовали то, 
что наша группа это именно $\SO(n)$, нам нужно было лишь то, чтобы 
группа была многообразием.

\begin{theorem}
Левоинвариантная мера Хаара на $\SO(n)$ единственна.
\end{theorem}

{\sc Доказательство.} Из предыдущих рассуждений следует, что мера Хаара единственна
в классе мер, задаваемых гладкими дифференциальными формами. Покажем, что
любая левоинвариантная мера Хаара $\nu$
задается  гладкой формой. Введем произвольным образом координаты в окрестности $\cO$ единицы
 в $\SO(n)$.
Это также нормализует меру Лебега $d \lambda(g)$ в этой окрестности. Выберем
неотрицательную
$C^\infty$-функцию $f_\epsilon$ на $\SO(n)$ с носителем в $\epsilon$-окрестности
единицы, такую, что $\int f_\epsilon(g)\,d\lambda(g)=1$. Рассмотрим  меру
$$
\int_\cO\bigl (f_\epsilon(h)\,\nu(h^{-1}g)\bigr) d\lambda(h)
$$
Так как $\nu(h^{-1}g)=\nu(g)$, этот интеграл равен $\nu(g)$. С другой стороны, мы получили
гладкую меру (убедитесь в этом), и тем самым мера $\nu(g)$ - гладкая.
\hfill $\square$

\sm

См. \cite[\S7]{Zhe1},  \cite[5.2]{Ros}.

\sm

{\bf\punct Группы Ли.} 

\begin{definition}
Группа Ли -- это $C^\infty$-гладкое многообразие $G$, снабженное групповой операцией умножения. При этом
произведение должно быть гладким отображением $G\times G\to G$, а взятие  обратного 
элемента -- гладким отображением $G\to G$.
\end{definition}

Разумеется, группы Ли являются топологическими группами. 
Теория групп Ли -- отдельная история, которой мы в этих записках,
по возможности, не касаемся. Отметим лишь, что класс групп Ли можно описать более осязаемым образом.

\begin{theorem}
Любая замкнутая подгруппа в группе $\GL(n,\R)$ является гладким многообразием%
\footnote{Более того, это подмногообразие вещественно аналитично.}, 
а тем самым и группой Ли.
\end{theorem}

\begin{theorem} {\rm (Адо)}
Для любой группы Ли $G$ существует гомоморфизм $\psi$ из $G$ в достаточно большую
группу $\GL(N,\R)$, который задает инъективное отображение окрестности единицы
группы $G$ в окрестность единицы группы $\GL(N,\R)$%
\footnote{Простое доказательство этой теоремы есть в \cite{NAdo}.}.
\end{theorem}

{\sc Замечание.}
Вложить группу Ли как замкнутое подмногообразие в группу $\GL(N,\R)$ иногда
 не удается. Бывают случаи, 
когда у всех гомоморфизмов $G\to \GL(n,\R)$ есть неустранимое ядро (оно является
некоторой дискретной подгруппой в центре группы $G$%
\footnote{Для тех, кто проходил топологию. Пространство
$\SL(n,\R)$ гомеоморфно произведению $\SO(n)$ на односвязное пространство
строго положительно определенных матриц (что вытекает из полярного разложения
матриц).
У группы $\SO(n)$ есть двулистное накрытие -- спинорная группа.
 Поэтому двулистное накрытие есть у $\SL(n,\R)$. Однако реализовать его как подгруппу в $\GL(N,\R)$
невозможно. Действительно (для тех кто проходил теорию представлений), 
все конечномерные представления алгебры Ли $\frs\frl(n)$ известны, и все они
 соответствуют однозначным представлениям группы $\SL(n,\R)$.}). 
 
 Кроме того, для многих групп Ли возможны
 <<неудачные вложения>> -- может возникнуть эффект типа <<обмотки тора>> -- группа
$G$ вкладывается в $\GL(N,\R)$, но образ вложения оказывается не замкнут
(понятно, что это можно проделать уже для группы $\R$).
\hfill $\lozenge$

\sm

{\sc Замечание.} Топологическая группа, являющаяся топологическим многообразием,
 автоматически
является группой Ли. Когда-то этот вопрос составлял содержание 5ой проблемы Гильберта.
 \hfill $\lozenge$

\begin{theorem}
На любой группе Ли существует единственная с точности до пропорциональности левоинвариантная
дифференциальная форма максимальной степени.
\end{theorem}

Доказательство было выше.

См. \cite[\S\S3-7]{Zhe1},\cite[гл. 5]{Zhe2},\cite[\S6]{Kir}, \cite[Ch.2, \S4.3]{Ros},
 \cite[Add. C]{Ner2}

\section[Существование меры Хаара]{Разные доказательства существования меры Хаара.} 

\COUNTERS

В этом параграфе приводятся четыре доказательства существования меры Хаара. Два из них годятся
для компактных групп, другие два -- в общем случае. Самое простое по идее, но утомительное
в деталях -- оригинальное доказательство Хаара. Доказательство А.Вейля -- его модификация,
в которой  рассуждения с мерами  заменяются на теорему Рисса--Маркова, \cite[IV.18]{RS1}.
Два доказательства для компактных групп  основаны на той же теореме и на доводах
типа <<теорем о неподвижной точке>>. Во всех доказательствах мы будем
полагать группы сепарабельными.

\sm

{\bf\punct  Последовательности Коши, равномерная непрерывность и т.п.}
Для топологических групп имеет смысл ряд понятий, которых нет для общих топологических
пространств. 

\sm

1) Последовательность $x_j\in G$ называется {\it последовательностью Коши}, если
$x_i^{-1} x_j$ стремится к единице, если $i$, $j\to\infty$.

\sm

2) Есть аналог понятия $\epsilon$-{\it сети}. А именно, пусть $U$ - подмножество в $G$,
а $\cO$ -- окрестность единицы в $G$. Мы скажем, что {\it $\cO$-сеть} - это покрытие множества
$U$ подмножествами вида $x_\alpha \cO$.

\sm

3) Функция на $f$ на $G$ {\it равномерно непрерывна}, если для для любого $\epsilon>0$
существует окрестность единицы $\cO$, такая, что при $x^{-1} y\in \cO$
выполнено $|f(x)-f(y)|<\epsilon$.

\sm

4) Семейство функций $f_\alpha$ на $G$ {\it равностепенно непрерывно},
 если для для любого $\epsilon>0$
существует окрестность единицы $\cO$, такая, что при $x^{-1} y\in \cO$
выполнено  $|f_\alpha(x)-f_\alpha(y)|<\epsilon$. Для компактных
групп остается в силе {\it теорема Асколи--Арцела}: подмножество  пространства непрерывных
функций компактно, если оно ограничено и равностепенно непрерывно.

\sm

Это означает, что на топологической группе определена {\it равномерная структура} в смысле 
Н.Бурбаки и А.Вейля. Впрочем, мы не будем пользоваться абстрактным  понятием равномерной структуры. 

Отметим, что $x_i^{-1} x_j$ можно заменить на $x_i x_j ^{-1}$,
$x^{-1} y$ на $x y^{-1}$, а сдвинутые окрестности $x\cO$ на $\cO x$. Это, вообще говоря,
 приведет к неэквивалентным понятиям ({\it левая и правая равномерные структуры}).
В случае
компактных групп этого различия не возникает.

\sm

См. \cite[гл. 2; гл.3, \S2 ]{Bou1}, \cite[2.6]{Kir}, \cite[\S4]{HR1}.

\sm

{\bf\punct Инвариантные метрики.}

\begin{theorem} {\rm(}Г.Биркгоф--Какутани{\rm)}
{\rm a)}
На топологической группе есть левоинвариантная метрика тогда и только тогда, когда на 
группе есть счетная база окрестностей  единицы.

\sm

{\rm b)} В случае компактных групп эта метрика может быть выбрана 
двусторонне инвариантной.

\end{theorem}

{\sc Набросок доказательства.} a) Мы выбираем последовательность $V_1$, $V_2$, $V_3$,
\dots окрестностей единицы такую, что $V_j^2\subset V_{j-1}$. Для любого
числа $r=2^{-k_1}+2^{-k_2}+\dots + 2^{-k_m}$, где $k_1<k_2<\dots<k_m$ положим 
$$
U_r:=V_{k_1}V_{k_2}\dots V_{k_m}
$$
Определим функцию $\psi$ на $G$ как
 супремум $r$ таких, что $g\notin U_r$. 

Инвариантная метрика определяется как 
$d(x,y)=\sup_{g\in G}|\phi(gx)-\phi(gy)|$.

\sm

b) Мы полагаем $d(x,y)=\max_{g\in G} d(xg,yg)$.
\hfill $\square$

\sm

Отметим, что явные выражения для левоинвариантных метрик используются 
довольно редко, однако всегда
приятно знать, что топология метризуема.

\sm

{\sc Задача.} На группе $\SL(2,\R)$ не существует двусторонне инвариантной метрики, совместимой с топологией. {\sc Указание.} Если бы такая метрика $d$ существовала, 
 она была бы удовлетворяла
условию $d(1,x)=d(a^{-1}x a)$.
\hfill $\lozenge$

\sm

См. \cite[\S8]{HR1}. 

\sm

{\bf \punct Единственность меры Хаара для компактных групп.}
Мы доказываем такое условное утверждение. {\it Если на компактной группе есть
вероятностная левоинвариантная мера Хаара, то она единственна, двусторонне
инвариантна, а также  инвариантна относительно отображения $g\mapsto g^{-1}$.}

Сначала заметим, что отображение $g\mapsto g^{-1}$ переводит левоинвариантную меру в правоинвариантную. Далее, возьмем правоинвариантную  меру Хаара $\mu$ и левоинвариантную
меру $\nu$. Положим, что эти меры - вероятностные. Мы утверждаем, что они совпадают.
Возьмем их свертку $\mu*\nu$, т.е., меру, определяемую из условия:
$$
\int_G f(z)\,d(\mu*\nu)(z)= \int_{G\times G} f(xy)\,d\mu(x)\,d\nu(y)
$$
для любой непрерывной функции $f$. В силу правой инвариантности $\mu$ этот
 интеграл  равен
 $$
 \int_{G} \Bigl( \int_{G} f(x)\,d\mu(x)\Bigr)\,d\nu(y) = \int_{G} f(x)\,d\mu(x).
 $$
 Т.е., $\mu*\nu=\mu$. Аналогично,  $\mu*\nu=\nu$.
 \hfill$\square$

\sm

{\bf\punct Первое доказательство существования в случае компактных групп (фон Нейман, 1934).}
Пусть $G$ -- компактная группа, пусть $C(G)$ -- пространство непрерывных
вещественных
функций $G$ с обычной нормой $\max|f|$. 
 Обозначим через $L_h$, $R_h$ операторы левого и правого сдвига
в $C(G)$, 
$$
L_h f(g)= f(hg), \qquad R_h f(g)=f(gh).
$$
Очевидно, что левые и правые сдвиги коммутируют,
$$
L_h R_q= R_q L_h.
$$
Мы построим линейный функционал $I(f)$ на
$C(G)$, инвариантный относительно сдвигов,
$$
I(f)=I(L_hf)=I(R_q f),
$$
 а потом сошлемся на теорему Рисса--Маркова.

Фиксируем функцию $f$ на группе и возьмем множество $S$ всех ее левых сдвигов.
По теореме Арцела, $S$ является компактом. Обозначим через$K_f=K_f^{left}$  замкнутую
выпуклую оболочку множества $S$. В силу той же теоремы Арцела, $K_f=K_f^{left}$ -- компакт.

Очевидно также, что для любой функции $\phi\in K_f$ выполнено $K_\phi\subset K_f$.

\begin{lemma}
Множество $K_f=K_j^{left}$ содержит постоянную функцию.
\end{lemma}

{\sc Доказательство.}
Рассмотрим cледующий непрерывный функционал на $C(G)$:
$$
Z\phi=\max(\phi) - \min(\phi).
$$
Пусть $\psi$ -- точка его минимума на $K_f$, допустим, что $\psi$
не является константой. Пусть $\Xi\subset K$ -- (компактное) множество точек, где $\psi$
достигает своего максимума $M$. Возьмем малое $\epsilon$.
Пусть $\cU$ -- множество, где $\psi(x)<\min (\psi)+\epsilon$. Рассмотрим конечное открытое покрытие множества
$\Xi$ сдвигами $g_j \cU$. Возьмем произвольную функцию вида
$$
\psi^*(g)= a\psi(h)+ \sum b_j \psi(g_jg), \qquad \text{где $a+\sum b_j=1$, $a>0$, $b_j>0$}
$$
Тогда $\psi^*\in K_\psi\subset K_f$. С другой стороны
$$
\max(\psi^*)<\max(\psi), \qquad \min(\psi^*)\ge\min (\psi).
$$
Таким образом, $Z(\psi^*)<Z(\psi)$, и мы приходим к противоречию.
\hfill $\square$

\sm 

Эта константа и является значением функционала $If$.

\sm

{\sc Задача.} Не ссылаясь на следующую лемму, докажите, что $K_f$ содержит единственную
константу. \hfill $\lozenge$

\begin{lemma}
Пусть константа $c$ содержится в $K_f^{left}$, а константа $c^\circ$ содержится в
$K^{right}_f$. Тогда $c=c^\circ$. В частности, $K_f^{left}$ содержит единственную
константу. 
\end{lemma}

{\sc Доказательство.} Пусть конечная выпуклая комбинация $\sum a_j f(h_j g)$ приближает
$c$ с точностью до $\epsilon$, а выпуклая комбинация $\sum b_i f(gq_i)$ приближает
$c^\circ$ с точностью до $\epsilon$. Рассмотрим операторы
$$
\cL:= \sum a_j L_{h_j}\qquad \cR:=\sum b_i R_{q_i}
.
$$
Тогда
$$
\cR \cL f= \cL\cR f.
$$
По построению, 
$$
c-\epsilon\le \cL f(g)\le c+\epsilon.
$$
Поэтому функция $\cR\cL f\in K^{left}_{\cL f}$ меняется в тех же пределах.
Аналогично, $c^\circ-\epsilon\le \cL\cR f\le c^\circ+\epsilon$.
Поэтому $|c-c^\circ|\le 2\epsilon$ для любого $\epsilon>0$. 
\hfill $\square$

\sm

\begin{lemma}
Функционал $If$ аддитивен: $I(f_1+f_2)=If_1+If_2$.
\end{lemma}

{\sc Доказательство.} Рассмотрим конечную выпуклую комбинацию
$\cL=\sum a_j L_j$ такую, что 
$$|\cL f_1- I(f_1)|\le\epsilon.$$
Отметим, что 
$$I (\cL f_2)= I(f_2),$$ 
так как 
$K^{left}_{\cL f_2}\subset K^{left}_{f_2}$.
Теперь рассмотрим выпуклую комбинацию $\cM$ левых сдвигов,
такую, что 
$$
|\cM \cL f_2 - I(f_2)|<\epsilon
$$
Условие
$$|\cM\cL f_1- I(f_1)|\le\epsilon.$$
выполняется автоматически. Поэтому
$$
|\cM\cL(f_1+f_2)-I(f_1)-I(f_2)|<2\epsilon
,$$
что завершает доказательство. \hfill $\square$

\sm

Таким образом, $I$ является непрерывным линейным функционалом
на $C(G)$, причем $f\ge 0$ влечет $If\ge 0$. В силу теоремы Рисса-Маркова
он задается положительной мерой на $G$. \hfill$\square$

\sm

См. \cite[Гл. 5, \S7]{Rud}.

\sm

{\bf\punct Другое доказательство существования  в случае компактных групп.}
Рассмотрим множество $\cK$ всех вероятностных мер на $G$, снабженное слабой топологией
($\mu_j\to \mu$, если для любой непрерывной функции $f$ последовательность
$\int f\,d\mu_j$ сходится к $\int f\,d\mu$). Известно, что $\cK$  компактно
(например, это следует из теоремы Банаха--Алаоглу о $*$-слабой компактности
шара в пространстве, сопряженном к банахову). Стоит заметить также, что множество $\cK$
метризуемо (докажите это).

Группа $G$ действует на $\cK$ сдвигами. Действие это непрерывно в следующем смысле.
Если $g_j$ стремится к $g$, то $g_j \mu$ стремится к $g\mu$ (докажите это, высказывание
 следует из равномерной непрерывности непрерывных функций). В частности, орбиты группы
  $G$ компактны.

 Мы хотим показать, что группа $G$ имеет на $\cK$ неподвижную точку.
 
 Рассмотрим  всевозможные выпуклые замкнутые подмножества в $\cK$, инвариантные
 относительно действия группы $G$.
 Множество $\frS$ таких множеств является частично упорядоченным
 множеством  относительно включения.
  Возьмем в нем минимальный элемент $L$. Возьмем
 в $L$ две различных точки, $\mu$, $\nu$. Возьмем орбиту $\cO$ точки $(\mu+\nu)/2$ 
 относительно группы $G$. Ее замкнутая выпуклая оболочка должна совпадать с $L$.
 Возьмем крайнюю точку $\xi$ множества $L$ (по теореме Крейна--Мильмана,
 такие точки существуют).
 В силу компактности множества $\cO$, мы имеем $\xi\in \cO$
 (см. \cite{Phe}, Следствие из теоремы 1.4). Поэтому $\xi$ имеет вид
 $(g\mu+g\nu)/2$ c $g\mu$, $g\nu\in L$. А поэтому наша точка -- не крайняя.
 
 Итак, если $L\in \frS$ состоит из более чем одной точки, то $L$ не минимально.
С другой стороны, ясно, что минимальные элементы в $\frS$ есть.
Действительно, возьмем какой-нибудь элемент $L_1\in \frS$. Возьмем $L_2\in \frS$,
строго содержащийся в $L_1$, и т.д. Таким образом, мы получаем убывающую
цепочку подмножеств $L_1\supset L_2\supset\dots$. В силу компактности,
$$
L_\infty:=\cap_{j=1}^\infty L_j\,\in \frS
$$
не пусто
Теперь мы уменьшим $L_\infty$ и т.д., пока мы не дойдем до минимального
элемента. Это довод трансфинитной индукции,
его оформление -- вопрос вкусовой. Проще всего сослаться на теорему
Хаусдорфа, о том, что любое линейно упорядоченное подмножество
частично упорядоченного множества содержится в максимальном линейно
упорядоченном подмножестве,
см. \cite[Добавл. A.1]{Rud}. \hfill$\square$

 \sm
 
 {\sc Задача.} Пусть компактная группа $G$  действует на компактном пространстве
 $M$. Покажите тем же способом что на $M$ есть инвариантная вероятностная мера. 
 \hfill $\lozenge$
 
 \sm
 
 {\bf\punct Единственность меры Хаара для локально компактных групп.}
 См. \cite{Bou1}, VII.8.
 Пусть $\mu$ -- некоторая левоинвариантная мера Хаара, $\nu$ -- некоторая правоинвариантная
 мера Хаара. 
  Пусть $f$, $g$ -- непрерывные функции на $G$ с компактным носителем.
 
 \sm

{\sc Задача.} a) Покажите, что функция 
$
h(x)=
\int_G f(x^{-1}y)\,d\nu(y)
$ 
непрерывна на $G$.

\sm

b) Покажите, что для любой ненулевой неотрицательной функции
$f$ выполнено $\int f\,d\mu>0$. \hfill $\lozenge$ 

\sm

Используя инвариантность мер, мы получаем
\begin{multline*}
\int f(x)\,d\mu(x)\cdot \int g(y)\,d\nu(y)= \iint f(x) g(y)\,d\nu(y)\,d\mu(x)
=\\=
\int f(x)\Bigl[\int g(y)\,d\nu(y)\Bigr]\,d\mu(x)=\int f(x)\Bigl[\int g(yx)\,d\nu(y)\Bigr]\,d\mu(x)
=\\=
\int \Bigl[\int  f(x) g(yx)\,d \mu(x)\Bigr]\,d \nu(y)= 
\int\Bigl[ \int  f(y^{-1}x) g(x) \,d\mu(x) \Bigr] \,d \nu(y)
=\\=
\int \Bigl[\int  f(y^{-1}x) \,d\nu(y)\Bigr]\,g(x) \,d\mu(x).
\end{multline*}
В левой и правой частях этого равенства стоят линейные функционалы от $g$, которые
задаются мерами
$$
\Bigl(\int f(x)\,d\mu(x)\Bigr)\cdot \,d\nu
\quad\text{и}\quad
\left[\int  f(y^{-1}x) \,d\nu(y)\right]\cdot \,d\mu.
$$
Поэтому мы имеем совпадение мер
$$
d\nu=\biggl[\biggl(\int f(z)\,d\mu(z)\biggr)^{-1} \int  f(y^{-1}x) \,d\nu(y)\biggr]\cdot d\mu
.$$
В квадратных скобках стоит непрерывная функция $S[\mu,\nu,f](x)$ от $x$, которая не может зависеть от $f$. Далее замечаем, что $\nu$ абсолютно непрерывна относительно
$\mu$, тот же довод показывает, что $\mu$ абсолютно непрерывна относительно
$\nu$, т.е., эти меры эквивалентны.

 Фиксируем  $\nu$ и ненулевую неотрицательную функцию $f$ 
 (от которой ничего не зависит).
 Тогда $\int  f(y^{-1}x) \,d\nu(y)$ однозначно определен, а поэтому
 $\mu$ восстанавливается однозначно по $\nu$, с точностью до
постоянного множителя $\int f(z)\,d\mu(z)$. Т.е., мера $\mu$
единственна с точностью до пропорциональности.
\hfill $\square$

\sm

См.   См. \cite[VII.8]{Bou1}, \cite[\S26]{Zhe2}, \cite[\S7]{Wei}, 
\cite[\S60]{Hal},
\cite{Nach}.

\sm

{\bf\punct Существование меры Хаара на локально компактных группах. Доказательство А.Вейля.}
Обозначим через $\cC(G)$ пространство не\-пре\-рыв\-ных вещественнозначных
 функций с компактным носителем
на $G$. Обозначим через $\cC(K)\subset \cC(G)$ пространство непрерывных функций
с носителем, содержащимся в компакте $K$.  Мы скажем, что последовательность $f_j\in\cC(G)$ сходится к $f$,
если носители всех функций $f_j$ содержатся в некотором компакте и $f_j\to f$ равномерно.
Через $\cC_+\subset \cC(G)$ мы обозначим конус неотрицательных функций. 
Мы докажем существование инвариантного относительно левых сдвигов
линейного функционала на $\cC(G)$, после чего останется сослаться на теорему Рисса--Маркова,
см. \cite{RS1}, IV.18.

Для ненулевой функции $\phi\in \cC_+$ определим функционал $I_\phi f$ на $\cC_+$ следующим образом.
Мы рассматриваем всевозможные конечные линейные комбинации%
\footnote{$L_{g}$ -- по-прежнему, оператор левого сдвига.} 
$\sum c_j L_{g_j} \phi$ с неотрицательными коэффициентами такие, что
\begin{equation}
f(x)\le \sum c_j L_{g_j} \phi(x).
\label{eq:est}
\end{equation}
Далее, положим
\begin{equation*}
I_\phi f=\inf \sum c_j,
\end{equation*}
где точная нижняя грань берется по всем оценкам (\ref{eq:est}).

\sm

{\sc Замечание.} Здесь уместно остановиться и понять, что это означает в случае $G=\R$.
Пусть носитель неотрицательной функции $\phi_\epsilon$  содержится в маленьком отрезке
 $[-\epsilon,\epsilon]$ и $\int \phi_\epsilon(x) \,dx=1$. 
 Тогда свертка
$$
f*\phi_\epsilon (x)=\int f(y)\phi_\epsilon(x-y)\, dy
$$
примерно равна $f(x)$. Далее заменяем интеграл на интегральную сумму,
и получаем выражение вида $f(a_j)(a_j-a_{j-1})\phi_\epsilon(x-a_j)$, 
примерно равное все той же
$f(x)$. Теперь нужно чуть-чуть увеличить коэффициенты и получить оценку
$f(x)\le \sum c_j \phi_\epsilon(x-a_j)$.  \hfill $\lozenge$

\sm

{\sc Задача.} Докажите аккуратно, что в случае $G=\R$ имеет место
$$
\qquad\qquad\qquad\qquad\qquad\qquad
\lim_{\epsilon\to 0}I_{\phi_\epsilon} f=\int f(x)\,dx.
\qquad\qquad\qquad\qquad\qquad\qquad
\lozenge
$$

\sm

Вернемся к доказательству теоремы. Понятно, что $I_\phi (af)=a I_\phi(f)$, для любого
числа $a>0$. Функционал $I_\phi$ также полуаддитивен,
$$
I_\phi(f_1+f_2)\le I_\phi f_1+ I_\phi f_1.
$$
Легко видеть, что для $\phi$, $\alpha$, $f\in \cC_+$
$$
I_\phi f\le I_\phi \alpha\cdot I_\alpha f
.$$
Переставив буквы, мы напишем также
\begin{equation}
I_\phi \alpha\le I_\phi f \cdot I_f \alpha.
\end{equation}
Теперь {\it мы фиксируем эталонную функцию $\alpha$ и будем нормировать
всевозможные $\phi$ условием $I_\phi \alpha=1$}. В силу двух последних
неравенств, $I_\phi f$ может меняться лишь в пределах
\begin{equation}
(I_f \alpha)^{-1}\le I_\phi f \le I_\alpha f.
\label{eq:vilka}
\end{equation}

\begin{lemma}
\label{l:weil-1.5}
Для любого компакта $K$ семейство функционалов $I_\phi$ равностепенно
непрерывно на $\cC_+(K)$.
\end{lemma}

{\sc Доказательство.} Возьмем функцию $h\in \cC_+(G)$, равную 1 на $K$.
Пусть $f_1-f_2\le\epsilon$. Тогда 
$$I_\phi f_1- I_\phi f_2\le \epsilon I_\phi h\le \epsilon I_\alpha h.$$
Далее, меняем местами $f_1$ и $f_2$.
\hfill $\square$

\begin{lemma}
\label{l:weil-2}
Для любых $f_1$, $f_2\in\cC_+(G)$ для любого $\epsilon>0$ существует окрестность
$U$ единицы, такая, что для любой функции $\phi$  с носителем в $U$
выполнено
\begin{equation} 
I_\phi f_1+ I_\phi f_2 \le  I_\phi (f_1+ f_2)+\epsilon
.\end{equation}
\label{eq:III}
\end{lemma}



{\sc Доказательство.} Возьмем функцию $h\in \cC_+$, равную 1 на
носителях $f_1$, $f_2$. Мы хотим показать, что для любого
$\lambda>0$ и любого $\epsilon>0$ следующая оценка
\begin{equation}
I_\phi f_1+ I_\phi f_2 \le  I_\phi (f_1+ f_2+\lambda h)+\epsilon
\label{eq:promezhut}
\end{equation}
выполнена для любых $\phi$ с достаточно малым носителем.
Обозначим
$$
\xi_1= \frac{f_1}{f_1+ f_2+\lambda h}, 
\qquad \xi_2= \frac{f_2}{f_1+ f_2+\lambda h}.
$$
Возьмем $\delta>0$ и окрестность $U$ единицы, такую, что
$$
|L_g \xi_1-\xi_1|<\delta,\qquad |L_g \xi_2-\xi_2|<\delta
$$
при $g\in U$. Пусть носитель $\phi$ содержится в $U$. Пусть 
$$
f_1+ f_2+\lambda h\le \sum c_j L_{g_j} \phi.
$$
 Тогда
\begin{multline*}
f_1(x)=(f_1(x)+ f_2(x)+\lambda h(x))\xi_1(x)\le \sum c_j [\xi_1(x) \phi (g_j x)]
\le\\ 
\le \sum c_j [(\xi_1(g_j)+\delta) \phi (g_j x)]
=
\sum c_j \Bigl(\frac{f_1(g_j)}{ (f_1(g_j)+ f_2(g_j)+\lambda h(g_j)}+\delta\Bigr)
  \phi (g_j x).
\end{multline*}
Аналогично выписывается неравенство для $f_2$. В итоге получаем
$$
I_\phi f_1+I_\phi f_2\le \sum c_j
\Bigl(\frac{f_1(g_j)+f_2(g_j) }{ (f_1(g_j)+ f_2(g_j)+\lambda h(g_j)}+2\delta\Bigr)
\le \sum c_j(1+2\delta)
.
$$
и мы получаем желаемую оценку (\ref{eq:promezhut}).

Теперь мы устремляем $\lambda$ к нулю, ссылаемся на лемму \ref{l:weil-1.5}
и получаем
(\ref{eq:III}).
\hfill $\square$

\sm

Завершим доказательство теоремы. Выберем возрастающую последовательность
компактов $K_1\subset K_2\subset\dots$, исчерпывающую  $G$. В каждом $\cC_+(K_m)$
выберем счетное плотное множество. Тем самым мы получим счетное плотное подмножество
$\Omega$ в $\cC_+(G)$. Возьмем фундаментальное семейство $U_1\supset U_2\supset\dots$ окрестностей единицы и последовательность функций $\phi_m$ с носителем в $U_m$.

Мы утверждаем, что из последовательности функционалов $I_{\phi_m}$
можно выбрать подпоследовательность, сходящуюся во всех точках $f^{(\alpha)}\in\Omega$.
Это делается обычной <<диагональной процедурой>>.
Ввиду (\ref{eq:vilka}), мы можем выбрать подпоследовательность, сходящуюся в точке
 $f^{(1)}$. Из нее мы выбираем подпоследовательность, сходящуюся в точке  $f^{(2)}$
и т.д. Мы получаем убывающую цепочку подпоследовательностей. Далее выбираем
по номеру $m_\alpha$ из $\alpha$-ой подпоследовательности так, что
$m_{\alpha} > m_{\alpha-1}$. В итоге мы получаем подпоследовательность
с желаемыми свойствами.

В силу равностепенной непрерывности, последовательность функционалов
$I_{\phi_{m_\alpha}}$ сходится на всем $\cC_+(G)$, а предельный функционал
$I$ непрерывен. В силу   (\ref{l:weil-2}), мы имеем
$$
I f_1+I f_2=I(f_1+f_2).
$$
Легко понять, что этот функционал продолжается по линейности
с $\cC_+(G)$ на все пространство $\cC(G)$. В силу
 (\ref{eq:vilka}), этот функционал положителен на $\cC_+$.
По теореме Рисса--Маркова он задается борелевской мерой.

\sm

См. \cite[гл. 7, \S1]{Bou2}, \cite[\S7]{Wei}, \cite[\S26]{Zhe2},
\cite{Nach}, \cite[\S15]{HR1}.

\sm

{\bf\punct Существование меры Хаара на локально компактных группах. Доказательство Хаара.}
Пусть $X$ -- локально компактное метризуемое сепарабельное пространство. Мы скажем,
что 
{\it объем%
\footnote{Англ. content.}}  (см. \cite[\S\S53-54]{Hal})  на $M$ -- это функция, которая каждому компактному множеству $K$ ставит в
 соответствие  число $\lambda(K)$ так, что
 
 \sm
 
 1) $\lambda(K)\ge 0$;
 
  \sm
 
 2) если $K\supset M$, то $\lambda(K)\ge \lambda(M)$;
 
  \sm
 
 3) для любых $K$, $L$ выполнено $\lambda(K\cup L)\le\lambda(K)+\lambda(L)$;
 
  \sm
 
 4) если $K$, $L$ не пересекаются, то $\lambda(K\cup L)=\lambda(K)+\lambda(L)$. 
 
  \sm

Функция объема, вообще говоря, до счетно аддитивной меры не продолжается. 

\sm

{\sc Пример.} Рассмотрим пространство, состоящее из сходящейся последовательности
 и ее предела. Положим, что конечные подмножества имеют меру ноль, а счетные замкнутые подмножества -- меру 1.
 \hfill $\lozenge$
 
\sm

Однако по функции объема можно следующим способом построить меру. Для любого
открытого множества $U$ положим
$$
\mu(U)=\sup_{\text{$K\subset U$, $K$ компактно}} \lambda(K).
$$
 По $\mu$ обычным способом строится внешняя мера, которая  оказывается
 счетно-аддитивной мерой на борелевской $\sigma$-алгебре.
  На компактных множествах она может   отличаться от функции объема,
$\mu(K)\ge\lambda(K)$.

\sm
 
Теперь возьмем локально компактную группу $G$. Пусть $U$ -- открытое ограниченное подмножество в $G$,
содержащее 1,   а $E$ --  ограниченное множество (т.е., множество с компактным замыканием). Рассмотрим всевозможные покрытия
$E$ левыми сдвигами $x_j U$ множества $U$.  Минимальное возможное число элементов покрытия
мы обозначим через $E:U$.

\sm

{\sc Пример.} В качестве  $G$  возьмем $\R^2$, а в качестве $U$ -- круг радиуса $r$.
Что получится? \hfill $\lozenge$ 

\sm

Величина $E:U$ как функция от $E$ обладает свойствами 1)-3) объема, а свойство аддитивности
4) выполнено если $K$, $L$ удовлетворяют $K\cap ULU^{-1}=\varnothing$.  

Фиксируем множество $V$, и для любого множества $U$ положим
$$
\lambda_U(E)=\frac{E:U}{V:U}.
$$
Заметим, что 
$$
E:U\le  (E:V)\cdot (V:U).
$$
Поэтому
\begin{equation}
O\le \lambda_U(E)\le E:V.
\label{eq:tih1}
\end{equation}
При этом
\begin{equation}
\lambda(V)=V:E.
\label{eq:tih2}
\end{equation}

Теперь надо применить теорему Тихонова. Рассмотрим множество $\Phi$
всех функций $f(E)$ на пространстве всех замкнутых подмножеств, удовлетворяющих
условиям (\ref{eq:tih1})--(\ref{eq:tih2}). Оно компактно. Для любого $U$ обозначим через
$\Delta(U)$ множество всех функций $\lambda_V$ c
$V\subset U$. Тогда для любого конечного набора $U_1$, \dots, $U_k$
выполнено
$$
\Delta(U_1)\cap\dots\cap \Delta(U_k)\supset \Delta(U_1\cap\dots\cap U_k)
,$$
в частности, пересечение непусто. Замыкания $\ov{\Delta(U)}$ множеств $\Delta(U)$
компактны (а их конечные пересечения, по-прежнему, не пусты)
 Поэтому есть точка, общая для всех множеств $\ov{\Delta(U_1)}$. Несложно
 показать, что это и будет инвариантным объемом.

\sm

См. \cite{Haa}, \cite[\S58]{Hal}, \cite[\S9.11]{Bog}.

\section[Инвариантные меры на однородных пространствах]%
{Модулярный характер. Инвариантные меры на однородных пространствах}

\COUNTERS

{\bf\punct Модулярный характер.}
Пусть $\nu$  -- левоинвариантная мера Хаара на группе $G$.
Легко видеть, что ее правый сдвиг $R_g\nu$ является левоинвариантной мерой.
Следовательно, эти меры различаются числовым множителем
$$R_g\nu=\Delta(g)\nu.$$
 C другой стороны, $R_{g_1}R_{g_2}=R_{g_1g_2}$.
Поэтому
$$
\Delta(g_1)\Delta(g_2)=\Delta(g_1 g_2),
$$
т.е., мы получаем гомоморфизм $\Delta_G$ из $G$
в мультипликативную группу $\R^\times_+$ положительных вещественных чисел, он называется
{\it модулярным характером}. 

В конкретных случаях модулярный характер легко вычисляется. Группа называется 
{\it унимодулярной} если он равен единице, т.е., если мера Хаара двусторонне инвариантна.

\sm

{\sc Задача.} Покажите, что мера $\Delta_G(g)^{-1}\nu$ правоинвариантна.
 Покажите, что она совпадает с образом $\nu$ при отображении $g\mapsto g^{-1}$.
\hfill $\lozenge$ 

\sm

{\sc Пример.} У группы $\SL(2,\R)$ гомоморфизмов $\R^\times_+$
нет. Поэтому $\SL(2,\R)$  унимодулярна.
 Легко проверить, что группа $\SL(n,\R)$ порождается подгруппами, изоморфными
 $\SL(2,\R)$ (мы вкладываем $\SL(2,\R)$ в качестве диагонального блока). Поэтому
 $\SL(n,\R)$ унимодулярна. Похожим образом показывается  унимодулярность
 $\SO(p,q)$, $\SL(n,\C)$, $\SO(n,\C)$.
\hfill $\lozenge$

\sm
 
 {\sc Задача.} Вычислите модулярный характер группы верхнетреугольных
 вещественных матриц $B_n$ (одно из возможных соображений: модулярный характер
 равен 1 на коммутанте группы (т.е. на группе строго верхнетреугольных матриц
 $T_n$). Поэтому достаточно вычислить его на диагональных матрицах. Другой рецепт
 предлагается в следующем пункте. \hfill $\lozenge$
 
 \sm

См. \cite[VII.1.4]{Bou2}, \cite[\S8]{Wei}, \cite[26.9]{Zhe2}.
 
 \sm
 
 {\bf\punct Вычисление модулярных характеров для матричных групп.}
 Рассмотрим замкнутую подгруппу $G$ в $\GL(n,\R)$. Через $\frg$ обозначим касательное пространство к $G$ в единице. Сопряжение $g\mapsto hgh^{-1}$ на группе индуцирует
 оператор $\Ad(g):\frg\to \frg$ в касательном пространстве. Берем левоинвариантную
 дифференциальную форму $\omega$. Мы сдвигаем форму левым сдвигом из единицы 
в $g$. Далее сдвигаем ее обратно правым сдвигом. Она должна умножиться на
$\Delta(g)$. С другой стороны, она умножается на $\det(\Ad(g))$. Итак
$$\Delta(g)=|\det(\Ad(g))|.$$

\sm

{\sc Задача.} Вычислите таким способом модулярный характер
для $\GL(n,\R)$, группы верхнетреугольных матриц, для группы строго верхнетреугольных матриц.
\hfill $\lozenge$

\sm 

{\bf \punct Разложения групп.} Следующее высказывание часто дает способ явно записывать
меру Хаара в <<системах координат>>.

\begin{proposition}
	Пусть $G$ -- локально компактная унимодулярная группа, $K$, $H$ -- ее замкнутые подгруппы.
	Допустим, что отображение $K\times H\to G$, заданное формулой
	$\Pi:(k,h)\mapsto kh$, является биекцией с точностью до почти всюду.
	Пусть $d\lambda_K$ -- левая мера Хаара на $K$, а $d\rho_H$ -- правая мера Хаара на
	$H$.  Тогда образ меры $d\lambda_K\times d\rho_H$ при отображении $\Pi$
	совпадает с мерой Хаара на $G$.
\end{proposition} 

{\sc Доказательство.} Рассмотрим прообраз меры Хаара при отображении $\Pi$.
Получится мера на $K\times H$, инвариантная относительно левых сдвигов на $K$ и правых
сдвигов на $H$. Пусть $H^\circ$ -- группа, антиизоморфная $H$. Это то же  множество
$H$, на котором введено умножение $g\circ h:=hg$. Тогда наша мера становится левоинвариантной
мерой на $K\times H^\circ$, т.е., произведением мер Хаара.
\hfill $\square$.

\sm

{\sc Примеры.}  a) Пусть $G=\GL(n,\R)$, $N_-$ -- подгруппа строго нижнетреунольных матриц,
$B_+$ -- подгруппа верхнетреугольных матриц. Почти любой элемент $g\in G$ представим
в виде произведения $g=nb$, где $n\in N_-$, $b\in B_+$, причем такое разложение
единственно.

\sm

b) Пусть $G=\GL(p+q,\R)$, $\cN$ -- группа блочных $p+q$-матриц
вида $\begin{pmatrix}
1&0\\C&1
\end{pmatrix}$,
а $\cB$ -- группа блочных матриц вида $\begin{pmatrix}
A&B\\0&D
\end{pmatrix}$. Почти любой элемент $g\in G$ представим
в виде произведения $g=nb$, где $n\in \cN$, $b\in \cB$, причем такое разложение
единственно.

\sm

c)  Возьмем в $G=\GL(n,\R)$ ортогональную подгруппу $\O(n)$ и строго верхнетреугольную подгруппу
$B_+$. Тогда любой элемент $g\in G$ однозначно представим в виде произведения
$g=ub$, где $u\in \O(n)$, $b\in B_+$.

\sm

d) Возьмем группу  аффинных изометричных преобразований $\R^n$.
Любой элемент этой группы однозначно представим в виде вращения,
оставляющего ноль на месте (т.е., элемента $\O(n)$) и сдвига.
\hfill $\lozenge$

\sm

{\sc Задача.} Докажите эти высказывания. \hfill $\lozenge$

 \sm
 
{\bf \punct Инвариантные меры на компактных пространствах.}

\begin{theorem}
\label{th:invariant-compact}
Пусть $K$ - компактная группа, $H$ -- замкнутая подгруппа.
Тогда на однородном пространстве $K/H$
 существует единственная $K$-инвариантная 
вероятностная мера.
\end{theorem} 
 
{\sc Доказательство.}
Рассмотрим отображение $K\to K/H$. Берем образ меры Хаара%
\footnote{Пусть на пространстве $M$ определена мера
	$\mu$. Пусть $\psi:M\to N$ -- отображение. Мера $\nu$ 
	на  $N$ определяется из условия $\nu(A):=\mu(\psi^{-1} A)$.}
Обратно, пусть $\nu$ -- вероятностная инвариантная мера на $K/H$,
$\kappa$ -- вероятностная мера Хаара на $K$.
Рассмотрим непрерывную функцию $f(x)$ на $K/H$.
Тогда 
$$
\int_{K/H} f(kx)\,d\nu(x)=\int_{K/H}f(x)\,d\nu(x).
$$
Интегрируем это равенство по $K$, в правой частми стоит константа, поэтому мы получаем
$$
\int_K\int_{K/H} f(kx)\,d\nu(x)\,d\kappa(k)=\int_{K/H}f(x)\,d\nu(x).
$$
Переставив пределы интегрирования в левой части, мы получим 
$$
 \int_{K/H}\Bigl[\int_K f(kx)\,d\kappa(k)\Bigr]\,d\nu(x).
$$
В квадратных скобках стоит $K$-инвариантная функция, т.е., константа. Мы интегрируем
постоянную функцию по вероятностной мере $\nu$ и получаем результат, независимый от
$\nu$. 
Поэтому $\nu$ единственна. 
\hfill $\square$

\sm

{\sc Задача.} Пусть $G$ -- локально компактная группа со счетной базой, $K$ 
-- компактная подгруппа. Тогда на $G/K$ существует $G$-инвариантная мера.
\hfill $\lozenge$

\sm

{\bf\punct Примеры компактных однородных пространств.%
\label{ss:examples-compact}}

\sm

1) {\it Сферы.} Сфера $S^{n-1}$ является $\SO(n)$-однородным пространством.
Инвариантная мера $\mu$ на ней нам хорошо известна. С другой стороны, фиксировав 
единичный вектор $e\in\R^n$, мы получаем отображение $g\mapsto ge$ из
$\SO(n)\to S^{n-1}$. Поэтому образ меры Хаара на $\SO(n)$ должен совпадать
с инвариантной мерой на $S^{n-1}$ (с точностью до нормировочного множителя, к примеру, мы можем
обе меры взять вероятностными).

Иными словами, {\it для любого измеримого подмножества $\Omega$
сферы мера множества элементов $g\in \SO(n)$, таких, что
$g e\in\Omega$, равна мере множества $\Omega$.}

\sm

{\sc Задача.} Рассмотрим стереографическую проекцию $S^{n-1}$ на
$\R^{n-1}$. Найдите образ инвариантной меры при этой проекции.
\hfill $\lozenge$

\sm

{\sc Задача.} Рассмотрим проекцию сферы на первую координатную ось. Найдите образ 
инвариантной меры.
Что будет при $n=2$? \hfill $\lozenge$

\sm

2) {\it Грассманианы.} Множество $\Gr_{p,q}$ всех $p$-мерных подпространств в
$\R^{p+q}$ является однородным пространством $\O(p+q)/\O(p)\times \O(q)$
(убедитесь в этом). Поэтому на  $\Gr_{p,q}$ существует $\O(p+q)$-инвариантная мера.
Явное выражение для меры см. в следующем параграфе. Множество $p$-мерных подпространств в 
$\C^{p+q}$  является однородным пространством $\U(p+q)/\U(p)\times \U(q)$.

\sm

3) {\it Пространство флагов.} Рассмотрим пространство $\Fl_n(\C)$, точками которого являются
возрастающие последовательности подпространств  (флаги) в $\C^n$.
$$
0=L_0\subset L_1\subset L_2\subset\dots\subset L_n=\C^n,\qquad \dim L_k=k.
$$
  Это однородное пространство $\U(n)/\U(1)^n$, где подгруппа 
  $\U(1)^n$ состоит из диагональных матриц. Аналогично пространство
  $\Fl_n(\R)$ флагов 
  в $\R^n$ является однородным пространством $\O(n)/\O(1)^n$. 
  Группа $\O(1)$ состоит из двух элементов $\pm1$.
  
\sm  
  
 4) Пусть $K$ -- компактная группа. Тогда
 $K\times K$ действует на $K$ левыми и правыми сдвигами,
 $$
 (h_1, h_2):\, g\mapsto h_1g h_2^{-1}.
 $$ 
 Стабилизатор точки 1 -- это диагональная подгруппа $diag(K)$,
 состоящая из пар $(h,h)\in K\times K$. Инвариантная мера, разумеется, является мерой Хаара.
 
\sm 
 
 5) {\it Многообразия Штифеля.} Это многообразие $\mathrm{St}(k,n)$,
  точками которого являются
 наборы из $k$ попарно ортогональных единичных векторов в $\R^n$.  Оно
 является однородным пространством $\O(n)/\O(k)$. Равносильно, многообразие
 Штифеля --  это, 
  множество изометричных вложений $\R^k\to\R^n$. На нем также действует группа
 $\O(k)\times\O(n)$. Кстати, какой в этом случае стабилизатор точки?
 
\sm 
 
 6) {\it Лагранжев грассманиан.} Рассмотрим пространство $\C^n$ со стандартным скалярным произведением.
 Рассмотрим множество вещественных $n$-мерных подпространств в $\C^n$,
 на которых мнимая часть скалярного произведения равна 0. Это пространство является
 однородным пространством $\U(n)/\O(n)$. См. \cite[\S3.3]{Ner2}.
 
 \sm

 {\bf\punct Инвариантные меры на однородных пространствах.}
Рассмотрим более общую ситуацию.

\begin{theorem}
Пусть $G$ -- локально компактная унимодулярная группа, $H$ -- ее замкнутая унимодулярная
подгруппа. Тогда на $G/H$ существует единственная с точностью до пропорциональности
инвариантная мера.
\end{theorem}

Стоит сразу отметить, что в случае групп Ли это утверждение тривиально,
см. следующий пункт.

\sm

{\sc Доказательство.} Пусть $\mu$, $\nu$ - меры Хаара на $G$ и $H$.
 Заметим, что на каждом классе смежности
$\xi=xH$ канонически определена мера $\nu_\xi$ как образ меры $\nu$ при отображении
$h\mapsto xh$.

Пусть $\cC(G)$, $\cC(G/H)$ -- пространства непрерывных функций с компактным
носителем на $G$ и $G/H$ соответственно.
 Рассмотрим оператор
$\Pi:\cC(G)\to \cC(G/H)$, заданный как усреднение по классам смежности
$\xi=xH$,
$$
\Pi f(\xi)=\int_\xi f(y)\,d\nu_\xi(y)=\int_H f(x q)\,d\nu(q).
$$

\begin{lemma}
{\rm a)}  Оператор $\Pi$ сюръективен. 

\sm

{\rm b)} Если $\Pi f=0$, то $\int_G f\, d\mu=0$.
\end{lemma}
   
   {\sc Доказательство леммы.} a)   Действительно, пусть $\phi$ --  непрерывная функция на
   $G/H$ с носителем $K$. Рассмотрим   компактное  множество $\wt K\subset G$,
   такое, что проекция $\wt K$ на $G/H$ содержит $K$. Рассмотрим функцию
   $f\in \cC(G)$, положительную на $\wt K$. Положим
   $$
   f^\circ(g)=\frac{f(g) \phi(gH)}{\Pi f}
   ,$$
 где $\phi$ рассматривается как функция на $G$. Тогда $\Pi f^\circ =\phi$.
 
\sm 
 
b) Пусть $\int f(x)\,d\nu_\xi=0$ для всех $\xi$. Пусть $F(x)$ -- непрерывная функция
на $G$ с компактным носителем. Тогда
\begin{multline*}
0=\int_G F(x) \int_H f(xq)\,d\nu(q)\,d\mu(x)
=\int_H \int_G F(x) f(xq)\,d\mu(x)\,d\nu(q)=\\=
\int_H \int_G F(xq^{-1}) f(x)\,d\mu(x)\,d\nu(q)
= \int_G\Bigl( \int_H  F(xq^{-1})\,d\nu(q)\Bigr)\, f(x)\,d\mu(x).
\end{multline*}
В силу сюръективности оператора $\Pi$, мы можем сделать
выражение в больших скобках равным 1 на носителе функции $f$.
\hfill $\square$

\sm

Вернемся к доказательству теоремы. Определим инвариантный
линейный функционал $I$ на $\cC(G/H)$ из условия
$$\int_G f\, d\mu= I (\Pi f).
$$
В силу утверждения b) леммы, он корректно определен,
в силу утверждения a) он определен на всем $\cC(G/H)$.
По теореме Рисса--Маркова, он задается положительной мерой.

С другой стороны, пусть $I$ -- инвариантный функционал на
$\cC(G/H)$. Тогда $I(\Pi f)$ -- левоинвариантный функционал на
$\cC(G)$, а такой функционал единственен с точностью до пропорциональности.
Поэтому единственен и $I$.
\hfill $\square$

\sm

{\sc Задача.} Где была использована унимодулярность групп?
\hfill $\lozenge$

\sm

{\sc Замечание.} Эта теорема верна в несколько большей общности:
инвариантная мера на $G/H$ существует, если для любого $h\in H$
выполнено $\Delta_G(h)=\Delta_H(h)$. Доказательство остается таким же.
\hfill $\lozenge$

\sm 
 
См. \cite[\S7.2]{Bou2}, \cite[\S9]{Wei} 
 
\sm 
 
 {\bf \punct Инвариантные меры на однородных многообразиях.}
 В случае групп Ли предыдущая теорема доказывается значительно проще.
 Пусть группа Ли $G$ транзитивно действует на гладком многообразии
 $M$. Пусть $H$ -- стабилизатор точки, так что $M\simeq G/H$.
Пусть $G$ -- группа Ли, $H$ -- ее замкнутая подгруппа. Тогда $G/H$
имеет структуру многообразия.

Пусть
 $\frg$ -- касательное пространство к $G$ в единице.
Отображение $x\mapsto g^{-1} x g$, где $x$, $g\in G$
индуцирует представление группы $G$ в $\frg$, обозначим его через
$\Ad_\frg(g)$. Аналогично определяется представление $\Ad_\frg(h)$
группы $H$. Заметим, что $h\in H$ действует и на $\frg$, и на $\frh$,
а потому группа  $H$ действует и на фактор-пространстве $\frg/\frh$. Обозначим
эти операторы через $\Ad_{\frg/\frh}(h)$. Очевидно,
$$
\det \Ad_{\frg/\frh}(h)=\frac{\Ad_\frg(h)}{\Ad_\frh(h)}.
$$
 
\begin{theorem}
Инвариантная форма объема $\omega$ на $G/H$ существует тогда и только тогда, когда
$\det \Ad_{\frg/\frh}(h)=1$ на $H$ тождественно.
\end{theorem}

{\sc Доказательство.} Допустим, инвариантная форма $\omega$
существует. Пусть $x_0$ - точка пространства $G/H$, которая
стабилизируется подгруппой $H$. Касательное пространство $T_{x_0}$ в этой точке
-- $\frg/\frh$. Значение формы  $\omega(x_0)$ на $T_{x_0}$
под действием $h\in H$ умножается на $\det \Ad_{\frg/\frh}(h)$, а потому этот
определитель равен 1.

Обратно, пусть $\det \Ad_{\frg/\frh}(h)=1$. Мы  задаем
форму в $x_0$ (такая форма одна с точностью до умножения на константу).
Для произвольного $y\in G/H$ выбираем элемент $g\in G$, такой, что
$gx_0=y$ и перенесем форму из $x_0$ в $y$. Результат этой операции однозначно
определен, так как  $gx_0=y$, $g'x_0=y$ влечет, что $g'=gh$ c $h\in H$.
\hfill $\square$

\sm

{\sc Замечание.} Чтобы существовала инвариантная мера на $G/H$, необходимо
и достаточно, чтобы $\det \Ad_{\frg/\frh}(h)=\pm 1$ для всех $h$.
Случаи с плюс-минусом в самом деле бывают, например на проективной плоскости
нет инвариантной формы объема (и вообще форм объема нет, по причине неориентируемости),
а  $\O(3)$-инвариантная мера есть.
\hfill $\lozenge$

\sm

{\sc  Замечание.} Если у групп $G$, $H$ нет нетривиальных гомомрфизмов в $\R_+^\times$,
то условие теоремы выполнено автоматически. 
\hfill $\lozenge$

\sm

{\sc  Замечание.}  Если группа $G$ унимодулярна, а подгруппа $H$ дискретна, 
то на $G/H$ есть инвариантная мера.
Пример: $\GL(2,\R)/\GL(2,\Z)$ -- пространство всех решеток в $\R^2$ (решетка в $\R^2$ -- 
замкнутая подгруппа, изоморфная $\Z^2$).
\hfill $\lozenge$

\sm

{\sc Задача.} a) На пространстве $\Fl_n(\C)$ флагов в $\C^n$, см. \ref{ss:examples-compact}, действует группа
$\GL(n,\C)\supset \U(n)$. Стабилизатор точки - это подгруппа $B_n$ верхнетреугольных матриц.
Покажите, что на пространстве флагов нет  $\GL(n,\C)$-ин\-ва\-ри\-ант\-ной меры.

\sm

b) Аналогично, на вещественном грассманиане  $\Gr_{p,q}(\R)$,  см. \ref{ss:examples-compact}, 
действует группа  $\GL(p+q,\R)$.  Покажите, что $\O(p+q)$-ин\-ва\-ри\-ант\-ная мера не
является $\GL(p+q,\R)$-инвариантной.
\hfill $\lozenge$

\sm
 
{\bf \punct Квазиинвариантные меры на однородных пространствах.} 

\begin{proposition}
Пусть $G$ --
 локально компактная группа со счетной базой. Пусть $H$ -- замкнутая подгруппа.
 Тогда на $G/H$ существует  мера, квазиинвариантная относительно $G$.
\end{proposition}

{\sc Доказательство.} Пусть $\nu$ -- левоинвариантая мера Хаара на $G$.
 Рассмотрим строго положительную интегрируемую функцию $h(g)$ на $G$.
 Рассмотрим меру $h\cdot \nu$ и ее образ $\zeta$ при отображении $\pi:G\to G/H$.
 Пусть $S\subset G/H$ имеет меру ноль. Тогда $\pi^{-1} S\subset G/H$
 имеет меру ноль. Поэтому $g \pi^{-1} S\subset G/H$ имеет меру ноль,
 поэтому $\pi g \pi^{-1} S=gS$ имеет меру ноль. 
 \hfill $\square$
 
 \sm 
 
 Отметим, что в конкретных случаях наличие квазиинвариантной
 меры обычно очевидно. Отметим без доказательства следующую теорему.
 
 \begin{theorem}
 	Квазиинвариантная мера на однородном пространстве локально компактной
 	группы единственна с точностью до умножения на положительную п.в. функцию.
 	\end{theorem} 
 
 См. \cite[VII.2.5]{Bou2}.

 \section{Инвариантные меры на грассманиане и ортогональной группе}
 
\COUNTERS 
 
 {\bf\punct Координаты на грассманиане.} Рассмотрим пространство
 $\R^{p+q}=\R^p\oplus \R^q$. Элементы группы
 $\GL(p+q,\R)$ мы будем записывать как
 блочные матрицы $\begin{pmatrix}\alpha&\beta\\\gamma&\delta\end{pmatrix}$
размера $p+q$. Мы будем считать, что эти матрицы действуют на 
$\R^{p+q}$ умножением справа на матрицы-строчки.

Обозначим через  
   $\Gr_{p,q}$ грассманиан всех $p$-мерных подпространств в
  $\R^{p+q}$. Через $\Mat_{p,q}$ мы обозначим пространство матриц
  размера $p\times q$.
  Любому элементу $T\in \Mat_{p,q}$ мы поставим в соответствие
  его график, т.е., множество векторов вида $(x,xT)\in \R^p\oplus \R^q$.
  Таким образом, мы получаем отображение $\Mat_{p,q}\to \Gr_{p,q}$.
  Его образ -- открытое плотное множество в $\Gr_{p,q}$. Дополнение
  до образа состоит из подпространств, имеющих ненулевое пересечение с $0\oplus \R^q$.
  Очевидно, что это множество меры нуль.

 Таким образом пространство $\Mat_{p,q}$
   может рассматриваться как всюду плотная карта на многообразии
$\Gr_{p,q}$.

\begin{proposition}
 Преобразованию грассманниана,
заданному матрицей 
$g=\begin{pmatrix}\alpha&\beta\\\gamma&\delta\end{pmatrix}$,
соответствует дробно-линейное преобразование пространства
$\Mat_{p,q}$, заданное формулой
\begin{equation}
z\mapsto T^{[g]}:=(\alpha+T\gamma)^{-1}(\beta+T\delta)
.
\label{eq:lin-frac}
\end{equation}
\end{proposition}

{\sc Доказательство.} Применим нашу матрицу к вектору $(x,xT)$,
$$
\begin{pmatrix}
x&xT
\end{pmatrix}
\begin{pmatrix}
\alpha&\beta\\\gamma&\delta
\end{pmatrix}=
\begin{pmatrix}
x(\alpha+T\gamma)&x(\beta+T\delta)
\end{pmatrix}
.$$
Полагая $y= x(\alpha+ T\gamma )$, мы получаем, что новое подпространство состоит
состоит из векторов $(y, y(a+Tc)^{-1}(\beta+T\delta))$, что и требовалось доказать.
\hfill $\square$

\sm

См. \cite[2.3.2, 2.3.3]{Ner2}

\sm

{\bf\punct Инвариантная мера.}

\begin{proposition}
Мера на $\Gr_{p,q}$, инвариантная относительно группы
$\O(p+q)$, в координатах $\Mat_{p,q}$ задается формулой
\begin{equation}
\det (1+TT^t)^{-(p+q)/2}\,
dT
.
\label{eq:grass-density}
\end{equation}
\end{proposition}

Утверждение вытекает из двух лемм, которые позволяют преобразовать выражение
(\ref{eq:grass-density}) под действием ортогональной группы.

\begin{lemma}
\label{l:jacobian}
Якобиан дробно-линейного преобразования
$$
T\mapsto T^{[g]}=(\alpha+T\gamma)^{-1}(\beta+T\delta)
$$
равен
$$
\det(\alpha+T\gamma)^{-p-q} \cdot
 \det \begin{pmatrix}
\alpha&\beta\\\gamma&\delta
\end{pmatrix}^p.
$$
\end{lemma}

\begin{lemma}
\label{l:1+TT}
Если $g\in \O(p+q)$,
то 
$$
\det \Bigl(1+T^{[g]} (T^{[g]})^t\Bigr)=\det(1+TT^t) \cdot|\det(\alpha+T\gamma)|^{-2}.
$$
\end{lemma}

Приступим к их доказательству.

\sm

См. \cite[\S2.11]{Ner2}.

\sm

{\bf\punct Вычисление якобиана.} 

\begin{lemma}
Дифференциал дробно-линейного преобразования {\rm(\ref{eq:lin-frac})}
равен
$$
 (\alpha+T\gamma)^{-1} \cdot dT (-\gamma T^{[g]}+\delta).
$$
\end{lemma}

{\sc Доказательство.} Нам достаточно разложить с точностью до $\epsilon$
выражение
$$
\bigl(\alpha+ (T+\epsilon\Delta)\gamma\bigr)^{-1}
\bigl(\beta+(T+\epsilon\Delta)\delta\bigr).
$$
Преобразуем сначала первый множитель 
\begin{multline*}
\bigl(\alpha+ T\gamma+\epsilon\Delta\gamma\bigr)^{-1}
=(\alpha+ T\gamma)^{-1}\bigl(1+\epsilon \Delta\gamma (\alpha+ T\gamma)^{-1}\bigr)^{-1}
=\\=
(\alpha+ T\gamma)^{-1}- \epsilon \cdot 
(\alpha+ T\gamma)^{-1} \Delta\gamma (\alpha+ T\gamma)^{-1}+o(\epsilon).
\end{multline*}
Остается перемножить два множителя
\begin{multline*}
\Bigl[(\alpha+ T\gamma)^{-1}- \epsilon \cdot 
(\alpha+ T\gamma)^{-1} \Delta\gamma (\alpha+ T\gamma)^{-1}+o(\epsilon)\Bigr]
\cdot
\Bigl[ (\beta+ T\delta)+\epsilon\Delta \delta \Bigr]
=\\=
(\alpha+ T\gamma)^{-1} (\beta+ T\delta)
+\epsilon (\alpha+ T\gamma)^{-1}\Delta 
\Bigl(-\gamma (\alpha+ T\gamma)^{-1}(\beta+ T\delta)+\delta  \Bigr)+
o(\epsilon),
\end{multline*}
что и требовалось доказать.
\hfill $\square$

\sm

Итак, дифференциал имеет вид
$$ A\ dz\, B$$
где $A$, $B$ -- явно выписанные матрицы. Очевидно (см. п.1.4), что якобиан
равен 
$$
(\det A)^q \det(B)^p= (\alpha+ T\gamma)^{-p}
\det (-\gamma T^{[g]}+d)^q.
$$
 Сейчас мы покажем, что
\begin{equation}
\det (-\gamma T^{[g]}+d)=\det (\alpha+ \gamma T)^{-1} \det \begin{pmatrix}
\alpha&\beta\\\gamma&\delta
\end{pmatrix}
.
\label{eq:gtd}
\end{equation}
Для этого мы воспользуемся формулой для определителя блочной матрицы,
которая важна и сама по себе.

\begin{theorem}
Пусть $\begin{pmatrix} A&B\\ C&D \end{pmatrix}$
-- блочная матрица размера $(p+q)\times (p+q)$.
Тогда
\begin{align}
\det
\begin{pmatrix} A&B\\ C&D \end{pmatrix}
&=\det A\cdot \det (D- CA^{-1} B)
\label{eq:block-det-1}
;\\
&=\det D\cdot \det (A-BD^{-1} C)
\label{eq:block-det-2}
.
\end{align}
\end{theorem}

Формально в равенстве (\ref{eq:block-det-1}) нужна невырожденность блока $A$.
Фактически в правой части стоит рациональное выражение с устранимой особенностью.

\sm

{\sc Доказательство.} Утверждение (\ref{eq:block-det-1}) вытекает из равенства
$$\qquad\qquad
\begin{pmatrix} A&B\\ C&D \end{pmatrix}
\begin{pmatrix} 1&-A^{-1}B\\ 0&1 \end{pmatrix}
=\begin{pmatrix} A&0\\ C&D- CA^{-1} B \end{pmatrix}
\qquad\qquad\square
$$

Теперь мы замечаем, что выражение $-\gamma (\alpha+ T\gamma)^{-1}(\beta+ T\delta)+\delta $
имеет форму $D- CA^{-1} B$. Поэтому
$$
\det\Bigl(-\gamma (\alpha+ T\gamma)^{-1}(\beta+ T\delta)+\delta\Bigr)=
\det (\alpha+ T\gamma)^{-1}\det
\begin{pmatrix}
\alpha+T\gamma& \beta+ T\delta\\
\gamma&\delta
\end{pmatrix}
.
$$
Что касается последней матрицы, то 
$$
\det\begin{pmatrix}
\alpha+T\gamma& \beta+ T\delta\\
\gamma&\delta
\end{pmatrix}=
\det\left[
\begin{pmatrix}
1&T\\0&1
\end{pmatrix}
\begin{pmatrix}
\alpha& \beta\\
\gamma&\delta
\end{pmatrix}\right]=\det \begin{pmatrix}
\alpha& \beta\\
\gamma&\delta
\end{pmatrix},
$$
что доказывает 
(\ref{eq:gtd}) и завершает вычисление якобиана.
\hfill $\square$

\sm

{\sc Задача.} Если $p=q$, то есть такой способ вычисления якобиана. Любое
дробно-линейное преобразование разлагается в произведение преобразований
вида 
$$
T\mapsto T+A,\qquad T\mapsto BTC,\qquad T\mapsto T^{-1}.
$$
Мы вычисляем якобианы этих преобразований, и проверяем, что
$J(g,T)=\det(\alpha+T\gamma)^{-1}$ удовлетворяет цепному правилу.
$$
\qquad\qquad\qquad\qquad\qquad
J(gh,T)=J(g,T) J(h,T^{[g]}).
\qquad\qquad\qquad\qquad\qquad\lozenge
$$

\sm

{\bf\punct Доказательство леммы \ref{l:1+TT}.}
\begin{multline*}
1+T^{[g]} (T^{[g]})^t=
1+(\alpha+T\gamma)^{-1}(\beta+T\delta) (\beta^t+\delta^tT^t)(\alpha^t+\gamma^tT^t)^{-1}
=\\=
(\alpha+T\gamma)^{-1}\Bigl[(\alpha+T\gamma)(\alpha^t+\gamma^tT^t)
+ (\beta+T\delta) (\beta^t+\delta^tT^t)\Bigr]
(\alpha^t+\gamma^tT^t)^{-1}.
\end{multline*}
В квадратных скобках мы имеем
$$
\Bigl[\dots\Bigr]=(\alpha \alpha^t+\beta\beta^t) + 
T(\gamma\alpha^t+\delta\beta^t)+(\alpha\gamma^t+\beta\delta^t)T^t+
T(\gamma\gamma^t+\delta\delta^t)T^t
.
$$ 
Теперь мы вспоминаем, что матрица 
$\begin{pmatrix}\alpha&\beta\\\gamma&\delta\end{pmatrix}$
ортогональна, т.е.
$$
\begin{pmatrix}\alpha&\beta\\\gamma&\delta\end{pmatrix}
\begin{pmatrix}\alpha&\beta\\\gamma&\delta\end{pmatrix}^t=
\begin{pmatrix}1&0\\0&1\end{pmatrix}=
\begin{pmatrix}
\alpha \alpha^t+\beta\beta^t&\alpha\gamma^t+\beta\delta^t\\
\gamma\alpha^t+\delta\beta^t&\gamma\gamma^t+\delta\delta^t
\end{pmatrix}.
$$
Поэтому мы получаем
$$\Bigl[\dots\Bigr]=1+TT^t,$$
а 
$$
1+T^{[g]}( T^{[g]})^t=
(\alpha+T\gamma)^{-1} \,(1+TT^t)\, (\alpha^t+\gamma^t T^t)^{-1}
,$$
что и влечет утверждение леммы.
\hfill $\square$

\sm

{\bf \punct Преобразование Кэли и мера Хаара на ортогональной группе.}
Преобразование Кэли матрицы $g$ задается формулой
$$
T:=(1+g)^{-1} (1-g) 
.
$$
Легко видеть, что это преобразование обратно самому себе, $g$
выражается через $T$ по формуле
$$
g:=(1+T)^{-1} (1-T) 
.
$$
Обозначим через $\ASymm_n$ пространство кососимметричных матриц размера $n$.

\begin{theorem}
{\rm a)} Преобразование Кэли переводит $\ASymm_n$ в $\SO(n)$.

\sm

{\rm b)} Пусть $g\in \SO(n)$, причем его собственные числа отличны от $-1$.
Тогда его преобразование Кэли содержится в $\Asymm_n$.
\end{theorem}

{\sc Доказательство.} a) Пусть $T=-T^t$. Тогда
$$
\left(\frac {1-T}{1+T}\right) \left(\frac {1-T}{1+T}\right)^t=
\left(\frac {1-T}{1+T}\right) \left(\frac {1+T}{1-T}\right)=1.
$$

b) Пусть $g\in \SO(n)$.  Тогда
$$
\left(\frac {1-g}{1+g}\right)+ \left(\frac {1-g}{1+g}\right)^t=
\left(\frac {1-g}{1+g}\right)+ \left(\frac {1-g^{-1}}{1+g^{-1}}\right)=
\left(\frac {1-g}{1+g}\right)+\left(\frac {g-1}{1+g}\right)=0.
$$

{\sc Задача.} Почему мы рассматриваем $g\in\SO(n)$, а не $g\in\O(n)$.
Что будет, если применить преобразование Кэли к ортогональной матрице
 с определителем, равным $-1$?
 \hfill $\lozenge$
 
 \sm
 
 Таким образом, мы можем рассматривать $\ASymm_n$ как карту на многообразии
 $\SO(n)$. Дополнение до этой карты -- гиперповерхность $\det(1+g)=0$.
 
 \sm

 \begin{proposition}
 \label{pr:8.4}
 Фиксируем $u$, $v\in \SO(n)$. Преобразование
 $$
 g\mapsto u^{-1} g v
 $$
в координатах $\ASymm_n$ задается дробно-линейным преобразованием
$z\mapsto (\alpha+T\gamma)^{-1} (\beta+T\delta)$, 
где
\begin{equation}
\begin{pmatrix}
\alpha&\beta\\ \gamma&\delta
\end{pmatrix}=
\frac12 
\begin{pmatrix}
u+v&u-v\\u-v&u+v
\end{pmatrix}
\label{eq:uv}
.\end{equation}
\end{proposition}

{\sc Доказательство.} Нам надо вычислить композицию трех отображений
$$
\ASymm_n\to\SO(n)\to\SO(n)\to \ASymm_n.
$$ 
Пусть $T\in \ASymm_n$. Он переходит
\begin{multline*}
T\,\mapsto\, (1+T)^{-1}(1-T)\,\mapsto\, u^{-1}(1+T)^{-1}(1-T)v\,\mapsto
\\
\mapsto\,
\Bigl(1+u^{-1}(1+T)^{-1}(1-T)v\Bigr)^{-1}\Bigl(1- u^{-1}(1+T)^{-1}(1-T)v\Bigr)
.
\end{multline*}
Преобразуем последнее выражение. Для этого умножим обе скобки {\it слева} на
$(1+T)u$. Получим
\begin{multline*}
\Bigl((1+T)u+(1-T)v\Bigr)^{-1} \Bigl((1+T)u-(1-T)v\Bigr)
=\\=
\Bigl((u+v)+T(u-v)\Bigr)^{-1} \Bigl((u-v)+T(u+v)\Bigr),
\end{multline*}
что и требовалось доказать. \hfill $\square$ 

\sm

Стоит заметить, что матрица (\ref{eq:uv}) является ортогональной. Это легко
проверить непосредственно умножением, но мы сошлемся на равенство
$$
\frac12 
\begin{pmatrix}
u+v&u-v\\u-v&u+v
\end{pmatrix}=J
\begin{pmatrix} 
u&0\\0&v
\end{pmatrix} J^{-1},\qquad\text{где}\,\,
J=
\frac 1{\sqrt 2}
\begin{pmatrix}
1&-1\\1&1
\end{pmatrix}
.
$$
Матрица $J$, как легко видеть, ортогональна.

\begin{theorem}
Мера Хаара в координатах $\ASymm_n$ задается формулой
$$
\det(1+TT^t)^{-(n-1)/2} dT.
$$
\end{theorem}

{\sc Доказательство.}
Мы докажем, что это выражение $\SO(n)$-инвариантно (без ссылки на теорему существования
меры Хаара). Как преобразуется
выражение $\det(1+TT^t)$ под действием ортогональной группы мы уже знаем,
см. лемму \ref{l:1+TT}. Нам надо вычислить как преобразуется
$dT$, т.е., вычислить якобиан нашего дробно-линейного преобразования на кососимметрических
матрицах. Дифференциал нам известен,
\begin{equation}
 (\alpha+T\gamma)^{-1} \cdot dT \cdot (-\gamma T^{[g]}+\delta)
 .
 \label{eq:differential}
\end{equation}

\begin{lemma}
$$(-\gamma T^{[g]}+\delta)^t= (\alpha+T\gamma)^{-1}.$$
\end{lemma}

{\sc Доказательство.} В принципе это можно доказать <<в лоб>>, используя ортогональность
 нашей матрицы. Мы пойдем другим путем. Заметим, что дифференциал
(\ref{eq:differential}) имеет вид $\Delta \mapsto A\Delta B$, где $\Delta$ --
кососимметрическая матрица, а $A$, $B$ -- фиксированные матрицы.
Это отображение должно переводить кососимметрические матрицы в кососимметрические.
Заметим, что матрица
$$A^{-1} A\Delta B A^{t-1}= \Delta B A^{t-1}$$
должна быть кососимметричной при любой кососимметричной матрице $\Delta$.
Но тогда $B A^{t-1}$ должна быть скалярной матрицей, $B=s\cdot A^{t}$. Для того,
чтобы вычислить этот скаляр, достаточно вычислить определитель матрицы
$B$. Но мы это уже делали, см. (\ref{eq:gtd}). Выходит, что $\det B=\det A$
и $s=\pm 1$.  Из соображений связности группы $\SO(n)$ мы имеем
$s=+1$. Впрочем, конец доказательства теоремы от этого не зависит.
\hfill $\square$ 

\sm

В силу леммы \ref{l:AXA}, якобиан нашего отображения равен 
$$\pm\det(\alpha+T\gamma)^{-(n-1)/2},$$
что нам и требовалось.
\hfill $\square$ 

\sm

См. \cite[\S3.1]{Hua}

\sm

{\sc  Задача.} Покажите, что преобразование Кэли переводит
унитарную группу $\U(n)$ в пространство антиэрмитовых матриц 
($T=-T^*$). Покажите, что мера Хаара задается формулой
$\det(1+TT^*)^{-n} dT$.
\hfill $\lozenge$ 

\sm

{\bf\punct Комментарии к вычислению.} Вычисление меры Хаара с помощью
формальных манипуляций
с дробно-линейными отображениями может показаться странным. Объясним, что произошло.

Рассмотрим прямую сумму $V\oplus W=\R^n\oplus \R^n$ двух $n$-мерных пространств.
Каждое из них снабдим  стандартным скалярным произведением.
Введем в
$V\oplus W$ симметричную билинейную форму $B$, задаваемую матрицей 
$\begin{pmatrix}1&0\\0&-1\end{pmatrix}$.
Обозначим через $\cL$ множество всех максимальных изотропных 
подпространств  в $V\oplus W$.
 Напомним, что подпространство
$L$ называется  {\it изотропным}, если форма $B$ равна 0 на $L$ (подробнее см. \cite[\S\S2.1-2.3]{Ner2}). Размерность
$L$ не может превосходить $n$, иначе $L\cap W\ne 0$, а форма
$B$ на нем одновременно равна $0$ и отрицательна. В силу того же довода
в случае $n$-мерного пространства $L\cap W= 0$. Поэтому $n$-мерное
изотропное подпространство
является графиком оператора $V\to W$.

\sm

{\sc Задача.} Покажите, что $n$-мерное подпространство в $V\oplus W$
является изотропным тогда и только тогда, когда оно является графиком ортогонального
оператора $V\to W$. \hfill $\lozenge$

\sm

Таким образом, пространство $\cL$ отождествляется с ортогональной группой $\O(n)$.
Кстати, отсюда следует, что $\cL$ состоит из двух компонент связности.

Рассмотрим в $V\oplus W$ графики единичного оператора и минус-единичного оператора 
$V\to W$. Обозначим их через $P$, $Q$. Эти подпространства изотропны, а форма
$B$ относительно разложения $P\oplus Q$ записывается матрицей 
$\begin{pmatrix}0&1\\1&0 \end{pmatrix}$.

\sm

{\sc Задача.}  График оператора $T:P\to Q$ является изотропным подпространством
тогда и только тогда, когда $T$ кососимметричен.\hfill $\lozenge$

\sm

Таким образом, одной и той же точке грассманиана $\cL$ соответствует ортогональная матрица
и кососимметрическая матрица. Это и есть преобразование Кэли. Становится
понятным и то, что техника для работы с грассманианом подходит и для $\SO(n)$.

\sm 

См. \cite[\S2.4]{Ner2}.


\newpage

\begin{center}
	\huge
	Глава B. Разное
\end{center}

\section{Зоопарк топологических групп}

\COUNTERS

Здесь мы продолжаем список примеров, начатый в п.\ref{ss:zoo-nach}. 
Параграф устроен как зоопарк, т.е. как набор изолированных клеток с примученными
зверьми. В последующих параграфах картинки этих зверей  используются
редко, в каждом случае дается точная ссылка на клетку.

\sm

{\bf \punct Алгебры Ли.} Алгебры Ли почти не используются в данных записках,
но добиваться полного их неупоминания не очень разумно (фактически они упоминаются в некоторых 
пунктах этого параграфа). Поэтому мы вкратце напомним
 простейшие определения и факты.

Вещественная (комплексная) {\it алгебра Ли} $\frg$ -- это {\it конечномерное} линейное пространство над $\R$
(соответственно, $\C$), снабженное билинейной операцией $\frg\times\frg\to\frg$,
 называемой {\it коммутатором} и обозначаемой
$[x,y]$, удовлетворяющей свойствам:
\begin{align*}
[\alpha x+\beta y,z]=\alpha [x,z]+\beta[y,z],\quad \text{где $\alpha$, $\beta\in\C$ (линейность)};\\
[x,y]=-[y,x] \qquad \text{(антикоммутативность)};\\
[[x,y],z]+[[y,z],x]+[[z,x],y]=0\,\,
\text{(тождество Якоби)}.
\end{align*}

{\it Соответствие между группами Ли и алгебрами Ли}. Из любой группы Ли изготавливается алгебра Ли по следующим рецептам.

\sm

1) Пусть $G$ -- замкнутая подгруппа в $\GL(n,\R)$ (напомним, что она автоматически является
подмногообразием). Пусть $\frg$ -- касательное пространство к $G$ в единице. 
Оказывается, что для любых
$x$, $y\in \frg$ выполнено $xy-yx\in \frg$. Коммутатор определяется как $[x,y]:=xy-yx$.

\sm

2) Рассмотрим общий случай. Пусть группа $G$ задана как многообразие с операцией умножения.
Алгебра Ли $\frg$, по-прежнему, является  касательным пространством  к $G$ в единице.
  Рассмотрим всевозможные непрерывные
гомоморфизмы $\R\to G$ (<<однопараметрические подгруппы>>), их можно рассматривать как параметризованные кривые в  $G$. Оказывается,
что любая такая однопараметрическая подгруппа однозначно определяется своим вектором $X$ скорости при $t=0$, и любой
вектор $x\in\frg$ является вектором скорости некоторой однопараметрической подгруппы.
Обозначим однопараметрическую подгруппу через $\exp(tX)$, в случае матричных групп это
в самом деле обычная матричная экспонента. Пусть $X$, $Y\in \frg$. Положим
\begin{equation}
\gamma(t,s)=\exp(tX)\, \exp(sY)\,\exp(-tX)\, \exp(-sY).
\label{eq:[]}
\end{equation}
Можно показать, что тейлоровское разложение $\gamma(t,s)$ в нуле начинается с
\begin{equation}
\gamma(t,s)= 2ts \sigma + o(t^2+s^2).
\label{eq:[]=1}
\end{equation}
Величина $\sigma\in\frg$ и есть $[X,Y]$.

\sm

{\sc Задача.} В случае матричной группы убедитесь, что тейлоровское разложение
(\ref{eq:[]}) имеет вид (\ref{eq:[]}), причем $\sigma=XY-YX$.
\hfill $\lozenge$

\sm 

{\sc Задача.} Описать алгебры Ли групп $\O(n)$, $\U(n)$, $\Sp(2n,\R)$, $\SL(n,\C)$.
\hfill $\lozenge$

\sm

Отображение $G\mapsto \frg$ является функториальным в следующем смысле: 
если $\theta:G_1\to G_2$ -- гомоморфизм групп Ли, то его дифференциал $d\theta:\frg_1\to \frg_2$
в единице является гомоморфизмом алгебр Ли.

\sm

Разным группам Ли может соответствовать одна и та же алгебра Ли. Например, возьмем
группу Ли $\SL(n,\C)$. У нее есть центр $\Z_n$, состоящий из скалярных матриц
$e^{2\pi i/n} \cdot 1$. Возьмем в центре любую подгруппу $C$ и рассмотрим фактор-группу
$\GL(n,\C)/C$.  Алгебры Ли всех таких групп одинаковы. Можно также взять произведение
$\SL(n,\C)$ с любой конечной группой $\Gamma$. Понятно, что алгебра Ли от этого не изменится.  Оказывается, что в общем случае
неоднозначность имеет примерно такой же вид. 

\sm

Для любой (конечномерной) вещественной алгебры Ли $\frg$ существует единственная связная 
односвязная группа Ли $\wt G$, чья алгебра Ли совпадает с $\frg$. Любая другая связная группа
Ли с алгеброй Ли $\frg$ изоморфна фактор-группе $\wt G$ по некоторой дискретной подгруппе,
содержащейся в центре группы $\wt G$.

Обратное соответствие тоже функториально, любому гомоморфизму
(конечномерных)  алгебр Ли $\pi:\frg_1\to\frg_2$ соответствует единственный гомоморфизм
$\wt G_1\to \wt G_2$ с данным дифференциалом $\pi$ в единице.

\sm

{\it Комплексные группы Ли и комплексификации.}  Любая комплексная алгебра Ли может
рассматриваться как вещественная алгебра Ли удвоенной размерности. Соответствующая группа Ли
обладает структурой комплексно аналитического многообразия (и наоборот,
комплексно аналитической группе Ли соответствует комплексная алгебра Ли).

Далее заметим, что 
  любой вещественной алгебре Ли $\frg$ мы можем поставить в соответствие
  ее {\it комплексификацию}
   $\frg_\C=\frg\otimes_\R \C$. Иными словами, мы выбираем базис $e_k$ в $\frg$
  и берем линейные комбинации $e_k$ с комплексными коэффициентами.
  Говорят также, что алгебра $\frg$ есть {\it вещественная форма} алгебры
  $\frg_\C$.
  
  Отметим, что вещественных форм у комплексной алгебры может быть много (см. следующий пункт),
  а может и не быть совсем.
  
  Очевидно, что гомоморфизм алгебр Ли $\frg\to\frh$ продолжается до гомоморфизма их комплексификаций
  $\frg_\C\to\frh_\C$.

\sm

\cite[главы 1-2]{Zhe1}, \cite[\S\S16, 34]{Zhe2}, \cite[гл.1]{Dix}, \cite[2.1.4., Add. C]{Ner2}, \cite[\SS 2.2-2.5]{Ros}, \cite{VO}, \cite[глава 10]{Pon}, \cite{Sophus}.

\sm

{\bf \punct Классические группы и простые группы Ли.} Понятно, что групп Ли очень много, 
однако среди них выделенную роль играют 10 серий так называемых {\it классических групп}.

\sm

{\it Группы вещественных и комплексных матриц.}
Во-первых, это полные линейные группы $\GL(n,\R)$, $\GL(n,\C)$, т.е., группы всех матриц 
размера $n\times n$ над вещественными и комплексными числами соответственно.

Далее это группы матриц, сохраняющих всевозможные <<формы>>, симметрические, кососимметрические
и эрмитовы. 

Напомним, что любая невырожденная  симметрическая билинейная форма на комплексном линейном 
пространстве приводится к виду 
$$
B(x,y)=\sum x_j y_j.
$$
Группа комплексных матриц, сохраняющих эту форму называется {\it комплексной ортогональной группой}
и обозначается через $\O(n,\C)$. Над вещественным полем невырожденная симметричная билинейная форма приводится к виду
$$
B(x,y)=\sum_{j\le p} x_j y_j- \sum_{ p< j \le p+q} x_j y_j.
$$
Группа сохраняющих ее линейных преобразований называется {\it псевдоортогональной группой}
и обозначается через $\O(p,q)$.

Перейдем к невырожденным кососимметричным формам. Такие формы бывают лишь на четномерных пространствах,
и они приводятся к виду
$$
B(x,y)=\sum_{j\le n} (x_j y_{j+n}-x_{j+n} y_j)
.
$$
Это дает нам {\it симплектические группы} $\Sp(2n,\R)$ и $\Sp(2n,\C)$.

Наконец, рассмотрим невырожденные эрмитовы формы на комплексных пространствах. Они
приводятся к виду
\begin{equation}
B(x,y)=\sum_{j\le p} x_j \ov y_j- \sum_{ p< j \le p+q} x_j \ov y_j
\label{eq:hermitian-form}
\end{equation}
Группа матриц, сохраняющих такую форму, называется {\it псевдоунитарной группой}
и обозначается через $\U(p,q)$.

\sm

{\it Группы кватернионных матриц.}
Оставшиеся 3 серии классических групп
 естественнее всего задавать кватернионными матрицами. Понятно, что
обратимые кватернионные матрицы образуют группу, она обозначается через $\GL(n,\H)$ 
(тело кватернионов обозначается через $\H$).

Перейдем к обсуждению форм над $\H$.
Алгебра кватернионов не коммутативна, поэтому говорить о линейных пространствах над $\H$ неаккуратно.
Давайте, для определенности, рассматривать <<левые модули>> над $\H$. Такой
модуль имеет вид $\H^n$, где $\H^n$ - пространство кватернионных векторов-строк
$(x_1,\dots,x_n)$. Такие строчки можно умножать на кватернионы слева,
$$
\lambda (x_1,\dots,x_n)= (\lambda  x_1,\dots,\lambda  x_n)
.$$
Аналог линейных операторов -- автоморфизмы левых модулей. Легко проверить, что они имеют вид
$x \mapsto xA$, где $A$ - кватернионная матрица размера $n$.

{\it Эрмитова форма} $\la x,y\ra$ на левом кватернионном модуле $V$ -- 
это отображение $V\times V\to \H$, удовлетворяющее условиям
$$
\la x,y+z\ra=\la x,y\ra+\la x,z\ra,\qquad \la x,y\ra=\ov{\la y,x\ra},\qquad
\la \lambda x, \mu y\ra= \lambda \la x,y\ra \ov \mu
,$$
где $x$, $y$, $z\in V$, а $\lambda$, $\mu\in \H$. Такие формы 
имеют вид
\begin{equation}
\la x, y\ra=\sum_{kl} x_k a_{kl} \ov y_l, \qquad
\text{где $a_{kl}=\ov a_{lk}$}.
\label{eq:quaternion-hermitian}
\end{equation}
Невырожденные эрмитовы формы%
\footnote{Невырожденность, как обычно, означает, что нет ненулевого вектора 
$a$, такого, что $\la a,y\ra=0$ для всех $y\in V$.}
 приводятся к виду
(\ref{eq:hermitian-form}). Группа операторов, сохраняющих такую эрмитову форму,
называется {\it псевдосимплектической группой} и обозначается через $\Sp(p,q)$.

Если $q=0$, мы пишем просто $\Sp(p)$. Это кватернионно унитарная группа (она также 
называется симплектической группой, хотя это и создает определенное
терминологическое несоответствие).

Наконец, над телом кватернионов бывают {\it антиэрмитовы формы}, в пре\-ды\-ду\-щем определении
мы должны заменить условие $\la x,y\ra=\ov{\la y,x\ra}$ на $\la x,y\ra=-\ov{\la y,x\ra}$,
а в выражении (\ref{eq:quaternion-hermitian}) наложить на коэффициенты условие
$a_{kl}=- \ov a_{lk}$. Невырожденная антиэрмитова форма на кватернионном
пространстве приводятся к виду
$$
P(x,y)=\sum_{l=1}^n x_l i\ov y_l
,$$
где $i$ -- мнимая единица. Группа операторов, сохраняющих такую форму, почему-то не имеет названия
и обозначается через $\SOS(2n)$.

Последние три серии групп несложно определить, не ссылаясь непосредственно на кватернионную алгебру.
Для этого заметим, что тело кватернионов изоморфно алгебре комплексных матриц порядка
2, имеющих вид $\begin{pmatrix} a&b\\-\ov a&\ov b \end{pmatrix}$. Поэтому группу $\GL(n,\H)$
можно определить как группу невырожденных комплексных блочных матриц размера
$n+n$, имеющих структуру
\begin{equation}
\begin{pmatrix}
\Phi&\Psi\\ -\ov \Psi&\ov \Phi
\end{pmatrix}
.
\label{eq:quaternionic-matrix}
\end{equation}
Группа $\SOS(2n)$ отождествляется с группой  матриц такой блочной структуры, удовлетворяющих
дополнительному условию
$$
g\begin{pmatrix}0&1\\1&0 \end{pmatrix}g^t=\begin{pmatrix}0&1\\1&0 \end{pmatrix}.
.$$
Это условие можно заменить на равносильное условие
$$
g \begin{pmatrix} i&0\\0&-i\end{pmatrix} g^*= \begin{pmatrix} i&0\\0&-i\end{pmatrix}.
$$
Т.е. 
$$
\SOS(2n)=\GL(n,\H)\cap \SO(2n,\C)=\GL(n,\H)\cap\U(n,n)=\SO(2n,\C)\cap \U(n,n).
$$

{\it Вариации.}
Отметим, что на некоторых из этих групп определитель не равен тождественно единице,
и поэтому он задает нетривиальный гомоморфизм в некоторую абелеву группу. Ядро 
этого гомоморфизма для групп типа $\GL(\cdot)$ обозначается%
\footnote{В кватернионном случае мы берем определитель матрицы (\ref{eq:quaternionic-matrix}),
он вещественен и положителен.} через $\SL(\cdot)$. В оставшихся  случаях перед названием
группы ставится буква $S$: $\SU(p,q)$, $\SO(p,q)$, $\SO(n,\C)$.

Фактор-группа классической группы по центру всегда обозначается добавлением приставки $P$:
$\PSL(n,\R)$ и т.д.

Алгебры Ли обозначаются теми же буквами, что и группы, только маленькими  готическими.

\sm

{\it Комплексификации классических алгебр.}
Операция комплексификации действует на классические алгебры Ли как:
\begin{align*}
&\frs\frl(n,\R) \to \frs\frl(n,\C), \quad \frs\frl(n,\C) \to \frs\frl(n,\C)\oplus \frs\frl(n,\C),
\quad \frs\frl(n,\H) \to \frs\frl(2n,\C),
\\
&\frs\fru(p,q)\to \frs\frl(p+q,\C),
\\
&\frs\fro(p,q)\to \frs\fro(p+q,\C),\quad \frs\fro(n,\C)\to  \frs\fro(n,\C)\oplus  \frs\fro(n,\C),
\quad  \frs\fro^*(2n)\to  \frs\fro(2n,\C)
\\
&\frs\frp(2n,\R)\to \frs\frp(2n,\C),\quad
\frs\frp(2n,\C)\to \frs\frp(2n,\C)\oplus \frs\frp(2n,\C),
\\
&\frs\frp(p,q)\to \frs\frp(2(p+q),\C).
\end{align*} 

{\it Простые алгебры Ли.} Напомним, что группа называется {\it простой}, если в ней нет
нетривиальных нормальных подгрупп. 

Алгебра Ли $\frg$ называется {\it простой}, 
если в ней нет нетривиальных
идеалов. Кроме того, одномерная алгебра Ли не считается простой.
На уровне соответствующей односвязной группы Ли эта простота означает отсутствие
замкнутых нормальных подгрупп размерности $>0$ (т.е. терминология чуть-чуть не
 соответствует терминологии общей теории групп).
 
 Следующие классические алгебры Ли просты (кроме достаточно очевидных исключений в малой размерности):
 \begin{align*}
 \frs\frl(n,\R),\quad  \frs\frl(n,\C), \quad  \frs\frl(n,\H), \quad \frs\fru(p,q),
 \\
 \frs\fro(p,q),\quad  \frs\fro(n,\C), \quad  \frs\fro^*(2n),\\
 \frs\frp(2n,\R),\quad  \frs\frp(2n,\C),\quad  \frs\frp(p,q).
 \end{align*}
 
 Оказывается, что простые алгебры Ли допускают полную классификацию 
 (Вильгельм Киллинг и Эли Картан, конец XIX века).
 
 \sm
 
 a) Есть три бесконечных серии%
 \footnote{На самом деле  ортогональные алгебры принято расщеплять на две серии, $\frs\fro(2n,\C)$, $\frs\fro(2n+1,\C)$,
 причины для этого становятся ясными, как только начинаешь работать с ними, см., например, ниже п.\ref{ss:zoo-comp}.}
 комплексных простых алгебр  Ли,
 $\frs\frl(n,\C)$, $\frs\fro(n,\C)$, $\frs\frp(2n,\C)$, а также 5 {\it особых алгебр},
 имеющих размерности 14, 52, 78, 133, 248, см. \cite[гл. 3]{VO}. 
 В принципе, особые алгебры могут быть описаны
 в терминах матриц над алгеброй Кэли, что, скорее является философским фактом, чем
 инструментом работы с ними. О разных моделях особых алгебр см. \cite[гл. 5]{VO}.
 
 \sm
 
 b) Есть 10 бесконечных серий классических вещественных алгебр Ли (которые перечислены выше)
  и 23 особых алгебры.  Особые алгебры включают 5 комплексных особых алгебр, а также их
  вещественные формы (см. \cite[гл. 4]{VO}).
  
 \sm
 
См. \cite[\S41]{Arn-mech}, \cite[гл. 14]{Zhe1}, \cite{Zhe2}, \cite[ch. 2-3, Add. C]{Ner2}, \cite[гл. 11]{Pon}, \cite[\S3.1]{Ros}, \cite{VO},
\cite{Ada}, \cite{Sophus}.

О структурной теории классических групп над общими полями и кольцами см. \cite{Die}, \cite{VaS}.
 
 \sm
 
 {\bf\punct Связные компактные группы Ли.} Классические компактные группы -  это $\SO(n)$, $\SU(n)$,
 $\Sp(n)$. Отметим, что все эти группы определяются единообразно как группы  вещественных,
 комплексных и кватернионных унитарных
 матриц. Кроме того, каждая из 5 особых комплексных групп Ли имеет единственную компактную форму.
 Разумеется, компактен <<тор>> $\R/\Z$. Оказывается, 
 что этот набор примеров является исчерпывающим. Точное утверждение формулируется так.  
 Алгебра Ли любой компактной группы Ли
 есть прямое произведение алгебр вида $\frs\fro(\cdot)$, $\frs\fru(\cdot)$,
  $\frs\frp(\cdot)$ и алгебры $\R^n$ с нулевым коммутатором.
  
  \sm

Отметим, что любая полупростая группа Ли содержит единственную с точностью
до сопряжения максимальную компактную подгруппу, например,
$$
\SL(n,\R)\supset\SO(n),\quad\O(n,\C)\supset\O(n),\quad \SL(n,\H)\supset \Sp(n),\quad \Sp(2n,\R)\supset \U(n).
$$

\sm

См. \cite{Zhe1}, \cite[\S45]{Zhe2}, \cite[гл. 11]{Pon}, \cite{Ada}, \cite{Sophus}.

\sm

{\bf \punct $p$-адические числа.}
Теория $p$-адических чисел -- отдельный большой предмет. Мы их упоминаем постольку,
поскольку не упомянуть их при осуждении топологических групп невозможно.

Пусть $p$ -- простое число.
 Представим произвольный ненулевой элемент 
$x\in\Q$ в виде $x=\frac {a}{b}\cdot p^{-k}$, где $a$, $b$ -- целые числа, взаимно простые с
$p$. Определим его $p$-адическую норму как $|x|_p=p^k$. Кроме того, положим $|0|_p=0$.
  Легко видеть, что
$$
|xy|_p=|x|_p \,|y|_p,\quad |x+y|_p\le \max ( |x|_p, |y|_p)
.
$$
Введем на $\Q$  расстояние $d(x,y)=|x-y|_p$, очевидно, оно удовлетворяет
ультраметрическому неравенству $d(x,z)\le d(x,y)+d(y,z)$. Обозначим через $\Q_p$
пополнение
$\Q$ по этой метрике. Легко видеть, что получается вполне несвязное топологическое
пространство, сложение и умножение (а также норма) продолжаются по непрерывности на $\Q_p$.

\sm

{\sc Задача.} Покажите, что $\Q_p$ отождествляется с множеством
$p$-ичных записей, бесконечных влево и конечных вправо,
$$
\dots a_3 a_2 a_1 a_0, a_{-1}\dots a_{-k}, 
$$
где $a_j=0$, $1$, $2$,  \dots $p-1$. Такой записи ставится в соответствие сходящийся ряд
$\sum_{j=-k}^\infty a_j p^j$. Сложение и умножение таких записей проводится <<столбиком>>.
Покажите, что обратный элемент однозначно определен. Покажите, что поле
$\Q_p$ локально компактно. \hfill $\lozenge$

\sm

Полученный объект называется полем $p$-адических чисел.
Числа $z$, у которых нет знаков после запятой (или, равносильно, $|z|\le 1$),
называются целыми $p$-адическими числами (мы их будем обозначать $\bbO_p$).

\sm

{\sc Задача.} Опишите
меру Хаара на аддитивной группе $\Q_p$; на мультипликативной группе $\Q_p^\times$.
\hfill $\lozenge$

\cite[гл. 2]{Ser1}, \cite{PR}, \cite[\S10.1]{Ner2}, \cite[\S II.1]{GGP}.

\sm

{\bf \punct   Матричные группы.} Разумеется, кроме матриц над $\R$, $\C$ и $\H$, можно рассматривать
группы  матриц над другими полями и кольцами. Упомянем некоторые случаи, в которых возникают
содержательные теории.

\sm

1){\it Матрицы над $p$-адическими числами.}
 К существенно новому миру приводит рассмотрение $p$-адических матричных групп, что лежит за
пределом настоящих заметок. Для тех, кто с $p$-адическими числами встречался раньше, предложим
несколько простых задач

\sm

{\sc Задача.} Убедитесь что группа $\GL(n,\Q_p)$ локально компактна. Опишите меру Хаара на ней.
Покажите, что группа $\GL(n)$ над кольцом целых $p$-адических чисел компактна.
\hfill $\lozenge$

\sm

2) {\it Адели и адельные группы.} Элемент кольца аделей $\bbA$ - это последовательность
$$
(z_2, z_3, z_5, z_7,\dots; z_\infty),
$$
где $z_p\in\Q_p$, $z_\infty\in\R$, причем $z_p$ является целым $p$-адическим
для всех $p$, кроме конечного числа.

\sm

{\sc Задача$^*$.} a) Покажите, что пространство аделей локально компактно. Опишите меру Хаара на аддитивной группе аделей.

\sm

b) Для каждого $z\in \Q$ рассмотрим адель $(z,z,\dots;z)$. Покажите, что образ $\Q$ --
дискретная подгруппа в аддитивной  группе $\bbA$. Покажите, что группа $\bbA/\Q$
компактна.

\sm

с) Рассмотрим группу $\A_{fin}$, состоящую из последовательностей 
$$(z_2, z_3, z_5, z_7,\dots)$$
без $z_\infty$, таких, что $z_p$ - целый для почти всех
$p$. Покажите, что образ $\Q$ всюду плотен в $\A_{fin}$.

\sm

d) Опишите множество всех гомоморфизмов $\Q$ в мультипликативную группу $\C$.
\hfill $\lozenge$

\sm

Матричные группы над кольцом аделей -- важный объект теории чисел (хотя почему это интересно, сразу не ясно).

\sm

3) {\it Формальные ряды}. Рассмотрим конечное поле $\bbF$ и кольцо формальных рядов 
$\bbF[t]$ над $\bbF$,
т.е., выражений вида $r(t)=\sum_{k\ge 0} a_k t^k$, где $a_k\in \bbF$. Введем в нем топологию, положив, что последовательность $r^{(j)}(t)$ сходится к $r(t)$, если для каждого $k$ мы имеем 
покоэффициентную стабилизацию, т.е.
$a_k^{(j)}=a_k$, начиная с некоторого номера.  Иными словами, для любого $M$
выполнено:
\begin{equation}
r^{(j)}(t)-r(t)=0 \,(\mathrm{mod}\,\, t^M r[t]), \,\, \text{начиная с некоторого номера $j$.}
\end{equation}
Аналогично, мы вводим кольцо формальных лорановских рядов  $\bbF[[t]]$, 
т.е., выражений вида $r(t)=\sum_{k\ge -N} a_k t^k$ (мы разрешаем слагаемые с отрицательными
степенями, но их число может быть лишь конечным). Последовательность $r_j(t)$ сходится
к $r(t)$, если существует единое $N$, такое, что $t^N r_j(t)\in \bbF[t]$
для всех $j$ и
имеет место покоэффициентная стабилизация.

\sm

{\sc Задача.} a) Покажите, что  $\bbF[[t]]$ является полем.

\sm

b) Покажите, что поля $\bbF[[t]]$ локально компактны, а операции сложения, умножения и обращения
непрерывны. Покажите, что $\bbF[t]$ -- компактное кольцо.
\hfill $\lozenge$

\sm

Это еще один родственник $\R$ и $\Q_p$.
Согласно теореме Понтрягина--Ко\-валь\-ско\-го,
\cite[\S 27]{Pon}, любое локально-компактное поле характеристики 0 является
или $\R$, или $\C$, или конечным расширением%
\footnote{Пусть поле $L$ содержит поле $K$. Тогда говорят, что $L$ -- {\it расширение
поля} $K$. Очевидно, $L$ является линейным пространством над $K$. Расширение называется
{\it конечным}, если $L$ конечномерно над $K$.} $p$-адического поля.
Любое локально-компактное поле характеристики $p$ является конечным расширением
поля $\bbF[[t]]$.

Разумеется, группы $\GL(n)$ над локально компактными полями локально компактны.

\sm

См. \cite[гл. II-III]{GGP},  \cite{PR}, \cite[\S27]{Pon}, \cite[Chap. 10-11]{Ner2}.

\sm

{\bf\punct Бесконечномерная унитарная группа.%
\label{ss:unitary-infty}} Пусть $H$ -- сепарабельное гильбертово пространство.
На алгебре $B(H)$ ограниченных операторов пространства $H$ есть несколько стандартных
топологий: равномерная, сильная, слабая (а также реже используемые ультрасильная, ультраслабая
и $*$-слабая). 

При ограничении этих топологий на полную унитарную группу остаются лишь равномерная топология
и слабая топология (а остальные упомянутые топологии совпадают со слабой).

\sm

а) {\it Равномерная топология.} Мы говорим, что последовательность $g$ сходится к
$g$, если $\|g_j-g\|\to 0$ (упражнение: проверьте, что эта топология не сепарабельна). 
Равномерная топология,  не очень естественна в следующем
смысле: гомоморфизмов из других топологических групп в $\bfU(\infty)$ оказывается
на удивление мало. 

\sm

{\sc Задача.} Пусть аддитивная группа $\R$ действует в $L^2(\R)$ сдвигами
$$T(s) f(x)=f(x+t).
$$
Покажите, что гомоморфизм
 $s\mapsto T(s)$ всюду разрывен.
\hfill $\lozenge$

\sm

Оказывается (проверка этого выходит за пределы данных заметок%
\footnote{У обеих групп классификация унитарных представлений хорошо известна,
их разрывность относительно равномерной топологии легко проверяется}.),
что нетривиальных гомоморфизмов $\SL(2,\R)$ в $\bfU(\infty)$ нет вообще,
а гомоморфизм из группы матриц  $\begin{pmatrix} a&b\\0&1 \end{pmatrix}$
является тривиальным на подгруппе матриц $\begin{pmatrix} 1&b\\0&1 \end{pmatrix}$.

\sm

{\sc Задача.}
Выведите отсюда аналогичные высказывания для $\SL(n,\R)$ и для группы
верхнетреугольных матриц.

\sm

b) Хорошей топологией на группе $\bfU(\infty)$ является 
{\it слабая операторная  топология}. Оказывается (проверьте эти утверждения),
что 

\sm

--- на $\bfU(\infty)$ слабая и сильная операторные топологии совпадают;

\sm

---  умножение на $\bfU(\infty)$ непрерывно по совокупности переменных;

\sm

---  слабая топология на $\bfU(\infty)$  метризуема и сепарабельна%
\footnote{Напомним, что на алгебре всех операторов слабая топология не является 
метризуемой и раздельно непрерывной.}.

\sm

Близкими родственниками бесконечномерной унитарной группы являются 
бесконечномерная ортогональная группа $\bfO(\infty)$ и бесконечномерная симплектическая (кватернионно-унитарная) группа $\mathbf{Sp}(\infty)$.

\sm

См.  \cite{Olsh-uni}, \cite{Tsa}, \cite[\S8.2]{Ner1}, \cite{Kui}, \cite[\S 4.5]{Pes}.

\sm

{\bf\punct Группы линейных операторов. Шейловская топология.%
\label{ss:shale}}
Напомним определение {\it шаттеновских классов}. Пусть $1\le p<\infty$. Пусть $A$
-- компактный оператор. Пусть $\lambda_1$, $\lambda_2$, \dots -- сингулярные
 числа оператора  $A$, т.е. ненулевые собственные числа оператора $\sqrt {A^*A}$.
Оператор $A$ содержится в {\it шаттеновском классе} $\frL_p(H)$
если $\sum_j \lambda_j^p<\infty$ (см. \cite[\S XI.9]{DS1}). Среди этих классов есть два выделенных и часто встречающихся:
операторы Гильберта--Шмидта ($p=2$) и ядерные операторы ($p=1$).
Класс $\frL_p(H)$ является идеалом в $B(H)$, он также является банаховым пространством
относительно нормы
$$
\|A\|_p=\Bigl(\sum_j \lambda_j^p \Bigr)^{1/p}
$$
(отметим, что сделанные выше высказывания не очевидны). Положим также, что класс
$\frL_\infty(H)$ -- это пространство компактных операторов, снабженное стандартной 
операторной нормой (ее можно записать как $\|A\|_\infty=\max \lambda_j$).

\sm

{\sc Пример.} Фиксируем $p\in [1,\infty]$. Рассмотрим группу всех обратимых операторов
вида $1+A$, где $A\in\frL_p$. Снабдим эту группу топологией, индуцированной
$\|\cdot\|_p$-нормой. Тогда получается полная сепарабельная топологическая группа.
\hfill $\lozenge$

\sm

Ввиду обилия операторных топологий можно напридумывать тьму различных операторных групп. 
Однако в сколько-либо содержательных построениях (разумеется, это функция
времени и вкусов комментаторов) участвуют не столь уж многие из них (полная унитарная группа
с ее родственниками оказываются интересными объектами, а  группы из последнего
примера, как будто, не очень). Естественные топологии, возникающие в теории 
унитарных представлений операторных групп, -- это топологии шейловского (D. Shale) типа.

\sm

{\it Шейловская топология. Пример}. Рассмотрим группу $\GLO$ обратимых операторов
в вещественном гильбертовом пространстве, представимых в виде
$U(1+T)$, где $U$ -- ортогональный оператор, а $T$ -- оператор Гильберта--Шмидта.
Эта группа возникает как группа симметрий  бесконечной гауссовой меры,
см. ниже \S 8. Рассмотрим полярное разложение $g=R S$ элемента группы
$\GLO$, здесь $R$ -- ортогональный оператор, а $S$ -- положительный оператор.
Заметим, что $S-1$ -- оператор Гильберта--Шмидта. Такой оператор $S$
однозначно представим в виде $S=\exp R$, где $R$ -- самосопряженный оператор Гильберта--Шмидта.
Как абстрактное множество $\GLO$ является произведением  ортогональной группы 
$\bfO(\infty)$ и пространства самосопряженных операторов Гильберта--Шмидта.
Теперь мы снабжаем ортогональную группу слабой топологией,  а второй сомножитель 
гильберт--шмидтовской топологией. Группу $\GLO$ мы снабжаем топологией произведения.

\sm

{\it  Другие группы с шейловской топологией.} 
a) Есть две группы такого типа, играющие выделенную роль в теории 
представлений и матфизике ({\it<<группы автоморфизмов канонических антикоммутационных
и коммутационных соотношений}>>).
Рассмотрим комплексное гильбертово пространство  $H$ со скалярным произведением
$\la v, w\ra$. Представим $\la v, w\ra$ в виде суммы действительной и мнимой
частей,
$$\la v, w\ra=(v, w)+i\{v,w\}.$$
Овеществим наше пространство $H$, то есть положим, что
$H_\R$ -- это тоже самое гильбертово пространство, только мы рассматриваемое
как вещественное. Тогда $(v,w)$ -- скалярное произведение в этом пространстве,
а $\{v,w\}$ -- кососимметрическая билинейная форма. Соответственно мы получаем
 бесконечномерную ортогональную группу $\bfO(\infty)$, состоящую из обратимых
вещественно-линейных
операторов, сохраняющих $(v,w)$, и бесконечномерную симплектическую группу 
$\mathbf{Sp}(\infty)$, состоящую из обратимых
вещественно-линейных операторов, сохраняющих $\{v,w\}$. При этом
$$
\bfO(2\infty)\cap \mathbf{Sp}(2\infty)=\bfU(\infty)
.$$

Теперь в каждой из групп $\bfO(2\infty)$, $\mathbf{Sp}(2\infty)$
возьмем подгруппу, состоящую из всех операторов,
отличающихся от унитарного оператора на гильберт-шмидтовскую поправку,
т.е., операторов, представимых в виде $U(1+T)$,  где $U$ -- унитарен,
а $T$ -- оператор Гильберта--Шмидта.
 Эти две группы являются естественными
 группами симметрий фермионного и бозонного   пространств Фока (Ф. А. Березин),
 см. \cite{Ber}, \cite[Chap. 4, 6]{Ner1}).




\sm

b) Как заметил Г.И.Ольшанский, есть 20 групп похожего вида, которые получаются так.
Мы берем пары $G(n)\supset K(n)$, где $G(n)$ -- классическая группа Ли, а $K(n)$ --
максимальная компактная подгруппа (например $\U(n,n)\supset \U(n)\times \U(n)$).
Берем аналогичные группы бесконечных матриц,
 из $G(\infty)\supset K(\infty)$, ограниченных в смысле $\ell_2$.
 Далее берем в $G(\infty)$ подгруппу, состоящую из всех матриц,
которые отличаются от матриц из $K(\infty)$ на гильберт--шмидтовскую поправку
(как было чуть выше для $\mathbf{GL}(\infty,\R)\supset \mathbf{O}(\infty,\R)$ и
$\mathbf{Sp}(2\infty,\R)\supset \mathbf{U}(\infty)$).

Аналогично мы можем брать пары 
$G^\circ(n)\supset K(n)$, где  $G^\circ$ -- классическая компактная группа Ли, а 
$K(n)$ -- симметрическая подгруппа%
\footnote{Формальное определение. Инволюцией в группе называется автоморфизм $\sigma$,
	такой, что $\sigma^2=1$. Симметрическая подгруппа -- множество неподвижных
	точек некоторой инволюции.  Примеры инволюций на классических группах
	-- комплексное сопряжение, отображение $g\mapsto g^{t-1}$, 
	сопряжение матрицей $\begin{pmatrix} 1&0\\0&-1\end{pmatrix}$.},
 и аналогичным образом переходим к пределу
(как это было чуть выше для  $\mathbf{O}(2\infty)\supset \mathbf{U}(\infty)$).

 Для таких
<<бесконечномерных классических групп>> есть хорошо развитая теория унитарных представлений,
см.  \cite{Olsh-GB}, \cite[гл. 9]{Ner1}, \cite{Pick1},  см. также ниже \S11.

\sm

См. \cite{Sha}, \cite[\S\S I.1, IV 3-4, VI.2, IX.1-2]{Ner1}.

\sm

{\bf\punct Группы 
	преобразований лебеговских пространств.%
	\label{ss:Ams}}
Пусть $M$ -- лебеговское пространство, эквивалентное отрезку $[0,1]$.
Рассмотрим группу $\Ams(M)$, состоящую из всех биективных с точностью до п.в.
отображений $M\to M$, сохраняющих меру $\mu$; два преобразования 
$g_1$, $g_2\in\Ams(M)$ считаются эквивалентными, если $g_1(m)=g_2(m)$ на множестве полной меры.
На группе $\Ams(M)$ вводится топология из следующего условия: последовательность
$g_j$ сходится к $g$, если для любых двух измеримых множеств
$A$, $B\subset M$ имеет место сходимость $\mu(g_j(A)\cap B)\to \mu(g(A)\cap B)$.

\sm

{\sc Задача.} Для любого преобразования $g\in \Ams(M)$ определен унитарный линейный оператор
$$
T(g)f(m)=f(g(m))
$$
в $L^2(M)$. Тем самым получается вложение $\Ams(M)\to \bfU(\infty)$.
 Покажите, что топология на $\Ams(M)$ индуцирована слабой операторной топологией на $\bfU(\infty)$.
Покажите, что подгруппа $\Ams(M)\subset \bfU(\infty)$ замкнута.
\hfill $\lozenge$

\sm

{\sc Задача.} Несепарабельная топология на $\Aut(M)$.

\sm

a) Введем на
$\Aut(M)$   метрику, положив, что $d(g,h)$ равно мере множества точек 
$m\in M$, таких, что $g(m)\ne h(m)$. Покажите, что эта метрика двустроннне инвариантна,
что полученная топология несепарабельна и совпадает с топологией, индуцированной равномерной
топологией на $\bfU(\infty)$.

\sm

b) Пусть $\SO(3)$ естественным образом действует на сфере, и соответственно в $L^2(S^2)$. Покажите, что топология
на $\SO(3)$, индуцированная  равномерной топологией на $\bfU(\infty)$, дискретна.
\hfill$\lozenge$

 \sm
 
 Далее, рассмотрим группу $\Gms(M)$ биекций п.в.  из $M$ в себя, переводящих
 меру $\mu$ в эквивалентную меру. Для таких отображений определена {\it производная Радона--Никодима}
$g'(m)$ из условия:
$$
\int_A g'(m)\,d\mu(m)=\mu(A)
$$
для любого измеримого множества $A\subset M$. Для $g\in \Gms(M)$ рассмотрим
унитарный оператор в $L^2(M)$, заданный формулой
$$
T(g)f(m)=f(g(m))\, g'(m)^{1/2}
.$$ 
Мы определяем топологию на $\Gms(M)$ как топологию, индуцированную
со слабой топологии на унитарной группе.

\sm

{\sc Замечание.} Это определение топологии равносильно такому более сложному, но легче применяемому. 
Пусть $A$, $B$ --
измеримые подмножества. Обозначим через $\kappa[g;A,B;t]$ распределение функции
$g'(m)$ на множестве $A\cap g^{-1}(B)$ (это мера на множестве положительных чисел, $t$
здесь обозначает координату на $\R^\times_+$). Последовательность $g_j$ сходится к $g$,
если для любых $A$, $B$ имеют место слабые сходимости мер 
$\kappa[g_j;A,B;t]\to \kappa[g;A,B;t]$, $t\cdot \kappa[g_j;A,B;t]\to t\cdot\kappa[g;A,B;t]$
на
$\R^\times_+$.
\hfill $\lozenge$

\sm

См. 
\cite[\S VIII.4, \S X.3]{Ner1}, \cite{Kech1}, \cite{Pes},  \cite{Ner-biinv}.
Мы вернемся к этим группам в \S12. Стоит отметить, что свойствам
индивидуальных элементов группы $\Ams$ посвящена огромная литература
(эргодическая теория).

\sm

{\bf\punct  Полная бесконечная симметрическая группа.%
\label{ss:S-infty}} Рассмотрим группу $\ov S(\infty)$
 всех перестановок натурального
ряда. Обозначим через $K_\alpha$ подгруппу в $S_\infty$, состоящую из перестановок,
оставляющих на месте точки $1$, $2$, \dots, $\alpha$.
Топология на $\ov S(\infty)$ вводится из условия: подгруппы $K_\alpha$ открыты и образуют фундаментальную
систему окрестностей единицы. Разумеется, эта топология вполне несвязна.

\sm

См. \cite{Cam}, \cite{Olsh-kiado}, \cite[\S\S VIII.1--2]{Ner1}, \cite{Pes}.

\sm

{\bf\punct Бисимметрическая группа.\label{ss:bisymmetric}} Есть содержательная теория унитарных представлений 
<<бесконечной симметрической группы>>, но это не теория представлений полной симметрической группы
$\ov S_\infty$ и не теория представлений дискретной группы $S_\infty$, состоящей из перестановок
с конечным носителем. Б\'ольшая часть работ по представлениям бесконечной симметрической группы
может быть отнесена к  теории представлений следующего объекта. {\it Бисимметрическая
группа} $\bbS(\infty)$ -- это подгруппа в $\ov S_\infty\times \ov S_\infty$, состоящая из пар
$(g,h)$, таких, что $gh^{-1}\in  S_\infty$. Мы не обсуждаем здесь, почему так получается
и ограничиваемся ссылками на \cite{Olsh-chip}, \cite{Ner-umn} (в последней статье большая
современная библиография). Ниже бисимметрическая группа появляется в
п.\ref{ss:bisymmetric-0.5} и в \S10.

\sm

{\bf\punct Группы петель и группы отображений.} Группа петель -- это группа гладких отображений
из окружности $S^1$ в компактную группу Ли $G$ (например, $G=\SO(n)$, $\U(n)$, $\Sp(n)$),
группа снабжается топологией равномерной сходимости всех производных.

Разумеется,  вместо окружности можно рассматривать  любое многообразие, а
вместо $G$  -- любую группу Ли. Однако, наиболее содержательная теория существует именно для групп петель.

\sm

См. \cite{PS}, \cite{Ism}, \cite[\S IX.7]{Ner1},  \cite{Ism-su2}.
 
\sm

Пусть $M$ -- пространство с непрывной вероятностной мерой.
Для некоторых групп%
\footnote{Нужны группы, обладающие нетривиальными действиями аффинными изометриями гильбертова 
пространства (нетривиальное аффиное действие - значит не сводящееся к линейному посредством 
смещения начала координат). Такие действия есть, например,
у группы $\SL(2,\R)$, $\SL(2,\C)$, у групп $\SO(1,n)$, $\SU(1,n)$, у упомянутых ниже групп
автоморфизмов деревьев Брюа--Титса. При наличии такого действия есть <<конструкция Араки>>, автоматически
изготовляющая представления групп $\cF(M,G)$.}
$G$ есть любопытная теория унитарных представлений групп $\cF(M,G)$ измеримых функций
$M\to G$. Мы ограничиваемся ссылками \cite{VGG}, \cite[\S X.1-2]{Ner1}.

\sm

{\bf\punct Классические группы диффеоморфизмов.} Известны следующие
естественные группы диффеоморфизмов связных гладких многообразий. 

a) Группа $\Diff(M)$ всех $C^\infty$-диффеоморфизмов многообразия $M$.

\sm

b) Фиксируем форму объема $\Omega$ на многообразии $M^n$ и возьмем группу $\SDiff(M^n,\Omega)$ диффеоморфизмов, сохраняющих 
форму
$\Omega$.

\sm

c) Рассмотрим $2n$-мерное многообразие.
 Фиксируем симплектическую форму $\omega$, т.е., дифференциальную 2-форму $\omega$, такую, что $d\omega=0$, а 
  $n$-ная внешняя степень $\omega^{\wedge n}=\omega\wedge\dots\wedge\omega$ не обращается в ноль ни в одной точке.
  
  Стандартный пример такой формы -- форма 
  $$\sum_{j=1}^n dp_j\wedge dq_j$$
  на $\R^{2n}$ с координатами $(p_1,\dots,p_n;q_1, \dots, q_n)$ (согласно теореме Дарбу,
  локально любая симплектическая форма приводится к такому виду).
 Далее рассматривается группа $\Symp(M,\omega)$ всех диффеоморфизмов, сохраняющих  форму $\omega$. Она называется
 группой {\it симплектических
 диффеоморфизмов} или {\it симплектоморфизмов}, см. \cite[\S40, гл.8]{Arn-mech}.
 
 \sm
 
 d) Рассмотрим $2n+1$-мерное многообразие и {\it контактную форму} $\lambda$ на нем.
Напомним, что это 1-форма, такая, что $\lambda\wedge (d\lambda)^{\wedge n}$
не обращается в ноль ни в одной точке.

Стандартный пример контактной формы -
форма
$$
\Lambda=dz+\sum_{j=1}^n x_j \,dy_j
$$
на $\R^{2n+1}$ с координатами $(x_1,\dots,x_n; y_1,\dots, y_n;z)$.
Локально любая контактная форма приводится к виду $h\cdot \Lambda$,
где $h$ -- не обращающаяся в ноль функция.

Далее рассматривается группа $\Cont(M,\lambda)$ всех диффеоморфизмов
которые переводят $\lambda$ в $\psi\cdot \lambda$, где $\psi$ -- функция на многообразии.
Получается группа {\it контактных диффеоморфизмов} или {\it контактоморфизмов},
см. \cite[Добавл.4]{Arn-mech}, \cite[\S 5.2]{ALV}. 

Это можно понимать еще так. Ядро формы $\lambda$ в фиксированной точке $x\in M$ -- 
это гиперплоскость в касательном пространстве к $M$ в точке $x$. Такая
гиперплоскость определена в каждой точке, и
таким образом на многообразии определено поле касательных гиперповерхностей (<<распределение>>
в смысле Фробениуса, \cite[\S 5.3]{ALV}). Контактоморфизмы -- это диффеоморфизмы, переводящие
это распределение в себя.

\sm

{\it Топологии.} В случае компактного многообразия $M$ на группе $\Diff(M)$ всех
$C^\infty$-диффеоморфизмов имеется единственная естественная топология -- топология равномерной
сходимости диффеоморфизмов и их производных.

В случае некомпактного многообразия  $M$ естественно возникают разные группы $C^\infty$-диффеоморфизмов.
На группе всех диффеоморфизмов естественная топология - топология равномерной сходимости 
диффеоморфизмов и всех их производных на каждом компактном подмножестве в  $M$.

Введем топологию на группе $C^\infty$-диффеоморфизмов $\Diff_{comp}(M)$ с компактным носителем.
Последовательность диффеоморфизмов $g_j$ сходится к $g$, если существует компактное множество
$K\subset M$, такое, что $g_j$ тривиальны вне $K$, а на $K$ имеет место сходимость
$g_j\to g$ равномерно со всеми производными. 

\sm

{\it  Инварианты форм.} У формы объема $\omega$ на компактном многообразии $M$ есть
единственный инвариант под действием группы всех диффеоморфизмов -- полный объем $\int_M\omega$
(M.~Морс).
Поэтому для данного компактного многообразия <<группа диффеоморфизмов, сохраняющих объем>>,
определена канонически с точностью до изоморфизма.

В случае симплектических форм это уже не так. Разумеется, симплектическая форма задает класс
вторых де Рамовских когомологий многообразия, и эти классы могут быть разными. 
Однако когомологичные формы могут быть неэквивалентными,  их различение является тонкой задачей
и входит в предмет симплектической топологии. Аналогична ситуация с контактными формами.

\sm

{\it Группы классов.} Вообще говоря, перечисленные группы несвязны.
Группой классов (mapping class group) называется фактор группы диффеоморфизмов (одного из типов
a-d) по компоненте связности. Группа классов диффеоморфизмов двумерной сферы $S^2$ совпадает
с $\Z_2$ (диффеоморфизмы или сохраняют ориентацию, или меняют ее).
Но у сферы $S^6$ группа классов, сохраняющих ориетацию, неожиданно оказывается группой $\Z_{28}$.
Если взять два 7-мерных шара $B_1$, $B_2$, взять диффеоморфизм их границ
$\partial B_1\to \partial B_2$, и склеить границы посредством диффеоморфизма, нетривиального в
$\Z_{28}$,
то получится сфера Милнора -- топологическая 7-мерная сфера, не диффеоморфная $S^7$.

Впрочем, группа классов диффеоморфизмов двумерной компактной поверхности 
(т.н. группа Тейхмюллера) уже является крайне нетривиальным объектом, много
изучавшимся. См. \cite{FM}, \cite{Iva}, \cite{Kash}

\sm

{\it Простота групп.} Теперь возьмем связные компоненты классических групп диффеоморфизмов.
В случае полной группы диффеоморфизмов $\Diff(M)$ (или $\Diff_{comp}(M)$)  мы получаем группу,
которая является простой как дискретная группа. Для разных многообразий 
группы $\Diff(M)$ не изоморфны, а любой автоморфизм $\Diff(M)$ 
является внутренним.

Аналогичное утверждение верно для группы диффеоморфизмов, сохраняющих объем.
В случае групп  симплектоморфизмов имеют место похожие утверждения, но они формулируются
несколько сложнее. А именно, компонента связности $\Symp(M,\omega)$ содержит <<большую>>
простую нормальную подгруппу. Фактор по ней в случае компактных многообразий -- абелева группа Ли
(она явно описыватся);
в случае некомпактных -- фактор является группой Ли $H$, которая или абелева,
или содержит одномерную подгруппу $R$, такую, что фактор $H/R$ -- абелев. При 
этом группа симплектоморфизмов помнит как многообразие так и симплектическую форму.
См.  \cite{Ban}, \cite{Ban0}.

\sm

 Встает вопрос, почему принято рассматривать именно такие группы, а не какие-либо другие 
 (например, почему бы не рассмотреть группу диффеоморфизмов,
 сохраняющих какую-нибудь хорошую
17-форму). Ответ на него дается теоремой Картана о примитивных псевдогруппах, 
обсуждаемой в п.\ref{ss:pseudogroups}.


\sm

{\bf \punct Алгебры Ли векторных полей.}
Мы вкратце обсудим алгебры Ли групп классических групп диффеоморфизмов (они понадобятся нам в 
пп. \ref{ss:evolution}, \ref{ss:hofer}, см. хорошее введедение в \cite[\S39]{Arn-mech}).
Для  простоты, рассмотрим случай $\R^n$.
Рассмотрим гладкое векторное поле $v=(v_1,\dots, v_n)$ и дифференциальное уравнение
$$
x'(t)= v(x(t)),\qquad g(0)=a.
$$
Определим $g_s(a):=x(s)$. Мы получаем отображение $g_s$ из $\R^n$ в себя, 
точнее диффеоморфизм из некоторой открытой области  $\Omega_t^+\subset \R^n$ 
в некоторую открытую область  $\Omega_t^-\subset \R^n$
(потому что решения дифференциального уравнения могут уходить на бесконечность
за конечное время, в случае компактных многообразий $g_s$ определен всюду). Если $t$ стремится к нулю, то области $\Omega_t^+$
и $\Omega_t^-$ исчерпывают $\R^n$. Эти отображения удовлетворяют свойству
$$
g_t\circ g_s= g_{t+s}.
$$
на естественной области определения $\Omega_{t+s}^+$ обеих частей.
Семейство отображений $t\mapsto g_t$ называется {\it потоком векторного поля $v$},
см. \cite[гл. 1]{Arn-ur}, \cite[\S1.1]{Ros}. Это отображение является аналогом экспонеты 
(как обычной экспоненты, 
так и экспоненты в группах Ли).

Как известно, векторному полю $v$ ставится в соответствие дифференциальный оператор
$$
D_v:=\sum v_j(x)\frac \partial{\partial x_j}.
$$
Формально
$$
\exp(t D_v) f(x)=f(g_t(x)).
$$
Коммутатор операторов дифференциальных операторов $D_u$, $D_v$
равен 
$$
\sum_j (D_u v_j- D_v u_j)\frac \partial{\partial x_j}
,$$
т.е., имеет вид $D_w$, где $w$ - вектор с компонентами $w_j=(D_u v_j- D_v u_j)$.
Этот вектор называется {\it коммутатором} $[u,v]$ векторных полей $u$, $v$
(геометрическую интерпретацию коммутатора см., например, в \cite[\S39]{Arn-mech}).
Эта операция удовлетворяет тождеству Якоби, и в итоге мы получаем
{\it алгебру Ли векторных полей}. Ее естественно считать алгеброй Ли группы
всех диффеоморфизмов.

Отметим, что полученная алгебра Ли - настоящая, а в каком смысле группа диффеоморфизмов
является группой Ли - вопрос отдельный. В любом случае,  на эвристическом уровне это так.

Теперь рассмотрим группу $\SDiff$ диффеоморфизмов $g$ пространства $\R^n$, сохраняющих стандартный объем. 
Это значит, что  якобиан отображения $g$ во всех точках равен 1, или, иными словами, матрица
Якоби $J(g,x)$ содержится в $\SL(n,\R)$ для всех $x$. Если $g_t$ -- поток,
состоящий из диффеоморфизмов, сохраняющих объем, то
$$
\lim_{\epsilon\to 0} \frac 1\epsilon \Bigl(J(g_\epsilon,x)-1\Bigr)\in \frs\frl(n,\R)
.
$$
Если $u$ -- векторное поле, порождающее поток,
то матрица, которая стоит под знаком предела, равна
$$
\begin{pmatrix}
 \frac {\partial u_1}{\partial x_1}& \frac {\partial u_1}{\partial x_2}&\dots\\
  \frac {\partial u_2}{\partial x_1}& \frac {\partial u_2}{\partial x_2}&\dots\\
  \vdots&\vdots&\ddots
\end{pmatrix},
$$
т.е., это линейная часть векторного поля в точке. След этой матрицы равен нулю,
т.е. 
$$
\div u=\sum \frac {\partial u_j}{\partial x_j}=0.
$$
Таким образом, алгебра Ли группы диффеоморфизмов, сохраняющих объем, состоит из бездивергентых
векторных полей (тут еще нужно обратное рассуждение, см. \cite[\S36A]{Arn-mech}).

Для группы симплектических диффеоморфизмов та же аргументация показывает, что
матрица Якоби должна лежать в $\Sp(2n,\R)$, а линейная часть 
векторного поля должна лежать в алгебре Ли $\frs\frp(2n,\R)$
(разные другие описания таких векторных полей,  см.  в \cite[\S40]{Arn-mech}).

\sm

См. \cite[\S36]{Arn-mech},  \cite{Ism}.

\sm

{\bf \punct Теорема Картана о примитивных псевдогруппах.%
\label{ss:pseudogroups}}
Пусть $M$ - гладкое вещественное или аналитическое комплексное многообразие.
{\it Локальный диффеоморфизм} $M$ это диффеоморфизм $\phi$ открытого подмножества
$U\subset M$ на открытое подмножество $V\subset M$.

\sm

{\it Псевдогруппа} $\Gamma$ -- это набор локальных диффеоморфизмов многообразия,
удовлетворяющий следующим условиям.

\sm

1) Если $\psi:U\to V$ содержится в $\Gamma$, а $W\subset U$ -- открытое подмножество,
то ограничение $\psi\bigr|_W$ содержится в $\Gamma$.

\sm

2) Если $\psi\in \Gamma$, то $\psi^{-1}\in\Gamma$.

\sm

3) Если $\phi:U\to V$, $\psi:V\to W$ содержатся в $\Gamma$, то $\psi\circ\phi\in \Gamma$.

\sm

4) Пусть $\psi:\cup U_j \to W$ -- локальный диффеоморфизм. Если 
все $\psi\bigr|_{U_j}$ содержатся в $\Gamma$, то $\psi\in\Gamma$.

\sm

{\sc Пример.}  Пусть $M=\C\cup\infty$ -- сфера Римана. Рассмотрим псевдогруппу всех
голоморфных локальных диффеоморфизмов $M$. Это бесконечномерный объект, в то время как 
группа всех голоморфных диффеоморфизмов $M$ изоморфна $\SL(2,\C)/\{\pm 1\}$
и состоит из дробно-линейных преобразований.  Т.е., локальных решений системы Коши-Римана 
очень много, а глобальных мало. Это показывает различие между групповой и псевдогрупповой точками зрениями.
\hfill $\lozenge$

\sm

Псевдогруппа называется {\it транзитивной},  если выполнено следующее условие:

\sm

(T) Для любых двух точек $x$, $y\in M$ существует $\psi\in \Gamma$, такой, что
$\psi x=y$.

\sm

Псевдогруппа называется {\it примитивной}, если выполнено условие

\sm

(P)  не существует $\Gamma$-инвариантного слоения%
\footnote{Слоение -- это разбиение многообразия на подмногообразия, которое локально
устроено как разбиение $\R^n$ на семейство параллельных аффинных подпространств.}
ни в какой области  $\U\subset M$.

\sm

Псевдогруппа называется {\it гладкой порядка $k$}, если выполнено условие:

\sm

($\mathrm S^1$) Пусть $\psi:U\to V$ -- локальный диффеоморфизм.
Пусть для любого $x\in U$ существует $\psi_x\in\Gamma$, такой, что
$\psi_x(x)=\psi(x)$ и все частные производные порядка $\le k$ у $\psi$
и $\psi_x$ в точке $x$ совпадают. Тогда $\psi\in \Gamma$.

\sm

{\sc Замечание.} Если псевдогруппа определяется как семейство диффеоморфизмов,
удовлетворяющих некоторой системе дифференциальных уравнений (как для классических групп диффеоморфизмов),
то последнее условие автоматически выполнено.

\sm

Классические группы диффеоморфизмов дают 4 различных типа псевдогрупп.
Теперь приведем еще два примера.

\sm

a) {\it Диффеоморфизмы с локально постоянным якобианом.}
Рассмотрим многообразие с атласом $W_j$. В каждой карте возьмем форму объема
$\omega_j$, определенную с точностью до постоянного множителя. Пусть на пересечениях карт
$W_i\cap W_j$ формы  $\omega_i$, $\omega_j$
отличаются постоянным множителем. Мы рассматриваем диффеоморфизмы, сохраняющие эту
структуру. Иными словами, рассмотрим произвольную пару  карт $W_k$, $W_l$. Мы требуем
чтобы прообраз формы $\omega_k$ при диффеоморфизме $g$ совпадал с $s\cdot \omega_l$
на $g^{-1} W_k\cap W_l$, где $s=s_{ij}\in \R$.

\sm

{\sc Пример.} Рассмотрим пространство $\R^k$ со стандартной формой объема. 
Рассмотрим фактор-пространство
$\R^k\setminus 0$ по отношению эквивалентности $x\sim 2x$.
Оно снабжено
 структурой 
только что описанного типа.
\hfill $\lozenge$

\sm

b) {\it Проективно симплектические диффеормозфизмы.}
Аналогично, рассмотрим многообразие с атласом $W_j$, в каждой карте возьмем симплектическую
форму, определенную с точностью до постоянного множителя. На пересечениях карт эти формы
совпадают с точностью до постоянного множителя. Далее берем группу (<<проективно симплектических>>) диффеоморфизмов, сохраняющих эту структуру.

\sm

{\sc Пример.} Рассмотрим пространство $\R^{2n}$ со стандартной симплектической формой
$\sum_{j=1}^n dx_j\wedge dx_{j+n}$. Рассмотрим фактор-пространство
$\R^{2n}\setminus 0$ по отношению эквивалентности $x\sim 2x$. Оно снабжено
проективно симплектической структурой.
\hfill $\lozenge$

\sm

{\it Классификационная теорема.} Любая примитивная транзитивная гладкая псевдогруппа диффеоморфизмов
локально изоморфна псевдогруппе одного из 12 следующих типов. Во-первых, это упомянутые 
выше 6 типов псевдогрупп: всех диффеоморфизмов; диффеоморфизмов, сохраняющих объем; диффеоморфизмов
с постоянным якобианом; симплектических диффеоморфизмов; проективно симплектических диффеоморфизмов;
контактных диффеоморфизмов. Во-вторых, это 6 аналогичных типов псевдогрупп
комплексно-аналитических диффеоморфизмов.

\sm

Это странная великая теорема, которая редко  используется напрямую, но неявно она маячит за 
описаниям групп симметрий разных геометрических структур: она показывает, что интересных возможностей здесь не так уж много.

\sm

{\sc Замечание.} Еще раз подчеркнем, что речь  идет именно о локальной классификации. 
\hfill $\lozenge$

\sm

{\sc Замечание.} Контактная псевдогруппа является гладкой порядка 2. 
Все остальные псевдогруппы Картана являются гладкими порядка 1.
\hfill $\lozenge$

\sm

{\sc Пример.} Локальный диффеоморфизм $\R^n$ называется {\it конформным}, если его дифференциал в
 каждой точке есть композиция гомотетии и ортогонального преобразования. Мы получаем
  псевдогруппу
конформных отображений. Легко проверить, что она транзитивна, примитивна и является гладкой
 порядка 1. Если $n\ne 2$, то она не совпадает ни с одной из бесконечномерных картановских
  псевдогрупп. Поэтому она конечномерна. Для читателей, знакомых со структурной теорией групп Ли: 
  выведите отсюда теорему Лиувилля об описании всех локально конформных преобразований $\R^n$.
  \hfill $\lozenge$
  
  \sm

См. \cite[\S\S 7.6-7.7]{ALV}, \cite{SS}.

\sm

{\bf \punct Геодезические потоки на группах диффеоморфизмов.%
\label{ss:evolution}}
Оказывается, что есть много эволюционных уравнений с частными производными,
которые интерпретируется как уравнения на геодезические на группах диффеоморфизмов.
Приведем самый старый и самый элементарный пример (В.И.Арнольд).

Рассматривается идеальная несжимаемая жидкость. Идеальность означает отсутствие вязкости
(внутреннего трения)
 и теплопроводности. Движение такой жидкости можно рассматривать как траекторию на группе
 $\SDiff(\R^3)$
 диффеоморфизмов, сохраняющих объем. А именно, мы смотрим, куда приходит каждая частица
 за время $t$, и получаем диффеоморфизм, зависящий от $t$. Ввиду того, что жидкость несжимаема,
 эти диффеоморфизмы сохраняют объем. Для простоты, положим, что на жидкость не действуют внешние
 силы. Движение такой жидкости описывается следующим {\it уравнением Эйлера.}
 Обозначим скорость жидкости в точке $x$ в момент времени $t$ через $v=v(t,x)$.
 Обозначим давление через $p$. Тогда
 \begin{align}
  \frac \partial {\partial t} v+ (v\cdot \nabla) v=-\grad p;
  \label{eq:eu1}
  \\
  \div v=0.
  \label{eq:eu2}
 \end{align}
Здесь
$$
v\cdot \nabla:= \sum v_j \frac \partial {\partial x_j}.
$$
Давление $v$ определено с точностью до аддитивной константы. В принципе, мы могли 
бы давление не вводить, а написать
$$
 \rot\Bigl(\frac \partial {\partial t} v+ (v\cdot \nabla) v\Bigr)=0
.$$
Сейчас мы покажем (на эвристическом уровне), что это уравнение
является уравнением на геодезические относительно правоинвариантной
$L^2$-метрики на группе $\SDiff(\R^3)$.

Обозначим через $g(t,x)$ положение частицы, находившейся
в начальный момент времени в  точке $x$, в момент времени $t$.
Тогда отображение $x\mapsto g(t,x)$ -- это диффеоморфизм, зависящий от времени.
Скорость частицы задается формулой
$$
v(t,x)=  \frac \partial {\partial t} g(t,x).
$$
Теперь вычисдяем ускорение частицы, т.е. считаем
$$
\frac d {d t} v(t,g(t,x))= \frac \partial {\partial t} v(t,y)\biggr|_{y=g(t,x)}+
\sum_j \frac{\partial v(t,y)} {\partial y}\biggr|_{y=g(t,x)}\cdot v(t,x)
$$
В итоге, уравнение (\ref{eq:eu1}) превращается в
\begin{equation}
\frac d {d t} v(t,g(t,x))=-\grad p
\label{eq:eu3}
\end{equation}
(собственно, это и есть вывод уравнений Эйлера).

Давайте рассматривать группу $\SDiff(\R^3)$ как <<подмногообразие>> в <<многообразии>> 
$\Diff(\R^3)$  всех диффеоморфизмов.  Касательное пространство к $\Diff(\R^3)$
мы отождествим с пространством всех векторных полей. Так как это можно сделать разными
способами, 
уточним последнее высказывание. Если $h_\epsilon(x)$ - траектория
в группе диффеоморфизмов, то касательный вектор к ней в точке $\epsilon=0$ --
это векторное поле, значение которого в точке $x$ равно
$v(x)=\frac d{d\epsilon}h_\epsilon(x)\Bigr|_{\epsilon=0}$.
Введем на  <<многообразии>> $\Diff(\R^3)$ риманову метрику, положив, что скалярное
произведение касательных векторов равно
$$
\la u, v\ra=\int_{\R^3} \sum_j u_j(x) v_j(x)\, dx.
$$

 Ограничим нашу риманову метрику на группу
$\SDiff(\R^3)$. Касательное пространство к $\SDiff(\R^3)$ состоит
из бездивергентных полей. Легко видеть, что полученная риманова метрика является правоинвариантной.
Рассмотрим траекторию на $\SDiff(\R^3)$, удовлетворяющую уравнению Эйлера.
Уравнение (\ref{eq:eu3}) показывает, что вектор ускорения траектории - это градиентное векторное
поле. Пусть $v$ бездивергентно. Тогда
$$
\la \grad p, v\ra=\int_{\R^3} 
\sum_j \frac{\partial p}{\partial x_j} v_j\,dx=
-\int_{\R^3} p\, \Bigl(\sum_j \frac{\partial v_j}{\partial x_j}\Bigr)\,dx
=0.
$$
Таким образом, ускорение траектории ортогонально подмногообразию. Это один
из вариантов определения геодезических на подмногообразии.

\sm

См. \cite{Arn1}, \cite[Добавл.2 ]{Arn-mech}, \cite{AKh}, \cite{Bre}, \cite{EM}, \cite{Khe}, \cite{HW},
\cite{MW}, \cite{Mich}.

\sm

{\bf \punct  Пространство всех подмногообразий.} 
Сначала рассмотрим пространство $\R^3$, снабженное стандартным объемом. 
Рассмотрим множество $\cL$ всех подмногообразий $\R^3$, эквивалентных окружности
  $S^1$ в $\R^3$. Мы будем рассматривать $\cL$ как <<бесконечномерное многообразие>>.
  Оказывается, что это многообразие является симплектическим.
  
  Касательный вектор к $\cL$ -- это сечение нормального расслоения к кривой. 
  Возьмем точку $\ell$ пространства $\cL$, т.е., вложенную окружность, возьмем ее параметризацию
  $\gamma(t)$, где $t\in [0,T]$. Возьмем два сечения нормального расслоения,
  $\xi_1(t)$, $\xi_2(t)$. По определению $\xi_j(t)$ - это вектор в $\R^3$, 
  определенный с точностью до эквивалентности 
  $$
  \xi_j(t)\sim \xi_j(t)+a(t)\gamma'(t),\qquad\text{где $a(t)\in\R$.}
  $$
  Рассмотрим величину 
  $$
  \omega_\ell(\xi_1,\xi_2):=\int_0^T \bigl(\gamma'(t),\xi_1(t), \xi_2(t)\bigr)
  \,dt
  $$
  где $\bigl(\dots\bigr)$ обозначает смешанное произведение трех векторов.
  Легко видеть, что это выражение не зависит от выбора параметризации
 $\gamma(t)$ и задает кососимметрическую билинейную форму на касательном пространстве.
 
 \sm
 
 {\sc Задача.} a) Покажите, что эта 2-форма замкнута. Покажите, что группа 
 $\SDiff(\R^3)$ действует на $\cL$ симплектоморфизмами, т.е., преобразованиями, сохраняющими эту
 симплектическую структуру.
 
 \sm
 
 b) Возьмем гамильтониан%
 \footnote{См. \cite[\S38]{Arn-mech}, формальное
 	определение см. ниже в следующем пункте.} на $\cL$, 
 равный евклидовой длине кривой. Рас\-cмот\-рим соответствующее
 гамильтоново векторное поле. Покажите, что соответстующее сечение нормального расслоения
 (касательный вектор к $\cL$) устроено так: в каждой точке берется бинормальный вектор длины,
 равной кривизне в точке (это т.н. binormal curvature flow; в идеальной гидродинамике, описываемой
 уравнением Эйлера, это соответствует перемещению $\delta$-образного вихря).
 \hfill $\lozenge$

 \sm
 
 Отметим, что точно так же вводится симплектическая структура на пространстве всех подмногообразий 
 коразмерности 2 в многообразии с фиксированной формой объема.
 
 \sm
 
 На самом деле, здесь довольно большой простор для фантазии.
 А именно возьмем $m$-мерное многообразие $M$ с фиксированной $k$-формой $\omega$.
 Возьмем пространство $\cL$ всех $n$-мерных подмногообразий в $M$.
 На нем канонически определена $(k+n-m)$-форма $\omega^\square$, 
 причем 
 $$d(\omega^\square)=d(\omega)^\square.
 $$
В частности, если $\omega$ была замкнута, то замкнута и $\omega^\square$.
Соответственно, группа диффеоморфизмов, сохраняющих форму
$\omega$, действует на $\cL$ преобразованиями, сохраняющими форму
$\omega^\square$.

\sm

См. \cite{HV}, \cite{Has}, \cite[\S3.5]{HW}, \cite{Ism}, \cite{MW}, \cite{Viz}.

\sm

{\bf\punct Группа гамильтоновых диффеоморфизмов и метрика Хофера.%
\label{ss:hofer}}
Рассмотрим симплектическое многообразие $M$ с 2-формой $\omega=\omega_x$
и $C^\infty$-функцию $H$ ({\it гамильтониан}) на $M$. В случае некомпактного многообразия мы предполагаем,
что $H$ имеет компактный носитель. 
Фиксируем точку $x\in M$ и рассмотрим дифференциал
$dH_x$ функции $H$
 в точке $x$. Напомним, что дифференциал функции, по определению, является линейным функционалом на
 касательном пространстве. В силу невырожденности билинейной формы $\omega$
 существует касательный вектор
 $\xi$, такой, что для любого касательного вектора 
 $v$ выполнено
 $$
 dH_x(v)=\omega_x(\xi,v).
 $$
 Для каждой точки $x$ мы имеем касательный вектор в этой точке, т.е.,
 мы получили векторное поле на $M$,
 обозначим его через $\sgrad H:=\xi$. Такие векторные поля называются {\it гамильтоновыми}.
 Если $M=\R^{2n}$ со стандартной формой $\sum dp_j\wedge dq_j$, то это векторное поле задается формулой
 $$
\Bigl(-\frac{\partial H}{\partial q_1}, \dots,- \frac{\partial H}{\partial q_n}, 
 \frac{\partial H}{\partial p_1} ,\dots,  \frac{\partial H}{\partial p_n}\Bigr)
 .
 $$
Гамильтоновы поля сохраняют симплектическую форму, см. \cite[\S\S 37-40]{Arn-mech},
но, вообще говоря, они образуют идеал в алгебре Ли всех симплектических векторных полей
(в случае односвязных многообразий любое симплектическое векторное поле гамильтоново.).

Пусть $g$ -- симплектоморфизм. Мы скажем, что $g$ {\it гамильтонов},
если существует  семейство $g_t$ диффеоморфизмов, гладко зависящих
от времени, такое, что $g_0$ -- единичный диффеоморфизм, а $g_T=g$,
и существует семейство гамильтонианов $H_t$, такое, что
\begin{equation}
\frac d{dt} g_t(x)= \sgrad H_t(g_t(x))
.
\label{eq:sgrad}
\end{equation}
Гамильтоновы диффеоморфизмы образуют подгруппу $\Ham(M,\omega)$ в компоненте связности группы симплектоморфизмов,
в случае односвязных многообразий эти две группы совпадают.

\sm

{\it Метрика Хофера.} На группе $\Ham(M,\omega)$ вводится двусторонне инвариантная метрика
$d(q,r)$ по следующему правилу. Рассматриваются всевозможные пути (гамильтоновы изотопии) $g_t$, $t\in [0,T]$,
в группе $\Ham(M,\omega)$,
удовлетворяющие (\ref{eq:sgrad}) и соединяющие $q$ и $r$ (т.е., $g_0=q$, $g_T=r$).
Определим длину пути $\ell(g_t)$ как
$$
\ell(g_t):=\int_0^T \bigl|\max_{x\in M} H_t(x) - \min_{x\in M} H_t(x)\bigr|\, dt.
$$
Расстояние Хофера $d(q,r)$ -- это точная нижняя грань длин всех путей, соединяющих
$q$ и $r$.
Здесь неожиданно сложным оказывается утверждение, что 
 $d(q,r)\ne0$ при $q\ne r$ (у этого высказывания, впервые полученного Хофером,
 известно несколько доказательств, и все они тяжелые).
 
 \sm
 
 {\sc Задача.} Проверьте остальные аксиомы метрического пространства и инвариантность метрики.
 \hfill $\lozenge$
 
 \sm
 
 {\sc Замечание.} Пусть, для определенности, многообразие некомпактно. Пусть 
 $p>1$. Для любых гамильтоновых диффеоморфизмов $q$, $r$ определим выражение
 $$
d_p(q,r):=\inf_{g_t}\int_0^T \biggl(\int_M  |H_t(x)|^p \,dx\biggr)^{1/p}\, dt 
,
$$
где точная нижняя грань берется по всем гамильтоновыи изотопиям, соединяющим
$q$ и $r$. Оказывается, что $d_p(q,r)$ равно нулю тождественно
 (см. доказательство в \cite[Следствие IV.8.21]{AKh}). \hfill $\lozenge$

\sm

Если принять за аксиомы некоторые сложно доказываемые утверждения, то метрику Хофера
$d(q,r)$ можно оценивать снизу (оценки сверху можно получать путем предъявления гамильтоновых изотопий).
Приведем  пример таких <<аксиом>>.

\sm

{\it Энергия смещения.} Пусть $A\subset M$ -- подмножество. 
Его энергия смещения $e(A)$ это точная нижняя грань $d(e,g)$
по всем $g$, таким, что $g(A)\cap A=\varnothing$.
Если таких $g$ нет, то мы полагаем $e(A)=+\infty$.
Имеет место следующее утверждение:

\sm

{\it Энергия смещения шара радиуса $r$ в стандартном $\R^{2n}$
равна $\pi r^2$.}

\sm

{\sc Задача.} Вывести отсюда, что шар радиуса $R\subset \R^{2n}$
нельзя перевести симплектоморфизмом внутрь цилиндра $p_1^2+q_1^2<r^2$ радиуса $r<R$.
\hfill $\lozenge$

 \sm
 
 Есть также оценки для энергии смещения гладких компактных лаграранжевых подмногообразий%
 \footnote{Напомним, что $n$-мерное подмногообразие $L$ в $2n$-мерном
 симплектическом многообразии $M$ называется лагранжевым, если в любой точке
 $x\in L$ касательное пространство к $L$ является $\omega_x$-изотропным. 
 Любая гладкая кривая в $\R^2$ -- лагранжево подмногообразие. Пример лагранжква многообразия в
 $\R^{2n}$ -- тор, заданный системой уравнений $p_j^2+q_j^2=1$.}. 
 Рассмотрим
 случай $\R^{2n}$. Возьмем 1-форму $\lambda=\sum_j p_j\,dq_j$, так что 
 $d\lambda=\omega$. Эта форма определяет класс первых вещественных когомологий
лагранжева подмногообразия $L$. Рассмотрим ее значения на всевозможных петлях в $L$.
 Получится подгруппа в аддитивной группе $\R$. Допустим, что она дискретна,
 пусть $a>0$ -- ее образующая. Тогда энергия смещения $L$ оценивается как
$e(L)\ge a/2$.

\sm

{\it  Геодезические.} Рассмотрим для простоты случай $M=\R^{2n}$.
Рассмотрим кривую $g_t$ в группе $\Ham$. Отнормируем ее <<скорость>> как 
$$\bigl|\max_{x\in M} H_t(x) - \min_{x\in M} H_t(x)\bigr|=1.$$
Мы скажем, что $g_t$ -- локально кратчайшая кривая,
если у любой точки $t$ есть окрестность такая, что для любых 
$a$, $b$ из этой окрестности $d(g_a,g_b)=|a-b|$.
Оказывается, что однопараметрические подгруппы в $\Ham$ 
(и их сдвиги) являются локально кратчайшими. Более того, кривая
$g_t$ является локально кратчайшей
тогда и только тогда, когда существуют точки $x_+$, $x_-\in M$,
такие, что для всех $t$
$$
\max_{x\in M} H_t(x)=H_t(x_+),\qquad \min_{x\in M} H_t(x)=H_t(x_-).
$$

\sm

См. \cite{McD}, \cite{Pol1}, \cite{Pol2}, \cite{Aud}.
 
 \sm
 
 {\bf\punct Геодезические на группах диффеоморфизмов и близость образов.}
 Есть такой тип задач прикладного характера. Рассмотрим, например, две замкнутые
 гладкие несамопересекающиеся кривые на плоскости. Как установить, похожи
 они или нет, и как установить меру их сходства? Один из способов
 ответа на этот вопрос -- такой. Мы вводим на группе диффеоморфизмов
 какую-нибудь подходящую метрику и вычисляем минимальную длину
 изотопии, переводящей одну кривую в другую. Мы ограничимся ссылками
 \cite{You}, \cite{MM}, \cite{DZ}.

 \sm
 
 {\bf\punct Группа диффеоморфизмов окружности.%
 	\label{ss:diff1}} Группа $\Diff_+(S^1)$ диффеоморфизмов окружности, сохраняющих
 ориентацию, играет выделенную роль в теории представлений.
 Причина этого, видимо, в том, что она допускает многочисленные вложения в классические группы
 с шейловской топологией (см. п.\ref{ss:shale}), что позволяет далее применять технологию <<вторичного квантования>>.
 Приведем пример вложения $\Diff(S^1)$ в группу $\GLO$ почти ортогональных операторов, 
 упомянутую в п.\ref{ss:shale}.
 
 Рассмотрим окружность с координатой $\phi\in \R/2\pi \Z$.
 Фиксируем $s$ из промежутка $(0,1)$. Рассмотрим пространство $C^\infty(S^1)$ вещественных
 гладких функций на окружности, введем в нем скалярное произведение 
 $$
 \la f_1, f_2\ra_s =\int_0^{2\pi} \int_0^{2\pi} \left|\sin\frac {\phi_1-\phi_2}2\right|^{s-1}
 f_1(\phi_1)\,f_2(\phi_2)\,d\phi_1\,d\phi_2.
 $$
 Оказывается, что это скалярное произведение при введенных выше ограничениях на $s$
 является положительно определенным. Обозначим через $V_s$ пополнение $C^\infty(S^1)$
 по этому скалярному произведению (по запасу функций это пополнение является соболевским
 пространством $H^{-s/2}$, но нам явный вид скалярного произведения тоже важен).
 Пусть  группа $\Diff_+(S^1)$ действует в $V_s$ по формуле
 $$
 T_s(q)\, f(\phi)=f(q(\phi))\, q'(\phi)^{(1-s)/2}.
 $$
 Оказывается, что операторы $T_s(q)$ лежат в группе $\GLO(V_s)$,
 тем самым мы получаем действия группы $\Diff_+(S^1)$ на пространстве с гауссовой мерой,
 см. ниже \S 8.
 
 \sm
 
См. \cite[Chapt. VI, \S IX.6]{Ner1}, \cite{PS}, \cite{Ner-Gamma}, \cite{Ner-MIRAN}.
 
 \sm
 
{\bf \punct Топологии на группах диффеоморфизмов.} Понятно, что таких топологий много.
Например, мы можем брать группы диффеоморфизмов гладкости $C^k$. Если многообразие не компактно,
то возникает много возможностей накладывать ограничения на диффеоморфизм на бесконечности
(вопрос тут скорее в том, с какими целями это делать). Мы приведем два  неочевидных примера.

\sm

1) {\it Соболевские диффеоморфизмы,} \cite{IKT}. Рассмотрим $n$-мерное многообразие. 
Пусть $s>n/2+1$ -- вещественное число.  Тогда множество диффеоморфизмов соболевской 
гладкости $s$ замкнуто относительно умножения и образует топологическую группу.

\sm 

2) {\it Классы Данжуа--Карлемана,} \cite{KMR}. 
Гладкая функция на области $\Omega\subset\R^n$ является вещественно аналитической, если
существуют константы
$C$ и $\rho$,  такие, что 
 любая частная производная $\partial^\alpha f/\partial x^\alpha$ 
 оценивается сверху как
$$
\Bigl|\frac{\partial^\alpha f}{\partial x^\alpha}\Bigr|\le C \rho^\alpha \alpha!
$$
Рассмотрим последовательность $M_j>0$, пусть $M_0=1$,
 и последовательность логарифмически выпукла, $M_k^2\le M_{k-1} M_{k+1}$. 
Класс Данжуа--Карлемана $C^M$ состоит из функций, для которых предыдущая оценка заменена на
$$
\Bigl|\frac{\partial^\alpha f}{\partial x^\alpha}\Bigr|\le C \rho^\alpha \alpha! \,
 M_{|\alpha|}.
$$
Класс $C^M$ называется {\it квазианалитическим}, если любая функция $f\in C^M$ однозначно определяется
своими тейлоровскими коэффициентами в точке. Оказывается, что
при $M_k=(\ln k)^k$ мы получаем квазианалитический класс. Критерий квазианалитичности:
$$
\sum_{k=1}^\infty \Bigl(\inf_{j\ge k}(j! M_j)^{1/j}\Bigr)^{-1}=\infty.
$$
Оказывается, что диффеоморфизмы квазианалитического класса $C^M$ образуют топологическую группу.

\sm

 {\bf\punct Индуктивные (прямые) пределы групп.} Рассмотрим цепочку вложений топологических групп
$$
\dots\longrightarrow G_{k} \longrightarrow G_{k+1}\longrightarrow G_{k+2} \longrightarrow\dots
,$$
Рассмотрим объединение
$$
\cG=\lim\limits_{\longrightarrow} G_k
$$
 всех групп этой цепочки (мы отождествляем группу $G_k$ с подгруппой в
$G_{k+1}$ для каждого $k$). Снабдим это объединение топологией, положив, что множество в
$\lim\limits_{\longrightarrow} G_k$ открыто, если его пересечение со всеми группами
$G_k$ открыто. Полученная топологическая группа называется 
{\it индуктивным} (или {\it прямым}) {\it пределом} групп $G_k$.

Последовательность $g_j$  в прямом пределе сходится, если все $g_j$ содержатся в одной
достаточно большой группе $G_k$ и последовательность сходится в этой группе.

\sm

{\sc Пример.}  Рассмотрим естественное вложение $\U(n)$ в $\U(n+1)$, матрица
$g\in \U(n)$ дописывается как $\begin{pmatrix}g&0\\0&1\end{pmatrix}$. Индуктивный
предел -- это группа всех бесконечных унитарных матриц $g$, таких, что
$g-1$ имеет лишь ненулевое число матричных элементов. 
\hfill $\lozenge$

\sm

{\sc Пример.}  Можно придумать много разных индуктивных пределов унитарных групп.
 Например, можно взять вложение $\U(2^n)$ в $\U(2^{n+1})$, заданное как 
$\begin{pmatrix}g&0\\0&g\end{pmatrix}$. Кажется, получаемые таким образом группы не 
очень интересны. \hfill $\lozenge$

\sm

См. \cite{Bou1}.

\sm
 
{\bf\punct  Проективные (обратные) пределы групп.%
	\label{ss:group-inverse-1}}
Рассмотрим теперь цепочку сюръективных гомоморфизмов
$$
\dots\xleftarrow{\psi_k} G_{k} \xleftarrow{\psi_{k+1}}
 G_{k+1}\xleftarrow{\psi_{k+2}} G_{k+2}
 \xleftarrow{\psi_{k+3}}\dots
$$
Рассмотрим множество 
$$\frG=\lim\limits_{\longleftarrow} G_k$$  
всевозможных последовательностей $g_k\in G_k$, таких, что
для всех $k$ выполнено $\psi_k(g_k)=g_{k-1}$. Это множество снабжено естественной
групповой
операцией: произведение последовательностей $\{h_k\}$ и $\{g_k\}$ равно
$\{h_kg_k\}$.

 Группа $\frG$ сюръективно отображается на любую группу $G_j$, 
 а именно последовательности $\{g_k\}$ ставится в соответствие ее элемент   $g_j$. 
 Обозначим это отображение через $\pi_j$.
 Очевидно $\psi_j\pi_j=\pi_{j-1}$.  Для любого открытого множества $\cV\subset G_j$ мы потребуем, 
 чтобы его прообраз $\pi^{-1}_j(\cV)$ был открыт в $\frG$, и чтобы всевозможные такие множества образовывали 
 базу топологии в $\frG$. Тогда $\frG$ становится топологической
 группой.

 \sm
 
 {\sc Задача.} a) Если все $G_k$ компактны, то и их проективный предел компактен.
 
 \sm
 
 b) Пусть $K$, $L$ -- компактные группы, $\psi:K\to L$ -- сюръективный гомоморфизм.
 Тогда образ меры Хаара на $K$ есть мера Хаара на $L$. \hfill $\lozenge$
 
 \sm
 
 Пусть все группы $G_k$ компактны. Нормируем меры Хаара, так, чтобы все меры были вероятностными.
 В силу последнего упражнения и теоремы Колмогорова
  определен проективный предел мер на
 $\frG$, он и будет мерой Хаара на $\frG$. 
 
\sm

Разумно расширить определение проективного предела. Рассмотрим счетное частично 
упорядоченное множество $A$, пусть для любых двух элементов $\alpha$, $\beta\in A$
существует элемент $\gamma$ такой, что $\gamma\ge \alpha$, $\gamma\ge\beta$. 
Рассмотрим семейство топологических групп $G_\alpha$, нумеруемых элементами
$\alpha\in A$. Пусть для любой пары $\alpha<\beta$ задан сюръективный гомоморфизм
$\psi_{\beta\alpha}:G_\beta\to G_\alpha$, причем при $\alpha<\beta<\gamma$
мы имеем $\psi_{\beta\alpha}\psi_{\gamma\beta}=\psi_{\gamma\alpha}$.
Элементом проективного предела $\frG$ являются всевозможные наборы 
$\{g_\alpha\in G_\alpha\}$, где $\alpha$ пробегает $A$, такие, что
$\psi_{\beta\alpha}g_\beta=g_\alpha$ для всех $\alpha<\beta$. Групповая структура
и топология вводятся так же, как и выше. Проективный предел компактных групп,
по-прежнему, компактен.

\sm

См. \cite[\S III.7]{Bou1}, \cite{RZ}.

\sm
 
 {\bf\punct Примеры проективных пределов.%
 	\label{ss:group-inverse-2}}

{\sc Задача.} Пусть $p$- простое число. Рассмотрим следующую цепочку отображений циклических групп
$$
\dots \longleftarrow \Z_{p^k} \longleftarrow \Z_{p^{k+1}} \longleftarrow \dots
$$
($Z_{p^k}=\Z_{p^{k+1}}/p^k \Z_{p^{k+1}}$). Опишите проективный предел.
Ответ здесь -- целые $p$-адические  числа.  \hfill $\lozenge$

\sm

{\sc Задача$^*$.}  Рассмотрим  такое частичное упорядочение на множестве натуральных чисел:
$\beta\ge\alpha$, если $\beta$ делится на $\alpha$. Каждому числу $\alpha$ ставится
в соответствие циклическая группа $\Z_\alpha$. Паре $\beta\ge\alpha$
ставится в соответствие ненулевой гомоморфизм $\Z_\beta\to \Z_\alpha$. Опишите проективный
предел. \hfill $\lozenge$

\sm

{\it Пополнения дискретных групп.} Пусть $\Gamma$ -- счетная дискретная
 группа. Рассмотрим  в ней убывающее  семейство нормальных подгрупп
 $$
 \Gamma= \Gamma_0\supset  \Gamma_1\supset  \Gamma_2\supset\dots
 $$
 Положим, что их пересечение равно 1. Рассмотрим цепочку гомоморфизмов
 $$
\dots \longleftarrow  \Gamma/ \Gamma_k  \longleftarrow  \Gamma/ \Gamma_{k+1} \longleftarrow \dots
 $$
 На каждом шагу этой цепочки мы имеем дискретную группу, но проективный предел
 $\ov \Gamma$ -- уже
 нетривиальная топологическая группа. Обозначим через $\pi_k$ проекцию $\ov\Gamma$
 на $\Gamma/\Gamma_k $.
 
 \sm
 
 {\sc Задача.} 
  Группа $\Gamma$ является плотной подгруппой в $\ov\Gamma$. Прообраз единицы при отображении
  $\pi_k$ является открыто-замкнутой подгруппой в $\ov\Gamma$, причем эта подгруппа
  совпадает с замыканием $\Gamma_k$. Группа $\ov\Gamma$ полна.
  \hfill $\lozenge$
  
  \sm
  
  Здесь возникает большой и разнообразный зоопарк конструкций. Отметим, что  в упражнении чуть выше
  мы брали $\Gamma=\Z$, $\Gamma_k=p^k\Z$, и ответ уже получается нетривиальным.
  
  Естественно эту конструкцию сразу обобщить. Рассмотрим  счетную дискретную группу
  $\Gamma$ и семейство нормальных подгрупп $\Gamma_\alpha$, где $\alpha$
  пробегает некоторое множество $A$. На множестве $A$ тогда возникает
  частичный порядок, $\alpha\le\beta$, если $\Gamma_\alpha\supset \Gamma_\beta$.
   Мы требуем, чтобы любое пересечение 
  $\Gamma_\alpha\cap \Gamma_\beta$ cодержало некоторую подгруппу $\Gamma_\delta$,
  и чтобы пересечение всех подгрупп $G_\alpha$ сводилось к единице. При
  этих условиях мы можем построить проективный предел
  $$
  \ov\Gamma=\lim\limits_{\longleftarrow} \Gamma/\Gamma_\alpha
  .$$

{\it Проконечные пополнения.} 
Пусть $\Gamma$ -- счетная группа, у которой пересечение
подгрупп конечного индекса тривиально. Тогда в качестве семейства $\Gamma_\alpha$
мы можем взять все подгруппы конечного индекса.

В этом случае проективный предел является компактной группой.

\sm

{\sc Замечание.} В свое время долго стояла т.н. конгруэнц-проблема об описании
проконечного пополнения групп $\GL(n,\Z)$ при $n\ge 3$. Ответ --  произведение 
групп $\GL(n)$ над целыми $p$-адическими числами по всем простым
$p$, см. \cite{BMS}. Для $\GL(2,\Z)$ пополнение значительно больше, чем это произведение,
см. \cite{RZ}.
\hfill$\lozenge$

\sm

 Группа называется {\it проконечной}, если она является проективным пределом
 конечных групп.

\sm

{\it Пронильпотентные пополнения.} Пусть $G$ -- группа, $H$ -- подгруппа.
Коммутантом $\{G,H\}$ называется подгруппа, порожденная всеми коммутаторами
$\{g,h\}=h^{-1}g^{-1}hg$, где $g$ пробегает $G$, а $h$ пробегает $H$.
Для счетной дискретной группы $\Gamma$ рассмотрим цепочку подгрупп
({\it нижний центральный ряд})
\begin{equation}
\Gamma_0=\Gamma,\Gamma_1=\{\Gamma,\Gamma_0\},\quad \Gamma_2=\{\Gamma,\Gamma_1\},\quad
\Gamma_3=\{\Gamma,\Gamma_2\},\quad \dots
\label{eq:central-series}
\end{equation}
Группа называется {\it нильпотентной}, если некоторая подгруппа
$\Gamma_j$ тривиальна. Группа называется {\it виртуально  нильпотентной}, если
пересечение $\cap_j\Gamma_j$ тривиально. Примеры виртуально нильпотентных групп:
свободная группа, группа чистых кос. Для таких групп определено пронильпотентное пополнение
$\lim\limits_{\longleftarrow} \Gamma/\Gamma_\alpha$.

\sm

{\it Группы Галуа и топология Крулля.} Пусть $F$ -- счетное (или конечное) поле.
Пусть $F$ - {\it алгебраическое расширение} $F$ (т.е., любой элемент из $K$ есть корень
многочлена с коэффициентами из $F$).   Говорят, что 
$K$ - {\it нормальное расширение} $F$, если всякий многочлен
с коэффициентами из $F$, имеющий корень в $K$, разлагается на линейные множители в поле
$K$. Пусть $K$ -- нормальное расширение $F$. {\it Группой Галуа}
 $G_{K/F}$ называется группа автоморфизмов поля $K$, сохраняющих
$F$. Пусть $L\subset K$ -- нормальное расширение поля $F$.  Тогда элементы
$G_{K/F}$ переводят $L$ в себя, а потому мы  имеем гомоморфизм $G_{K/F}\to G_{L/F}$.

 Рассмотрим счетное поле $K$ (например, $\Q$).
Рассмотрим его алгебраическое замыкание $\ov K$. Рассмотрим группу Галуа $\Gal_{\ov K/K}$
 поля $K$, т.е. 
группу автоморфизмов поля $\ov K$, тривиальных на $K$. 
Очевидно, что группа Галуа $\Gal_{\ov K/K}$ переводит любое нормальное расширение 
$L$ в себя. Тем самым, мы получаем гомоморфизм
$\Gal_{\ov K/ K}\to \Gal_{L/K}$. Рассмотрим все конечные нормальные расширения $L_\alpha$ поля
$F$, содержащиеся в $K$. Если $L_\alpha\subset L_\beta$, то мы имеем отображение
$$
G_{L_\beta/F}\to G_{L_\alpha/F}
,$$  
а следовательно, определен обратный предел $\lim\limits_{\longleftarrow} G_{L_\alpha/F}$.

\sm

{\sc Задача.} Покажите, что естественное отображение 
$G_{K/F}\to \lim\limits_{\longleftarrow} G_{L_\alpha/F}$ является биекцией. 
Таким образом, группа Галуа является компактной группой (а построенная топология  на ней
называется <<топологией Крулля>>). \hfill $\lozenge$

\sm

Оказывается также, что любая проконечная группа является группой Галуа.

\sm 

См. \cite{Ser2}, \cite{RZ}, \cite{Bab}.

\sm 

{\bf\punct Автоморфизмы и иерархоморфизмы деревьев.\label{ss:hier}}

{\it Автоморфизмы.} Дерево -- это связный граф без циклов.
{\it Дерево Брюа--Титса} $T_n$ -- это бесконечное дерево, у которого из каждой вершины исходит
$n+1$ ребро. Автоморфизм дерева -- это биекция множества вершин на себя,
переводящая соседние вершины в соседние. Топология на группе $\Aut(T_n)$ автоморфизмов дерева вводится из
условия: стабилизаторы конечных поддеревьев являются открытыми подгруппами. 
Группы $\Aut(T_n)$ -- довольно простой и приятный пример локально компактных групп.
Интересно, что эти группы с точки зрения теории представлений оказываются
родственниками $\SL(2,\R)$, при этом сами деревья Брюа--Титса имитируют
плоскость Лобачевского (П.Картье). Отметим, что при простом $n=p$ группа
$\Aut(T_p)$ содержит $p$-адическую $\SL(2,\Q_p)$.

\sm

См. \cite{Ser-trees}, \cite{Olsh-tree}, \cite{FTT}, \cite[\S10.4]{Ner2}.

\sm

{\it Иерархоморфизмы.} Понятие {\it границы дерева} $\Abs(J_n)$ достаточно очевидно. Дадим формальное определение.
Назовем {\it путем} в дереве
 последовательность вершин $a_1$, $a_2$, $a_3$, \dots, такую, что
$a_j$ и $a_{j+1}$
-- разные концы одного ребра. Два пути $a_j$, $b_j$ эквивалентны,
если существует $k\in\Z$ такое, что, начиная с некоторого номера,
 выполнено $a_{j+k}$. Точкой границы назовем класс эквивалентных путей.
 
 Разрежем дерево посередине какого-нибудь  ребра.
 Два получившихся поддерева назовем {\it ветками}.
 Часть границы, примыкающую к ветке, мы будем считать открыто-замкнутым
 множеством. Тогда граница дерева становится вполне несвязным компактом.

 {\it Веером} назовем конечный набор веток с попарно не пересекающимися
 границами, причем объединение этих границ покрывает всю  границу
 дерева. Пусть даны два веера с одним и тем же количеством ветвей.
 Отобразим один веер на другой с сохранением структуры дерева
 в каждой ветви.
 Это отображение индуцирует гомеоморфизм границы дерева. Такие гомеоморфизмы
 мы будем называть
  {\it иерархоморфизмами}. Легко видеть, что иерархоморфизмы дерева Брюа--Титса образуют группу,
  обозначим ее через $\Hier(T_n)$.

 Эта группа локально-компактна. А именно, на $\Aut(T_n)$
 мы вводим ее естественную топологию и полагаем, что
 $\Aut(T_n)$ открыта в $\Hier(T_n)$. Классы смежности
 $\Hier(T_n)/\Aut(T_n)$ образуют счетное множество.
 Каждый класс смежности отождествляется с 
 $\Aut(T_n)$, а на $\Hier(T_n)$ мы вводим топологию несвязного объединения.

\sm

Хотя эта группа локально компактна, с точки зрения конструкций представлений
она похожа на группу диффеоморфизмов окружности. Кстати, при простом $n=p$
она содержит группу локально-аналитических диффеоморфизмов $p$-адической проективной
прямой.

\sm

{\it Группы Р.Томпсона%
	\footnote{Группа Ричарда Томпсона. Не путать с группой Джона Томпсона}.}
 Уложим дерево Брюа-- Титса $T_2$ на плоскость. Тогда на границе
возникает циклический порядок. Группа Томпсона определяется как группа всех иерархоморфизмов,
сохраняющих этот порядок. Это -- странная дискретная группа, допускающая несколько
внешне непохожих друг на друга реализаций.

\sm

1. {Кусочно-линейные отображения.}
 Рассмотрим окружность $\R/\Z$. Рассмотрим группу $\mathrm{PL}_2$ кусочно линейных монотонно возрастающих
 гомеоморфизмов этой окружности в себя, удовлетворяющих следующим двум свойствам:
 
 \sm 
 
  a) все точки излома двоично-рациональны (т.е., знаменатели координат имеют вид $2^k$);
  
  \sm
  
  b) угловые коэффициенты линейных участков графика имеют вид $2^m$, где $m\in\Z$.

\sm

{\sc Задача.} Покажите, что  группа $\mathrm{PL}_2$ совпадает с группой Томпсона.
Во-первых, надо нарисовать дерево. Его вершины нумеруются отрезками
вида $[a2^{-k}, (a+1)2^{-k}$, где $a<2^k$ -- целые числа. Две вершины соединяются ребром,
если один отрезок является половинкой другого. В итоге получается ветка
дерева $T_2$. Очевидно, что наша группа действует иерархоморфизмами
этой ветки. Убедитесь, что для дерева $T_2$ ветка и все дерево иерархоморфны.
\hfill $\lozenge$

\sm

2. Кусочно дробно-линейные отображения. Рассмотрим
вещественную проективную прямую $\R\P^1=\R\cup \infty$.
Рассмотрим группу $\mathrm{PPSL}(2)$, состоящую из кусочно
дробно-линейных
$$
x\mapsto (ax+b)/(cx+d),\qquad \text{где $a$, $b$, $c$, $d\in \Z$, $ad-bc=1$}  
$$
отображений,
причем в точках склейки они имеют гладкость $C^1$. Отметим, что
точки склейки автоматически являются рациональными. 
На самом деле, тут присутствует картинка из учебников комплексного анализа.
\begin{figure}
	\epsfbox{haar.2}\quad \epsfbox{haar.5}
\caption{На левой половине рисунка изображена стандартная картинка на плоскости Лобачевского,
	инвариантная относительно $\PSL(2,\Z)$. На правой половине --
	две копии плоскости Лобачевского и с идеальным $n$-угольником на каждой копии, стороны 
	$n$-угольников выбраны из дуг левой картинки. Мы берем отображение сегментов, сохраняющее
	циклический порядок, на каждом сегменте мы должны перевести картинку из дуг в картинку из дуг,
	и этим отображение однозначно определено.%
	\label{fig:ideal}}
\end{figure}
Мы рассматриваем верхнюю полуплоскость $\Im z>0$ комплексной плоскости
как модель плоскости Лобачевского, прямые соответствуют 
окружностям, перпендикулярным вещественной оси. 
Отражения относительно <<прямых>> - это инверсии.
Мы рисуем идеальный треугольник, ограниченный лучами
$\Re z=0$, $\Re z=1$ и полуокружностью $|z-1/2|=1/2$.
Далее мы отражаем его относительно сторон,
 новые треугольники
отражаем относительно их сторон, и т.д.
Отметим, что эти отражения бывают двух типов. Отражение относительно вертикальной прямой
$\Re z=n>0$ добавляет к картинке 
вершину $n+1$ (а если брать прямую $\Re z=n<0$ - то $n+1$). Отражение относительно
дуги окружности с концами $\frac ab$, $\frac cd$ (дроби предполагаются несократимыми)
 добавляет вершину 
$\frac{a+b}{c+d} $.

По множеству треугольников очевидным образом строится дерево. Вершины соответствуют
треyгольникам, две вершины соединены ребром, если треугольники соседние.
Группа $\mathrm{PPSL}(2)$ действует иерархоморфизмами этой картинки, см. \ref{fig:ideal}.

\sm

{\it Функция Минковского $?^{-1}(x)$.}
Интересно, что группы 
$\mathrm{PL}_2$ и $\mathrm{PPSL}(2)$ сопряжены в группе гомеоморфизмов окружности.
Соответствующее отображение $?(x)$ осуществляет биекцию между множеством двоично рациональных
точек  окружности $\R/\Z$ и множеством рациональных точек проективной прямой
$\R\cup \infty$. Правило отождествления ясно  из структуры дерева,
$$
?^{-1}\Bigl(\frac k{2^n}\Bigr)=\frac ab,\quad ?^{-1}\Bigl(\frac {k+1}{2^n}\Bigr)=\frac cd
\quad
\Rightarrow\quad ?^{-1}\Bigl( \frac{2k+1}{2^{n+1}}\Bigr)=
\frac{a+b}{c+d}
$$
(дроби опять предполагаются несократимыми).
Отметим, что эта странная монотонная функция имеет производную, равную 
нулю или бесконечности  всюду, так же как и обратная к ней функция $?(y)$. При этом $?^{-1}(\cdot)$
является биекцией мкжду множеством рациональных чисел и множеством квадратичных иррациональностей.

\sm

{\it Почти ортогональные действия групп иерархоморфизмов.}
Фиксируем  вершину $v$ дерева $T_n$  и возьмем единственную
 вероятностную меру $d\omega(x)$ на границе, инвариантную относительно
стабилизатора точки $v$. А именно, к вершине $v$ примыкает $n+1$
ветка, тем самым граница разбивается на $n+1$ подмножество
меры $1/(n+1)$, каждое из этих подмножеств разбивается на $n$ 
одинаковых подмножеств, и т.д. Введем на границе расстояние,
для этого соединим путем $w_j$, $j\in \Z$, две точки границы $x$, $y$.
А именно, для дюбого $j$ вершины 
$w_j$ и $w_{j+1}$ должны быть соседними, $w_j$ попарно различны,
$\lim_{j\to-\infty}=x$, $\lim_{j\to+\infty}=y$. Мы положим, что
$\rho(x,y)$ равно минимальному расстоянию от $w_j$ до
начальной вершины $v$.  Фиксируем $s\in (0,n)$ и введем на пространстве
вещественных
кусочно постоянных функций на $\Abs(T)$ скалярное произведение
по формуле
$$
\la f,g\ra_s:=
\int_{\Abs(T_p)} \int_{\Abs(T_p)}
\rho(x,y)^{-n+s}
f(x) g(y)\, d\omega(x)\, d\omega(y).
$$
Это скалярное произведение оказывается положительно определенным.
Обозначим пополнение по этому скалярному прозведению через $V_s$.
Определим действие группы $\Hier(T_p)$ в $V_s$ по формуле
$$
T_s(q) f(x) = f(q(x)) \,q'(x)^{(n-s)/2}
,$$
где $q'(x)$ -- производная Радона--Никодима преобразования
$q(x)$ (кстати, как ее определить в геометрических терминах?).
Как и в случае группы диффеоморфизмов окружности, мы получаем
$T_s(q)\in \GLO(V_s)$, более того, в нашем случае оператор 
$T_s(q)$ отличается от ортогонального оператора на поправку
конечного ранга.

\sm

См. \cite{Ner-derevo}, \cite{Ner-trees}, \cite{GS}, \cite{CFP}, \cite{Kap},
\cite{Imb}.

\sm

{\bf \punct Формальные нильпотентные группы%
	\label{ss:mal}} (для тех, кто знаком с группами Ли).
Рассмотрим бесконечномерную $\Z_+$-градуированную алгебру Ли
 $\frg=\oplus_{\alpha=1}^\infty\frg_\alpha$ с конечномерными однородными подпространствами
$\frg_\alpha$. Обозначим через $\frg$ пространство формальных рядов
$S=\sum_{\alpha=1}^\infty X_\alpha$, где $X_\alpha\in \frg_\alpha$.

Рассмотрим универсальную  обертывающую алгебру $U(\frg)$ алгебры Ли $\frg$
(см. \cite[гл. 2]{Dix-U}, \cite[\S19]{Zhe2}). Она является
градуированной ассоциативной алгеброй,  $U(\frg)=\oplus_\alpha U(\frg)_\alpha$
(с конечномерными однородными слагаемыми). Пополним ее, положив, что $\ov U(\frg)$ --
это алгебра формальных рядов вида $\sum_{\alpha=1}^\infty z_\alpha$,
где $z_\alpha\in U(\frg)_\alpha$. 

Мы имеем корректно определенное отображение
$X\mapsto e^X$ из $\ov\frg$ в $\ov U(\frg)$. Обозначим его образ
через  $\frG$. Применяя формулу Кэмпбелла--Хаусдорфа (\cite[\S21]{Zhe2}, \cite{Reu}), мы обнаруживаем,
что  $\frG$ замкнута относительно умножения. Таким образом мы получили группу.

Формула Кэмпбелла-- Хаусдорфа замечательна степенью своей неудобности.
Группа $\frG$ может быть описана следующим способом, который в конкретных случаях бывает более
приятным.

Рассмотрим тензорное произведение $\ov U(\frg)\otimes \ov U(\frg)$ и
коумножение (см. \cite[\S 2.7]{Dix-U}).
$$
\Delta \ov U(\frg)\,\to\, \ov U(\frg)\otimes \ov U(\frg).
$$
Напомним, что это гомоморфизм алгебр, который определен на элементах алгебры $\frg$ (которые
отождествлены с образующими обертывающей
алгебры) как
$$
\Delta(X)= X\otimes 1+1\otimes X
.$$ 

Оказывается, что группа $\frG$ состоит из элементов алгебры $\ov U(\frg)$,
удовлетворяющих условию
\begin{equation}
\Delta(z)=z\otimes z
\label{eq:z-z}
.
\end{equation}
Содержательная ситуация возникает уже для свободных алгебр Ли.

\sm

{\it Группы Ли, соответствующие свободным алгебрам Ли.}  Рассмотрим свободную алгебру Ли 
с $n$ образующими. Опишем соответствующую группу Ли.

 Известно, что  обертывающая алгебра свободной алгебры Ли
-- это свободная ассоциативная алгебра  $\mathrm{Free}_n$ с $n$ образующими,
т.е., алгебра $\cF_n$ некоммутативных многочленов от $n$ букв
$a_1$, \dots, $a_n$. Через  $\ov\cF_n$ мы обозначим алгебру формальных
рядов от тех же образующих. Коумножение 
$\Delta:\mathrm{Free}_n\to \mathrm{Free}_n\otimes \mathrm{Free}_n$ достаточно определить на образующих,
$$
\Delta(a_j)=a_j\otimes 1+1\otimes a_j
.$$
Пусть $u$, $v$ -- конечные слова, составленные из букв $a_1$, \dots, $a_n$.
Мы скажем, что слово $w$ -- {\it тасовка} $u$ и $v$, если буквы слова $w$
можно покрасить в два цвета, красный и синий, так, что слово,
составленное из синих букв, совпадет с $u$, а слово, составленное из красных
букв, совпадет с $v$. Например, тасовками слов $ab$ и $cd$ являются
$abcd$, $acbd$, $acdb$, $cdab$, $cadb$, $cabd$.

Рассмотрим произведение $w=a_{i_1}\dots a_{i_k}$. Тогда
$$
\Delta(w)=\sum_{\text{$u$, $v$, такие, что $w$ -- тасовка $u$ и $v$.}} u\otimes v
.$$
Пусть $\sum_w c_w w$ -- элемент группы Ли свободной алгебры $\mathrm{Free}_n$.
Тогда уравнение (\ref{eq:z-z}) переписывается как следующая система уравнений на  коэффициенты
$c_w$:
\begin{equation}
c_u c_v=\sum_{\text{по всем тасовкам $w$ слов $u$, $v$}} c_w.
\label{eq:ccc}
\end{equation}
Пространство решений этой системы (то есть наша группа) допускает параметризацию.
А именно, коэффициенты $c_w$ перед словами Линдона%
\footnote{{\it Слово Линдона}
	-- это слово $w$, у которого для любого его представления 
	$w=uv$ выполнено  $w<v$ в смысле лексикографического порядка..}
 $w$ можно выбирать
произвольно, а все остальные коэффициенты уже однозначно (один за другим)
находятся из уравнений (\ref{eq:ccc}). См. \cite{Reu}.

\sm

{\it Мальцевские пополнения.}
Пусть $\Gamma$ дискретная нильпотентная группа без кручения. Согласно
А.И.Мальцеву (см. \cite{Mal}, \cite[17.3]{KarM}), в $\Gamma$ можно выбрать конечный набор 
образующих $g_1$,  \dots, $g_r$ так, что 

\sm 

1) для любого $p$, подгруппа $\Gamma^p$, порожденная $g_p$, $g_{p+1}$, $g_r$,
является нормальной, причем $\Gamma_n$ лежит в центре $\Gamma$, 
группа
$\Gamma_{r-1}/\Gamma_{r}$ лежит в центре $\Gamma/\Gamma_r$, и.т.д.
($\Gamma_{p}/\Gamma_{p+1}$ лежит в центре $\Gamma/\Gamma_{p+1}$);

\sm 

2) любой элемент группы однозначно представим в виде
$g^\bfn:=g_1^{m_1}\dots g_n^{n_r}$; тем самым $\Gamma$ отождествляется 
как множество с $\Z^r$;

\sm

3) в формуле для умножения 
\begin{equation}
g^\bfm g^\bfk=g^{\bfk(\bfm,\bfn)}
\label{eq:ggg}
\end{equation}
зависимость $\bfk(\bfm,\bfn)$ является полиномиальной.

\sm

Теперь в формуле (\ref{eq:ggg}) мы можем положить степени $m_j$, $n_j$
призвольными вешественными числами, и что дает нам умножение на 
на $\R^r$. Получается некоторая нильпотентная группа Ли $G=G(\Gamma)$,
а также ее алгебра Ли $\frg(\Gamma)$.
При этом однородное пространство $G/\Gamma$ компактно,
а сама конструкция функториальна, т.е., любой гомомрфизм
$\Gamma_1\to \Gamma_2$ продолжается до гомоморфизма
$G(\Gamma_1)\to G(\Gamma_2)$.

Пусть теперь группа виртуально нильпотентна, пусть
$\Gamma_p$ -- ее нижний центральный ряд, см. (\ref{eq:central-series}).
Тогда есть цепочка сюръективных гомомoрфизмов 
$$
\dots\longleftarrow \Gamma/\Gamma_l \longleftarrow \Gamma/\Gamma_{l+1}\longleftarrow\dots
,$$
а следовательно, и цепочки  гомомoрфизмов
\begin{align*}
\dots\longleftarrow G(\Gamma/\Gamma_l) \longleftarrow G(\Gamma/\Gamma_{l+1})\longleftarrow\dots
\\
\dots\longleftarrow \frg(\Gamma/\Gamma_l) \longleftarrow \frg(\Gamma/\Gamma_{l+1})\longleftarrow\dots
\end{align*}
У этих цепочек есть проективные пределы, бесконечномерная алгебра Ли $\frg(\Gamma)$
и бесконечномерная группа $G(\Gamma)$. 

Мальцевское пополнение свободной группы обсуждалось выше. Кажется, сейчас самый интересный
объект возникающего зоопарка связан с группами чистых кос. Соответствующая алгебра Ли
была построена Т.~Коно \cite{Koh}, это алгебра Ли с образующими
$x_{ij}$, где $1\le i,j\le n$, причем $i\ne j$, а $x_{ij}=x_{ji}$
(то есть, всего образующих $n(n-1)/2$). Соотношения имеют вид
\begin{align*}
[x_{ij},x_{ik}+x_{jk}]=0\quad \text{или $[x_{ij},x_{ij}+x_{ik}+x_{jk}]=0$}
\\
[x_{ij},x_{kl}]=0,\quad\text{если $i$, $j$, $k$, $l$ попарно различны}.
\end{align*}
Алгебра Коно имеет неожиданно много интересных действий.
Наример, она действует линейными операторами в 
в тензорных произведениях $n$ представлений полупростой группы
(конструкция В.~Г.~Книжника--А.~Б.~Замолодчикова), в пространствах нулевого
веса представлений алгебры Ли $\frs\frl(n)$; она действует векторными полями на произведениях
двумерных сфер (А.~А.~Клячко), коприсоединенных орбитах $\GL(n,\C)$ (см. ссылки в \cite{Ner-Mal}).

\sm

См. \cite{Reu}, \cite{Koh}, \cite{Ner-Mal}, \cite{Bab}.

\sm

{\bf\punct Свертки мер.%
\label{ss:svertki}} Напомним определение {\it свертки мер}.
Пусть $G$ - группа, $\mu$, $\nu$ -- меры (положительные или знаконеопределенные) с компактным
носителем. Свертка мер $\mu*\nu$ определяется из условия: для любой непрерывной функции
$f$ выполнено
$$
\int_G f(g)\,d \mu*\nu(g)=\int_G\int_G f(g_1 g_2)\,d\mu(g_1)\, d\nu(g_2)
$$
(см. \cite[10.1]{Kir}, \cite[\S19-20]{HR1}). Эта операция линейна по каждому аргументу и ассоциативна.
На локально компактной группе свертка раздельно непрерывна относительно слабой сходимости мер.
Отметим, что при свертке дельта-мер мы имеем $\delta_a*\delta_b=\delta_{a+b}$.
Далее эта операция продолжается по линейности и непрерывности. В случае 
конечной группы мы получаем групповую алгебру, в случае групп $\R$, $\Z$
мы получаем обычную свертку, известную из курса анализа.

Если $\rho$ -- унитарное представление группы $G$, то для мер с компактным носителем мы можем
определить операторы 
$$
\rho(\mu):=\int \rho(g)\,d\mu(g).
$$
Это задает представление нашей сверточной алгебры,
$$
\rho(\mu_1+\mu_2)=\rho(\mu_1)+\rho(\mu_2),\qquad \rho(\mu_1*\mu_2)=\rho(\mu_1)\rho(\mu_2).
$$

Полугруппа вероятностных мер на $\R$ с операцией свертки -- стандартный объект
теории вероятностей и случайных процессов (безгранично делимые
распределения, случайные процессы с независимыми
приращениями, случайные блуждания, разложения вероятностных
законов и т.д). Полугруппы вероятностных мер на локально компактных
группах много изучались, мы ограничимся здесь лишь ссылками
\cite{Hey}, \cite{Breu}, \cite{App}, \cite{Boug}. 
Мы же сейчас рассмотрим некоторые подполугруппы сверточных полугрупп.

\sm

{\bf\punct Гипергруппы.
\label{ss:hypergroups}}
 {\it Пример.}. Рассмотрим
плоскость $\R^2$ и семейство концентрических окружностей с центром в нуле.
Обозначим через $p_a$ равномерную вероятностную меру, сосредоточенную на окружности
радиуса $a$.  Свертка  $p_a*p_b$ не имеет вида $p_c$, она  размазана по кольцу $|a-b|\le r \le a+b$
и является радиально симметричной. Поэтому она должна представляться как комбинация мер
$p_c$.

\sm

{\sc Задача.} Покажите, что 
$$
p_a*p_b=\int_{|a-b|}^{a+b} h_{a,b}(c) p_c\, dc,
$$
где 
$$
\qquad
h_{a,b}(c)=\frac 2\pi \frac {c}{\bigl[(a+b+c)(a+b-c)(a+c-b)(b+c-a)\bigr]^{1/2}}.
\qquad\qquad\lozenge
$$
Таким образом мы получаем операцию, которая из пары $\delta$-мер на полупрямой
изготовляет некоторую меру на полупрямой. 

\sm

{\it Гипергруппы.} Рассмотрим локально компактное метризуемое пространство $V$. 
Обозначим через $\cM(V)$ множество вероятностных мер на $V$ с компактным носителем.
Возьмем отображение  $V\times V\to \cM(V)$, точки $a$ пространства $V$  
мы отождествляем с дельта-мерами $\delta_a$. Отображение мы обозначим через
$$
\delta_a\star \delta_b\mapsto \mu_{a,b}.
$$
Естественно потребовать от этой операции ассоциативности,
$(\mu_a\star \mu_b)\star \mu_c=\mu_a\star (\mu_b\star \mu_c)$, 
т.е.,
$$
\int_V \mu_{v,c}\,d\mu_{a,b}(v)=\int_V \mu_{a,v}\,d\mu_{b,c}(v) 
$$
Кроме того, чтобы проимитировать другие аксиомы группы,
обычно требуют наличие единицы $\delta_e$, такой, чтобы
$\delta_e\star \delta_a=\delta_a\star \delta_e=\delta_a$
и инволюции $\dagger:V\to V$ такой, что 
$(\delta_a\star\delta_b)^\dagger= \delta_b^\dagger\star\delta_a^\dagger$
(это аналог взятия обратного элемента).
Наконец, мы требуем, чтобы носитель $\delta_a\star\delta_b$ 
содержал $e$ тогда и только тогда, когда $\delta_a=\delta_b^\dagger$.
Нужно еще наложить какие-то условия непрерывности, которые могут 
варьироваться (в литературе варьируется и предшествующая часть определения),
тут важны не детали определения, а то, что интересные и много изучавшиеся гипергруппы в самом деле
существуют.

Теперь мы приведем несколько массовых примеров гипергрупп.

\sm

{\it Гипергруппы классов сопряженности.}
Пусть $G$ -- локально компактная группа, $K$ -- ее компактная подгруппа. Рассмотрим пространство
классов сопряженности $G//K$ группы $G$ относительно подгруппы $K$. Каждый класс сопряженности
является орбитой компактной группы, а поэтому на нем есть единственная
вероятностная мера, инвариантная относительно $K$. Таким образом, каждому классу сопряженности $\gamma$ ставится в соответствие
дельта-образная мера $\delta_\gamma$, сосредоточенная на классе $\gamma$.  Очевидно, 
свертка мер вида $\delta_\gamma$ определяет 
структуру гипергруппы на $G//K$.

\sm

{\it Двойные классы смежности.} Пусть $G$ -- локально компактная  группа, $K$ -- ее компактная подгруппа. 
Рассмотрим пространство двойных классов смежности $K\setminus G/K$, т.е.,
фактор-множество группы $G$ по отношению эквивалентности
$$g\sim k_1 g k_2,\qquad k_1,k_2\in K.$$
Каждый двойной класс смежности является орбитой компактной группы $K\times K$,
а поэтому на нем есть единственная инвариантная вероятностная мера. Мы рассматриваем 
$\delta$-образные инвариантные меры, сосредоточенные на двойных классах смежности
и свертки таких мер.

\sm

{\it Орбиты группы автоморфизмов.} Пусть компактная группа $L$ действует автоморфизмами
локально компактной  группы $G$. Мы снова берем инвариантые вероятностные меры, сосредоточенные на орбитах 
группы $L$ и их свертки.
Пример с окружностями на плоскости относится к этому типу. 

\sm

{\sc Замечание.} 1) Пусть $G$ -- группа, а $K$ -- подгруппа. Рассмотрим произведение
$G\times K$ и <<диагональное>> вложение $i:K\to G\times K$, заданное формулой  
$k\mapsto (k,k)$. Легко видеть, что классы сопряженности $G//K$ находятся
во взаимно однозначном соответствии с двойными классами смежности
$i(K)\setminus (G\times K)/i(K)$.

\sm

2) Пусть группа $L$ действует на группе $G$ автоморфизмами. Тогда определено
полупрямое произведение $L\ltimes G$ (см., например, \cite[2.4]{Kir}). Легко видеть, что орбиты 
$L$ на $G$ находятся во взаимно однозначном соответствии с двойными классами
смежности $L\subset (L\ltimes G)/L$.

Поэтому все упомянутые выше типы гипергрупп сводятся к гипергруппам двойных классов смежности.
\hfill $\lozenge$

\sm

{\sc Замечание.} Здесь возможен двойственный язык, и он более употребителен. Мы берем пространство
$C(K\setminus G/K)$ непрерывных
функций с компактным носителем, биинвариантных относительно $K$, т.е., функций
$f$, таких, что 
$$
f(k_1 g k_2)=f(g),\qquad \text{при $k_1$, $k_2\in K$.}
$$
Рассматриваем это пространство как алгебру, где в качестве умножения взята свертка.
Понятно, что, по-существу, эта алгебраическая структура эквивалентна структуре
гипергруппы.\hfill $\lozenge$

\sm

Приведем три примера гипергрупп, которые заодно показывают содержательность этого понятия.

\sm

{\it Алгебры Гекке.} См. \cite{Iwa}. Рассмотрим группу $G=\GL(n,\bbF_q)$
матриц размера $n$ над конечным полем $\bbF_q$ из 
$q$ элементов. Пусть $B$ -- подгруппа верхнетреугольных матриц. 
Пусть $S_n\subset \GL(n,\bbF_q)$ -- группа матриц перестановок.
Любой двойной класс $B\setminus G/B$ смежности имеет вид $B w B$, где
$w\in S_n$, причем при разных $w$ получаются разные
классы (докажите это; это то, что на мехмате проходят по алгебре на первой лекции
первого курса).

Рассмотрим сверточную алгебру $\cH_q=C(B\setminus G/B)$ ({\it алгебру Гекке},
изобретенную и названную так Н.~Ивахори). Ее базисные
элементы нумеруются подстановками $w\in S_n$. Обозначим через $T_j$
элементы, отвечающие простым транспозициям $(j,j+1)$. Оказывается,
что алгебра Гекке -- это алгебра с образующими $T_j$ и
соотношениями 
\begin{align*}
 T_j^2=q+(q-1)T_j;
 \\ 
 T_j T_{j+1} T_j= T_{j+1} T_jT_{j+1};
 \\
 T_kT_l=T_l T_k\qquad \text{при $l\ne k\pm 1$.}
\end{align*}

Интересно, что при подстановке $q=1$ в эти формулы мы получаем соотношения для элементарных транспозиций 
$(j,j+1)$ в симметрической группе: первое соотношение превращается в $T_j^2=1$, 
а второе в $(T_j T_{j+1})^3=1$. Таким образом, алгебра $\cH_1$ -- это групповая алгебра
симметрической группы. Отметим, что алгебра Гекке $\cH_q$ определена при любом
$q\in\C$. При всех $q$, кроме корней из 1, эти алгебры изоморфны групповой алгебре симметрической группы.

\sm

{\it <<Гипергруппа обобщенного сдвига>>.} См. \cite{Koo}.
Пусть $\K$ -- одно из тел, 
$\R$, $\C$ или кватернионы $\H$. Рассмотрим пространство $\K^{n+1}$, 
снабженное эрмитовой формой с одним минусом и $n$ плюсами, ее матрица равна
$$\begin{pmatrix}
  -1&0&\dots\\
  0&1&\dots\\
  \vdots&\vdots&\ddots
 \end{pmatrix}
$$
(в случае $\K=\R$ эрмитовость означает симметричность). Рассмотрим группу $G=\U(1,n;\K)$,
состоящую из всех матриц, сохраняющих эту форму, это соответственно группа $\O(1,n)$,
$\U(1,n)$ или  $\Sp(1,n)$. В $G$ рассмотрим максимальную компактную подгруппу 
$\K= \U(1;\K)\times\U(n,\K)$ она состоит из матриц блочных матриц
размера $1+n$ вида $\begin{pmatrix}
                                                          a&0\\0& D
                                                         \end{pmatrix}$,
 где $a$ -- элемент тела с нормой 1, а $D$ -- унитарная матрица.
 Соответственно, однородное пространство $\U(1,n;\K)/ \U(1;\K)\times\U(n,\K)$
 - это пространство Лобачевского (при $\K=\R$), комплексное или кватернионное
 гиперболическое пространство. 
 
 Мы рассматриваем гипергруппу
 двойных классов смежности $K\setminus G/K$. Любой такой класс смежности имеет вид
 $K h_t K$, где 
 $$
 h_t=\begin{pmatrix}
      \ch t&\sh t&0&\dots\\
      \sh t&\ch t &0&\dots\\
      0&0&1&\dots \\
      \vdots&\vdots&\vdots&\ddots
     \end{pmatrix}, \qquad \text{где $t\ge 0$.}
 $$
 Таким образом, множество $K\setminus G/K$ отождествляется с полупрямой $t\ge 0$.
 Введем обозначения 
 $$
 d:=\dim \K,\qquad \alpha=nd/2-1,\qquad \beta= d/2-1.
 $$
 Пусть 
 $$_2F_1[a,b,c;x]=\sum_{j=0}^\infty \frac{(a)_j(b)_j}{(c)_j j!} x^j$$
 обозначает гипергеометрическую функция Гаусса. Тогда
 $$
 \delta_t\star \delta_s=\int_{|t-s|}^{t+s} L(t,s,u)\delta_u \,du
 ,$$
 где ядро $L(t,s,u)$ задается формулой
 \begin{multline*}
  L(t,s,u)=\frac{2^{-d(n+1)+2} (\ch s\ch t\ch u)^{\alpha-\beta-1}}
  {\pi^{1/2} \Gamma(\alpha+1/2) (\sh s\sh t\sh u)^{2\alpha} }
  \times\\ \times
  (1-B^2)^{\alpha-1/2}
 {}_2F_1\bigl[\alpha+\beta, \alpha-\beta; \alpha+1/2;\tfrac 12(1-B) \bigr]
, \end{multline*}
а
$$
B:=\frac{\ch^2 t+\ch^2 s+ \ch^2 u-1} {2\ch s\ch t \ch u}
.$$
Отметим, что выражение для $L(s,t,u)$
 симметрично по переменным $s$, $t$, $u$, в частности наша гипергруппа
коммутативна. 

Эти формулы можно интерпетировать так. Возьмем на гиперболическом пространстве
$G/K$ два интегральных оператора с ядрами $P(d(x,y))$, $Q(d(x,y))$,  зависящими
лишь от расстояния $d(x,y)$ между точками $x$, $y$. Тогда их произведение будет интегральным оператором
с ядром $R(d(x,y))$, задаваемым формулой
$$
R(u)=\int K(t,s,u) P(t) Q(s)\,dt\,ds.
$$

Интересно, что как и в случае алгебр Гекке, у нас получается семейство
операций, зависящее от непрерывных параметров (можно брать любые $\alpha\ge \beta\ge 1/2$).
Таким образом, гипергруппы $K\setminus G/K$ интерполируются по размерности тела 
$\K$ и по размерности пространства.

\sm

{\it Орбитальные меры.} Рассмотрим пространство эрмитовых матриц размера $n$, на котором действует 
естественным образом группа $\U(n)$. Орбиты $\cO_\Lambda$ группы $\U(n)$
нумеруются наборами $\Lambda=(\lambda_1,\dots,\lambda_n)$
собственных чисел матрицы. Рассмотрим инвариантную вероятностную меру $\delta_\Lambda$,
сосредоточенную на такой орбите. Мы снова получаем гипергруппу. Встает вопрос
о вычислении сверток $\delta_\Lambda*\delta_M$. В принципе, для этих сверток может быть
написана явная формула в элементарных функциях (см. \cite{DRW}), но выражение оказывается знакопеременной суммой
большого числа слагаемых, так что даже не удается понять, какой будет носитель свертки.
Что касается описание носителя свертки, то это буквально совпадает c вопросом
об описании множества, которое может пробегать спектр суммы двух матриц с данными спектрами
(это долго стоявшая проблема, решенная A.~A.~Клячко в 1993г., см. ответ, например, в 
\cite[Add. to Ch.2]{Ner2}).

\sm

Отметим, что гипергруппы двойных классов смежности -- структуры довольно сложные,
и перечислим случаи, когда они фактически применяются. Мы приводим
лишь по типовому примеру (дальше надо иметь в виду всевозможные полупростые группы).

\sm

1a) $G=\GL(n,\R)$, $K=\O(n)$.

\sm

1b) $G=\U(n)$, $K=\O(n)$.

\sm

1c) $G=\O(n)\ltimes \Symm(n)$, $K=\O(n)$, здесь $\Symm(n)$ -- пространство симметрических матриц.

Все гипергруппы 1a)-1c) коммутативны (И.~М.~Гельфанд).

\sm

2) Алгебры Гекке, обсуждавшиеся выше.

\sm

3a) Пусть $G=\GL(n,\Q_p)$, $K=\GL(n,\bbO+p)$, где $\bbO_p$ -- целые $p$-адические числа.
Гипергруппа $K\setminus G/K$ коммутативна.

\sm

3b) Пусть снова $G=\GL(n,\Q_p)$, а $K$ теперь подгруппа Ивахори (см. \cite[\S 10.5]{Ner2}), состоящая из матриц,
у которых элементы выше диагонали лежат в $\bbO_p$, элементы ниже диагонали содержатся
в $p\cdot \bbO_p$, а диагональные элементы содержатся в $\bbO_p^\times:=\bbO_p\setminus (p\cdot \bbO_p)$.
Алгебры $K$-биинвариантных непрерывных функций на $G$ -- это т.н. аффинные алгебры Гекке.

\sm

См. \cite{Buch}, \cite{DRW}, \cite{Iwa}, \cite{IM},  \cite{Jew}, \cite[Sect.3, 7]{Koo}, 
\cite{Mac}, \cite{Rosl}.

\sm

{\bf\punct Представления гипергрупп.%
\label{ss:repres-hyper}} Пусть $G$ -- локально компактная группа, $K$ -- ее компактная подгруппа.
Пусть $\rho$ -- унитарное представление группы   $G$ в гильбертовом пространстве 
$H$ (см. следующий параграф). Рассмотрим подпространство $V$, 
состоящее из всех $K$-неподвижных векторов. Через 
$P$ обозначим ортогональный проектор на $V$. 

\sm

{\sc Задача.} Обозначим через $\delta_K$ вероятностную меру Хаара на $K$, рассматриваемую как 
 меру на $G\supset K$. Рассмотрим оператор 
 $$
 \rho(\delta_K):=\int_K \rho(k)\,dk
 $$
 Покажите, что этот оператор совпадает с проектором $P$.
 \hfill $\lozenge$

\sm

Теперь заметим, что дельта-функция двойного класса смежности $\frg=KgK$ (см. предыдущий пункт)
может быть записана в виде свертки: 
$$
\delta_{KgK}=\delta_K*\delta_g*\delta_K
$$
Поэтому
$$
\rho\bigl(\delta_{KgK}\bigr)=\rho(\delta_{K}) \rho(\delta_{g}) \rho(\delta_{K})= P\rho(g) P.
$$
Таким образом, эти операторы фактически действуют в пространстве $V$.
Более того, мы получам представление гипергруппы двойных классов
смежности в пространстве $V$ в следующем смысле: если $\delta_\frg*\delta_\frh=\mu_{\frg,\frh}$,
то
$$
\rho(\delta_\frg*\delta_\frh)=\rho(\mu_{\frg,\frh}).
$$

С другой стороны, если представление $\rho$ неприводимо, а $V\ne 0$, 
то $\rho$ однозначно определяется представлением гипергруппы в $V$,
что делает гипергруппы инструментом работы с представлениями.
Мы не будем этого касаться подробнее.

\sm

{\bf\punct Операторные узлы.}
Мы продолжаем сюжет с гипергруппами, см. п.\ref{ss:hypergroups}.
Рассмотрим группу $\U(n+m)$ и ее подгруппу $\U(n)$.
 Введем на наших группах
метрику, определяемую евклидовой операторной нормой. Рассмотрим 
гипергруппу классов сопряженности $\U(n+m)//\U(m)$ (такие классы называются
<<{\it операторными узлами}>>%
\footnote{Английские термины: colligation, node.}). Оказывается, что
при $m\to\infty$ эта гипергруппа вырождается в полугруппу с настоящим умножением.
Сначала опишем это умножение. Рассмотрим группу $\U(n+\infty)$, 
которая состоит из блочных финитных унитарных матриц размера 
$(n+\infty)\times (n+\infty)$, через $\U(\infty)$ 
мы обозначим ее подгруппу, состоящую из матриц вида 
$\begin{pmatrix} 1&0\\0&u \end{pmatrix}$.
Рассмотрим пространство $\Gamma$  классов сопряженности 
$\U(n+\infty)//\U(\infty)$. Введем на нем {\it умножение операторных узлов} 
по формуле 
\begin{equation}
\begin{pmatrix} a&b\\ c&d \end{pmatrix} \circ
\begin{pmatrix} p&q\\ r&t \end{pmatrix}=
\begin{pmatrix}
 a&b&0\\
 c&d&0\\
 0&0&1
\end{pmatrix}
\begin{pmatrix}
 p&0&q\\
 0&1&0\\
 r&0&t
\end{pmatrix}
=\begin{pmatrix}
  ap&b&aq\\
  cp&d&cq\\
  r&0&t
 \end{pmatrix}
 .
 \label{eq:circ}
\end{equation}
Новая матрица имеет размер $n+\infty+\infty$, который выглядит отличным
от начального размера $n+\infty$. Мы, однако, заметим, что $\infty+\infty=\infty$.
Точнее мы установим какую-нибудь биекцую между двумя возникающими счетными множествами.

\sm

{\sc Задача.}
a) Покажите, что мы получили операцию на классах сопряженности, т.е.,
результат $\circ$-умножения (как класс сопряженности) зависит лишь 
от классов сопряженности начальных матриц, а не от выбора их представителей.

\sm

b) Покажите, что эта операция ассоциативна. \hfill $\lozenge$

\sm

Теперь объясним, в каком смысле построенный объект с умножением является
пределом конечномерных гипергрупп. Фиксируем $k$. Возьмем две матрицы 
$$
g_1=\begin{pmatrix} a&b\\ c&d \end{pmatrix},\qquad
g_2=\begin{pmatrix} p&q\\ r&t \end{pmatrix}\in \U(n+k)
.
$$
Возьмем большое $m=k+k+L$
и достроим наши матрицы до блочных матриц размера $n+k+k+L$,
$$
g_1^L=\begin{pmatrix} a&b&0&0\\ c&d&0&0\\ 0&0&1&0\\ 0&0&0&1  \end{pmatrix},\qquad
g_2^L=\begin{pmatrix} p&q&0&0\\ r&t&0&0\\ 0&0&1&0 \\ 0&0&0&1 \end{pmatrix}\in \U(n+k+k+L).
$$
Определим матрицу
$$
g_1
\circledcirc
g_2:=
\begin{pmatrix} a&b&0&0\\ c&d&0&0\\ 0&0&1&0\\ 0&0&0&1  \end{pmatrix}
\begin{pmatrix} p&0&q&0\\0&1&0&0 \\ r&0&t&0\\ 0&0&0&1 \end{pmatrix} \in \U(n+k+k+L)
$$
Рассмотрим соответствующие вероятностные инвариантные меры $\phi^L$, $\psi^L$ 
на классах сопряженности $\U(n+k+L)//\U(k+L)$. 
Обозначим  $\epsilon$-окрестность
класса сопряженности 
$g_1\circledcirc g_2$ через $\cN_\epsilon$. Оказывается, что
при больших $L$ свертка $\phi^L*\psi^L$ почти сконцетрирована
в $\cN_\epsilon$. Более точно, для любого $\epsilon>0$ для любого 
$\delta>0$ существует $L_0$, такое,  что для любого
$L\ge L_0$ мера $\phi^L*\psi^L$ множества  $\cN_\epsilon$
больше $1-\delta$. 

\sm

Отметим, что $\circ$-умножение, по разным причинам (не тем, которые мы обсуждаем), 
появилось в спектральной теории несамосопряженных
операторов (М.~С.~Лившиц) и в теории систем.

\sm

Далее, оказывается, что эта операция <<материализуется>> 
следующим образом. Для матрицы 
$g=\begin{pmatrix}  a&b\\c&d \end{pmatrix}$
мы, следуя Лившицу, пишем {\it характеристическую функцию,}
$$
\chi_g(\lambda)=a+\lambda b(1-\lambda d)^{-1} c
, $$
где $\lambda\in \C$. Мы можем также положить $\chi(\infty)=-d^{-1}$.
Мы получаем рациональную функцию на $\C$ со значениями в пространстве матриц 
$n\times n$. Очевидно, что эта функция голоморфна внутри круга $|z|=1$
(т.к. $\|d\|\le 1$).

\sm

{\sc Задача.} a) Проверьте, что для матриц
$g_1$, $g_2$, лежащих в одном классе 
$\U(n+\infty)//\U(\infty)$, характеристические 
функции совпадают.

\sm

b) Проверьте, что 
$$
\chi_{g\circ h}(\lambda)=\chi_{g} (\lambda) \chi_h (\lambda)
.$$

c) Покажите, что при $|\lambda|=1$ значение $\chi_g(\lambda)$
содержится в $\U(n)$. См. доказательство частного случая этого утверждения
в п.\ref{th:livshits}.

\sm

d) Если $|\lambda\|\le 1$, то $\|\chi_g(\lambda)\|\le 1$.
Если $|\lambda|\ge 1$, то $\|\chi_g(\lambda)\|$ является растягивающим
оператором, т.е. $\|\chi_g(\lambda)v\|\ge \|v\|$ для всех векторов $v\in\C^n$.
\hfill $\lozenge$

\sm

Дальше встает вопрос о том, восстанавливается ли класс сопряженности
по характеристической функции. Оказывается, что для унитарных
матриц общего положения ответ утвердителен, но возможны вырождения, которые
устроены следующим образом. Рассмотрим всевозможные унитарные матрицы вида
\begin{equation}
\begin{pmatrix}
 a&b&0\\ c&d&0\\ 0&0&e
\end{pmatrix},
\label{eq:nechist}
\end{equation}
Мы имеем дело с классами сопряженности, поэтому блок $e$
можно считать диагональной матрицой.
Понятно, что характеристическая функция вообще не видит блока $e$
(убедитесь в этом).
Оказывается, что любой класс сопряженности имеет представитель вида 
(\ref{eq:nechist}), такой, что у блока $d$ нет собственных
чисел, лежащих на окружности $|\lambda|=1$. 

Полный набор данных, достаточных для восстановления
операторного узла, это характеристическая функция 
и набор собственных чисел матрицы $d$, лежащих на единичной
окружности.

\sm

{\sc Замечание.} Это можно сказать иначе, 
Введем еще один инвариант,
а именно функцию 
$$
\theta_g(\lambda):=\det(1-\lambda d).
$$
Отметим, что эта функция однозначно определяется множеством
$\cM_g$ своих нулей.
С другой стороны, полюса характеристической характеристической
функции -- это в точности нули функции $\theta_g(\lambda)$ в области
$|\lambda|>1$. А нулей $\theta_g(\lambda)$, лежащих на окружности,
характеристическая функция не видит. Теперь осталось понять, что
операторный узел восстанавливается по характеристической функции
и $\theta_g$. Это несложное упражнение по теории инвариантов,
что, впрочем, лежит за пределами этих записок.
\hfill$\lozenge$

\sm

См. \cite{Bro}, \cite[Chapter 19]{Dym}, \cite[Добавление E]{Ner1}, \cite{Ner-hyper}.

\sm

{\bf\punct Шлейфы.} Оказывается, что феномен вырождения гипергрупп (впервые отмеченный Г.~И.~Ольшанским
в 1980) является весьма общим и порождает большое количество необычных алгебраических
структур. Приведем один пример. 

\sm

{\it Пример шлейфа.}
Рассмотрим ту же группу 
$G=\U(\infty)$ унитарных блочных финитных матриц размера $\infty$
и в ней подгруппы $K[\alpha]$, состоящие из блочных $\alpha+(\infty-\alpha)$-матриц вида 
$\begin{pmatrix} 1&0\\0&u \end{pmatrix},$
где $u$-вещественная оргогональная
матрица. Рассмотрим пространства двойных классов смежности
$$K[\alpha]\setminus G/K[\beta],$$
Мы будем записывать  элементы такого пространства как блочные матрицы
$
\begin{pmatrix}
 a&b\\c&d
\end{pmatrix}
$
размера 
$$\bigl(\alpha+(\infty-\alpha)\bigr)\times \bigl(\beta+(\infty-\beta)\bigr),$$
определенные с точностью до эквивалентности
$$
\begin{pmatrix}
 a&b\\c&d
\end{pmatrix}
\sim  \begin{pmatrix} 1&0\\0&u \end{pmatrix}\begin{pmatrix}
 a&b\\c&d
\end{pmatrix} \begin{pmatrix} 1&0\\0&v \end{pmatrix}.
$$
Мы определим отображение 
$$
K[\alpha]\setminus G/K[\beta]\,\times\, K[\beta]\setminus G/K[\gamma]\,
\to\, K[\alpha]\setminus G/K[\gamma]
$$
по той же формуле (\ref{eq:circ}).
 Это корректно определенное ассоциативное
умножение.

Чтобы избавиться от ссылок на гипергруппы и одновременно
не превратить появление подобных операций в черную магию, мы определим это 
умножение другим способом (и  этот способ является 
универсальным общим приемом).

Возьмем последовательность матриц 
$$
\theta_j[\beta]=
\begin{pmatrix}1_\beta&0&0&0\\
 0&0&1_j&0\\
 0&1_j&0&0\\
 0&0&0&1_\infty
\end{pmatrix}\in K[\beta]
.$$
Рассмотрим последовательность двойных классов смежности
\begin{equation} 
\frr_j:=K[\alpha]\cdot
\begin{pmatrix}
 a&b\\c&d
\end{pmatrix}\,\theta_j[\beta]\,
\begin{pmatrix}
 p&q\\r&t
\end{pmatrix}\cdot K[\gamma]
\label{eq:1000}
.\end{equation}

{\sc  Задача.} Убедитесь, что эта последовательность стабилизируется с некоторого
номера и ее предел равен $\circ$-произведению двойных классов смежности.
\hfill $\lozenge$

\sm

В итоге мы получаем некоторую алгебраическую структуру,
а именно категорию (определений категорий см., например, в \cite[\S 2.5]{Ner1}).
Объекты этой категории -- неотрицательные целые числа,
а множество морфизмов из $\beta$ в $\alpha$ -- это пространство $K[\alpha]\setminus G/K[\beta]$

Мы назовем эту категорию {\it шлейфом пары} $\U(\infty)\supset \O(\infty)$.

\sm

{\it Представления шлейфа.}
Теперь рассмотрим унитарное представление $\rho$ группы $\U(\infty)$ в гильбертовом пространстве 
$H$. Обозначим через $H[\alpha]$  пространство $K[\alpha]$-неподвижных векторов, а через
$P[\alpha]$ -- проектор на $H[\alpha]$. Для
$g\in \U(\infty)$ мы рассмотрим оператор 
$$\ov \rho(g):H[\beta]\to H[\alpha], $$
заданный формулой,
\begin{equation}
\ov \rho(g):=P[\alpha] \rho(g)\Bigr|_{H[\beta]}
.
\label{eq:ov-rho}
\end{equation}
Мы получаем оператор, который зависит лишь от двойного класса смежности $\frg$,
содержащего $g$. Действительно, возьмем векторы
$\xi\in H[\alpha]$, $\eta\in H[\beta]$, возьмем $k_1\in K[\alpha]$,
$k_2\in K[\beta]$ и рассмотрим 
<<матричный элемент>>
\begin{multline*}
\la \xi, \ov\rho(k_1 g k_2) \eta\ra=\\=
\la \xi, \rho(k_1 g k_2) \eta\ra= \la \xi, \rho(k_1) \rho( g) \rho( k_2) \eta\ra
=\la \rho(k_1)^{-1} \xi,  \rho( g) \rho( k_2) \eta\ra=\la \xi, \rho(g) \eta\ra
=\\=
\la \xi,\ov \rho(g) \eta\ra
.
\end{multline*}
Мы видим, что матричные элементы операторов $\ov\rho(k_1 g k_2)$ и $\ov \rho( g)$
равны, а тем самым операторы совпадают. Таким образом,
мы можем писать оператор $\ov\rho(\frg)$, где $\frg$ -- двойной класс смежности
(сравните с п.\ref{ss:repres-hyper}).

Теперь возьмем два класса смежности $\frg\in\K[\alpha]\setminus G/K[\beta]$,
$\frh\in \K[\beta]\setminus G/K[\gamma]$. Оказывается, что
верна <<{\it теорема мультипликативности}>>:
$$
\rho(\frg\circ \frh)=\rho(\frg)\,\rho(\frh).
$$
Иными словами, мы получили представление категории двойных классов смежности
(определение представлений категорий см., например, в 
\cite[\S 2.5]{Ner1}).
Сравните это утверждение с рассмотрениями  п.\ref{ss:repres-hyper}).

\sm

{\it Характеристические функции.}
Объясним, как материализовать $\circ$-умножение. Для простоты, положим 
$\alpha=\beta=\gamma$ (это не очень существенно).
Рассмотрим двойной класс смежности $\frg$,
содержащий $g=\begin{pmatrix} a&b\\c&d\end{pmatrix}$.
Фиксируем $\lambda\in \C$
и запишем следующее
уравнение (это что-то вроде извращенного варианта уравнения на собственные числа)
$$
\begin{pmatrix}
q_+\\ x\\ q_- \\ \lambda x
\end{pmatrix}
=
\begin{pmatrix}
 \begin{pmatrix} a&b\\c&d\end{pmatrix}&0\\
 0& \begin{pmatrix} a&b\\c&d\end{pmatrix}^{t-1}
\end{pmatrix}
\begin{pmatrix}
p_+\\ \lambda  x\\ p_- \\ x
\end{pmatrix}
.$$
Мы рассматриваем это равенство как систему линейных соотношений на переменные
$p$, $q$, $x$. Исключим из этой системы $x$ и получим
некоторую зависимость
$$
\begin{pmatrix}
 q_+\\q_-
\end{pmatrix}
=\chi_\frg (\lambda) \begin{pmatrix}
 p_+\\p_-
\end{pmatrix}
,$$
где <<{\it характеристическая функция}>>
$\chi_\frg(\lambda)$ -- рациональная 
функция на $\C$, принимающая значения в пространстве  матриц размера 
$2\alpha\times 2\alpha$.
Тогда, как и в случае операторных узлов, $\circ$-умножение переходит
в поточечное  умножение  характеристических функций. 

Как и в случае операторных узлов, эта функция удовлетворяет некоторым дополнительным
условиям, они несколько сложней, чем случае операторных узлов,
и мы не будем их здесь обсуждать, подробности
см. в \cite[\S IX.4]{Ner1}).

\sm

Отметим, что первый пример такого рода умножения двойных
классов смежности  был обнаружен P.~C.~Исмагиловым
в 60х годах. Явление это весьма общее. Например,
ниже в \S\S 8-11 у нас появляются различные большие группы симметрий мер,
пп.\ref{ss:GLO}, \ref{ss:poisson-quasiinvariant}, \ref{ss:bisymmetric}, \ref{ss:UUU},
для всех этих групп естественным образом возникают  подобные структуры.

\sm

См. \cite{Ner1}, \cite{Ner-umn}, \cite{Ner-colligations}, \cite{Ner-spheric},
\cite{Olsh-GB}, \cite{Olsh-chip}, \cite{Ner-buildings}.

\sm

{\bf \punct Шлейф бисимметрической группы и чипы.%
	\label{ss:bisymmetric-0.5}}
Пусть $G$ -- бисимметрическая группа, пусть
$K$ -- ее диагональ. Напомним, что
$K$ состоит из пар $(g,g)$, где
$g$ пробегает полную бесконечную симметрическую группу. Рассмотрим в
$K$ подгруппу $K[\alpha]$, определенную  выше
(мы берем перестановки $g$, фиксирующие
1, 2, \dots, $\alpha$). Дадим комбинаторное описание двойных классов смежности
$K[\alpha]\setminus G/K[\beta]$.

Элементы самой группы $G$ мы будем представлять себе
как диаграммы вида
$$
\epsfbox{bisym.1}
$$ 
Две копии натурального ряда мы обозначим через $\{1,2,3,\dots\} $  и
$\{1',2',3',\dots\} $. Две элемента $g$, $g'\in \ov S_\infty$ 
переставляют каждый свою копию. Левая и правая перестановки почти симметричны 
друг другу относительно вертикальной пунктирной оси за исключением конечного
участка вблизи самой пунктирной оси.

Теперь опишем двойные классы смежности. Соединим дугами пары
кружочков
$k$, $k'$ верхнего ряда с $k>\alpha$ и пары $l$, $l'$
нижнего ряда с $l>\beta$. 
\begin{figure}
$$\epsfbox{bisym.2}$$
$$\qquad\qquad\epsfbox{bisym.3}$$
\caption{Построение чипа (на дугах чипа,
	 не помеченных числами, по умолчанию предполагается 0).\label{fig:chip}}
\end{figure}
Припишем каждой горизонтальной дуге длину $1/2$,
каждой вертикальной дуге длину 0. Забудем о верхних
кружочках с номерами $>\alpha$ и нижних кружочках
с номерами $>\beta$.
Мы получаем картинку ({\it чип}),
устроенную следующим образом.
В верхнем ряду стоят кружочки с метками
$1$, $2$, \dots, $\alpha$ 
и $1'$, $2'$, \dots, $\alpha'$.
В нижнем ряду стоят кружочки с метками
$1$, $2$, \dots, $\beta$ 
и $1'$, $2'$, \dots, $\beta'$.
Картинка разделена на две половинки 
вертикальной линией.
Информация содержится в дугах.
Они бывают трех типов.

\sm

1. <<Вертикальные дуги>>. Они идут от верхних 
кружочков к нижним, и оба конца дуги содержатся
либо в левой, либо в правой половине картинки.
Длина дуги является целым числом (возможно, нулем).

\sm

2. <<Горизонтальные дуги>>. Они соединяют два
кружочка верхнего ряда, либо   два кружочка нижнего ряда. Кружочки обязательно находятся  по разные стороны пунктирной линии. Длина горизонтальной дуги
-- число вида $k+1/2$, где $k\ge0$ -- целое число.

\sm

3. Циклы. Длина цикла -- целое число $\ge 1$.
Число циклов счетно, причем почти все циклы имеют
длину 1. 

\sm

Циклы длины 1 можно выкинуть, после чего мы
получаем финитную картинку.

\sm

{\sc Задача.} a) Убедитесь, что чипы находятся во взаимно однозначном соответствии с двойными классами
смежности. Кстати, где мы использовали, что пара
подстановок лежит в бисимметрической группе?

\sm

b) Придумайте определение этого умножения в духе (\ref{eq:1000}).
\hfill $\lozenge$

\sm

Умножение чипов состоит в склейке картинок (и забывании точек слейки),
см. Рис. \ref{fig:product-chips}. Длины дуг складываются.
Циклы переходят без изменений в новый чип (конечно, при склейке могут образовываться
дополнительные циклы).

\begin{figure}
	$$\epsfbox{bisym.4}$$
	\caption{Умножение чипов.\label{fig:product-chips}}
\end{figure}

\sm

Так же, как в случае классических групп, умножение чипов
является пределом сверток, и, что более существенно,
для чипов верна теорема мультипликативности (см. предыдущий пункт).

\sm

См. \cite{Olsh-chip}, \cite{Ner-umn}, \cite{Olsh-semi}, \cite{Ner-hyper}.

\sm

{\bf \punct Супергруппы} (для тех, кто знаком с группами Ли).

 {\it Супералгебры Ли}. Супералгеброй Ли  называется следующая структура.
Рассмотривается линейное пространство $\frg$, разложенное в прямую сумму
$\frg=\frg_{\ov 0}\oplus \frg_{\ov 1}$, элементы этих подпространств мы будем
называть однородными элементами и обозначать их степени, т.е., 0 или 1, через $\deg (x)$.
Вводится билинейная операция $\frg \times \frg \to \frg$ (суперкоммутатор $[x,y]$),
причем 

\sm

1. Эта операция $\Z_2$-градуирована, т.е., для однородных элементов $x$, $y$
выполнено
$$\deg[x,y]=\deg x+ \deg y,$$

2. Операция <<суперантикоммутативна>>, т.е.,
для однородных элементов
$$
[x,y]=(-1)^{\deg x\cdot \deg y}[x,y].
$$

3. Операция удовлетворяет супертождеству Якоби,
$$
[x,[y,z]]=[[x,y],z]+ (-1)^{\deg x \cdot \deg y}[y,[x,z]].
$$

Таким образом, четная часть $\frg_{\ov o}$
является алгеброй Ли, нечетная часть $\frg_{\ov 1}$
является $\frg_{\ov1}$-модулем; кроме того есть коммутативное отображение
$\frg_{\ov 1}\times \frg_{\ov 1}\to \frg_{\ov 0}$, на которое супертождество
Якоби налагает дополнительные условия.

\sm

{\sc Примеры.}
a) {\it Алгебра $\frg\frl(p|q)$.} Элементами алгебры являются комплексные блочные
 матрицы
$X:=\begin{pmatrix}
a&b\\c&d
\end{pmatrix}
$
размера $(p+q)\times (p+q)$. Матрицы вида 
$\begin{pmatrix}
	a&0\\0&d
\end{pmatrix}
$ имеют степень $\ov 0$, матрицы $\begin{pmatrix}
	0&b\\c&0
\end{pmatrix}$   -- степень $\ov 1$. Если $X$, $Y$ --
однородные элементы, то их суперкоммутатор определяется
через матричное умножение как
$$
[X,Y]=XY-(-1)^{\deg X\cdot \deg Y}YX.
$$ 

b) {\it Алгебра супер-Вирасоро.} Рассмотрим алгебру $s-\Vir$,
с базисом, составленным из четных элементов $L_n$,  где $n$ пробегает $\Z$,
нечетных элементов
$J_\alpha$, где $\alpha$ пробегает $\Z$, и четного центрального элемента
$c$ (<<центральный>> означает, что $[c,L_n]=0$, $[c,J_\alpha]=0$
для всех $n$, $\alpha$), соотношения коммутации имеют вид
\begin{align}
[L_n,L_m]&=L_{n+m}+\frac{n^3-n}8\cdot \delta_{m+n,0}\cdot c;
\label{eq:super-vir1}\\
[L_m, J_\alpha]&= \bigl(\frac m2-\alpha\bigr)J_{\alpha+m};
\label{eq:super-vir2}\\
[J_\alpha,J_\beta]&=2L_{\alpha+\beta}+\frac 12 \bigl(\alpha^2-1/4\bigr)\cdot\delta_{\alpha+\beta,0} \cdot c.
\label{eq:super-vir3}
\end{align} 
Соотношения (\ref{eq:super-vir1}) для $L_n$ -- это обычные соотношения в алгебре Вирасоро,
см. ниже (\ref{eq:vir}), генераторы $L_n$ соответствуют векторным полям
$e^{in\phi}\frac\partial{i\partial \phi}$ на окружности. Генераторы
$J_\alpha$ сответствуют плотностям $e^{in\alpha}\,(d\phi)^{-1/2}$
веса $-1/2$, формула
(\ref{eq:super-vir2}) соответствует действию векторных полей на плотностях.
Наконец (\ref{eq:super-vir3}), если отбросить центральную 
поправку, соответствует умножению плотностей веса $-1/2$. 
\hfill $\lozenge$

\sm

{\it Грассманизация. Алгебры Ли.}
Согласно Ф.~А.~Березину и Г.~И.~Кацу, из супералгебры Ли может быть произведено
семейство алгебр Ли и   соответствующих групп по следующему рецепту. Рассмотрим набор антикоммутирующих переменных
$\theta_1$, $\theta_2$, \dots (набор может быть конечным или бесконечным),
$$\theta_i\theta_j=-\theta_j\theta_i,\qquad \theta_j^2=0.$$
Обозначим через $\Lambda$ алгебру, порожденную $\theta_j$,
т.е., грассманову алгебру (см., например, \cite[\S II.1]{Ner1}).
 Разложим $\Lambda=\Lambda_0\oplus \Lambda_1$, где подпространство 
 $\Lambda_0$ порождено одночленами четной степени над $\theta_j$,
 а подпространство $\Lambda_1$ одночленами нечетной степени.
 Пусть $\frg=\frg_{\ov 0}\oplus \frg_{\ov 1}$ -- супералгебра Ли.
 Пусть $x_k$
-- базис в $\frg_{\ov0}$, а $y_l$ -- базис в $\frg_{\ov1}$. Тогда $\frg(\Lambda)$ 
это алгебра, состоящая из линейных комбинаций вида
\begin{equation}
\sum_k a_k(\theta) x_k+ \sum_l b_l(\theta) y_l, 
\qquad\text{где $a_k(\theta)\in \Lambda_0$, $b_l(\theta)\in \Lambda_{1}$}
,
\label{eq:sum-super}
\end{equation}
 мы считаем, что множители, зависящие от $\theta$, выносятся
 из под знака коммутатора естественным образом,
\begin{align*}
 [a_k x_k,b_l y_l]=a_k b_l [x_k,y_l],\quad[a_kx_k,a_m x_m]=a_k a_m[x_k,x_m],
 \\
 \quad [b_l y_l,b_j y_j]=b_l b_j [y_l,y_j]
 ,\end{align*}
  а дальше
 продолжаем операцию $[\cdot,\cdot]$ по аддитивности.
 
 \sm
 
 {\sc Задача.} Покажите, что операция  $[\cdot,\cdot]$ на $\frg(\Lambda)$
 антикоммутативна и удовлетворяет тождеству Якоби. Иными словами, мы получили
 алгебру Ли. \hfill $\lozenge$

\sm

{\it Грассманизация. Супергруппы.}
Аналогично мы можем строить группы. Например, для алгебры $\frg\frl(p|q)$ мы рассмотрим
группу $\GL(p|q,\Lambda)$, состоящую из обратимых блочных матриц $\begin{pmatrix}
A&B\\C&D
\end{pmatrix}$ размера $p+q$, причем блоки $A$, $D$, составлены из элементов
$\Lambda_0$, а $B$, $D$ -- из элементов $\Lambda_1$ (убедитесь, что получилась
группа).

\sm

{\sc Задача.} a) Когда матрица $\begin{pmatrix}
A&B\\C&D
\end{pmatrix}$ обратима?

\sm

b) Опишите соответствующую алгебру Ли $\frg\frl(p|q;\Lambda)$.
\hfill $\lozenge$

\sm

Конечномерным представлением супералгебры $\frg\frl(p|q)$ естественно считать гомоморфизм
$\frg\frl(p|q)$  в некоторую супералгебру $\frg\frl(M|N)$. Такой гомоморфизм
порождает гомоморфизм алгебр Ли $\frg\frl(p|q;\Lambda)\to \frg\frl(M|N;\Lambda)$
и гомоморфизм групп $\GL(p|q;\Lambda)\to \GL(M|N;\Lambda)$. 

\sm 

{\sc Замечание.}
Подчеркнем, что
 отнюдь не все гомоморфизмы абстрактных алгебр Ли $\frg\frl(p|q;\Lambda)\to \frg\frl(M|N;\Lambda)$
и групп Ли $\GL(p|q;\Lambda)\to \GL(M|N;\Lambda)$ получаются таким 
образом. Группа $\GL(p|q;\Lambda)$ является чем-то вроде <<многообразия>>
над кольцом $\Lambda$ (как бывают многообразия над $\R$ и $\C$),
и мы рассматриваем лишь гомоморфизмы, согласованые  со структурой
<<многообразия>>.
\hfill $\lozenge$

\sm
 
Объясним, как строить группы по произвольным супералгебрам
$\frg=\frg_{\ov 0}\oplus \frg_{\ov1}$. Мы хотим построить <<группу Ли>>,
соответствующую алгебре Ли $\frg\frl(\Lambda)$.
 Через $\Lambda_0^+$ мы обозначим подпространство
в $\Lambda_0$, порожденное четными одночленами степени $>0$.
Рассмотрим алгебру $\frg^+(\Lambda)$, состоящую из выражений
(\ref{eq:sum-super}), где $a_k(\theta)\in\Lambda_0^+$, а $b_l(\theta)\in \Lambda_1$
по-прежнему. Как абстрактная алгебра Ли $\frg\frl(\Lambda)$
изоморфна полупрямому произведению $\frg_0$ и $\frg\frl^+(\Lambda)$.

 Сначала определим вспомогательную
 группу $G^+(\Lambda)$. Она состоит
из формальных выражений $\exp(X)$, где $X$  пробегает $\frg^+(\Lambda)$,
а умножение задается формулой
$$
\exp(X)\cdot \exp(Y)=\exp\Bigl(\ln(e^X e^Y)\Bigr), 
$$
где $\ln(e^X e^Y)$ мы вычисляем по формуле Кэмпбела--Хаусдорфа.
Фактически ряд Кэмпбела--Хаусдорфа всегда будет получаться конечным
из-за того, что $\theta_j^2=0$. Теперь {\it предположим,
что алгебре Ли $\frg_{\ov 0}$ соответствует некоторая группа $G^0$}
и что определено присоединенное действие $g: X\mapsto g^{-1} Xg$
группы $G^0$ на $\frg_{\ov 1}$
(если $\frg$ конечномерна, то это выполнено автоматически,
нам хочется иметь в виду также примеры типа супералгебры Вирасоро).
Тогда  $G^0$ действует и на $G^+(\Lambda)$ автоморфизмами
$$
g: \, \exp(X)\mapsto \exp(g^{-1} Xg),\qquad g\in G^0
.$$
Мы определим  
$G(\Lambda)$ как полупрямое произведение $G^0$ и 
$G^+(\Lambda)$.
Таким образом, $G(\Lambda)$ определена как абстрактная группа
(еще раз напомним, что мы рассмариваем  лишь очень специальные 
гомоморфизмы $G_1(\Lambda)\to G_2(\Lambda)$.

\sm

Для супералгебр Ли и супергрупп есть интересная теория представлений, в чем-то
похожая, а в чем-то нет, на обычную теорию представлений.

\sm

См. \cite{Ber-super}, \cite{Witt}, \cite{Serg}, \cite{dict}, \cite{CCTV}, \cite{Ner-super}.

\sm

{\bf\punct Полнота и пополнение.} Положим, для простоты, что группа $G$ имеет счетную базу
окрестностей в единице. Тогда на ней есть левоинвариантная метрика $d_l(x,y)$, совместимая 
с топологией. Можно рассмотреть пополнение по этой метрике. На полученный объект 
по непрерывности продолжается умножение, однако обратных элементов, вообще говоря, не будет.
Получится полугруппа с непрерывным умножением. Если мы возьмем правоинвариантную метрику
$d_r(x,y)=d_l(x^{-1},y^{-1})$,
то в пополнении получится другая полугруппа.

Назовем {\it лево-правой последовательностью Коши} последовательность $g_j$,
такую, что $g_k g_l^{-1}\to 1$ и $g_k^{-1} g_l\to 1$ при $k$, $l\to\infty$. 
Топологическая группа $G$ (удовлетворяющая первой аксиоме счетности) называется {\it полной} 
(в смысле Райкова%
\footnote{Н.~Бурбаки \cite[\S IV.3]{Bou1} дает другое (неэквивалентное) определение:
группа полна, если в ней сходятся левые
последовательности Коши и сходятся правые последовательности Коши. Для примера,
полная унитарная группа (со слабой топологией) и полная симметрическая группа не полны по Бурбаки.
Разумеется, предлагая определения, и Д.~А.~Райков, и Н.~Бурбаки не ограничивались классом
метризуемых групп.}),
если любая лево-правая последовательность Коши сходится. Равносильным
образом, группа полна относительно (не инвариантной) метрики $d_{lr}(x,y):=d_l(x,y)+d_r(x,y)$.
Если группа не полна, то к ней можно добавить пределы лево-правых последовательностей Коши.
Тогда получится полная группа, которая называется {\it пополнением} группы $G$. 

%

\sm 

{\sc Задача.} Рассмотрим полугруппы $\ov G_l$ и $\ov G_r$, получающиеся
из $G$ пополнением по лево- и право- инвариантным метрикам соответственно.
Покажите, что эти полугруппы антиизоморфны в следующем смысле. Отображение
$g\mapsto g^{-1}$ продолжается до непрерывной биекции  $s:\ov G_l\to\ov G_r$,
причем $s(xy)=s(y)s(x)$. \hfill $\lozenge$

\sm 

{\sc Задача.} a) Полная симметрическая группа $\bfS_\infty$ полна. 
Ее пополнение по левоинвариантной метрике -- полугруппа всех инъективных отображений $\N\to\N$.

\sm

b) Унитарная группа $\bfU(\infty)$ полна. Ее пополнение по левоинвариантной метрике 
-- полугруппа изометрий, т.е., операторов, удовлетворяющих условию
$V^*V=1$. \hfill $\lozenge$

\sm

{\bf\punct Польские группы.} Напомним, что топологическое пространство $X$ называется
{\it польским,} если на нем существует метрика, совместимая с топологией, превращающая
$X$ в полное сепарабельное метрическое пространство.

\sm

{\sc Задача.} Открытое подмножество польского пространства является польским пространством;
$G_\delta$-подмножество%
\footnote{Подмножество в метрическом пространстве называется $G_\delta$-подмножеством, если
оно представимо в виде счетного пересечения открытых множеств. Пример: множество иррациональных
чисел.}
 польского пространства является польским пространством.
 \hfill $\lozenge$

\sm

Топологическая группа называется {\it польской}, если она является польским топологическим
пространством. В этом случае она автоматически полна (в смысле Райкова).

Можно еще сказать, что польская группа -- это полная метризуемая сепарабельная группа.

Почти все группы, перечисленные выше -- польские. Исключение составляют несепарабельные
 группы (полная унитарная группа  и
 группа $\Ams(M)$, снабженные несепарабельными топологиями.
 Неполна также группа гамильтоновых диффеоморфизмов относительно метрики Хофера
(и хороших описаний ее пополнения не известно).

\sm

{\it Примеры польских групп.} a) Локально компактные группы со счетной базой топологии.
 В частности, связные группы Ли (или группы Ли с конечным или счетным числом компонент 
 связности).
 
 \sm

b) Пусть $X$ -- полное сепарабельное метрическое пространство. Пусть $\mathrm{Iso}(X)$
-- группа его изометрий. Мы скажем, что последовательность изометрий $g_j$ сходится к
$g$, если для любого $x\in X$ выполнено  $g_j x\to gx$.
 Тогда $\mathrm{Iso}(X)$ становится польской группой.
 
 \sm
 
c) Группа изометрий сепарабельного банахова пространства (в, частности, полная унитарная 
группа).

\sm

d) Группа гомеоморфизмов компактного метрического пространства: $g_j$ сходится
к $g$, если $\max_{x\in X} d\bigl(g_j x,gx)\to 0$.

\sm

e) Аддитивная группа любого пространства Фреше.

\sm

f) Полная симметрическая группа с топологией, введенной выше.

\sm

g) Группы $\Ams(M)$, $\Gms(M)$ преобразований пространства с мерой, снабженные сепарабельной топологией.

\sm

h) Замкнутая подгруппа польской группы является польской группой. Подгруппа, являющаяся 
$G_\delta$-множеством в польской группе, является польской группой (она автоматически замкнута).

\sm

e) Фактор-группа польской группы по замкнутой подгруппе -- польская группа.

\sm

g) Произведение счетного семейства польских групп -- польская группа. Поэтому
проективный предел (счетной цепочки) польских групп -- польская группа.

\sm

Существует содержательная абстрактная   теория польских групп. Хотя она скорее относится
к дескриптивной теории множеств, чем к анализу или алгебре, аналитик может найти в ней немало интересного.

\sm

{\it Автоматическая
непрерывность гомоморфизмов.}
  По-видимому, для польских групп имеет место такой философский факт: если гомоморфизм из польской группы $H$
   в польскую группу $G$ определен без использования аксиомы выбора, то он непрерывен.
    Одна из формальных теорем на этот счет звучит так:
    
\sm    
    
 Пусть $G$, $H$ -- польские группы.    Рассмотрим на $H$ сигма-алгебру, 
    порожденную открытыми множествами и множествами первой бэровской категории. Пусть 
    гомоморфизм $H\to G$ измерим
относительно этой $\sigma$-алгебры. Тогда он непрерывен (B.~J.~Pettis).

В частности, если на одной и той же группе введены две польских топологии, и порождаемые ими
борелевские $\sigma$-алгебры совпадают, то топологии совпадают.

Для многих групп $H$ известна автоматическая непрерывность гомоморфизмов в польские группы $G$
(без каких-либо дополнительных условий), в частности, это так для полной унитарной группы,
 полной симметрической группы, группы $\Ams(M)$, групп гомеоморфизмов отрезка
  и канторовского множества.

\sm

{\sc Задача.} Придумайте примеры разрывных гомоморфизмов топологических групп, отказавшись
от условия полноты; от условия сепарабельности. Докажите существование разрывного
гомоморфизма $\R\to\R$.
\hfill $\lozenge$

\sm

{\it Универсальные группы.} a) <<Гильбертов куб>> это счетное произведение 
$[-1,1]^\infty$ отрезков $[-1,1]$. 
Любая польская группа изоморфна некоторой замкнутой подгруппе в
группе гомеоморфизмов гильбертова куба (В.~В.~Успенский).

b) Любая польская группа изоморфна подгруппе группы Урысона, см. п.\ref{ss:uryson}.

\sm

c) Пусть $G$ -- польская группа, у которой есть фундаментальная система окрестностей единицы,
состоящая из открытых подгрупп. Тогда $G$ изоморфна замкнутой подгруппе в полной симметрической
группе.

\sm

{\it Группы, не имеющие унитарных представлений.}
Стоит заметить, что бесконечномерная унитарная группа универсальной группой в этом смысле
не является: бывают польские группы, которые в нее не запихиваются.
Примеры: 

\sm

a) группа гомеоморфизмов отрезка (М.~Г.~Мегрелишвили). 

\sm

b) группа  операторов в гильбертовом пространстве
вида $1+K$, где $K$  компактен, с топологией, индуцированной операторной нормой (Н.~И.~Нессонов).

\sm

3) В аддитивной группе $H$ сепарабельного гильбертова пространства можно выбрать дискретную подгруппу
$L$, такую, что $H/L$ не имеет нетривиальных унитарных представлений (W.Banaszczyk).

\sm

См. \cite{BK}, \cite{Kech}, \cite{Pes},  \cite{RosSo},  \cite{Tsa}.

\sm

{\bf\punct Диффеология.} На бесконечномерных группах и на бесконечномерных <<многообразиях>>
можно вводить дифференцируемую структуру, не вводя даже топологии. 

Рассмотрим множество $X$ и набор $\cD=\cD(X)$ отображений из открытых подмножеств
$U$ пространств $\R^k$ в $X$ ({\it параметризации}), причем этот набор удовлетворяет следующим аксиомам:

\sm

1) Любое отображение открытого множества в точку лежит в $\cD$.

\sm

2) Пусть $U\subset \R^k$, $V\subset\R^n$ -- области,  а $F:V\to U$ -- гладкое отображение.
Тогда для любой параметризации $P:U\to X$ отображение $P\circ F$
является параматризацией.

\sm

3) Пусть $R:U\to X$ - отображение, такое, что у любой точки $a\in U$
есть окрестность $V_a$, такая, что ограничение $R$ на $U_a$
содержится в $\cD$. Тогда $U$ содержится в $\cD$.

\sm

Такой набор   $\cD$ называется {\it диффеологией}. Неформально говоря, мы декларирум, какие
отображения  из конечномерных областей в $X$ мы считаем дифференцируемыми.

Если $f$ -- функция на $X$, мы называем ее  {\it гладкой}, если все композиции $f\circ P$,
где $P$ пробегает $\cD$, являются гладкими.

Рассмотрим более общую ситуацию, пусть 
$X$, $Y$ -- диффеологические пространства. Отображение $F:X\to Y$ называется {\it гладким},
если $P\in \cD(X)$ влечет $F\circ P\in\cD(Y)$. 

\sm

{\it Диффеологической группой} мы назовем группу $G$, снабженную  диффеологией $\cD$,  
так, что для любых двух параметризаций $P_1$, $P_2:U\to G$,
отображение $P_1\cdot P_2:U\to G$ является параметризацией,
кроме того, взятие обратного являются гладким отображением.

\sm

{\sc Замечание.} Определение диффеологии, конечно, является чрезвычайно  общим, и под него могут подпадать
довольно экзотичные объекты. Например, мы можем положить, что отображение
$U\to \R/\Q$ является гладким, если оно представимо как композиция гладкого отображения
$U\to \R$ и проекции $\R\to \R/\Q$. В итоге мы получаем диффеологию на $\R/\Q$.
\hfill $\lozenge$

\sm

Так или иначе, в качестве минимального средства для введения гладкой структуры 
этот подход небезынтересен. После небольшой модификации можно
определить и комплексно-аналитические структуры.

\sm

См. \cite{Sou}, \cite{ology}

\sm

{\bf\punct Немного монстров.%
	\label{ss:bohr}} В следующих конструкциях присутствует аксиома выбора.

\sm

 a) Пусть $K$ -- компактная группа, пусть $S$ - несчетное множество. 
 Рассмотрим прямое произведение $K^S$ группы $K$ самой на себя в количестве
$S$ штук. Более аккуратно, рассмотрим множество всех функций  $I\to K$. Согласно теореме 
Тихонова, полученная группа $K^S$ компактна (разумеется, она не удовлетворяет первой аксиоме счетности).

\sm

b) {\it Боровский} (H.~Bohr) {\it компакт.} Пусть, в обозначениях предыдущего пункта,
$I=\R$, а $K$ -- окружность $e^{i\phi}$. Рассмотрим аддитивную группу $\R$.
Каждому ее элементу $x$ поставим в соответствие функцию $f_x(s):I\to K$, заданную формулой
$f_x(s)=e^{isx}$. Замкнем это семейство функций  в $K^S$. В силу компактности
$K^S$, полученное замыкание (обозначим его через $\cB$) является компактной группой. Она
обладает разными замечательными свойствами.

\sm

{\sc Задача.} 1) Функция на $\R$ называется {\it почти периодической},
если она равномерно приближается функциями вида $\sum_{k=1}^l c_k \exp(ia_kx)$.
Равносильное определение --
функция $f(x)$ почти периодическая, если множество функций $f_t(x)=f(x+t)$
предкомпактно в равномерной топологии. Покажите, что почти периодические функции -- 
это в точности функции, полученные ограничением непрерывных функций на $\cB$ на подгруппу $R$.

\sm

2) Группа $\cB$ компактна, поэтому на ней есть мера Хаара, а значит и пространство $L^2$.
Что соответствует интегралу по мере Хаара на языке почти периодических функций? 
Что соответствует скалярному произведению в $L^2$? Функции $e_t(x)=e^{itx}$, очевидно, допускают непрерывное продолжение с 
$\R$ на $\cB$. Покажите, что эти продолжения образуют ортонормированный базис в $L^2(\cB)$.

\sm

3)$^*$ Любой гомоморфизм $\R$  в любую компактную группу  продолжается по непрерывности
до гомоморфизма $\cB\to G$.
\hfill $\lozenge$

\sm

Некоторым недостатком этого замечательного объекта является то, что ни одной точки из
$\cB\setminus \R$ построить нельзя (что является типичным результатом применения аксиомы выбора; кстати, где она была применена?).

\sm

с) Произносимое выше в этом параграфе слово <<счетный>> во многих случаях
 может быть устранено: определение обратного предела,
проконечного пополнения, счетность поля в определении группы Галуа. На формальном уровне
все выживает. Но это может приводить приводит к неприятностям типа описанной в предыдущем абзаце%
\footnote{Стоит заметит, что независимо от аксиомы выбора, почти периодические
функии являются небезынтересным элементом анализа. Известно, например, что кольца Сатурна напоминают канторовское
множество (между каждыми двумя дырками есть дырка, самые большие дырки названы в честь их
первооткрывателей, Кассини и др.). Одно из существующих объяснений \cite{AvSy} --
что эта структура связана с почти периодическим потенциальным полем, 
создаваемым спутниками Сатурна.}.

\sm

См. \cite{Bohr}, \cite{Cor}, \cite[\S 16]{Dix}.

\sm

{\bf\punct  Группа Урысона.%
\label{ss:uryson}} Опишем еще одного монстра, но в другом духе.
{\it Пространство Урысона $\bbU$} -- это полное сепарабельное метрическое пространство,
удовлетворяющее следующему свойству:

\sm

--- Пусть $A$ -- конечное метрическое пространство, а $B$ -- некоторое метрическое
пространство, полученное из $A$ добавлением одной точки. Тогда любое изометрическое вложение
$A\to\bbU$ продолжается до изометрического вложения $B\to \bbU$.

Пространство, обладающее этим свойством, единственно с точностью до изометрии
и обладает рядом удивительных свойств:

\sm

--- Если $A$ -- {\it конечное} подмножество $\bbU$, тогда любое изометричное вложение
$A\to \bbU$ продолжается до обратимой изометрии $\bbU\to\bbU$. То же утверждение
остается в силе для любых компактных подмножеств.

\sm

--- Любое сепарабельное метрическое пространство $M$ вкладывается в $\bbU$
изометрично. Более того, вложить его можно так, что любая изометрия 
$M$ будет продолжаться до изометрии $\bbU$.

\sm

Этот объект был введен в посмертной публикации П.~С.~Урысона. Доказательство
Урысона, в определенном смысле, конструктивно, в том же смысле, как конструктивно
доказательство Кантора существования трансцендентных чисел. Сейчас известно 
несколько других доказательств, но сколько-либо удобной явной конструкции
пространства Урысона не известно. 

Оказывается, однако, что пространство Урысона в некотором смысле является 
метрическим пространством <<общего положения>>. А именно, рассмотрим 
пространство всевозможных метрик $\cM$ на множестве $\N$. Для 
каждой метрики рассмотрим пополнение $\N$ по этой метрике. 
Оказывается, что для всюду плотного $G_\delta$-множества 
пространства $\cM$ полученное пополнение будет пространством Урысона.

\sm

Соответственно, возникает и огромная группа всех изометрий
пространства $\bbU$. 

См. \cite[Ch.5]{Pes}, \cite{Ury}, \cite{Ver-ury}.

\sm

{\bf\punct Мантии.} 
{\it Вопрос Ольшансткого.} Обозначим через $\B(H)$ полугруппу операторов 
в бесконечномерном гильбертовом пространстве $H$. Через $\bfU(H)$
мы обозначим группу унитарных операторов в $H$. Снабдим 
$\B(H)$ слабой операторной топологией

\sm

{\sc Задача.} a) Полугруппа $\B(H)$ метризуема, сепарабельна и компактна.

\sm

b) Умножение в $\B(H)$ раздельно непрерывно.

\sm

c) Группа $\bfU(H)$ плотна в $\B(H)$.
\hfill $\lozenge$

\sm

Рассмотрим унитарное представление $\rho$ группы $G$
в гильбертовом пространстве $H$. Рассмотрим множество операторов
$\rho(G)\subset \bfU(H)$ и замкнем его в $\B(H)$. Замыкание 
$\ov{\rho(G)}$ в $\B(H)$
 замкнуто относительно умножения. Действительно, пусть
последовательности $\rho(g_j)$, $\rho(h_j)$ сходятся
к $A$ и $B$ соответственно. Тогда
\begin{multline*}
\lim_{k\to\infty}
\Bigl[\lim_{l\to\infty}
\rho(g_k h_l)\Bigr]=
\lim_{k\to\infty}
\Bigl[\,\lim_{l\to\infty}
\rho(g_k) \rho(h_l)\Bigr]=\\=
\lim_{k\to\infty} \rho(g_k)
\Bigl[\,\lim_{l\to\infty}
 \rho(h_l)\Bigr]= \lim_{k\to\infty}( \rho(g_k) B)
 = \Bigl[\lim_{k\to\infty} \rho(g_k)\Bigr] B=AB.
\end{multline*}
Итак, $AB$ содержится в $\ov{\rho(G)}$, а поэтому
$\ov{\rho(G)}$ является компактной полугруппой.
Мы назовем ее мантией группы $G$.

\sm

Если $G$ -- полупростая группа Ли с конечным центром, то 
эта конструкция дает одноточечную компактификацию
группы, см. \cite{HM}. Вообще в случае групп Ли
ответы получаются довольно простые и
не очень интересные. Однако в случае бесконечномерных групп
одноточечная компактифицация невозможна, и вообще компактификация
вынужденным образом должна быть непохожей на группу $G$.
Вопрос в том, можно ли такие объекты описать?

\sm

{\it Примеры.} a) Любое унитарное представление
представление полной унитарной группы $\bfU(\infty)$
продолжается по непрерывности до представления 
полугруппы $\bfB(\ell_2)$ всех сжимающих операторов.
Тем самым   $\bfB(\ell_2)$  является мантией $\bfU(\infty)$.
См. \cite{Olsh-uni}, \cite[\S VIII.3]{Ner1}.

\sm

b) Рассмотрим полугруппу $\Sigma$ всех
бесконечных матриц, состоящих из нулей и единиц, причем в каждой строке
и каждом столбце стоит не более одной единицы.
Очевидно, $\Sigma\subset\B(\ell_2)$. Снабдим $\Sigma$
индуцированной слабой топологией.

\sm

{\sc Задача.} Покажите, что полугруппа $\Sigma$
компактна, а симметрическая группа $S(\infty)$
плотна в $\Sigma$. \hfill $\lozenge$

\sm

Оказывается, что любое представление полной симметрической
группы $\ov S_\infty$ продолжается по непрерывности
до представления $\Sigma$. По существу, это утверждение равносильно
теореме А.~Либермана о классификации унитарных представлений $\ov S_\infty$, см. \cite{Olsh-kiado},
\cite[\S VIII.1-2]{Ner1}.

\sm

c)  Для группы $\Ams([0,1])$ есть полугруппа $\Pol$, обладающая похожими
свойствами, она обсуждается ниже в \S12. 

Обычно ситуация более сложна

\sm

d) Рассмотрим группу $\GLO(\ell_2)$, многократно обсуждавшуюся выше.
Рассмотрим группу $G^\circ=\GLO(\ell_2\oplus \ell_2)$,
состоящую из блочных операторов
\begin{equation}
\begin{pmatrix}
a&b\\c&d
\end{pmatrix},
\label{eq:abcd}
\end{equation}
 и в ней подгруппу
$K$, состоящую из матриц вида 
$\begin{pmatrix}
1&0\\0&u
\end{pmatrix}$,
 где $u$ -- элемент полной вещественной ортогональной группы.
Рассмотрим двойные классы смежности
$K\setminus G^\circ/K$ с естественной фактор-топологией
и введем на этом множестве умножение по правилу
(\ref{eq:circ}). Оказывается, что любое унитарное представление
группы $\GLO(\ell_2)$  продолжается по непрерывности
до представления полугруппы $K\setminus G^\circ/K$,
см. \cite[п.1.11]{Ner-spheric}.

Неизвестно, вся ли это мантия, но можно показать, что
остающаяся неопределенность относится к описанию точных
 условий на матрицу
(\ref{eq:abcd}). 

\sm

e) В последнем примере бросается в глаза сходство со шлейфовыми конструкциями.
Оно не случайно. Дело в том, что  полугруппа $\B(\ell_2)$
содержит нулевой элемент, который в унитарных
представлениях $\bfO(\infty)$ действует как проектор
на подпространство $\bfO(\infty)$-неподвижных векторов
(все неприводимые представления $\bfO(\infty)$
реализуются в тензорах над стандартным представлением,
а любое унитарное представление является прямой суммой
неприводимых, отсюда сделанное замечание легко следует).
Поэтому проекторы $P[\alpha]$ из формулы (\ref{eq:ov-rho})
содержатся в мантии группы. А поэтому в мантии
содержатся и операторы 
$$
\wt\rho(g)=P[\alpha]\,\rho(g)\,P[\beta]
,$$  
которые почти не отличаются от операторов
$\ov\rho(g)$. А именно, если записать оператор
$\wt\rho(g)$  как блочный оператор
$$
H[\beta]\oplus H[\beta]^\bot \,\to\, H[\alpha]\oplus H[\alpha]^\bot,
$$ 
мы получим матрицу
$$
\begin{pmatrix}
\ov \rho (g)&0\\ 0 &\rho(g)
\end{pmatrix}.
$$
Выходит, что шлейф фактически лежит в мантии, а мантия является
<<бесконечным элементом>> шлейфа. В существующих на сегодняшний
день теориях представлений бесконечномерных классических групп
и бесконечных симметрических групп это так и есть
(а описания шлейфов, в общем, известны). Но надо
иметь в виду, что конструкция с умножением двойных
классов смежности работает не для всех бесконечномерных
групп, в то время как для определения мантии нужны лишь унитарные
представления.

\sm

{\it Универсализация мантии.}
В существующих теориях представлений
бесконечномерных групп  замыкания в большинстве случаев 
или описаны или в принципе понятны (об одном исключении, см. следующий пункт)
Как зависят эти объекты зависят  от представления, и существует ли универсальная мантия?

В принципе можно предложить такую универсальную конструкцию.
Пусть $\rho_\alpha$ -- набор всех попарно неэквивалентных унитарных представлений
группы $G$,
пусть $H_\alpha$ -- соответствующие гильбертовы пространства.
Рассмотрим все вложения 
$$\rho_\alpha:G\to \bfU(H_\alpha)\to \B(H_\alpha)$$
Возьмем диагональное вложение $G$ в (тихоновское) прямое произведение всех $\B(H_\alpha)$,
и замкнем образ. Получается компактная полугруппа, которая действует в любом
унитарном представлении группы $G$.

Отметим, что эта постановка
вопроса не лишена опасностей, что показывает уже случай $G=\R$ (см. п.\ref{ss:bohr}),  
 в действительности их больше, и на данный момент универсальную мантию
 лучше понимать как эвристический объект.
 
 \sm

См. \cite{Olsh-semi}, \cite{Ner1}.

\sm

{\bf\punct О замыканиях группы диффеоморфизмов окружности.}  Вернемся к группе $\Diff_+(S^1)$
диффеоморфизмов окружности, сохраняющих ориентацию. Ее алгебра Ли -- алгебра
$\Vect(S^1)$ векторных
полей на окружности. В ее комплексификации $\Vect(S^1)_\C$ удобно взять  базис
$L_n=e^{in\phi}/i\partial\phi$, соотношения коммутации принимают вид
$$
[L_n,L_m]=(m-n)L_{m+n}.
$$
С точки зрения теории представлений более интересна {\it алгебра Вирасоро}
$\Vir$. К базису
добавляются еще один элемент $\zeta$, a соотношения коммутации меняются на
\begin{equation}
[L_n,\zeta]=0,\qquad [L_n,L_m]=(m-n)L_{m+n}+\frac{n^3-n}{12} \zeta.
\label{eq:vir}
\end{equation}
Тогда $\zeta$ становится центральным элементом, а фактор-алгебра
$\Vir/\C\zeta$ совпадает с $\Vect(S^1)_\C$.
Важную роль играют представления алгебры Вирасоро со старшим весом,
т.е., представления, в которых есть порождающий вектор
$v$, такой, что
\begin{align*}
L_n v=0\,\, \text{при $n>0$};
\\
L_0 v=hv,\qquad \zeta v= cv
\end{align*}
для некоторых постоянных $h$, $c$. Такие представления 
интегрируются до проективных представлений группы $\Diff_+(S^1)$,
$$
\rho(g_1)\rho(g_2)= \sigma(g_1,g_2)\rho(g_1g_2).
$$   

\sm

A) {\it Тривиальные представления.} Реализуем окружность как  $\R/\Z$ или
как отрезок $[0,1]$ со склеенными концами. Фиксируем $s\in \R$ и рассмотрим
тривиальное представление $\rho_s$ группы диффеоморфизмов окружности в пространстве  $L^2[0,1]$,
заданное 
 формулой
\begin{equation}
\rho_s(q)f(x)= f(q(x)) \,q'(x)^{1/2+is}.
\label{eq:trivial}
\end{equation}
Понятно, что в слабом замыкании $\ov{\rho_s(\Diff)}$ лежат похожие отображения
с более общими монотонными $q(x)$. Мы опустим обсуждение этого и заметим,
что замыкание содержит много операторов типа умножения на функцию,
$$
A f(x)= f(x) a(x).
$$
Они получаются пределами последовательностей
$\rho_s(q_n)$, где $q_n(x)$ стремится $x$ в топологии
равномерной сходимости, а $q_n'(x)$  быстро осциллирует
(вроде $q_n(x)=x+(\sin nx)/n$).
Опишем возможные предельные операторы.
Для диффеоморфизма $q$ рассмотрим отображение $[0,1]\to [0,1]\times (0,\infty)$,
заданное формулой
$x\mapsto (x,q'(x))$. Рассмотрим образ $\mu[q]$ меры Лебега при этом отображении.
Это вероятностная мера на $[0,1]\times [0,\infty]$,
причем ее проекция на$[0,1$] совпадает с мерой Лебега.
\begin{figure}
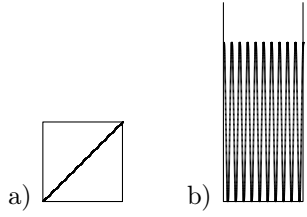

a)  \epsfbox{haar.4}
\qquad
b) \epsfbox{haar.3}
\caption{a) Функция $q(x)=x+ \frac 1n \sin nx$.\, b) $q'(x)$.}
\end{figure}
Такие меры удобно задавать функциями $x\mapsto \nu_x$ из $[0,1]$ в полугруппу
вероятностых мер на $\R_+^\times$.

\sm

{\sc Задача.} a) Покажите, что для любой функции $x\mapsto\nu_x$, удовлетворяющей условию
$$
\int_{\R_+^\times} t\,d\nu_x(t)=1,
$$
существует последовательность $q_n$ диффеоморфизмов, такая, что
$q_n(x)$ cходится к $x$, а последовательность мер
$\mu[q_n]$ cходится слабо к $\nu$.

\sm

b) Соответствующие операторы $\rho_s(q_n)$
слабо сходятся к оператору умножения на функцию,
$$
A_{s,\nu} f(x)= f(x) \cdot \int_0^\infty t^{1/2+is} d\nu_x(t).
$$

c) Умножение на какие функции можно получить таким способом?

\sm

d) Более утомительный вопрос: опишите слабое замыкание группы диффеоморфизмов.
\hfill $\lozenge$

\sm

B) {\it Полугруппа трубок.} 
Пусть $\ov\C=\C\cup\infty$ обозначает сферу Римана,
Пусть $D_+$, $D_-\subset \ov\C$ -- соответственно диски
$|z|\le 1$, $|z|\ge 1$, пусть $D_\pm^\circ$ -- 
их внутренности, пусть $S^1$ обозначает круг $|z|=1$.
 Отображение $f:D_\pm^\circ\to \ov \C$
называется {\it однолистным}, если $f$ голоморфно и
является вложением.
 Мы скажем, что отображение
$D_\pm\to \ov\C$ однолистно вплоть до границы,
если оно однолистно во внутренности и  продолжается
до гладкого вложения всего круга с невырожденным
дифференциалом.

Рассмотрим сохраняющий ориентацию диффеоморфизм 
$q:S^1\to S^1$. Склеим $D_+$ c $D_-$, отождествляя точки
$e^{i\phi}\in D_+$  с точками $ q(e^{i\phi})\in D_-$.
После склейки получится двумерное многообразие, 
гомеоморфное сфере. Можно показать, что комплексные структуры
тоже склеиваются, и мы таким образом получаем сферу Римана
(эта конструкция называется {\it сваркой}).

\begin{figure}
$$
\epsfbox{welding.2}\quad
\epsfbox{welding.1}
$$
a)\qquad\qquad\qquad\qquad\qquad\qquad\qquad\qquad\qquad b)
\caption{a) Сварка.   b) Элемент полугруппы $\Gamma$.%
	\label{fig:tubes}}
\end{figure}

\sm

{\sc Замечание.} Сварка определена при довольно
слабых ограничениях на отображение $q$. Оно должно быть
гомеоморфизмом и должно удовлетворять условию
квазисимметричности. А именно, функция $f$ вещественной
переменной называется {\it квазисимметричной}, если
$$
\frac {|x-y|}{|x-z|}\le t\quad\text{влечет} \quad
\frac {|f(x)-f(y)|}{|f(x)-f(z)|}\le \eta(t),
$$
где $\eta(t)$ -- некоторое возрастающее непрерывное отображение 
$[0,\infty]\to [0,\infty]$, $\eta(0)=0$. В частности, липшицевы
отображения квазисимметричны.
\hfill $\lozenge$

\sm

Итак, у нас получилась сфера Римана и два однолистных вплоть до 
границы отображения $D_+\to \ov\C$, $D_-\to \ov\C$.

Обратно, пусть даны  однолистные вплоть до границы отображения
$p_+:D_+\to \ov\C$, $p_-:D_-\to \ov\C$, причем их образы покрывают
сферу, а образы внутренностей $p_+(D_+^\circ)$, $p_-(D_-^\circ)$
 не пересекаются. Тогда мы имеем отображение $q:=p_+^{-1}\circ p_-$,
 определенное на окружности $|z|=1$. 
 
 \sm

Определим {\it полугруппу трубок} $\Gamma$. Ее элемент -- это следующий набор данных:

\sm

--- сфера Римана $\ov\C$;

\sm 

--- два однолистных вплоть до границы   отображения $p_+:D_+ \to \ov\C$,
$p_-:D_- \to \ov\C$
с непересекающимися образами.

\sm

Два элемента $(\C,p_+,p_-)\in \Gamma$, $(\C,p'_+,p'_-)\in \Gamma$
считаются одинаковыми, если существует биголомормофное отображение
$w:\ov\C\to\ov\C$, такое, что $p'_\pm=p_\pm\circ w$. 

Определим произведение элементов $(\ov\C,p_+,p_-)$, $(\ov\C,r_+,r_-)\in\Gamma$.
Для этого рассмотрим области $\ov\C\setminus p_-(D_-^0)$, $\ov\C\setminus r_+(D_+^0)$,
это два диска с параметризованными границей. 
Склеим их, отождествляя пары точек $p_-(e^{i\phi})$ и $q_+(e^{i\phi})$,
и получим одномерное комплексное многообразие $W$, гомеоморфное
сфере, т.е., опять сферу Римана. Кроме того, оно помнит
отображения $p_+$ и $q_-$, то есть мы опять получили элемент полугруппы
$\Gamma$.

Теперь определим сходимость последовательности диффеоморфизмов
$g^{(j)}$.  Применим к ним сварочную конструкцию и получим последовательность
$(\ov\C,p_+^{(j)}), p_-^{(j)})$. 
Мы скажем, что последовательность диффеоморфизмов сходится к элементу
$(\ov\C,q_+,q_-)$
полугруппы $\Gamma$,  если выбрав подходящим образом координаты
на каждой из сфер Римана (они могут различаться на дробно-линейное преобразование),
 мы можем получить равномерную на компактных подмножествах
 $D_\pm^\circ$ сходимость $p_+^{(j)}\to q_+$, $p_-^{(j)}\to q_-$.
 
 \begin{figure}
 $$	\epsfbox{welding.3}
 	\qquad\qquad \epsfbox{welding.4}$$
 	\caption{Последовательность сварок и ее предел в $\Gamma$.
 		\label{fig:carath}}
 \end{figure}
 
 Пример семейства диффеоморфизмов, сходящегося к элементу полугруппы $\Gamma$
 изображен на Рис. \ref{fig:carath}. Сфера Римана изображена как плоскость,
 ноль находится в центре окружностей,
 $\infty$ находится на бесконечности. Жорданов контур разделяет сферу на две области.
 Мы берем отображаем круг $D_-$ во внутреннюю область, полагая $p_+(0)=0$,
 $p_+'(0)>0$. Круг $D_-$ мы отображаем на внешнюю область, полагая
 $p_-(\infty)=\infty$, и коэффициент Лорана $c_{-1}$ в разложении 
 $$
 p_-(z)=c_{-1}z+c_0+c_1z^{-1}+\dots
 $$
 положителен. Далее устремляем частоту зубчатки  к бесконечности,
  применяем теорему Каратеодори 
 о сходимости однолистных функций, см. \cite[\S II.5]{Gol} и получаем картинку
 справа.
 
 Любое представление группы $\Diff_+(S^1)$ со старшим весом продолжается
 по голоморфности до представления полугруппы $\Gamma$
 (см. \cite{Ner-Gamma}). То, что 
 это продолжение непрерывно в смысле введенной сходимости, строго не доказано.
 Однако в упомянутой статье есть явные задания таких представлений интегральными операторами,
 и, по крайней мере, поточеченая сходимость ядер имеет место.

\sm

{\it Интерполяция.} Как будто, две обсуждавшихся выше полугруппы,
связанных с группой $\Diff_+(S^1)$, выглядят совершенно различно.
В п.\ref{ss:diff1} было описана конструкция семейства представлений
группы
$\Diff_+(S^1)$, зависящая от параметра $s\in(0,1)$.
Для этого семейства описание мантий неизвестно. Однако у  
этого семейства есть понятные пределы при $s\to 0$ и $s\to 1$.
При $s=0$ мы получаем сумму симметрических степеней тривиального представления
$\rho_0$ (см. \ref{eq:trivial}), это ясно из
разложения (\ref{eq:fock-orthogonal}), обсуждаемого ниже.
При $s=1$ мы получаем представление вида $W\otimes W^*$, 
где $W$ -- некоторое представление группы  $\Diff_+(S^1)$
со старшим весом (обсуждение этого -- за пределами настоящих записок).
Так или иначе, мы получаем  на одном конце полугруппу первого типа, на другом -- второго.
О том, что происходит в промежутке -- неизвестно.  Те же самые вопросы 
возникают в связи с разными другими сериями представлений
$\Diff_+(S^1)$, описанных в \cite[\S IX.6]{Ner1}.

\sm

См. \cite{Ner-Gamma}, \cite[гл. VII, \S IX.6]{Ner1}, \cite{Ner-cantor}.

\sm

{\bf \punct Большие группы и инвариантные меры.}
В связи с многочисленными примерами <<больших>> топологических групп встает вопрос об аналогах мер Хаара. На этот вопрос есть
два ответа, отрицательный и положительный.

\sm

{\it Отрицательный ответ} был получен А.~Вейлем и звучит так.
Любая топологическая группа, на которой существует левоинвариантная мера, является плотной
подгруппой в локально компактной группе \cite[Добавл. 1]{Wei}.

\sm

{\it Положительный ответ} состоит в том, что, несмотря на теорему Вейля, есть много
интересных мер, связанных с бесконечномерным
группами и похожих на меру Хаара, но не удовлетворяющих каким-нибудь формальным условиям.

\sm

1. Ниже в \S\S10-11 мы обсуждаем аналоги меры Хаара для бесконечной симметрической группы
и для бесконечномерной унитарной группы. Эти объекты  получаются обратными
пределами групп%
\footnote{Но не в категории групп.}, они не являются группами,
зато бесконечная симметрическая и унитарная группы действуют на них левыми и правыми
умножениями, при этом получаются интересные регулярные представление.  В принципе,
эти объекты можно рассматривать как аналоги грассманианов, что в любом случае верно,
потому что у грассманианов есть похожие пределы (D.~Pickrell). Интересные
инвариантные меры известны и для грассманианов над конечными и $p$-адическими полями
(см.\cite{Ner-p}, \cite{Ner-p-adic}).

\sm

2. С другой стороны, есть довольно много мер на топологических группах, квазиинвариантных
относительно всюду плотной подгруппы. Простой пример дает теорема Камерона--Мартина,
см. ниже п.\ref{ss:cameron}. Есть также следующая конструкция Шавгулидзе
\cite{Shav}. Рассмотрим
группу $\Diff^2_{+}[0,1]$, диффеоморфизмов $\phi$ отрезка $[0,1]$ гладкости $C^2$, у которых
$\phi(0)=0$, $\phi''(0)=1$. Рассмотрим отображение, 
$$
F:\phi\mapsto\frac {\phi''}{\phi'}
,$$
которое ставит в соответствие диффеоморфизму 
непрерывную функцию $F$ с $F(0)=0$. Введем на пространстве таких функций меру Винера
(см., например, ниже п.\ref{ss:functional-spaces}). Перенесем эту меру на $\Diff^2_{.}[0,1]$
посредством этого отображения. Оказывается, что получается мера на 
$\Diff^2_{.}[0,1]$, квазиинвариантная относительно преобразований
$\phi\mapsto \alpha\circ\phi$, где $\alpha$-диффеоморфизм гладкости
$C^3$. 

Семейства мер, сосредоточенных на диффеоморфизмах окружности
разных гёльдеровских классов, квазиинвариантных относительно 
групп более гладких диффеоморфизмов, обсуждаются в \cite{Ner-MIRAN}.

Впрочем, продолжить такого рода конструкции до интересного гармонического анализа
на группах пока не удалось.

\sm

3. Есть много конструкций мер, инвариантных или квазиинвариантных относительно
действия  бесконечномерных групп. Мы в \S\S 8-9
рассматриваем два примера, наиболее простых и общезначимых. См. также ссылки
в \S\S 8-11, старый обзор \cite{Ner-MIRAN} и 
\cite{Ism}, \cite{Pick-mac}, \cite{Ald}, \cite{Ner-cantor}, \cite{Ism-su2},
\cite{Alb}.

\sm

{\sc Задача.} Докажите, что на гильбертовом пространстве
не существует конечной борелевской меры, инвариантной
относительно группы всех вращений (такая мера, во всех
отношениях хорошая и вполне полезная, была, однако, 
построена И.Сигалом, см. ниже \S 8; впрочем, она сосредоточена
вне гильбертова пространства).\hfill $\lozenge$

\section{Преобразование Фурье на компактных группах.}

\COUNTERS

Этот параграф отличается от прочих свои абстракционизмом. 
Наша цель -- доказать  теорему
\ref{th:fourier}. Формулировка и, особенно, доказательство включают
представления компактных групп. Мы формально манипулируем с ними, не обсуждая, что 
это значит в конкретных случаях. Для конечных групп текст остается осмысленным,
некоторые доказательства в этом случае существенно упрощаются.

\sm

{\bf \punct Унитарные представления.} Пусть $G$ -- топологическая группа, а $H$
-- гильбертово пространство (бесконечномерное или конечномерное; {\it мы  считаем гильбертово пространство
 сепарабельным
по определению}). 
Унитарное представление $\rho$ группы $G$ это непрерывный гомоморфизм
из $G$ в унитарную группу пространства $H$, снабженную слабой операторной топологией
(см.  п.\ref{ss:unitary-infty}).

Условие непрерывности можно сформулировать так: для любых векторов
$v$, $w\in H$ функция $\la \rho(g)v, w\ra$ ({\it матричный элемент}) непрерывна
на $G$.

\sm

{\sc Задача.}   a) Покажите, что непрерывность достаточно проверить в единице группы.

\sm

b) Пусть $e_j$ - ортонормированный базис в $H$. Покажите, что непрерывность матричных элементов
$\la\rho(g)e_k, e_l\ra$ влечет непрерывность представления  $\rho$.
\hfill $\lozenge$

\sm

{\sc Замечание.}  На самом деле непрерывность конкретных представлений проверять не надо. Оказывается,
что из измеримости матричных элементов следует их непрерывность (см. \cite[\S29, следствие 2]{Nai}).
Если бы мы построили
разрывное представление, то вместе с ним построили бы и неизмеримую функцию.
\hfill $\lozenge$

 \sm
 
 {\bf\punct Примеры унитарных представлений.} Приведем несколько примеров (точнее, массовых
 конструкций) унитарных представлений. 
 
 \sm
 
 1) Пусть группа $G$ действует (справа) на пространстве $M$ с мерой $\mu$  преобразованиями, сохраняющими
 меру%
 \footnote{Формально от действия необходимо потребовать непрерывности, т.е. для любых
 измеримых подмножеств $A$, $B\subset M$, функция $g\mapsto \mu(Ag\cap B)$ должна быть непрерывна на группе
 $G$. Равносильно, гомоморфизм $G\to\Ams(M)$ должен быть непрерывен (см. п.\ref{ss:Ams}).
 Фактически следить за этим обычно не надо, см. замечание в конце предыдущего пункта.}. 
 Тогда  $G$ действует в $L^2(M,\mu)$ унитарными преобразованиями
 $$
 T(g) f(x)=f(xg).
 $$
 
 2) Пусть группа $G$ действует (справа) на пространстве $M$ с мерой $\mu$  преобразованиями,
 оставляющими меру квазиинвариантной. Обозначим через $g'(x)$ производную Радона--Никодима преобразования%
 \footnote{См. п.\ref{ss:Ams}, если речь идет о действии на $\R^n$, то это якобиан.} $g$ в точке $x$.
 Для любого $s\in \R$ определено унитарное
 представление $G$ в $L^2(M,\mu)$ по формуле
 $$
 T_s(g)f(x)=f(xg) g'(x)^{1/2+is}.
 $$
 
 {\scЗадача.}
 Проверьте унитарность этих операторов. Проверьте, что мы получили представление группы $G$.
 \hfill $\lozenge$
 
 \sm
 
 3) В случае локально-компактных групп мы, в частности, можем рассматривать $L^2$
 на группе и всевозможных однородных пространствах.

 Представление локально компактной группы в $L^2(G)$ левыми сдвигами
 $$
 L_h f(g)=f(h^{-1}g)
 $$
 называется {\it леворегулярным представлением}. Аналогично определяется
 {\it праворегулярное представление}, оно действует операторами
 $$R_h f(g)=f(gh).$$
 
 \sm
 
 4) Пусть, $G$, по прежнему, действует на пространстве с мерой, преобразованиями,
 оставляющими меру квазиинвариантной. Пусть $c(x,g)$ -- комплекснозначная
  функция на $G\times X$, удовлетворяющая тождеству 
  (цепному правилу)
  \begin{equation}
  c(x,g_1g_2)=c(xg_1,g_2) \, c(x,g_1)
  \label{eq:cocycle}
  \end{equation} 
и $|c(g,x)|=1$. Тогда формула
$$
T(g)f(x)=f(xg)\, g'(x)^{1/2}\, c(x,g)
$$
задает унитарное представление группы (проверьте это).  
 
\sm 
 
{\bf \punct Несколько определений.} Пусть $\rho$ -- унитарное представление
группы $G$ в пространстве $H$. Пусть $V\subset H$ -- замкнутое 
подпространство, инвариантное относительно всех операторов $\rho(g)$.
 Представление $G$ в $V$ называется {\it подпредставлением}
представления $\rho$. Представление называется {\it неприводимым}, если у него
нет $G$-инвариантных подпространств, отличных от $0$ и $H$. В противном
случае оно называется {\it приводимым}. 

Очевидно, что для любого подпредставления $V$ его ортогональное дополнение
$V^\bot$ тоже
является подпредставлением. 

Группа $G$ действует также в пространстве $H'$, двойственном к
$H$ по формуле
$$
\rho'(g)\ell (v)= \ell(\rho(g^{-1}) v), \qquad v\in H,\, \ell\in H'
.$$
Это представление называется {\it двойственным}, или {\it контраградиентным}, 
или {\it комплексно сопряженным}. По теореме Рисса--Фишера $H'$ отождествляется
с $H$, но посредством антилинейного оператора. Поэтому двойственное представление
можно определить так: мы меняем структуру линейного пространства в $H$,
положив, что умножение на $i$ есть умножение на $-i$, а операторы представления
оставляем теми же. Если же мы задаем представление посредством матриц в некотором
 ортонормальном
базисе, то переход к двойственному представлению соответствует поэлементному
 комплексному сопряжению матриц%
 \footnote{Операция поэлементного сопряжения матричных элементов оператора зависит от выбора базиса.
 	Однако, получающееся при этом представление будет одним и тем же, с точностью до эквивалентности.}.
 
 Если $\rho_1$, $\rho_2$ -- представления $G$ в $H_1$ и $H_2$ соответственно, то
 операторы $\rho_1(g)\oplus\rho_2(g)$ определяют представление
 $G$ в $H_1\oplus H_2$. Оно называется {\it прямой суммой} двух представлений.
Аналогично определяются конечные и счетные прямые суммы. 
Операторы $\rho_1(g)\otimes\rho_2(g)$ действуют в тензорном произведении $H_1\otimes H_2$.
Это определяет {\it тензорное произведение} представлений группы $G$.

Пусть $\rho_1$, $\rho_2$ -- представления групп $G_1$ и $G_2$. Тогда
мы имеем представление $\rho_1(g_1)\otimes \rho_2(g_2)$ группы $G_1\times G_2$.

\sm

Представления  $\rho_1$, $\rho_2$ группы $G$ в гильбертовых пространствх
$H_1$, $H_2$ называются {\it эквивалентными}, если между ними есть
унитарный  оператор $J:H_1\to H_2$, такой, что
$$
\rho_2(g) =J\rho_1(g) J^{-1}.
$$

Определим морфизм представлений.  Пусть $\rho_1$, $\rho_2$ -- представления  группы $G$
 в $H_1$ и $H_2$ соответственно. {\it Сплетающий оператор} $J:H_1\to H_2$ -- это ограниченный 
линейный оператор, удовлетворяющий тождеству
$$
\rho_2(g) J=J\rho_1(g).
$$

{\bf \punct Лемма Шура.}

\begin{lemma}
Пусть $H_1\to H_2$ -- сплетающий оператор. Тогда $\ker J$ -- подпредставление
в $H_1$, а замыкание $\im J$ -- подпредставление в $H_2$.
\end{lemma}
 
{\sc Доказательство.} Пусть $Jv=0$. Тогда 
$$
J\rho_1(g) v=\rho_2(g)Jv=0
,$$ 
 т.е., $\rho_1(g) v$ тоже принадлежит $\ker J$.
 
 Пусть $w\in \Im J$, $w=Jv$. Тогда
 $$
 \rho(g_2)w= \rho(g_2)Jv=J\rho_1(g) v\in \Im J.
 $$
Поэтому подпространство $\im J$ инвариантно, а поэтому
инвариантно и его замыкание. \hfill $\square$

\begin{theorem}
{\rm (лемма Шура)}

{\rm a)} Пусть $\rho_1$, $\rho_2$ -- неэквивалентные неприводимые унитарные
 представления группы $G$.
Тогда любой сплетающий оператор между ними равен 0.

\sm

{\rm b)}  Пусть $\rho$ -- неприводимое унитарное представление. Тогда любой оператор,
сплетающий его с собой, имеет вид $\lambda\cdot 1$, где $\lambda$ -- скаляр.
\end{theorem} 

{\sc Доказательство.} Пусть $\rho$ -- представление, $J$ --
{\it самосопряженный} оператор, сплетающий
$\rho$ с самим собой. Возьмем подмножество $\Omega$ в  спектре $J$
и соответствующий спектральный проектор $P_\Omega$ (см. \cite[\S VII.3]{RS1}). Он коммутирует
с любым оператором, коммутирующим с $J$, а поэтому является сплетающим
для $\rho$. Образ $P_\Omega$ является подпредставлением. Поэтому $P_\Omega$
равен нулю или единице. Поэтому спектр $J$ состоит из одной точки, т.е.,
$J=\lambda\cdot 1$.

Пусть теперь $S$ -- произвольный оператор, сплетающий $\rho$  собой.
Тогда операторы $S+S^*$ и $i(S-S^*)$ -- тоже сплетающие, они самосопряжены,
а поэтому являются скалярными. Поэтому и $S$ -- скалярный оператор.
Утверждение b) доказано. 

Докажем a). Оператор $J^*$ является сплетающим, а поэтому  $J^*J$  сплетает
$\rho_1$  с самим собой.  Так как $J^*J$ самосопряжен, то $J^*J=s\cdot 1$, где $s$
--
неотрицательный скаляр.
Аналогично $JJ^*=t\cdot 1$. Если $s$ или $t=0$, то $J=0$. Далее
$$
s\cdot J^*=J^*J J^*=t\cdot J^*
.$$
Поэтому $s=t$, а, следовательно, $s^{-1/2} J$ -- унитарный оператор.
Но тогда $\rho_1$, $\rho_2$ эквивалентны. 
\hfill $\square$

\sm

{\sc Задача.} Где в доказательстве была использована унитарность представления?
\hfill $\lozenge$

\sm

{\sc Задача.} Пусть $\rho_1$, $\rho_2$, $\rho_3$ -- попарно неэквивалентные представления.
Опишите все сплетающие операторы из $\rho_1\oplus \rho_1$ в следующие представления: $\rho_2\oplus \rho_3$; 
в $\rho_1\oplus \rho_1$; в
$\rho_1\oplus\rho_2\oplus \rho_3$; $\rho_1\oplus\rho_1\oplus \rho_1$.
\hfill $\lozenge$

\sm

{\bf\punct Вложение представлений в пространство функций на группе.%
\label{ss:J}}
Рассмотрим  унитарное представление группы $G$ в пространстве $H$.
Фиксируем единичный вектор $\xi\in H$. Для любого $v\in H$ рассмотрим функцию $F_v(g)$
на $G$, заданную формулой 
$$
F_v(g)=\la v, \rho(g^{-1}) \xi\ra.
$$
Мы получили линейный оператор $J$ из $H$ в пространство непрерывных
функций на $G$, именно $Jv(g)=F_v(g)$.

\begin{proposition} Оператор $J:H\to C(G)$ непрерывен
и коммутирует с действием группы
$$J (\rho(h) v) (g)= Jv(gh).$$ 
\end{proposition}

{\sc Доказательство.}
$$
J(\rho(h) v) (g)= \la \rho(h) v, \rho(g^{-1}) \xi\ra=
\la  v, \rho(h^{-1}) \rho(g^{-1}) \xi\ra=
\la  v, \rho\bigl((gh)^{-1}\bigr) \xi\ra
$$
Если $\|v-w\|\le \epsilon$, то 
$$|Jv(g)-Jw(g)|= \la v-w,\rho(g^{-1}) \xi\ra  <\epsilon,$$
что влечет непрерывность.
\hfill $\square$ 

\sm

В пространстве непрерывных функций на произвольной группе скалярного произведения нет.
Однако в случае компактных групп $C(G)\subset L^2(G)$,
поэтому мы получаем непрерывное отображение $H\to L^2(G)$. 

В случае
неприводимого представления $\rho$, мы, очевидно, получаем вложение $H\to L^2(G)$.
Далее, оператор $J^*J$ является сплетающим оператором в $H$, и,
следовательно, скалярен. Мы получаем следующее утверждение.

\begin{proposition}
Любое неприводимое представление компактной группы изометрично вкладывается в регулярное
представление.
\end{proposition}

\sm

{\bf \punct Соотношения ортогональности.}
{\it Ниже группа $G$ будет предполагаться компактной}. 
Пусть $\rho$ -- унитарное представление группы $G$ в пространстве $H$.
Для любых двух векторов $v$, $w\in H$ определим функцию на группе
({\it матричный элемент})
по формуле
$$
r_{v,w}(g)=\la \rho(g) v,w\ra
.
$$

\begin{theorem} Пусть $G$ -- компактная группа. Пусть мера Хаара
всей группы $G$ выбрана  вероятностной.

\sm

{\rm a)} Если $\rho_1$ и $\rho_2$ неэквивалентны,
то любой матричный элемент представления $\rho_1$ ортогонален
в $L^2(G)$ любому матричному элементу $\rho_2$.

\sm

{\rm b)} Матричные элементы одного неприводимого представления $\rho$ удовлетворяют соотношению
\begin{equation}
\int_G
r_{v_1,w_1}(g)\,\ov {r_{v_2,w_2}(g)}\, dg
=\frac{1}{\dim\rho}\la v_1, v_2 \ra\cdot \ov{\la w_1, w_2 \ra}
\label{eq:soot-ortog}
\end{equation}

 {\rm c)} Любое неприводимое представление компактной группы конечномерно.
\end{theorem}

{\sc Доказательство.}
Пусть даны два неприводимых представления
$\rho_1$, $\rho_2$. Рассмотрим скалярное произведение их
матричных элементов
\begin{equation}
S(v_1,w_1,v_2,w_2):=
\int_G
r^1_{v_1,w_1}(g)\,\ov {r^2_{v_2,w_2}(g)}\, dg
.
\end{equation}
 Выражение $S(v_1,w_1,v_2,w_2)$  линейно по $v_1$, антилинейно по $w_1$,
антилинейно по $v_2$, линейно по $w_2$. Далее,
\begin{multline*}
\int_G
r_{v_1,w_1}(g)\,\ov {r_{v_2,w_2}(g)}\, dg=
\int_G
r_{v_1,w_1}(gh)\,\ov {r_{v_2,w_2}(gh)}\, dg
=\\=
\int_G
\la \rho(gh) v,w\ra\cdot\ov{\la \rho(gh) v,w\ra}\,dg=
S(\rho_1(h)v_1,w_1,\rho_2(h)v_2,w_2)
.
\end{multline*}
Таким образом,
\begin{equation}
S(v_1,w_1,v_2,w_2)=S(\rho_1(h)v_1,w_1,\rho_2(h)v_2,w_2)
\label{eq:S1}
\end{equation}
Аналогично, рассматривая левые сдвиги, мы получаем
\begin{equation}
S(v_1,w_1,v_2,w_2)=S(v_1,\rho_1(h')w_1,v_2,\rho_2(h')w_2)
\label{eq:S2}
\end{equation}
В первом равенстве (\ref{eq:S1}) фиксируем $w_1$, $w_2$. Получается полуторалинейная форма
от $v_1$, $v_2$. Любая такая форма имеет вид $\la v_1,A v_2\ra$
для некоторого ограниченного оператора $A$. При этом
$$\la  v_1,A v_2\ra=\la \rho_1(h) v_1,A \rho_2(h) v_2\ra=
 \la  v_1,\rho_1(h)^{-1} A \rho_2(h) v_2\ra $$
 Поэтому, оператор $A$ -- сплетаюший.
 
 Если $\rho_1$ и $\rho_2$ не эквивалентны,
 то $A=0$ и $S=0$, что доказывает утверждение a).

Если $\rho_1$ и $\rho_2$ эквивалентны, то  $A$ -- скалярный оператор
$$
A=c(w_1,w_2)\cdot 1,
$$ причем скаляр $c(w_1,w_2)$ зависит от $w_1$ антилинейно, а 
от $w_2$ линейно. Поэтому $c(w_1,w_2)$ имеет вид
$$
c(w_1,w_2)=\la w_2, Bw_1\ra.
,$$
где $B$ -- линейный оператор. В силу (\ref{eq:S2}) оператор 
 $B$ -- тоже сплетающий, а поэтому скалярен.
Итак,
\begin{equation}
S=\sigma\cdot\la v_1,v_2\ra\cdot \ov{\la w_1,w_2\ra}
,
\label{eq:S=}
\end{equation}
 где $\sigma$ -- некоторая положительная константа. Сейчас мы ее вычислим.
 
 \sm

Выберем в пространстве представления ортонормированную  систему $e_1$, \dots, $e_N$.
Пусть $v$ -- единичный вектор. Тогда 
$$
\sum_{j\le N} |\la\rho(g) v,e_j\ra|^2\le 1.
$$
Проинтегрируем это неравенство по группе. Используя (\ref{eq:S=}), мы получаем
$$
\int_G |\la\rho(g) v,e_j\ra|^2\,dg=\sigma\cdot \la \rho(g) v,\rho(g) v\ra\cdot \ov{\la e_j,e_j\ra}=\sigma\cdot 1\cdot 1,
$$
и
получаем
$N\cdot \sigma \le 1$. Поэтому $N \le \sigma^{-1}$, т.е., $\rho$ конечномерно.
Утверждение c) доказано.

\sm

 Теперь берем 
ортонормированный базис $e_j$. Интегрируя равенство
$$\sum |\la\rho(g) v,e_j\ra|^2=1$$
по группе, получаем
$\dim \rho\cdot \sigma=1$, что доказывает утверждение b).
\hfill $\square$

\sm

\begin{theorem}
\label{eq:matrix-total}
Матричные элементы неприводимых представлений компактной группы
$G$ образуют тотальную систему
в $L^2(G)$.
\end{theorem}

Мы докажем это утверждение в общем случае чуть позже в п.\ref{ss:proof-total}.

\sm

{\sc Доказательство для $G=\U(n)$.}
Обозначим через $\cL$ линейную оболочку всех матричных элементов неприводимых представлений.
 Легко
видеть, что $\cL$ замкнута относительно умножения. Действительно,
произведение двух матричных элементов
$$
\la \rho_1(g)v_1,w_1\ra\cdot \la \rho_2(g)v_2,w_2\ra=
\bigl\la (\rho_1(g)\otimes \rho_2(g))\,(v_1\otimes v_2),w_1\otimes w_2\bigr\ra
$$ 
является матричным элементом тензорного произведения. Раскладывая тензорное произведение
на неприводимые представления, мы получаем, что его матричные элементы лежат в $\cL$.

Далее возьмем тавтологическое представление $g\mapsto g$ группы $\U(n)$ и сопряженное
представление $g\mapsto \ov g$. Их матричные элементы разделяют точки $\U(n)$, а поэтому,
по теореме Стоуна--Вейерштрасса породенная ими  алгебра функций плотна в $C(G)$. Поэтому она плотна и в
$L^2(G)$.
\hfill $\square$

\sm

Отметим, что тот же довод остается в силе для групп $\O(n)$ и $\Sp(n)$,
и вообще для любой группы, имеющей точное%
\footnote{Представление называется {\it точным}, если для любого
$g\ne 1$ оператор $\rho(g)$ отличен от единицы.}
конечномерное представление (или хотя бы набор представлений $\rho_j$, разделяющих точки группы,
т.е., для любых $g_1\ne g_2$ существует $\rho_\alpha$, такое что $\rho_\alpha(g_1)\ne \rho_\alpha(g_2)$).

Теорема \ref{eq:matrix-total} влечет такое следствие.

\begin{corollary}
В каждом унитарном неприводимом представлении $\rho^{(\alpha)}$ выберем
по ортонормированному базису $e_k^{(\alpha)}$. Тогда набор функций
$$\la \rho(g) e_k^{(\alpha)}, e_l^{(\alpha)}\ra$$
образует ортонормированный базис в $L^2(G)$. 
\end{corollary}


{\bf \punct Разложимость представлений компактной группы.}

\begin{theorem}
Любое унитарное представление компактной группы разлагается в прямую сумму конечномерных
представлений.
\end{theorem}

На первый взгляд кажется, что это сразу следует из приводимости
бесконечномерных представлений. Мы берем подпредставление $V$,
раскладываем пространство как $V\oplus V^\bot$, в каждом из слагаемых берем
по подпредставлению и т.д. Однако, может оказаться, что на каждом шагу 
все слагаемые будут бесконечномерны, и процесс тогда будет продолжаться 
до бесконечности. 

\sm

{\sc Доказательство.} Пусть $\rho$ -- унитарное представление в пространстве $H$.
Возьмем вектор $\xi\in H$ и построим по нему оператор $J:H\to L^2(G)$ как в п.\ref{ss:J}.
Рассмотрим сопряженный оператор $J: L^2(G)\to H$. В силу последнего
следствия, пространство $L^2(G)$ разлагается
в прямую сумму неприводимых конечномерных подпредставлений. Рассмотрим подпредставление $V$, на котором
$J^*$ является ненулевым. Тогда образ $J^*V$ будет неприводимым подпредставлением в $H$.

Теперь  мы берем неприводимое подпредставление $W_1$ в $H$, рассматриваем ортогональное дополнение
и выбираем в нем неприводимое подпредставление и т.д. Ссылаясь на трансфинитную индукцию
или какие-либо ее варианты (типа леммы Цорна),
мы получаем искомое утверждение. Разумеется, трансфинитная индукция тут не является необходимой.
Вопрос о том, как ее избежать, мы оставляем любителям в качестве упражнения.
\hfill $\square$
  
  \sm
  
  \begin{theorem}
  	\label{th:unique-decomposition}
  	 Разложение унитарного представления $R$ компактной группы
  	 на  неприводимые подпредставления единственно в следующем смысле:
  	 если $R\simeq \oplus \mu_k\simeq \oplus \nu_l$,
  	 то слагаемые в этих разложениях эквивалентны с точностью до перестановки
  	 номеров.
  \end{theorem}

  Мы не утверждаем, что прямые слагаемые в этих разложениях совпадают, и это, вообще говоря, неверно.
  
  Разложим представление $R$ на неприводимые $R=\oplus \mu_k$, и для
  каждого типа неприводимых представлений $\rho_\alpha$ возьмем
  сумму $R_\alpha$ всех слагаемых, изоморфных $\rho_\alpha$. Мы получаем
  {\it разложение $R$ на изотипические слагаемые},
  \begin{equation}
  R=\bigoplus_{\alpha} R_\alpha.
  \label{eq:isotip}
  \end{equation}
  
  \begin{theorem}
  	\label{th:isotip}
  	Разложение на изотипические слагаемые канонично в том смысле, что 
  	каждое слагаемое $R_\alpha$ однозначно определяется самим представлением $R$.
  	\end{theorem} 
  	
  	{\sc Доказательство теоремы \ref{th:isotip}.} Пусть
  	\begin{equation}R=\bigoplus R_\alpha=\bigoplus R'_\alpha
  		\label{eq:two-isotipic}
  	\end{equation}
  	 -- два разложения.
Пусть $\Pi_\alpha$ -- ортогональные проекторы на $R_\alpha$, $\Pi'_\alpha$ -- ортогональные
 проекторы на $R'_\alpha$ (эти операторы являются сплетающими).  	Рассмотрим единичный оператор $I:R\mapsto R$, он разумеется, тоже сплетающий. Запишем его в блочном виде
относительно разложений (\ref{eq:two-isotipic}). Эти блоки имеют вид
$\Pi'_\beta \cdot I\cdot  \Pi_\alpha=\Pi'_\beta \Pi_\alpha$.

Далее заметим, что при $\alpha\ne\beta$ выполнено $\Pi'_\beta \Pi_\alpha=0$ просто потому,
что нет ненулевых сплетающих операторов между $R_\alpha$ и $R'_\beta$. Действительно, разложим $R_\alpha=\oplus V_j$,
$R'_\beta=\oplus W_k$  в прямую сумму неприводимых представлений, в каждом $V_j$ реализовано
представление $\rho_\alpha$, в каждом $W_j$ - представление $\rho_\beta$. Возьмем
соответствующее блочное разложение оператора $\Pi'_\beta \Pi_\alpha$. Тогда каждый
блок является сплетающим оператором, и в силу леммы Шура он равен нулю. 

Итак, лишь блоки $\Pi'_\alpha \Pi_\alpha$ оператора $I$  могут быть ненулевыми.
В силу обратимости оператора $I$, он должен устанавливать биекцию между
$R_\alpha$ и $R'_\alpha$. Т.е.,$R_\alpha=R'_\alpha$. \hfill $\square$

\sm

{\sc Доказательство теоремы \ref{th:unique-decomposition}.}
Изотипическое разложение (\ref{eq:isotip}) канонично, а число неприводимых слагаемых,
на которые разлагается каждое $R_\alpha$ равно $\dim R_\alpha/\dim \rho_\alpha$.
\hfill $\square$

  \sm
  
  {\bf \punct Свертка.} Напомним определение {\it групповой алгебры} конечной группы $\Gamma$.
  Для каждого элемента $\gamma\in \Gamma$ мы определяем формальный базисный элемент
  $e_\gamma$. Групповая алгебра $\C[\Gamma]$ состоит из линейных комбинаций
  $\sum_{\gamma\in\Gamma} c_\gamma e_\gamma$, где $c_\gamma\in\C$. Умножение
  определяется из условия
  $e_\gamma e_\mu=e_{\gamma\mu}$, или
  $$
  \sum a_\gamma e_\gamma\cdot \sum b_\mu e_\mu= \sum_{\gamma\in\Gamma}
  \Bigl(\sum_{\nu\in \Gamma} a_{\gamma\nu^{-1}} b_\nu\Bigr) e_\gamma.
  $$
Введем аналогичную операцию ({\it свёртку}) на локально компактной группе
(чтобы не путаться, {\it положим что группа унимодулярна}).

Пусть $f_1$, $f_2$ -- функции на локально компактной группе $G$. Их
{\it свёртка} $f_1*f_2$ определяется как
$$
f_1*f_2(g)=\int_G f_1(gh^{-1})\,f_2(h)\,dh=\int_G f_1(h)\, f_2(h^{-1}g)\,dh.
$$

{\sc Задача.} a) Покажите, что свертка двух непрерывных функций с компактным
носителем -- непрерывная функция с компактным носителем.

\sm

b) Покажите, что свертка двух функций из $L^1(G)$ лежит в $L^1(G)$.

\sm

c) Покажите, что свертка -- ассоциативная операция  на $L^1(G)$.
\hfill $\lozenge$

\sm

{\sc Задача.} Пусть $U_j\subset G$ -- фундаментальная система окрестностей единицы.
Пусть $\phi_j$ -- неотрицательная непрерывная функция, носитель которой содержится в $U_j$,
такая, что $\int_G \phi_j(g)\,dg=1$. Тогда для любой непрерывной функции
$F$  последовательность $F*\phi_j$ сходится к $F$ в $C(G)$. Для любой функции
из $F\in L^1$  последовательность $F*\phi_j$ сходится к $F$ в $L^1(G)$;
аналогичное высказывание верно и для $L^2$. \hfill $\lozenge$

\sm

{\bf\punct Представления сверточной алгебры.}
Пусть $\rho(g)$ -- унитарное представление локально компактной унимодулярной группы $G$
в гильбертовом пространстве $H$. Для $F\in L^1(g)$ определим оператор
$$
\rho(F):=\int_G F(g)\rho(g)\, dg.
$$

Интеграл можно понимать так: это оператор $Q$, такой, что
  для любых $v$, $w\in H$ выполнено
$$
\la Q v, w\ra = \int_G F(g)\la\rho(g)v,w\ra\, dg.
$$  
Правая часть этого равенства определяет билинейную форму от $v$, $w$, 
которая, как легко видеть, является непрерывной.
Поэтому  она представима в виде $\la Q v, w\ra$.

\sm

Интеграл можно понимать и так: мы выбираем ортнормальный базис в $H$ и интегрируем
каждый матричный элемент
$$
\la Q e_k, e_l\ra = \int_G F(g)\la\rho(g)e_k,e_l\ra\, dg.
$$

\begin{theorem} {\rm a)} $\rho(F_1)\,\rho(F_2)=\rho (F_1*F_2)$.

\sm

{\rm b)} $\rho(F)^*=\rho(\check F)$, где $\check F(g)=\ov{F(g^{-1}}$.

\sm

{\rm с)} Отображение $F\mapsto \rho(F)$ непрерывно как отображение из
$L^1(G)$ в алгебру линейных операторов, снабженную равномерной топологией. 
\end{theorem}

{\sc Задача.} a) Проверьте эти утверждения. 

\sm

b) Докажите обратное утверждение. Пусть есть отображение $\mu$ из сверточной
алгебры $L^1(G)$ в алгебру операторов, удовлетворяющее перечисленным свойствам.
Пусть множество векторов вида $\mu(F)v$, где $v$ пробегает $H$, а $F$ пробегает
$L^1(G)$, тотально. Тогда $\mu$ получается из унитарного представления группы $G$.
См. \cite[10.2]{Kir}, \cite[\S29]{Nai}.
{\sc Указание.} Унитарный оператор $\rho(h)$ строится из условия
$\rho(h)\mu(F)v=\mu(L_h F)v$. Он существует, потому что 
$\la \mu(L_h F_1)v,\mu(L_h F_2)w\ra=\la \mu(F_1)v,\mu(F_2)w\ra$.
Далее остается проверить, что представление сверточной алгебры, определяемое
представлением $\rho(h)$, совпадает с $\mu$.
\hfill $\lozenge$.

\sm

{\bf\punct Доказательство теоремы \ref{eq:matrix-total}.%
	\label{ss:proof-total}}
Пусть теперь $L$ -- левое регулярное представление компактной группы.
Пусть $F$ непрерывна.
Тогда
$$
L(F)f(g)=\int_G F(h)\,f(gh)\,dh=
\int_G F(g^{-1}q)\,f(q)\,dq.
$$ 
Таким образом, $L(F)$ -- интегральный оператор с непрерывным ядром. Поэтому 
$L(F)$ -- оператор Гильберта--Шмидта. 

\sm

Функция $F$ называется {\it центральной}, если она постоянна на классах сопряженности,
$F(hgh^{-1})=F(g)$.

\sm

{\sc Задача.} Для любой центральной функции $F$ и любого $p\in G$ 
$$
L(F)L(p)=L(p) L(F).
$$
Кроме того, для центральной функции $F$ и любой функции $\Psi$
выполнено $F*\Psi=\Psi*F$. \hfill $\lozenge$

\sm

Пусть теперь $F$ -- центральная функция, удовлетворяющая 
\begin{equation}
F(g^{-1})=\ov{F(g)}.
\label{eq:F-self}
\end{equation}
Тогда оператор $L(F)$  самосопряжен.
 Его собственные
подпространства должны быть инвариантны относительно относительно операторов $L(h)$.
Действительно, если $L(F) f=s f$, то 
$$L(F)L(h)f =L(h) L(F)f= s\cdot L(h)f.$$
Возьмем набор непрерывных центральных функций $F_j$.
Общие собственные подпространства операторов $L(F_j)$ тоже $L(h)$-инвариантны. 
Заметим, что все такие подпространства конечномерны, кроме, быть может,
$\cap_j\ker(F_j)$.

 Теперь возьмем фундаментальную систему окрестностей единицы
$U_1\supset U_2\supset\dots$ и возьмем последовательность непрерывных центральных функций
$F_j$, удовлетворяющих (\ref{eq:F-self}), таких, что носитель $F_j$ лежит в $U_j$ и $\int F_j=1$
(такие функции легко строятся с помощью двусторонне инвариантной метрики).
Тогда $L(F_j)$ слабо сходится к единичному оператору, а поэтому $\cap_j\ker(F_j)=0$.
Мы расщепили $L^2(G)$ в прямую сумму конечномерных инвариантных подпространств.

Эти подпространства, в свою очередь, раскладываются в прямую сумму неприводимых 
подпредставлений. Остается доказать следующую лемму.

\begin{lemma}
Пусть $V$ -- неприводимое подпредставление  левого регулярного представления $G$.
Тогда любой элемент  $\phi\in V$ является матричным элементом некоторого неприводимого представления
$G$.
\end{lemma}

{\sc Доказательство.} Нам нужно предъявить представление и два вектора в нем.
Представление -- это само $V$. Один вектор -- функция $\phi$.

Пусть $f$ пробегает $V$. Рассмотрим линейный функционал $\ell_g(f)= f(g)$.
Он должен представляться в виде скалярного произведения 
$$
 f(g)=\ell_g(f)=\la f, \psi_g\ra
.$$
Понятно, что $L(h)\psi_g= \psi_{h^{-1}g}$
В качестве второго вектора в $V$ мы выбираем $\psi_e$, где e -- единица группы.
\hfill $\square$

\sm

{\bf \punct Преобразование Фурье на компактной группе $G$.} 
Пусть $V$ -- конечномерное евклидово пространство. Обозначим через
$\Hom(V)$ алгебру линейных операторов в $V$. 
Введем в $\Hom(V)$ скалярное произведение,
$$
\la\la A, B\ra\ra= \tr AB^*.
$$
Если выбрать в $V$ ортонормированный базис, то стандартные
матричные единицы $E_{ij}$ образуют ортонормированный базис в пространстве
$\Hom(V)$.

Пусть дано унитарное представление $\rho$ группы $G$ в $V$.
Тогда в $\Hom(V)$  действует группа
$G\times G$ левыми и правыми умножениями
$$
(g_1, g_2): A\mapsto g_1^{-1} A g_2.
$$

Обозначим через $\wh G$ множество всех
унитарных неприводимых представлений компактной группы  $G$, определенных с точностью до эквивалентности.
Если $G$ имеет счетную базу, то это множество счетно (в случае конечной группы оно конечно),
действительно, в этом случае $L^2(G)$ сепарабельно, а все неприводимые представления
содержатся в $L^2(G)$.

Пусть $\rho^{(\alpha)}$ пробегает $\wh G$, пусть $H^{(\alpha)}$ -- соответствующие
евклидовы  пространства. 
Рассмотрим линейное  пространство $\cL(\wh G)$, состоящее из всех последовательностей операторов
$$
B^{(\alpha)}\in \Hom(H^{(\alpha)}).
$$ 
Рассмотрим в $\cL(\wh G)$ подпространство $\cL^2(\wh G)$, состоящее из последовательностей
$B^{(\alpha)}$,
удовлетворяющих
$$
\sum_\alpha  \dim(\alpha) \tr B^{(\alpha)} (B^{(\alpha)})^*<\infty,
$$
где $\dim(\alpha):=\dim H^{(\alpha)}$. Введем в $\cL^2(\wh G)$
скалярное произведение по формуле
$$
\la B, C  \ra=\sum_\alpha  \dim(\alpha)  \tr B^{(\alpha)} (C^{(\alpha)})^*.
$$

{\it Преобразование Фурье на компактной группе} $G$ -- это отображение
$$L^1(G)\to \cL(\wh G),$$
заданное формулой
$$
\cF f(\alpha)=\rho^{(\alpha)}(f)=\int_G f(g)\, \rho^{(\alpha)}(g) \,dg.
$$

\begin{theorem}
\label{th:fourier}
{\rm a)} Левый  сдвиг $L_h$ на группе 
переходит в  оператор левого  поточечного
умножения на операторнозначную функции $\rho^{(\alpha)}(h)^{-1}$,
$$
\cF(L_h f)(\alpha)=\rho^{(\alpha)}(h)^{-1} \cF f(\alpha),\qquad 
$$
Аналогично, правый сдвиг переходит в
$$
\cF(R_h f)(\alpha)= \cF f(\alpha) \rho^{(\alpha)}(h)
.
$$

\sm

{\rm b)} Свертка переходит в поточечное умножение,
$$\cF(f_1*f_2)=\cF(f_1)*\cF(f_2).$$

\sm

{\rm c)} Обозначим $\check f(g):=\ov{f(g^{-1})}$. Тогда 
$(\cF (\ov f))^{(\alpha)}=(\cF f^{(\alpha)})^*$.

\sm

{\rm d)} Отображение $\cF$ является унитарным оператором $L^2(G)\to \cL^2(\wh G)$.

\sm

{\rm e)} Обратное отображение $\cL^2(\wh G)\to L^2(G)$
задается формулой
\begin{equation}
f(g)=\sum_\alpha d(\alpha) \tr B^{(\alpha)} \rho_\alpha(g)^{-1}.
\label{eq:cF-1}
\end{equation}

{\rm f)} Отображение $\cF :L^1(G)\to \cL(\wh G)$ инъективно.

\sm

{\rm g)} Образ  центральной функции является последовательностью,
состоящей из скалярных операторов.
\end{theorem}

{\sc Задача.} Пусть $G$ -- это окружность $\R/2\pi/Z$.
Убедитесь, что сформулированная теорема дает стандартный
список свойств разложения функции в ряд Фурье.
Что получится, если $G=\Z_n$ -- конечная циклическая группа?
\hfill $\lozenge$

\sm

{\sc Доказательство.}
Все это было уже доказано выше. Отображение (\ref{eq:cF-1}) является изометричным вложением
в силу соотношений ортогональности. В силу теоремы \ref{eq:matrix-total} его
образ --все  $L^2(G)$, т.е., оператор унитарен. Поэтому обратный оператор
совпадает с сопряженным, а сопряженный оператор к (\ref{eq:cF-1})  -- это
$\cF$.
\hfill $\square$

\sm

{\bf \punct Характеры.} {\it Характером} конечномерного представления $\rho(g)$
называется функция
$$
\chi_\rho(g)=\tr \rho(g)=\sum \la \rho(g) e_j, e_j\ra
,$$
где $e_j$ -- базис в пространстве представления.

Для компактной группы $G$ мы будем  обозначать  через $\chi^{(\alpha)}$
характеры представлений $\rho^{(\alpha)}$

\begin{theorem}
\label{th:characters}
Пусть $G$ -- компактная группа.

\sm 

{\rm a)} Характеры образуют ортонормированный базис в пространстве
центральных функций.

\sm

{\rm b)} Свертки характеров вычисляются как
\begin{align}
\chi^{(\alpha)}*\chi^{(\alpha)}=\frac 1{\dim \rho^{(\alpha)}} \chi^{(\alpha)};
\label{eq:chichi1}
\\
\qquad \chi^{(\alpha)}*\chi^{(\beta)}=0\quad \text{ при $\alpha\ne\beta$.}
\label{eq:chichi2}
\end{align}
Более точно, для любого матричного элемента $R$ представления $\rho^{(\beta)}$
при $\alpha\ne\beta$ выполнено $\chi^{(\alpha)}*R=0$. Для любого матричного
элемента $R$  представления $\rho^{(\alpha)}$ выполнено
$\chi^{(\alpha)}*R= (\dim \rho^{(\alpha)})^{-1} R$.

\sm

{\rm c)} Преобразование Фурье характера $\chi_\alpha$ равно
$$
  (\cF\chi^{(\alpha)})^{\beta}=\begin{cases}
0,& \qquad \text {при $\beta\ne\alpha$};
\\
\frac 1{\dim \rho^{(\alpha)} }, & \qquad \text {при $\beta=\alpha$};
\end{cases}
$$

{\rm d)} Для любого унитарного представления $\mu$ группы $G$ оператор
$$ (\dim \rho^{(\alpha)})^{-1} \mu(\chi^{(\alpha)})$$
является проектором на сумму
всех подпредставлений  $\mu$, изоморфных $\rho^{(\alpha)}$.

\sm

{\rm e)} Пусть $\nu$ -- конечномерное представление $G$. Тогда 
кратность вхождения $\rho^{(\alpha)}$ в $\nu$ равна $\la\nu, \chi^{(\alpha)} \ra$.
\end{theorem}

{\sc Доказательство.} 
Равенство
\begin{equation*} 
 \la \chi^{(\alpha)}, \chi^{(\beta)}\ra= \delta_{\alpha\beta}
\label{eq:chi-chi}
\end{equation*}
 следует из соотношений ортогональности.

Вычислим свертку двух матричных элементов 
$$
\int_G \la \rho_1(h^{-1}g)v_1,w_1\ra
\cdot
\la \rho_2(h)v_1,w_2\ra\,dh=
\int_G \la \rho_1(g)v_1,\rho(h)w_1\ra
\cdot
\la \rho_2(h)v_1,w_2\ra\,dh .
$$
Это скалярное произведение двух матричных элементов. Если $\rho_1$, $\rho_2$ не эквивалентны,
то мы получаем ноль. Если $\rho_1=\rho_2=\rho$, то мы получаем в правой части матричный
элемент
$$
\la v_1, v_2\ra\cdot \la \rho(h)w_1,w_2\ra
.$$
Отсюда легко следует утверждение b) (проверьте это; выведите также эту формулу 
используя преобразование Фурье).

\sm

Теперь вычислим $\rho^{(\alpha)}(\chi^{(\beta)})$.
Функция $\chi^{(\beta)}$ центральна, поэтому наш оператор является
скалярным. Достаточно посчитать его след
$$
\tr\rho^{(\alpha)}(\chi^{(\beta)})=\tr\int_G \chi^{(\beta)}(g) \rho^{(\alpha)}(g)\,dg
= \int_G \tr \rho^{(\beta)}(g)\cdot \tr \rho^{(\alpha)}(g)\,dg=\delta_{\alpha\beta}
.
$$

Это доказывает с). Но тогда функции $\cF (\chi^{(\alpha)})$ образуют ортнормированный
базис в образе пространства центральных функций. Переходя в $L^2(G)$, мы получаем утверждение
a).

Утверждение d) следует из c). Утверждение e) вытекает из a).
\hfill $\square$

\sm

{\bf \punct Лево-право регулярное представление.} 

\begin{theorem}
Пусть $G$, $H$ -- компактные группы. Любое неприводимое представление
группы $G\times H$ имеет вид $\mu\otimes \nu$, где $\mu$, $\nu$ 
-- неприводимые представления $G$ и $H$ соответственно.
\end{theorem}

Сначала напоминание из линейной алгебры. Пусть $V$, $W$ -- конечномерные линейные
пространства. Пусть $\Hom(V,W)$ обозначает множество линейных операторов из
$V$ в $W$. Пусть $V'$ -- двойственное пространство. Есть канонический изоморфизм
\begin{equation}
\Hom(V,W)\simeq V'\otimes W.
\label{eq:hom-otimes}
\end{equation}
А именно, пусть $\ell \in V'$ -- линейный функционал, $w\in W$.
Тогда вектору $\ell\otimes w\in V'\otimes W$ соответствует следующий оператор
ранга 1:
$$
A v= \ell(v)w.
$$
Далее соответствие продолжается по линейности.
Если мы фиксируем базисы $p_i\in V$, $q_j\in W$, то тем самым пространства отождествляются
со своими двойственными; вектору $p_i\otimes q_j$ соответствует оператор,
переводящий $p_i$ в $q_j$, все остальные базисные векторы переходят в ноль.
Можно еще сказать, что и $\Hom(V,W)$  и $V'\otimes W$ отождествляются с пространством
матриц, и мы ставим в соответствие матрице ее саму. Важно, однако, подчеркнуть, что наш
изоморфизм не зависит от выбора базисов.

\sm

Еще одно замечание. Рассмотрим представление
 $\rho$ компактной группы $G$ в пространстве $W$. 
  Рассмотрим  прямую сумму $k$ копий этого представления. Тогда пространство
  представления отождествляется с $\C^k\otimes W$, обозначим $\C^k$ через $V$.

\sm  
  
  {\sc Задача.} Покажите, что операторы вида $A\otimes 1$ в $V\otimes W$
    являются сплетающими, и все сплетающие операторы имеют такой вид.
\hfill $\lozenge$  

\sm

{\sc Доказательство теоремы.}
Рассмотрим  неприводимое представление группы $G\times H $ в пространстве
$S$.  Ограничим его на подгруппу $H$, и разложим в сумму изотипических слагаемых,
$$
 S=\bigoplus_{k=1}^N (V_k\otimes W_k)
$$
где в каждом $W_k$ реализуется неприводимое представление группы $H$
(в разных $W_k$ разные), а $V_k$ группа $H$ действует тривиально.
 
 Группа $G$ 
должна действовать $H$-сплетающими операторами, они могут иметь лишь вид
$$
\oplus_{k=1}^N (A_k(g)\otimes 1), \qquad \text{где $A_k(g)\in\Hom(W_k)$}
.$$
Подпространства $V_k\otimes W_k$ будут $G\times H$-инвариантны. Поэтому
такое слагаемое может быть лишь одно, $N=1$. Далее, операторы 
$A(g)$ должны образовывать представление группы $G$, и оно тоже
должно быть неприводимо.
\hfill $\square$

\sm

\begin{theorem}
\label{th:regular}
Лево-право-регулярное представление группы $G\times G$ в $L^2(G)$ 
разлагается в прямую сумму
$$
\bigoplus_{\alpha\in \wh G} \left(\rho^{(\alpha)}\right)'\otimes \rho^{(\alpha)}
.$$
\end{theorem}

{\sc Доказательство.} Согласно, теореме \ref{th:fourier}.d, наше представление
разлагается в прямую сумму
$$
\bigoplus_{\alpha\in \wh G} \Hom(H^{(\alpha)}, H^{(\alpha)})
.$$
Далее, мы используем изоморфизм \ref{eq:hom-otimes}.
\hfill $\square$

\sm

{\bf \punct Случай конечных групп.}

\begin{theorem}
Пусть группа $G$ конечна.

\sm

{\rm a)} Число попарно различных неприводимых представлений группы
$G$ равно числу классов сопряженности в группе $G$.

\sm

{\rm b)} Пусть $k_1$, \dots, $k_p$ -- размерности попарно различных
неприводимых представлений группы $G$. Тогда 
$$
k_1^2+\dots+k_p^2=\text{порядку группы $G$}
$$
\end{theorem} 

{\sc Доказательство.}  a) 
В силу  теоремы \ref{th:characters}.a число характеров равно числу классов сопряженности.

b) Следует из теоремы \ref{th:regular}. \hfill $\square$

\sm

{\sc Замечание.} Кроме того, размерности неприводимых представлений являются делителями
порядка группы. Доказательство этого утверждения находится несколько в стороне от темы
этих записок, см. \cite[\S11.1]{Kir}
\hfill$\lozenge$

\sm

{\sc Задача.}
a) Опишите все неприводимые представления группы движений правильного 7-угольника.

\sm 

b) Для групп собственных движений правильного тетраэдра,
куба, икосаэдра опишите все неприводимые представления. Покажите, что вторая группа
-- симметрическая группа $S_4$, а первая и третья -- знакопеременные группы
$A_4$ и $A_5$.

\sm

c) Покажите (исходя из утверждений этого пункта), что любая группа порядка
$p^2$, где $p$ -- простое число, абелева. Какие есть группы порядка 
$p^3$ и какие у них неприводимые представления?
\hfill $\square$

\sm

{\bf\punct Квазирегулярные представления компактной группы.}
Пусть, по-прежнему, $\rho^{(\alpha)}\in \wh G$, $H^{(\alpha)}$ -- пространство представления
$\rho^{(\alpha)}$. 

\begin{theorem}
Пусть $G$ -- компактная группа, $K$ -- ее замкнутая подгруппа.
Для $\rho^{(\alpha)}\in \wh G$ обозначим  через $\bigl(H^{(\alpha)}\bigr)^K$
подпространство векторов в $H^{(\alpha)}$, неподвижных относительно $K$.
Тогда представление $G$ в $L^2(G/K)$ изоморфно
$$
\bigoplus_{\alpha\in \wh G} \dim \bigl(H^{(\alpha)}\bigr)^K \cdot \rho^{(\alpha)}.
$$

\end{theorem}

{\sc Доказательство.}
Пространство $L^2(G/K)$ отождествляется с подпространством в 
$L^2(G)$, правоинвариантных относительно $K$, 
$$
f(gh)=f(g)
.$$
Применим к этому подпространству преобразование Фурье.
Получим
$$
\bigoplus_{\alpha\in \wh G} \Hom(H^{(\alpha)},\bigl(H^{(\alpha)}\bigr)^K )
$$
Надо еще заметить, что для $\rho^{(\alpha)}$ и $\bigl(\rho^{(\alpha)}\bigr)'$
размерности пространств $K$-неподвижных векторов совпадают.
\hfill $\square$

\sm

Литература к параграфу: \cite[8.1, 9.2, 11.1, 12.1-12.2]{Kir},
\cite[гл. 7]{HR2},  \cite[гл. 3]{Ada}, \cite[\S\S 27-30]{Zhe1}, \cite[\S28]{Zhe2}, \cite{Serre-finite}.

\section{Распределение собственных чисел}

\COUNTERS

{\bf\punct Формула Г.~Вейля.}
Рассмотрим группу $\U(n)$. Любой ее элемент сопряжением приводится к диагональному виду,
с числами вида $e^{i\phi_1}$, \dots, $e^{i\phi_n}$ по диагонали.
Числа эти определены с точностью до перестановок. Положим, для определенности,
что
\begin{equation}
0\le\phi_1\le \phi_2\le\dots\le\phi_n<2\pi.
\end{equation}
Обозначим через $\Lambda\subset\R^n$ множество всевозможных таких наборов.
Рассмотрим отображение 
\begin{equation} 
\Xi:\U(n)\to\Lambda,
\label{eq:Xi}
\end{equation}
 которое каждой матрице $g$ ставит в соответствие
 набор $\phi_j$. 

\begin{theorem}
\label{th:wey}
Образ вероятностной меры Хаара при отображении $\Xi$ имеет вид
\begin{equation}
\frac 1{(2\pi)^n } \prod_{k<l} \left|e^{i\phi_k}-e^{i\phi_l}\right|^2 \prod_k d\phi_k
=\frac {2^{n(n-1)/2}}{(2\pi)^n }  \prod_{k<l} \left|\sin\frac{\phi_k-\phi_l}2\right|^2 \prod_k d\phi_k
.
\label{eq:weyl}
\end{equation}
\end{theorem}

В частности, пусть $H$ -- центральная функция на $\U(n)$, т.е. функция, зависящая
лишь от набора собственных значений $e^{i\phi_1}$, \dots, $e^{i\phi_n}$ матрицы $g$,
$$
H(g)=h\left( e^{i\phi_1}, \dots, e^{i\phi_n}\right).
$$
Тогда
\begin{equation}
\int\limits_{\U(n)} H(g)\, dg:=\frac 1{(2\pi)^n n!}
\int\limits_0^{2\pi}\dots \int\limits_0^{2\pi}
h\left( e^{i\phi_1}, \dots, e^{i\phi_n}\right)
\prod_{k<l} \left|e^{i\phi_k}-e^{i\phi_l}\right|^2 \prod_k d\phi_k
.
\label{eq:wey-int}
\end{equation}

Стоит иметь в виду следующую формулу с определителем Вандермонда
\begin{equation}
\prod_{k>l} \left(e^{i\phi_k}-e^{i\phi_l}\right)=
\det\begin{pmatrix}
e^{i(n-1)\phi_1}&e^{i(n-1)\phi_2}&\dots& e^{i(n-1)\phi_n}\\
e^{i(n-2)\phi_1}&e^{i(n-2)\phi_2}&\dots& e^{i(n-2)\phi_n}\\
\vdots&\vdots&\ddots&\vdots\\
1&1&\dots&1
\end{pmatrix}
.
\label{eq:vander}
\end{equation}

{\bf\punct Доказательство.}
{\it Обозначения}. Через  $E_{kl}$ обозначим матричные единицы, т.е. матрицы,
у которых на месте $kl$ стоит 1, а остальные матричные элементы равны нулю. Единичную
матрицу мы   будем обозначать через $e$. Касательное пространство к многообразию
$M$ в точке $x$ мы будем обозначать через $T_x M$.

Через $\Fl_n$ мы обозначим пространство (многообразие), точкой которого является
упорядоченная
$n$-ка $\ell:=(\ell_1,\dots,\ell_n)$ попарно ортогональных прямых  в $\C^n$.
 Очевидно, $\Fl_n$ является однородным
 $\U(n)$-пространством, стабилизатор точки изоморфен  $n$-мерному тору
$\T^n$. Если в качестве точки выбран набор координатных осей 
$$\frl:=\{\C e_1, \dots,\C e_n\},$$
 то стабилизатор 
-- это в точности группа $\Delta$ диагональных матриц. Т.е, 
$$\Fl_n\simeq\U(n)/\Delta.$$
Каждой $n$-ке $\ell$ мы поставим в соответствие унитарную матрицу $h$,
такую, что $h\frl=\ell$. Элемент $h$ определен с точностью до домножения
справа на элемент из $\Delta$. 
В силу теоремы \ref{th:invariant-compact},
это пространство%
\footnote{Отметим, что это однородное пространство совпадает
с пространством флагов, обсуждавшимся в \S3.}
снабжено инвариантной мерой, мы будем обозначать ее
через $d\ell$.

\sm 

{\it  Разделение координат.} Теперь заметим, что множество унитарных
 матриц, у которых хотя бы два
 собственных числа совпадают,
 имеет меру ноль%
 \footnote{{\sc Задача.} Докажите, что эта алгебраическая поверхность имеет  коразмерность
 	 3 в $\U(n)$ (ниже, при доказательстве теоремы,
 	нам это уточнение не нужно)}%
. Обозначим
через $\U'$ множество матриц, у которых собственные числа различны. Через $\Lambda'$ обозначим
симплекс
$$ 
0\le\phi_1< \phi_2<\dots<\phi_n<2\pi
.
$$
Элемент симплекса $(\phi_1, \dots,\phi_n)\in \Lambda'$ мы будем обозначать через $\phi$.
 Той же буквой мы будем
обозначать диагональную  матрицу с собственными числами $\phi_j$.
Через $d\phi$ обозначим меру Лебега на $\Lambda'$.

Рассмотрим естественное отображение 
$$S:\Fl_n\times\Lambda'\to \U'.$$
А именно, паре $(\ell,\phi)$ мы ставим в соответствие линейное преобразование $\in\U(n)$,
 у которого каждая прямая $\ell_k$
является собственной с собственным значением $e^{i\phi_k}$.

Очевидно, мы получаем биекцию, причем она коммутирует с действием группы
$\U(n)$. 
На $\Fl_n$ группа $\U(n)$ действует вращением $n$-ок, на $\Lambda'$ она 
действует тривиально, на самой себе $\U(n)$ действует сопряжениями $g\mapsto hgh^{-1}$.

Если элементу $\ell \in\Fl_n$ мы ставим в соответствие $h\in \U(n)$ как выше,
a $\phi \in \Lambda'$,
то
$$ 
S(h,\phi)=h \exp(i\phi)\,h^{-1}.
$$

Мера Хаара на $\U(n)$ инвариантна относительно сопряжений. Ее образ относительно
отображения $S^{-1}$ должен быть $\U(n)$-инвариантной мерой на  $\Fl_n\times\Lambda'$.
В силу $\U(n)$-инвариантности, эта мера должна иметь вид
$$
Q(\phi)\,d\ell\, d\phi 
,$$
где 
 $Q(\phi)$ -- некоторая функция,
ее нам и надо вычислить. Можно это сказать другими словами. На $\U(n)$ есть инвариантная дифференциальная форма старшей степени. Ее прообраз при отображении $S$ является формой старшей степени на 
$\Fl_n\times \Lambda'$. С другой стороны, на $\Fl_n\times \Lambda'$ есть форма 
$d\ell\, d\phi$. Нам нужно вычислить отношение этих двух форм. 
В силу инвариантности, этот множитель достаточно вычислить  при
$\ell=\frl$. 
Вопрос фактически состоит в вычислении якобиана 
(определителя дифференциала),
 его достаточно посчитать 
в точках вида $(\frl,\phi)\in \Fl_n\times\Lambda'$ (где $\phi$ 
пробегает $\Lambda'$). Отображение $S$ индуцирует отображение
касательных пространств%
\begin{equation}
T_\frl\, \Fl_n\,\times\, T_{\phi}\,\Lambda'\,\to\, T_{\exp(i\phi)}\,\U(n)
.
\label{eq:TTT1}
\end{equation}

{\it Описание  касательных пространств}.

\begin{lemma}
\label{l:lie-u}
Касательное пространство $T_e\,\U(n)$ к $\U(n)$ в единице состоит из антиэрмитовых
{\rm($X^*=-X$)} матриц.
\end{lemma}

{\sc Доказательство.} Пусть
$X$ -- касательный вектор к $\U(n)$ в единице. Выпустим кривую $\gamma(\epsilon)$
в направлении этого вектора. Тогда
$$
1=\gamma(\epsilon)^*\gamma(\epsilon)=(1+\epsilon X+O(\epsilon^2))^* (1+\epsilon X+O(\epsilon^2))=1+\epsilon (X^*+X)+O(\epsilon^2),
$$
и, тем самым, $X^*+X=0$. Обратно, если $X$ антиэрмитова, то в качестве кривой можно взять
$\gamma(\epsilon)=\exp(\epsilon X)$. Тогда
$$
\qquad\qquad\qquad\qquad
\gamma(\epsilon)^*\gamma(\epsilon)= \exp(-\epsilon X)\exp(\epsilon X)=1.
\qquad\qquad\qquad\qquad
\square
$$

{\sc Задача.} Условие унитарности матрицы $g$, т.е. $g^*g-1=0$, состоит в занулении
$n^2$ вещественных функций. Покажите, непосредственно, что их градиенты
в единице независимы. Покажите их независимость в остальных точках группы.
\hfill $\lozenge$

\sm

Выберем в $T_e\,\U(n)$ базис $E_{kl}-E_{lk}$, $i(E_{kl}+E_{lk})$,
$i E_{kk}$.

\sm

Так как $\Fl_n\simeq \U(n)/\Delta$, касательное пространство 
$T_\frl\,\Fl(n)$ отождествляется с факторпространством
$T_e\,\U(n)/T_e\Delta$, т.е., с пространством антиэрмитовых матриц,
профакторизованному по пространству диагональных матриц.
В этом пространстве мы выбираем базис $E_{kl}-E_{lk}$, $i(E_{kl}+E_{lk})$.

Область $\Lambda'$ содержится в $\R^n$, поэтому касательное пространство
к ней отождествляется с $\R^n$. Так как $\Lambda'$ нумерует диагональные
матрицы, удобно считать, что базисом в ней являются матрицы $E_{kk}$.

{\it Вычисление якобиана.}
Вернемся к отображению  (\ref{eq:TTT1}). В правой части стоит пространство,
зависящее от точки $\phi$. Естественно рассмотреть левый сдвиг 
$g\mapsto \exp(-i\phi) g$, который сдвинет $e^{i\phi}$ в единицу.
Теперь мы имеем отображение
\begin{equation}
(h,\phi)\mapsto \exp(-i\phi) h \exp(i\phi) h^{-1}
\label{eq:eheh}
,\end{equation}
которое индуцирует линейное отображение касательных пространств
$$
T_\frl\, \Fl_n\,\times\, T_{\phi}\,\Lambda'\,\to\, T_{e}\,\U(n).
$$
Нам надо посчитать его определитель. Можно взять выписанные выше базисы
и записать наше отображение через эти базисы (оно окажется составленным из
блоков размера $1\times 1$ и $2\times 2$). Однако чуть приятнее перейти 
к комплексификациям этих пространств (определитель от этого не изменится),
\begin{equation}
(T_\frl\, \Fl_n)_\C \,\times\, (T_{\phi}\,\Lambda')_\C \,\to\, (T_{e}\,\U(n))_\C.
\label{eq:TC}
\end{equation}

 Рассмотрим касательный вектор
$$v:=\sum_{k\ne l} t_{kl} E_{kl}\in (T_\frl\, \Fl_n)_\C $$
 и касательный вектор 
$$u=\sum  \sigma_k E_{kk}\in (T_{\phi}\,\Lambda')_\C.$$
 Посмотрим, что с ними случится
при отображении (\ref{eq:eheh}). Иными словами, нам надо вычислить с точностью
до $O(\epsilon^2)$ выражение
\begin{multline*}
 \exp(-i\phi) (1+\epsilon v)\exp(i\phi)
(1+i\epsilon u) (1-\epsilon v)
=\\=
1+\epsilon \bigl(-v+ \exp(-i\phi)v\exp(i\phi)  +iu\bigr)+O(\epsilon^2)
\end{multline*}
Сопряжение диагональной матрицей легко считается, и мы получаем
$$-v+ \exp(-i\phi)v\exp(i\phi)  +iu=
\sum_{k\ne l} (e^{i(\phi_l-\phi_k)}-1)t_{kl}E_{kl}+ \sum_k \sigma_k E_{kk}
$$
Таким образом, матрица оператора (\ref{eq:TC}) оказывается диагональной. На диагонали стоят всевозможные
числа $e^{i(\phi_l-\phi_k)}-1$ и еще $n$ единиц. Перемножая эти числа, мы получаем
желаемый результат.

\sm

См. \cite[\S72]{Zhe1},  \cite[\S 7.4]{Wey}, \cite[\S5.3]{Ros}, \cite[\S 3.2]{Hua}.

\sm

{\bf\punct Вычисление нормировочной константы.} Нужно еще вычислить
числовой множитель в формуле (\ref{eq:vander}). Для этого нужно посчитать
следующий интеграл по тору
$$
\frac 1{(2\pi)^n}
\int_0^{2\pi}\dots \int_0^{2\pi} \prod_{k<l} \left|e^{i\phi_k}-e^{i\phi_l}\right|^2 \prod_k d\phi_k.
$$
Под интегралом стоит скалярный квадрат тригонометрического многочлена. Поэтому интеграл
есть сумма квадратов его коэффициентов. Выражение для 
$\prod_{k<l} (e^{i\phi_k}-e^{i\phi_l})$ дается формулой (\ref{eq:vander}).
Мы видим, что все коэффициенты многочлена равны $\pm 1$, а число слагаемых 
равно $n!$.

\sm

{\bf\punct Функции Шура.}
Пусть $\alpha_1\ge \alpha_2\ge\dots\ge \alpha_n$ -- целые числа. Положим
$$
\Delta_\alpha(z):=\det \{ z_j^{\alpha_k+n-k}\}:=
\det\begin{pmatrix}
z_1^{\alpha_1+n-1}&z_2^{\alpha_1+n-1}&\dots&z_n^{\alpha_1+n-1}\\
z_1^{\alpha_2+n-2}&z_2^{\alpha_2+n-2}&\dots&z_n^{\alpha_2+n-2}\\
\vdots&\vdots&\ddots&\vdots\\
z_1^{\alpha_n}&z_2^{\alpha_n}&\dots&z_n^{\alpha_n}
\end{pmatrix}
.
$$
Отметим, что $\Delta_0(z)$ -- это определитель Вандермонда,
$$\Delta_0(z)=\prod_{k<l} (z_k-z_l).$$
{\it Функция Шура $s_\alpha(z)$} определяется формулой
$$
s_\alpha(z)=\frac{\Delta_\alpha(z)}{\Delta_0(z)}
.
$$
Отметим, что 
\begin{equation}
s_{\alpha_1+1,\dots, \alpha_n+1}(z)=z_1\dots z_n\cdot s_{\alpha_1,\dots \alpha_n}(z).
\label{eq:shur-shift}
\end{equation}

 \begin{lemma}
 Функция Шура $s_\alpha(z)$ симметрична по переменным
 $z_1$,\dots, $z_n$. При $\alpha_n\ge 0$ она является однородным
 многочленом от  $z_1$,\dots, $z_n$ степени $\sum \alpha_j$.
 \end{lemma}
 
 {\sc Доказательство.} Вычтем второй столбец в $\Delta_\alpha$ из первого.
 Все выражения
$z_1^{\alpha_k+n-k}-z_2^{\alpha_k+n-k}$ делятся на $(z_1-z_2)$. Поэтому
многочлен
$\Delta_\alpha(z)$ делится на $(z_1-z_2)$. Аналогично, он делится на
$z_k-z_l$, а поэтому и на $\Delta_0(z)$.
 
 Что касается симметричности, то и числитель, и знаменатель антисимметричны,
 а поэтому отношение симметрично.
 \hfill $\square$

\sm

Рассмотрим $n$-мерный тор $\T^n$
$$
|z_1|=|z_2|=\dots=|z_n|=1.
$$
Обозначим через $L^2_{symm}(\T^n)$ пространство {\it симметричных} функций на $\T^n$
со скалярным произведением
$$
\la f, g\ra=\frac1{(2\pi)^n n!}
\int_{\T^n} f(z)\,\ov {g(z)}\, |\Delta_0(z)|^2\, \prod d\phi_j,\qquad z_j=e^{i\phi_j} .
$$

\begin{theorem}
Функции Шура образуют ортнормированный базис в пространстве $L^2_{symm}(\T^n)$.
\end{theorem}

{\sc Доказательство.}
$$
\la s_\alpha(z),s_\beta(z)\ra=
\!\!\!\int \frac{\Delta_\alpha(z)}{\Delta_0(z)}\cdot 
\ov{\biggl(\frac{\Delta_\beta(z)}{\Delta_0(z)}\biggr)}
\Delta_0(z) \ov{\Delta_0(z)}\prod d\phi_j
=\!\!\!\int \Delta_\alpha(z)\ov{\Delta_\beta(z)}\prod d\phi_j
.
$$
Теперь мы можем переформулировать утверждение в виде:
функции $ \Delta_\alpha(z)$ образуют ортонормированный базис в пространстве
кососимметрических функций на торе со стандартным $L^2$-скалярным произведением.
Но это очевидно: функции $  \Delta_\alpha(z)$ получаются кососимметризацией
базисных тригонометрических одночленов 
$$
\qquad\qquad\qquad\qquad
z_1^{\alpha_1}\dots z_n^{\alpha_n}=e^{i\alpha_1\phi_1}\dots e^{i\alpha_n\phi_n}.
\qquad\qquad\qquad\qquad\square
$$

{\bf\punct Характеры неприводимых представлений $\U(n)$.}
Характер -- функция от класса сопряженности, а поэтому характер унитарной
группы фактически есть функция от собственных значений 
$e^{i\phi_1}$, \dots, $e^{i\phi_1}$. Поэтому характер может рассматриваться как 
симметрическая функция на торе $\T^n$.  В силу формулы
(\ref{eq:wey-int}) характеры неприводимых представлений
 должны образовывать ортонормированный базис
в $L^2_{symm}(\T^n)$.

\begin{theorem}
\label{th:characters-u}
Характеры неприводимых представлений группы $\U(n)$ -- это в точности функции Шура.
\end{theorem} 

В следующим пункте мы приведем доказательство этой теоремы, впрочем оно, в определенном смысле,
является недостойным.

\sm

{\bf\punct Доказательство теоремы \ref{th:characters-u}.}

\sm

{\it Некоторые свойства функций Шура.}

\sm

{\sc Задача.} Покажите, что элементарные симметрические функции
$$\sigma_p(z):=\sum_{i_1<\dots < i_p} z_{i_1} z_{i_2}\dots z_{i_p}$$
являются функциями Шура. В этом случае $\alpha$ состоит из $p$ единиц, а 
дальше -- нули.
\hfill $\lozenge$

\sm

 Рассмотрим пространство $V_N$
симметрических
однородных многочленов степени $N$ от формальных переменных $z_1$,\dots, $z_n$.
Рассмотрим в $V_N$ (неортогональный) базис из т.н. <<мономиальных симметрических функций>>
(см. \cite{Mac})
$$
m_\alpha (z)=z_{1}^{\alpha_1}\dots z_{n}^{\alpha_n} +
\left\{\begin{matrix}\text{сумма всех {\it различных} мономов,}\\
\text{полученных перестановками $z_j$}
\end{matrix} \right\}
,$$
где 
\begin{equation}
\alpha_1\ge\dots\ge\alpha_n\ge 0,\quad \text{а $\sum\alpha_j=N$}
\end{equation}
 Упорядочим элементы базиса
в лексикографическом порядке. Самый большой элемент это 
$$m_{N,0,\dots,0}=\sum z_j^N.$$
Самый маленький имеет вид 
\begin{equation}
(z_1\dots z_n)^k \sigma_p(z),\qquad \text{где $N=kn+p$}.
\label{eq:minimal}
\end{equation}

Обозначим через $V_N[\alpha]$ линейную оболочку всех $m_\beta$ c $\beta\preccurlyeq\alpha$,
а через $V_N^\circ[\alpha]$ линейную оболочку всех $m_\beta$ c $\beta\prec\alpha$.

\sm

{\sc Задача.} a) Покажите, что 
$$
s_\alpha(z)=m_\alpha(z)+\sum_{\beta\prec \alpha} c_{\alpha\beta} m_\beta(z)
,$$
где $c_{\alpha\beta}$ -- некоторые числа (нам не интересные).

\sm

b) Выведите отсюда, что функции $s_\alpha(z)$ с $\sum \alpha_j=N$ образуют
базис в $V_N$, причем  $s_\alpha(z)$ перпендикулярна к  $V_N^\circ[\alpha]$
\hfill $\lozenge$

\sm

{\it Некоторые  представления группы $\U(n)$}. Рассмотрим
$p$-ые внешние степени $\lambda_p$ тавтологического $n$-мерного представления $g\mapsto g$ 
группы $\U(n)$. Это представления в пространстве поливекторов вида
$$
\sum_{i_1<\dots < i_p} c_{i_1,\dots, i_p} e_{i_p}\wedge \dots \wedge e_{i_p}.
$$

{\sc Задача.} Убедитесь, что характер этого представления
- это элементарная симметрическая функция 
$\sigma_p(z)$.
\hfill $\lozenge$

\sm

Пусть $\alpha_n\ge 0$. Рассмотрим тензорное произведение
\begin{equation}
\bigotimes_{p=1}^n \lambda_p^{\otimes(\alpha_p-\alpha_{p+1})},
\label{eq:otimesotimes}
\end{equation}
где $\alpha_{n+1}:=0$. Его характер равен
$$
\sigma_\alpha(z):=\prod_p \sigma_p(z)^{\alpha_p-\alpha_{p+1}}
.$$

{\sc Задача.} Убедитесь, что
$$
\sigma_\alpha(z)=
m_\alpha(z)+\sum_{\beta\prec \alpha} b_{\alpha\beta} m_\beta(z)
,$$
где $b_{\alpha\beta}$ -- некоторые числа.

\sm

{\it Доказательство.}
Дальше мы доказываем  индукцией по $\alpha$ то, что функции $s_\alpha(z)$
фиксированной степени $N$ являются неприводимыми характерами.
Для самого малого $\alpha$ это так, см (\ref{eq:minimal}), 
(\ref{eq:shur-shift}). Пусть мы доказали это утверждение для всех
$\alpha\prec\gamma$. Тогда
$$
\kappa(z):=
\sigma_\gamma(z)-\sum_{\alpha\prec\gamma} \la \sigma_\alpha,s_\alpha\ra\cdot s_\alpha(z)
$$
является характером (на уровне представлений это означает выбрасывание из представления всех
неприводимых компонент с характерами $s_\alpha$, $\alpha\prec\gamma$).

Далее $\kappa(z)$ ортогональна $V_N^\circ[\gamma]$ и содержится в $V_N[\gamma]$.
Поэтому $\kappa(z)$ пропорциональна $s_\gamma(z)$. Но  $\kappa(z)$ имеет вид
$$
\kappa(z)=m_\gamma(z)+\sum_{\alpha\prec \gamma} a_{\gamma\alpha} m_\alpha(z).
$$
А поэтому $\kappa(z)=s_\gamma(z)$, и тем самым $s_\gamma(z)$ является характером.
Так как $\la s_\gamma,s_\gamma\ra=1$, этот характер  неприводим.

\sm

Итак, $s_\alpha(z)$ является характером при $\alpha_n\ge 0$. Это ограничение
не существенно, так как мы можем умножать представление $\rho(g)$ на $\det(g)^k$,
тогда характер умножается на $z_1\dots z_p$, а индекс функции Шура сдвигается на
$(k,\dots,k)$, см. (\ref{eq:shur-shift}). Это завершает доказательство теоремы \ref{th:characters-u}.

\sm

{\sc Замечание.} Заодно мы получили классификацию неприводимых представлений группы
$\U(n)$. Они нумеруется сигнатурами $\alpha$ как функции Шура. Соответствующее 
представление $\rho_\alpha$ содержится в (\ref{eq:otimesotimes}) и является 
единственным слагаемым, которое не содержится в таких же тензорных произведениях
 с меньшими номерами $\alpha$.
Другой вопрос, что это описание не всегда удобно.
\hfill $\square$

\sm

\cite[гл. XI]{Zhe1}, \cite[гл. VII]{Wey}, \cite[\S 6.4]{Ros},  \cite[гл. 7]{Ada},
\cite[\S I.3]{Mac}.

\sm

{\bf\punct Зоопарк.%
\label{ss:zoo-comp}} Здесь мы без доказательства приводим несколько формул  в духе
(\ref{eq:weyl}). 

\sm

{\it Классические компактные группы.}
Элемент группы $\SO(2n,\R)$ имеет собственные числа вида $e^{\pm i\phi_k}$, где $k=1$, \dots,
$n$. Плотность  распределения собственных чисел равна
$$
C\cdot \prod_{1\le k<l\le n} \sin^2\left(\frac{\phi_k-\phi_l}2  \right)  \sin^2\left(\frac{\phi_k+\phi_l}2  \right).
$$

В случае групп $\SO(2n+1)$ собственные числа имеют вид  $e^{\pm i\phi_1}$, \dots,
$e^{\pm i\phi_n}$, 1. Плотность распределения равна
$$
C\cdot \prod_{1\le k<l\le n} \sin^2\left(\frac{\phi_k-\phi_l}2  \right)  \sin^2\left(\frac{\phi_k+\phi_l}2 
\right)\prod_{1\le k\le n} \sin^2\left(\frac {\phi_k}2\right).
$$

Наконец, в случае симплектической группы $\Sp(2n)$ набор собственных чисел имеет
вид  $e^{\pm i\phi_1}$, \dots,
$e^{\pm i\phi_n}$, а плотность их распределения
$$
C\cdot \prod_{1\le k<l\le n} \sin^2\left(\frac{\phi_k-\phi_l}2  \right) 
\sin^2\left(\frac{\phi_k+\phi_l}2  \right)\prod_{1\le k\le n} \sin^2 {\phi_k}.
$$

{\it Пространство эрмитовых матриц.} Пусть  $\K$ обозначает $\R$, $\C$ или тело кватернионов
$\H$. Обозначим 
$$d=\dim \K.$$
 Через $\Herm_n(\K)$ обозначим пространство эрмитовых матриц размера $n$ с коэффициентами
из $\K$. Такие матрицы приводятся унитарными (над $\K$) преобразованиями к диагональному виду.
Распределение их собственных чисел имеет вид
$$
C(\K)\cdot \prod_{1\le k<l\le n}|\lambda_k-\lambda_l|^{d}.
$$

Все эти формулы доказываются одинаково. Мы не будем здесь стремиться к полноте,
в следующем пункте мы приведем без доказательства более изощренный результат. 

\sm

{\sc Задача.} Выведите какую-нибудь из этих шести формул.
\hfill $\lozenge$

\sm 

См. \cite[гл. 7]{Wey}, \cite[\S 123]{Zhe1}, \cite[глава 3]{Hua}, \cite{Met}.

\sm 

{\bf\punct Эллиптические координаты.}
Это фольклорное название следующей конструкции.
Для $X\in \Herm_n(\K)$ через $[X]_p$ означим левый верхний уголок матрицы
$X$ размера $p\times p$. Обозначим через 
$$
\lambda_{p1}\le \lambda_{p2}\le\dots\le \lambda_{pp}
$$
собственные числа этого уголка. Согласно теореме Рэлея, собственные
числа уголков $[X]_p$ и $[X]_{p+1}$ чередуются,
$$
\lambda_{(p+1)1}\le \lambda_{p1}\le \lambda_{(p+1)2}\le \lambda_{p2}\le \lambda_{(p+1)3}
\le
\dots\le \lambda_{pp} \le \lambda_{(p+1)(p+1)}.
$$

Рассмотрим отображение $\Sigma$, которое каждой матрице $X\in\Herm_n(\K)$
ставит в соответствие набор всех собственных чисел всех уголков $[X]_p$.
Обозначим через $\cR$ множество всех наборов чисел $\lambda_{kj}$,
где $1\le k\le n$, $1\le j\le k$, удовлетворяющих условиям чередования.
 
 \begin{theorem}
 Образ меры Лебега при отображении $\Sigma:\Herm_n(\K)\to \cR_n$
 равен 
\begin{multline*}
d\Phi(\lambda)=
\frac{\pi^{n(n-1)d/2}}
{\Gamma^{n(n-1)/2} (d/2)} \cdot \frac
{\prod\limits_{2\le j\le n} \quad
\prod\limits_{1\le \alpha\le j-1,\,\,\,
                               1\le p\le j}
|\lambda_{(j-1)\alpha}-\lambda_{jp}|^{d/2-1} }
{\prod\limits_{2\le j\le n-1}\quad
\prod\limits_{1\le \alpha<\beta\le j}
  (\lambda_{j\beta}-\lambda_{j\alpha})^{d-2}}
\times \\ \times
\prod\limits_{1\le p<q\le n}
(\lambda_{nq} -\lambda_{np})  \prod\limits_{1\le j\le n}\quad \prod_{1\le\alpha\le j}
  d\lambda_{j\alpha}
.\end{multline*}
 \end{theorem}

 Матрица $X$ не восстанавливается по собственным числам $\lambda_{pj}$.
 Однако можно ввести недостающие координаты следующим образом.
 Обозначим через $e_k$ стандартный базис в $\R^n$.
 Предположим, что собственные числа каждой  матрицы $[X]_p$
 попарно различны (это выполнено всюду за исключением подмногообразия коразмерности
$1+d$). Рассмотрим набор единичных собственных векторов 
$v_{pj}\in \K^p$, отвечающих собственным числам $\lambda_{pj}$.
Отнормируем их условием 
$$
\la v_{pj}, e_1\ra>0.
$$
Далее введем набор чисел
$$
u_{pj}:=\frac{\la v_{pj},e_{p+1}\ra}{|\la v_{pj},e_{p+1}\ra|}\in \K
.
$$
По определению, $|u_{pj}|=1$. Т.е., мы получаем точку сферы $S^{d-1}$.
 В случае $\K=\R$ мы имеем нульмерную сферу $u_{pj}=\pm 1$.
 Таким образом, мы получаем определенное почти всюду 
 (за исключением набора подмногообразий коразмерности $d$ и $1+d$)
 отображение
 $$
\wt \Sigma: \Herm_n(\K) \to \cR_n\times \left(S^{d-1}\right)^{(n-1)(n-2)/2}.
 $$
 Предыдущая теорема уточняется следующим образом.
 
 \begin{theorem}
 Отображение $\wt \Sigma$ взаимно однозначно с точностью до п.в.
 Образ меры Лебега при отображении $\wt\Sigma$ есть произведение
 меры $d\Phi(\lambda)$ и равномерной вероятностной меры на произведении сфер. 
 \end{theorem}
 
См. \cite{Ner-triangle}.

\section{Симметрии гауссовой меры}

\COUNTERS

{\bf \punct Гауссова мера.} Рассмотрим прямую $\R$ с гауссовой мерой
$$(2\pi)^{-1/2}e^{-x^2/2} \,dx$$
 (мера всей прямой равна 1). Рассмотрим счетное произведение
этого пространства самого на себя,
$$
(\R^\infty,\mu):=\left(\R,\frac1{\sqrt{2\pi}}e^{-x^2/2} \,dx\right)
\times 
\left(\R,\frac1{\sqrt{2\pi}}e^{-x^2/2} \,dx\right)\times\dots
$$

Точки этого пространства мы будем обозначать как
$x=(x_1,x_2,\dots)$.

Множество называется $A\subset \R^\infty$ называется {\it цилиндрическим},
если оно имеет вид 
$$A=A_1\times \dots \times A_k\times \R\times \R\times\dots$$
Для цилиндрического множества мера $\mu(A)$ равна произведению мер 
множеств $A_j$. По определению счетного произведения мер, цилиндрические множества порождают
$\sigma$-алгебру всех измеримых множеств.

Функция $f$ называется {\it цилиндрической}, если она зависит лишь от первых нескольких
переменных, $f(x)=f(x_1,\dots,x_k)$. Для цилиндрических функций мы имеем
$$
\int\limits_{\R^\infty} f(x_1,\dots, x_k)\,d\mu(x)
=\frac{1}{(2\pi)^{k/2}} \int\limits_{\R^k} f(x_1,\dots, x_k)\,
e^{-(x_1^2+\dots+x_k^2)/2}
\,dx_1\dots dx_k.
$$

Отметим следующие простые равенства
\begin{align}
&\frac1{\sqrt{2\pi}}\int_{-\infty}^\infty x^{2k} e^{-x^2/2}\,dx=(2n-1)!!
\\
&\frac1{\sqrt{2\pi}}\int_{-\infty}^\infty e^{zx} e^{-x^2/2}\,dx=e^{z^2/2},\quad
 \text{где $z\in \C$};
\\
&\int_{\R^\infty} f_1(x_1)\,f_2(x_2)\,\dots f_n(x_n)\,d\mu(x)
=\prod_{j=1}^n \frac1{\sqrt{2\pi}}\int_{-\infty}^\infty f_j(x)\, e^{-x^2/2}\,dx
.
\end{align}

{\sc Задача.} Вспомните определение независимых случайных
величин и убедитесь, что функции вида $f_1(x_1)$, $f_2(x_2)$, \dots
независимы. \hfill $\lozenge$

\sm

{\bf\punct Некоторые множества полной меры.}

\sm

{\sc Задача.} Мера множества ограниченных последовательностей равна 0.
В частности, равна 0 мера пространства $\ell_2$.
\hfill $\lozenge$

\begin{proposition}
\label{pr:beppo}
Пусть $\lambda_k>0$, $\sum \lambda_k<\infty$. Тогда ряд 
$\sum \lambda_k x_k^2$ сходится почти всюду.
\end{proposition}   

{\sc Доказательство.} Мы применяем теорему Беппо Леви о монотонной сходимости.
\hfill $\square$

\sm

Мы видим, что мера сосредоточена на последовательностях, которые не очень быстро растут.

\begin{theorem}
\label{th:sqrt-ln}
Для почти всех $x\in\R^n$ выполнено
$$
\limsup_k \frac{|x_k|}{\sqrt {2\ln k}}=1
.
$$
\end{theorem}

Для доказательства нам будет нужно следующее утверждение
 (<<{\it лемма  Бореля--Кантелли}>>).
 
 \begin{theorem}
 Пусть $(X,\xi)$ -- пространство с вероятностной мерой.
 Пусть $A_1$, $A_2$, \dots -- независимые события.
 Вероятность того, что точка содержится  в бесконечном числе  множеств
 $A_j$, равна нулю,
 если  $\sum \xi(A_j)<\infty$, и единице, если $\sum \xi(A_j)=\infty$.
 \end{theorem}

См., например, \cite[VI.10]{Shir}.

\sm

{\sc Доказательство теоремы \ref{th:sqrt-ln}.}
Нам нужно доказать, что
\begin{equation}
\sum_{k=1}^\infty \int_{(1+\epsilon)\sqrt{2\ln k}}^\infty e^{-x^2/2}\,dx<\infty,
\qquad 
\sum_{k=1}^\infty \int_{(1-\epsilon)\sqrt{2\ln k}}^\infty e^{-x^2/2}\,dx=\infty
\label{eq:for-BC}
.\end{equation}
Оценим
\begin{multline*}
J(u):=
\int_u^\infty e^{-x^2/2}\,dx= -\int_u^\infty \frac 1x e^{-x^2/2}\,d(e^{-x^2/2})
=\\=-\frac 1x e^{-x^2/2}\Bigr|_u^\infty - \int_u^\infty \frac 1{x^2} e^{-x^2/2}\,dx
.
\end{multline*}
Последнее слагаемое
 меньше, чем $u^{-2} J(u)$. Следовательно,
$$
J(u)=\frac 1u e^{-u^2/2} \bigl(1+O(u^{-2})\bigr)
\qquad \text{при $u\to\infty$.}
$$
Поэтому
\begin{multline*}
\int_{(1\pm\epsilon)\sqrt{2\ln k}}^\infty e^{-x^2/2}\,dx\sim
\frac 1{(1\pm\epsilon)\sqrt{2\ln k}} e^{-(1\pm \epsilon)^2 (2\ln k)/2}=\\=
\frac 1{(1\pm\epsilon)\sqrt{2\ln k}} k^{-(1\pm \epsilon)^2}
\end{multline*}
Теперь утверждения
(\ref{eq:for-BC}) становятся очевидными.
\hfill $\square$

\sm

{\bf\punct Линейные функционалы.}

\begin{theorem}
\label{th:lin-fun}
Пусть $\sum|a_j|^2<\infty$. Тогда ряд 
$$\sum a_j x_j$$
сходится почти всюду.
\end{theorem}

Заметим, что если все $a_j$ - ненулевые, то всюду он сходиться не может,
просто потому, что, выбирая $x_j$ подходящим образом, мы можем добиться
того, чтобы слагаемые стремились к $\infty$.

С другой стороны, сходимость ряда в смысле $L^2$ очевидна. Функции $x_j$
образуют ортонормированную систему и $\sum|a_j|^2<\infty$ -- это
просто условие сходимости ряда по ортонормированной системе.

\sm

Сходимость ряда почти всюду следует из следующей теоремы Колмогорова--Хинчина.

\begin{theorem}
Пусть $\xi_j$ последовательность независимых случайных величин
 с нулевыми математическими
ожиданиями. Пусть ряд дисперсий сходится, $\sum_j D\xi_j<\infty$.
Тогда ряд $\sum \xi_j$ сходится почти всюду.
\end{theorem}

См. \cite[\S IV.2]{Shir}.

 Далее заметим , что  для $z=(z_1,z_2,\dots)\in \ell_2$ определена
 функция
 $$
\theta_z(x)= \exp \Bigl(\sum z_j x_j\Bigr)
 .$$

 \sm

{\sc Задача.} Покажите, что эта функция интегрируема, причем,
\begin{equation}
\int  \exp \Bigl(\sum z_j x_j\Bigr)\,d\mu(x)=\exp\Bigl(\tfrac 12\sum z_j^2\Bigr).
\label{eq:exp-z-x}
\end{equation}

{\bf \punct Сдвиги. Теорема Камерона--Мартина.%
	\label{ss:cameron}} 

\begin{theorem}
Пусть $b=(b_1,b_2,\dots)\in\ell_2$, $b_j\in\R$.
Тогда сдвиг $T_b:x\mapsto x+b$ переводит меру
$\mu$ в эквивалентную ей меру, причем производная Радона--Никодима равна
\begin{equation}
\gamma(x,b)=\exp\bigl(-\sum b_j x_j-\tfrac12\sum b_j^2\bigr).
\end{equation}
\end{theorem}

Для доказательства утверждений о квазиинвариантости мер мы будем использовать
следующую лемму (которую оставляем в качестве упражнения).

\begin{lemma}
\label{l:Radon-Nyk}
Пусть $X$ -- пространство с вероятностной мерой $\xi$. Пусть $\sigma$-алгебра
пространства $X$ порождена некоторым полукольцом $\Sigma$. Пусть
$g:X\mapsto X$ -- измеримое отражение, а $\psi$ -- неотрицательная функция, такая, что
для любого измеримого $A\in \Sigma$ выполнено
$$
\xi(A)=\int_A \psi(x)\,d\xi(x).
$$ 
Тогда $g$ оставляет меру квазиинвариантной, а $\psi$ -- производная Радона--Никодима 
$g$.
\end{lemma}

{\sc Доказательство теоремы.} Прежде всего заметим, что
функция $\gamma(x,b)$ положительна, интегрируема, причем ее интеграл равен
1. Нам достаточно показать,
что для измеримых множеств $A$ выполнено
$$
\int_{A}  \gamma(x,b)\, d\mu(x)= \mu(T_bA)=\int_{b+A}\,\,d \mu(x).
$$ 
Это, в свою очередь, достаточно проверить для цилиндрических множеств,
а для них равенство превращается в тавтологию.
\hfill $\square$

\sm

Заметим, что $\R^\infty$ является абелевой группой, а мера $\mu$ оказывается
квазиинвариантной относительно плотной подгруппы. Еще заметим, что выражение для производной 
Радона--Никодима получается формальным сокращением формально бессмысленного выражения
$$
\frac{e^{-\sum(x_j+b_j)^2/2}\,\prod d(x_j+b_j)}
{e^{-\sum x_j^2/2}\,\prod dx_j}.
$$

{\bf\punct Вращения. Теорема И.Сигала.} Пусть $S$ -- бесконечная ортогональная матрица.
Применим ее к вектору-столбцу $x\in\R^{\infty}$,
$$
Sx:=\begin{pmatrix}
s_{11}&s_{12}&\dots\\
s_{21}&s_{22}&\dots\\
\vdots&\vdots&\ddots
\end{pmatrix}
\begin{pmatrix}
x_1\\ x_2\\ \vdots
\end{pmatrix}
=
\begin{pmatrix}
s_{11}x_1+s_{12}x_2+\dots\\
s_{21}x_1+s_{22}x_2+\dots\\
\vdots
\end{pmatrix}.
$$
Заметим, что для каждого $\alpha$ выполнено
$\sum_j s_{\alpha j}^2=1$, поэтому, в силу теоремы
\ref{th:lin-fun} ряды в правой части сходятся почти всюду,
т.е. для фиксированной матрицы $S$ определено почти всюду 
 отображение $x\mapsto Sx$ из $\R^\infty$ в  $\R^\infty$.
 
 \begin{theorem}
 \label{th:segal}
{\rm a)} Преобразование $x\mapsto Sx$ сохраняет меру $\mu$.

\sm

{\rm b)} Пусть $S_1$, $S_2$ -- ортогональные матрицы. Тогда для почти 
всех $x\in\R^{\infty}$ выполнено 
$$
(S_1 S_2) x= S_1(S_2x).
$$
 \end{theorem}
 
 Прямое доказательство есть в \cite{ShF}. В этих записках мы предпочтем обходной путь.
 
 \sm
 
 \begin{proposition}
Пусть $X$ --  лебеговское  пространство $X$ с непрерывной вероятностной мерой $\nu$. 
 
 \sm
 
 {\rm a)} Пусть $U$ -- ортогональный оператор в вещественном $L^2(X)$,
  переводящий неотрицательные функции в неотрицательные,
 а функцию $f(x)=1$ в себя. Тогда существует преобразование
 $S:X\to X$, переводящее меру $\nu$ в себя, такое, что
 $$
 U f(x)= f(Sx).
 $$
 
 {\rm b)} Если мы опустим условие $U1=1$, то существует преобразование
 $S:X\mapsto X$, оставляющее меру квазиинвариантной, такое, что
 $$
 Uf(x)= f(Sx)\,S'(x)^{1/2}.
 $$
  \end{proposition}

{\sc Доказательство предложения.} a) Рассмотрим множество $\cK$ функций $h(x)$, таких, что
$0\le h(x)\le 1$.  Оператор $U$ переводит это множество в себя,
поскольку и $h$, и $1-h$ должны переходить в неотрицательные функции. Множество
$\cK$ компактно в слабой топологии гильбертова пространства. Крайние точки
$\cK$ -- это характеристические функции $\chi_A$ измеримых  множеств 
$A\subset X$ (множество $A$ определено своей характеристической функцией
$\chi_A(x)\in L^2$ с точностью до множества меры ноль).
Таким образом
$$
U\chi_A=\chi_B
$$
(разумеетеся, $A$ и $B$ определены с точностью до множества меры ноль). В силу унитарности
$U$, мы имеем 
$$\nu(B)=\|\chi_B\|^2=\|\chi_B\|^2=\nu(A).$$
Обозначим $B=\sigma(A)$.

\sm

{\sc Задача.} Покажите, что операция $\sigma$ согласована с
пересечением множеств, $\sigma(A_1\cap A_2)=\sigma(A_1)\cap \sigma(A_2)$,
а также c объединением и взятием разности. 
\hfill $\lozenge$

\sm

Известно, что любая такая операция индуцируется преобразованием $X\to X$, сохраняющим меру,
см., например, \cite{Bog}, теорема  9.5.1.

\sm

b) Положим $\psi:=U1$. Отождествим $X$  с отрезком $[0,1]$.
Положим 
$$
q(x)=\int_0^x \psi^2(t)\, dt.
$$ 
Тогда унитарный оператор
$$
V f(x)=f(q(x))\,q'(x)^{1/2}
$$
переводит $1$ в $\psi$,
а, следовательно, $V^{-1} U$ переводит 1 в 1. Мы свели утверждение
b) к утверждению a).

\sm

Обозначим через $\O_{fin}(\infty)$ группу финитных ортогональных матриц $S$,
т.е. ортогональных матриц таких, что $S-1$ имеет лишь конечное число ненулевых
элементов. Группу $\bfO(\infty)$ всех ортогональных операторов мы снабдим слабой операторной топологией
(тогда группа $\O_{fin}(\infty)$ плотна в $\bfO(\infty)$).

Далее заметим, что для матриц $S$ из $\O_{fin}(\infty)$ утверждение теоремы
\ref{th:segal} очевидно.

\begin{lemma}
Представление группы $\O_{fin}(\infty)$ в $L^2(\R^\infty)$ продолжается
по непрерывности до представления группы $\bfO(\infty)$.
\end{lemma}

{\sc Доказательство леммы.} Рассмотрим следующую систему функций 
$\theta_z\in L^2(\R^\infty)$:
\begin{equation}
\theta_z(x)=\exp(\sum z_j x_j),\qquad \text{где $z\in \ell_2(\C).$} 
\label{eq:theta-z}
\end{equation}
Эта система тотальна в $L^2(\R^\infty)$.
Кроме того,
$$
\la \theta_z,\theta_u\ra=\exp\Bigl(\sum z_j \ov u_j+\tfrac12 \sum z_j^2+ \tfrac12 \sum \ov u_j^2\Bigr).
$$
Для $S\in \bfO(\infty)$ мы имеем
$$
\la \theta_{zS},\theta_{uS}\ra=\la \theta_z,\theta_u\ra.
$$
Поэтому существует унитарный оператор 
$R(S)$ такой, что для всех $z$ выполнено
$$
R(S)\theta_z=\theta_{zS}
.$$

{\sc Задача.} a) Убедитесь в правильности последней логической связки.

\sm

b) Покажите, что $R(S)$ слабо непрерывно зависит от $S$.
\hfill $\lozenge$

\sm

Это завершает доказательство леммы. \hfill $\square$

\sm

{\sc Доказательство теоремы \ref{th:segal}.}
Для финитных ортогональных матриц $S$ преобразования
$$ 
R(S)f(x^t)=f(Sx^t)
$$
переводят положительные функции в положительные, а 1 в 1. 
В силу слабой непрерывности, это выполнено и для
$S\in\bfO(\infty)$. Поэтому $R(S)$ индуцируется преобразованием 
$$
x\mapsto \tau(x)=\bigl(\tau_1(x_1,x_2,\dots),\tau_2(x_1,x_2,\dots),\dots\bigr)
, $$
сохраняющим меру. Для сохранения меры необходимо, чтобы
$ \tau_j$ были независимыми одинаково распределенными гауссовыми случайными величинами.
В частности, они образуют ортонормированную систему в $L^2(\R^\infty)$. 
Так как наше преобразование переводит функции $\theta_z$ в $\theta_{zS}$,
для любого $z\in\ell_2$ для почти всех $x$ выполнено
$$
z\tau(x)^t=(zS) x^t\quad\text{или} \quad \sum_j z_j \tau_j(x)=\sum_k\Bigl(\sum_j z_j s_{jk}\Bigr) x_k
$$
Раскрытие скобок и перегруппировки слагаемых с точки зрения сходимости почти всюду
здесь потребовали бы обоснования. Однако, мы можем рассматривать это
  равенство как равенство в $L^2$. Так как $x_k$ -- ортонормированная система в
  $L^2(\R^\infty)$, $S$ -- ограниченная матрица, а $z\in\ell_2$,
  мы можем изменить порядок суммирования
  $$
  \sum_k\Bigl(\sum_j z_j s_{jk}\Bigr) x_k=\sum_j z_j \Bigl(\sum_k  s_{jk}x_k\Bigr)
  $$
 Следовательно,
$$
z\tau(x)^t=z(Sx^t)
,$$  
и это $L^2$-равенство выполнено для всех $z$. Полагая $z=(1,0,\dots)$,
мы получаем $\tau_1(x)=\sum_j s_{1j}x_j$.
\hfill $\square$

\sm

{\bf\punct Диагональные преобразования.}

\begin{theorem}
Пусть $\lambda_k>-1$, $\sum \lambda_k^2<\infty$. 
Тогда мера $\mu$ квазиинвариантна относительно преобразований
$$
(x_1,x_2,\dots)\mapsto \bigl( (1+\lambda_1)x_1, (1+\lambda_2)x_2,\dots\bigr)
,$$
причем производная Радона--Никодима равна
\begin{equation}
\prod_{k=1}^\infty \exp\left\{-\tfrac 12 (\lambda_k^2+2\lambda_k)x_k^2\right\}\, (1+\lambda_k)
.
\label{eq:lambda-prod}
\end{equation}
\end{theorem}

{\sc Замечание.}
Производная Радона--Никодима совпадает с формально
вычисленным   отношением
\begin{equation}
\frac {\prod e^{-(1+\lambda_k)^2x_k^2/2}\prod d(1+\lambda_k)x_k}
{\prod e^{-x_k^2/2} dx_k}
.\end{equation}
Но ситуация менее безобидна, чем кажется на первый взгляд. Например, перепишем
(\ref{eq:lambda-prod})  в виде
$$
\prod_{k=1}^\infty \exp\left\{-\tfrac 12 (\lambda_k^2+2\lambda_k)x_k^2\right\}\,
\prod_{k=1}^\infty (1+\lambda_k).
$$
Тогда для сходимости выражения $\prod (1+\lambda_k)$ необходимо условие $\sum \lambda_j<\infty$,
это условие более ограничительно, чем  $\sum \lambda_j^2<\infty$.
\hfill $\lozenge$

\sm

{\sc Доказательство.} Сначала покажем, что произведение (\ref{eq:lambda-prod})
сходится почти всюду. Это тоже самое, что сходимость почти всюду ряда из 
логарифмов
$$
\sum_k\left(-\tfrac 12 (\lambda_k^2+2\lambda_k)x_k^2+\ln(1+\lambda_k)\right).
$$
Имея в виду тейлоровское разложение $\ln(1+\lambda_k)$, мы запишем
$$\ln(1+\lambda_k)=\lambda_k-\gamma_k.$$
Величина $\gamma_k$ положительна, и при достаточно малых
$\lambda_k$ оценивается сверху как  $\gamma_k\le \lambda_k^2$. Наш ряд преобразуется к виду
$$
\sum_k\left(-\tfrac 12 \lambda_k^2 x_k^2- \lambda_k (x_k^2-1)- \gamma_k x_k^2
\right).
$$
Это выражение есть сумма трех сходящихся рядов. Первый и третий ряд сходятся почти всюду
абсолютно по теореме Беппо Леви. Далее, $(x_k^2-1)$ -- одинаково распределенные независимые
случайные величины с нулевым средним. Поэтому ряд
$\sum \lambda_k (x_k^2-1)$ сходится почти всюду по теореме Колмогорова--Хинчина.
Сходимость произведения почти всюду доказана.

Обозначим $k$-ый множитель произведения (\ref{eq:lambda-prod})
  через $\Phi_k$. По построению,
  \begin{multline*}
  \frac 1{\sqrt{2\pi}}\int_\R\Phi_k(x)\, e^{-x^2/2}\,dx
  = \frac 1{\sqrt{2\pi}}\int_\R
  \frac{e^{-(1+\lambda_k)^2x^2/2}(1+\lambda_k)}{e^{-x^2/2}}
  \, e^{-x^2/2}\,dx
    =\\=
  \frac 1{\sqrt{2\pi}}\int_\R e^{-(1+\lambda_k)^2x^2/2} \,d (1+\lambda_k)x
  = \frac 1{\sqrt{2\pi}}\int_\R  e^{-x^2/2}\,dx=1.
  \end{multline*} 
  Аналогично, для любого измеримого подмножества $A\subset \R$ выполнено
  \begin{equation}
  \frac 1{\sqrt{2\pi}}\int_{(1+\lambda_k)A}\Phi_k(x)\, e^{-x^2/2}\,dx
  =\frac 1{\sqrt{2\pi}}\int_A  e^{-x^2/2}\,dx.
  \label{eq:intA}
  \end{equation}
  
  \begin{lemma}
  \label{l:diag-L1}
  Последовательность $\prod_{k=1}^n\Phi_k(x_k)$ сходится в смысле $L^1(\R^\infty)$.
  \end{lemma}
  
{\sc Вывод теоремы из леммы \ref{l:diag-L1}.}  В силу леммы, интеграл от  $\prod_{k=1}^\infty\Phi_k(x_k)$
равен 1. Равен 1 и интеграл $\prod_{k=m+1}^\infty\Phi_k(x_k)$ по 
пространству $\R^\infty$ последовательностей $(x_m,x_{m+1}, \dots)$. 

  Пусть $A_1$, \dots, $A_m\subset \R$ -- измеримые множества.
  В силу (\ref{eq:intA}),
\begin{multline*}
\mu( (1+\lambda_1)A_1\times\dots\times (1+\lambda_m) A_m)
\prod _{k=1}^m \int_{A_k} \Phi_k(x_k)\,dx_k=
=\\=
\int\limits_{A_1\times\dots\times A_m}\prod_{k=1}^m\Phi_k(x_k)e^{-\frac12\sum x_k^2}\,dx_1\dots dx_k\cdot 1=\\=
\int\limits_{A_1\times\dots\times A_m\times\R\times \R\times\dots} \prod_{k=1}^m\Phi_k(x_k)
\prod_{k=m+1}^\infty\Phi_k(x_k)\,d\mu
.
\end{multline*} 
Теперь мы можем сослаться на лемму  \ref{l:Radon-Nyk}.
\hfill $\square$

\sm

Лемма 
\ref{l:diag-L1} непосредственно вытекает из следующей леммы.

\begin{lemma}
 \label{l:diag-L2}
  Последовательность $\prod_{k=1}^n\Phi_k(x_k)$ сходится в смысле $L^2(\R^\infty)$.
\end{lemma}

{\sc Доказательство леммы \ref{l:diag-L2}.} Мы покажем (<<в лоб>>), что последовательность
 частичных произведений
фундаментальна в $L^2(\R^{\infty})$. Вычислим
\begin{multline}
\label{eq:long-Phi}
\|\prod_{k=1}^n\Phi_k(x_k)-\prod_{k=1}^{n+m}\Phi_{k}(x_k)\|^2=
\\=
(2\pi)^{-\frac{n+m}2}\int_{\R^{n+m}} \prod_{k=1}^n\Phi_k(x_k)^2\,
\left\{\prod_{k=n+1}^{n+m}\Phi_k(x_k)-1  \right\}^2\,
e^{-\sum_{k=1}^{n+m} x_k^2/2}\prod_{k=1}^{n+m}dx_k
=\\=
(2\pi)^{-\frac{n}2}\int_{\R^{n}}
\prod_{k=1}^n\Phi_k(x_k)^2 e^{-\sum_{k=1}^{n} x_k^2/2}\prod_{k=1}^{n}dx_k
\times\\\times
(2\pi)^{-\frac{m}2}\int_{\R^{m}}
\left(\prod_{k=n+1}^{n+m}\Phi_k(x_k)^2-1  \right)\,
e^{-\sum_{k=n+1}^{n+m} x_k^2/2} \prod_{k=n+1}^{n+m}dx_k
,\end{multline}
после раскрытия больших фигурных скобок мы использовали (\ref{eq:intA}). Обозначим
$$
I_k:=\frac 1{\sqrt {2\pi}}\int \Phi_k(x_k)^2 e^{{-x_k}^2/2}\,dx_k
$$
Выражение 
(\ref{eq:long-Phi}) тогда записывается как
\begin{equation}
\prod_{k=1}^n I_k \cdot \left( \prod_{k=n+1}^{n+m} I_k-1 \right).
\label{eq:IIII}
\end{equation}
С другой стороны, $I_k$ -- гауссовский интеграл, он равен
$$
I_k=\frac{(1+\lambda_k)^2}{\sqrt{1+4\lambda_k+\lambda_k^2}}=1+O(\lambda_k)^2
$$
В силу этой асимптотики, первый множитель (\ref{eq:IIII}) имеет конечный
 предел при $n\to\infty$. Выражение
 $\prod_{k=n+1}^{n+m} I_k$ является хвостом сходящегося произведения, а поэтому стремится к 1.
 Поэтому (\ref{eq:long-Phi}) стремится к нулю при $n\to\infty$.
 \hfill $\square$

\sm

{\bf \punct Линейные преобразования.%
\label{ss:GLO}} Обозначим через $\GLO(H)$ группу
обратимых линейных преобразований $A$ вещественного гильбертова пространства $H$,
представимых в виде
 $$A=U(1+T),$$
 где $U$ -- ортогональный оператор, а $T$ -- оператор Гильберта--Шмидта.
 
\sm 
 
 {\sc Задача.}  Убедитесь, что это группа. 
 \hfill $\lozenge$
 
\begin{lemma}
Любой элемент $A$ группы $\GLO(\ell_2)$ представим в виде
$A=V_1 (1+\Lambda) V_2$, где $V_1$, $V_2$ ортогональны, а $\Lambda$ -- диагональная
гильберт-шми\-дтов\-ская матрица
с собственными числами $>-1$. 
\end{lemma}

{\sc Доказательство.}
Рассмотрим полярное разложение оператора $A\in \GLO(\ell_2)$, $A=U P$, где $P$ самосопряжен
и положителен, а $U$ ортогонален. Тогда
$$
P^2=P^*P=(1+T^*)U^*U(1+T)=(1+T^*)(1+T)=1+T^*+T+T^*T
$$
Оператор $T^*+T+T^*T$ является гильберт-шмидтовским. Приведем его к диагональному виду
$T^*+T+T^*T= W D W^{-1}$. Тогда $P=U \sqrt{1+ D} W^{-1}$, а $A=VW \sqrt{1+ D} W^{-1}$.
\hfill $\square$

\sm

В качестве следствия получаем такую теорему
(она обычно называется теоремой Фельдмана--Гаека).

\begin{theorem}
Группа $\GLO(\ell_2)$ действует на  $L^2(\R^\infty)$ преобразованиями
$x\mapsto Ax$, оставляющими меру $\mu$ квазиинвариантной.
\end{theorem}

{\sc Доказательство.} Мы разлагаем $A$ как в лемме, $A=V_1 (1+\Lambda) V_2$
и рассматриваем композицию преобразований
$$
x\mapsto V_2 x\mapsto (1+\Lambda) V_2x\mapsto V_1 (1+\Lambda) V_2x
.
$$
Понятно, что полученное преобразование оставляет меру квазиинвариантной. Нужно лишь
проверить, что для почти любого вектора $x\in \R^\infty$ выполнено
\begin{equation}
 V_1 \bigl((1+\Lambda) (V_2x)\bigr) = \bigl( V_1 (1+\Lambda) V_2\bigr)x.
\label{eq:(((())))}
\end{equation}
Обе части равенства имеют смысл п.в. С другой стороны, $x_j$ мы можем рассматривать
как элементы гильбертова пространства $L^2(\R^\infty)$. Тогда все ряды,
возникающие в вычислениях (\ref{eq:(((())))}) являются абсолютно сходящимися
в смысле $L^2$ и мы их
можем перегруппировывать
\hfill $\square$

\sm

{\sc Замечание.}
Производную Радона--Никодима преобразования $A=U(1+T)$, где $T$ --ядерный,
 можно получить спомощью формальной манипуляции 
\begin{multline*}
\frac{e^{-\la Ax,Ax \ra/2 } \,d(Ax)}
{e^{-x^2/2}\,dx}= \frac{e^{-\la U(1+T)x,U(1+T)x \ra/2 } \,d\bigl(U(1+T)x\bigr)}
{e^{-x^2/2}\,dx}
=\\=
e^{-\bigl\la (T+T^*+T^*T)x,x \bigr\ra/2} |\det\bigl(U(1+T)\bigr)|
.
\end{multline*}
Пусть оператор $T$ ядерный, тогда множители полученного выражения имеют смысл.
Несложно показать, что первый множитель -- это функция, лежащая в $L^1(\R^\infty)$
(для этого достаточно привести квадратичную форму 
$\bigl\la (T+T^*+T^*T)x,x \bigr\ra$ к главным осям и проинтегрировать ее). 
Далее, напомним, что для матриц вида $1+H$, где $H$ -- ядерная, 
определитель $\det(1+H)$ корректно определен. Мы полагаем
\begin{equation}
|\det\bigl(U(1+T)\bigr)|=\sqrt{\det \bigl((1+T)^*(1+T)\bigr)}.
\label{eq:det-GLO}
\end{equation}
Подчеркнем, что имеет смысл только модуль определителя. Знак определителя 
в данном случае
определить невозможно.

\sm

{\sc Задача.} Рассмотрим группу обратимых операторов в $\ell_2$,
представимых в виде $U(1+T)$, где $T$ -- ядерный оператор. Покажите, что функция
$|\det(\cdot)|$, определенная формулой (\ref{eq:det-GLO}),
удовлетворяет тождеству
$$
\qquad\qquad\qquad\qquad|\det(A)|\cdot |\det(B)|=|\det(AB)|.\qquad\qquad\qquad\qquad
\lozenge
$$

\sm

{\bf\punct Неисправимость действия.%
\label{ss:individ}} Еще раз напомним, что наши преобразования были
определены почти всюду, почти всюду были выполнены  и равенства
$(AB)x=A(Bx)$. Назовем действие группы $G$ на пространстве с мерой
$(X,\xi)$ преобразованиями, сохраняющими меру, {\it индивидуальным},
если преобразования $x\mapsto gx$ являются измеримыми по Борелю,
 определены для всех $x\in X$ и
равенство $g_1(g_2 x)=(g_1g_2)x$ выполнено для всех $x$. Оказывается, что наше
действие группы $\bfO(\infty)$ не может быть исправлено и превращено в индивидуальное
действие. Более того, есть такая теорема \cite{Gla}.

\begin{theorem} 
{\rm a)} Группа $\bfO(\infty)$ всех ортогональных матриц
не имеет нетривиальных индивидуальных действий на лебеговских пространствах с мерой.

\sm

{\rm b)} 
Рассмотрим группу $\O_2(\infty)$ ортогональных операторов вида $1+T$, где $T$
-- оператор Гильберта--Шмидта. 
Группа $\O_2(\infty)$ не имеет нетривиальных индивидуальных действий на лебеговских пространствах с мерой.
\end{theorem}

\begin{corollary}
 Полная бесконечномерная унитарная и полная бесконечномерная
симплектическая группа, а также группа $\Ams$ преобразований пространства с мерой
{\rm(}см. п.{\rm\ref{ss:Ams})}
не имеют нетривиальных индивидуальных действий на лебеговских пространствах с мерой.
\end{corollary}

В самом деле, эти группы содержат подгруппы, изоморфные $\bfO(\infty)$.
 
Не имеет индивидуальных инвариантных действий и группа измеримых отображений
отрезка в окружность $\R/2\pi \Z$.
 
\sm

{\bf\punct Функциональные пространства.%
	\label{ss:functional-spaces}} Рассмотрим абстрактное гильбертово
пространство $H$, рассмотрим базис $e_j$ и пространство $\wt{ (H,e_j)}$ ({\it сигаловское
расширение} пространства $H$), состоящее
из формальных выражений $\sum x_j e_j$, где $x_j\in \R$ произвольны.
Это пространство фактически является пространством $\R^\infty$, введем на нем меру, как
и выше. Формально пространство $\wt{ (H,e_j)}$ зависит от выбора базиса.
Однако, два пространства $\wt{ (H,e_j)}$, $\wt{ (H,e_j')}$ с мерой
канонически изоморфны, изоморфизм отождествляется ортогональным преобразованием
переводящий один базис в другой. Т.е., мы получаем канонически определенное
пространство с мерой $\wt H$.

Оказывается, что для конкретных функциональных пространств $H$ сигаловские
 расширения $\wt H$ допускают 
прозрачные описания. Сначала напомним некоторые стандартные факты из анализа.

\sm

{\sc Задача.} a) Пусть $\chi$ -- обобщенная функция на окружности $S^1=\R/2\pi \Z$.
 Покажите, что ее коэффициенты Фурье
 $$
 c_n=\frac 1{2\pi} \int_{S^1}\chi(\phi)\, e^{-in\phi}\,d\phi
 $$
имеют не более чем полиномиальный рост, т.е., существует $k$, такое, что
$c_n=O(n^k)$ при $n\to \pm\infty$. 

\sm

b) Для любой последовательности
$c_n$ полиномиального роста ряд
$\chi(\phi)=\sum c_n e^{-in\phi}$ сходится в смысле обобщенных функций.

\sm

c) Покажите также, что эти две операции взаимно обратны.
\hfill $\lozenge$

\sm

Напомним, что в связи с этим вводится шкала соболевских пространств $H^s(S^1)$,
где $s\in \R$.
Пространство $H^s(S^1)$ состоит из  функций,
у которых коэффициенты Фурье удовлетворяют условию
$$
\sum c_n^2 (1+n^2)^s<\infty.
$$
При $s\ge 0$ элементы этих пространств являются настоящими функциями, а при
$s<0$ -- обобщенными функциями. Отметим, что $H^n$, где $n\in\N$ состоит
из функций, $n$-ная производная которых лежит в $L^2$. А $n$-ная производная
функции из $L^2$ лежит в $H^{-n}$.  

\sm

{\it Белый шум.} Рассмотрим вещественное пространство $L^2[0,\pi]$ на отрезке.  
Сигаловское расширение (в данном случае оно называется белым шумом) состоит из рядов
$$
\sum_{n>0}c_n \sin n\phi.
$$
Для произвольной последовательности $c_n$ эта сумма не имеет смысла.
Но в силу предложения \ref{pr:beppo}, ряд
$$
\sum c_n^2 (1+n^2)^{-1/2-\epsilon}
$$
сходится почти всюду на $\wt{L^2[0,\pi]}$. Выбрасывая множество меры ноль, мы 
можем положить, что точки сигаловского расширения --  обобщенные функции
соболевского класса $H^{-1/2-\epsilon}$.  

\sm

{\it Броуновское движение {\rm(}мера Винера{\rm)}.}
Возьмем пространство вещественных функций $W$ на отрезке $[0,1]$,
таких, что
\begin{equation}
\int_0^1 f'(x)^2 \,dx <\infty,\qquad f(0)=0.
\label{eq:broun1}
\end{equation}
Введем в нем скалярное произведение
\begin{equation}
\la f_1,f_2\ra=\int_0^1 f_1'(x)\,f_2'(x)\,dx
.
\label{eq:broun2}
\end{equation}
Мерой Винера называется стандартная гауссова мера на $\wt W$.
В качестве ортонормального базиса в $W$ можно взять систему
$$
\frac{1}{\sqrt {2 \pi n}}\sin(2\pi n x).
$$ 

{\sc Задача.} a) Пространство $\wt W$ можно отождествить с соболевским пространством
$H^{1/2-\epsilon}$. Оператор дифференцирования $\frac d{dx}$ корректно 
определен как отображение $\wt W\to \wt{L^2}$ и является изоморфизмом пространств с мерой.

\sm

b) Пользуясь упомянутым базисом, покажите, что для фиксированного $x_0$
для почти всех $f\in \wt W$ определено значение в точке $f(x_0)$. 
\hfill $\lozenge$

\sm

Дальнейший розыгрыш этого подхода есть в \cite{ShF}, \S4.
Мы, для разнообразия, обсудим  другой путь.  Напомним сначала определение
ортогональной 
{\it системы Хаара}. Разобьем отрезок $[0,1]$ на $2^n$ равных кусков
$\left[\frac{a}{2^n},\frac{a+1}{2^n}\right]$. Рассмотрим функцию
$\chi$, которая на левой половинке этого отрезка равна $2^{n/2}$,
на правой половинке равна $-2^{n/2}$, а вне отрезка равна 0 (в итоге
$L^2$-норма $\chi$ равна 1). Система Хаара -- это набор всех таких функций
по всем уровням $n\ge 0$ и всем отрезкам $\left[\frac{a}{2^n},\frac{a+1}{2^n}\right]$.
Этот набор функций образует ортонормированный базис в $L^2[0,1]$. 
Занумеруем функции Хаара, перебирая слева направо отрезки уровней $n=0$, 1, 2, 3, \dots.

Первообразные $S_k(x)$ от функций Хаара образуют ортонормированный базис в
$W$. Графики этих функций изображены на Рис. \ref{fig:haar}.

\begin{figure}
\epsfbox{haar.1}  
\caption{Элемент ортогонального базиса в $H^1[0,1]$. \label{fig:haar}}
\end{figure}

Определим пространство $\wt W$ как пространство сумм рядов $\sum c_k S_k(x)$,
где $c_k$ -- независимые гауссовы переменные  распределенные с плотностью $e^{-t^2/2}$,
т.е., $c=(c_1,c_2,\dots)$ -- точка нашего пространства $\R^\infty$

\begin{theorem}
\label{th:brown-continuous}
Почти все элементы пространства $\wt W$ являются непрерывными функциями. 
\end{theorem}  

В силу теоремы \ref{th:sqrt-ln}, мы можем положить
$|c_k|\le C\sqrt {\ln k}$.
Теорема \ref{th:brown-continuous} следует из следующей леммы.

\begin{lemma}
Пусть $|c_k|\le C\sqrt {\ln k}$. Тогда ряд $\sum c_k S_k(x)$ сходится
равномерно на отрезке.
\end{lemma}

{\sc Доказательство леммы.} Рассмотрим сумму
$$
F_n(x)=\sum_{k=2^n}^{2^{n+1}-1} c_k S_k(x)
$$
(суммирование взято по всем функциям $S_k$ уровня $n$). Носители слагаемых 
не перекрываются, поэтому на каждом отрезке $[a 2^{-n},(a+1)2^{-n}]$
сумма равна $c_k S_k(x)$, где $S_k(x)$ -- базисная функция, соответствующая
этому отрезку. В частности, 
$$
\max| F_n(x)|= 2^{-n/2-1} \max_{2^n\le k \le 2^{n+1}-1} |c_k|.
$$
С другой стороны, $c_k$ растет примерно как $\sqrt{2\ln k}$.
Таким образом, ряд $\sum_n F_n(x)$ имеет сходящуюся мажоранту.
\hfill $\square$

\sm

{\sc Задача.} Докажите, что есть константы
$C_1$, $C_2$, такие, что для почти каждого элемента $\wt W$
для всех $x_1$ $x_2\in [0,1]$ выполнено следующее условие гёльдеровости
$$
\qquad \qquad\qquad
|f(x_1)-f(x_2)|\le |x_1-x_2|^{1/2} \bigl(C_1+C_2\bigl|\ln|x_1-x_2|\bigr|\bigr)
.
\qquad \qquad\qquad
\lozenge
$$ 

Теория броуновского движения излагается во многих книгах, мы же на этом месте остановимся.

\sm

{\bf\punct Примеры действий групп на пространствах с гауссовой мерой.}
 Пусть $G=\Diff[0,1]$ -- группа
группа $C^\infty$-диффеоморфизмов отрезка $[0,1]$, сохраняющих ориентацию. 

Группа $G$ действует ортогональными преобразованиями вещественного пространства
$L^2[0,1]$ по формуле
$$
T_0 (q) f(x)=f(q(x))\,q'(x)^{1/2}
,$$
где $q\in G$.
Поэтому она действует на белом шуме $\wt {L^2}$ преобразованиями, сохраняющими меру.

Более интересный пример. Рассмотрим
пространство $W$ функций с $f'\in L^2[0,1]$ из предыдущего
пункта. Опеределим в нем действие той же группы $G$, заданное формулой
\begin{equation}
T_1(q) f(x)=f(q(x))\,q'(x)^{-1/2}
.
\label{eq:T1}
\end{equation}

\begin{theorem}
Операторы $T_1(q)$ содержатся в группе $\GLO(W)$. Следовательно преобразования
$T_1(q)$ оставляют меру Винера квазиинвариантной.
\end{theorem} 

{\sc Задача.}
Оператор $A$ содержится в $\GLO(H)$ тогда и только тогда, когда 
$A^*A-1$ является оператором Гильберта--Шмидта (рассмотрите полярное разложение $A$).
\hfill $\lozenge$

\sm

{\sc Доказательство теоремы.} Мы покажем, что оператор $T_1(q)^*T_1(q)-1$
является ядерным (фактически, доказательство дает класс Шаттена $\cL^{1/2+\epsilon}$).
Вычислим соответствующую квадратичную форму,
\begin{multline*}
J(f_1,f_2):=\bigl\la(T_1(q)^*T_1(q)-1)f_1,f_2\bigr\ra_W=
\bigl\la T_1(q) f_1, T_1(q) f_2\bigr\ra_W-\la f_1,f_2\ra_W 
=\\=
\int_0^1 \left[f_1(q(x))\,q'(x)^{-1/2}\right]'\cdot \left[f_2(q(x))\,q'(x)^{-1/2}\right]'\,dx
-\int_0^1 f_1'(x)\, f_2'(x)\,dx.
\end{multline*}
После выполнения дифференцирования в первом интеграле получится четыре слагаемых.
Одно из них имеет вид
$$
\int_0^1 f_1'(q(x))\,f_2'(q(x)) q'(x)\,dx=\int_0^1 f_1'(y)\, f_2'(y)\,dy.
$$
Здесь мы сделали замену 
$$y=q(x),\quad x=r(y),$$ где $r(y)$ -- обратная функция. Это слагаемое
сокращается. Три  других
слагаемых дадут нам
\begin{multline*}
J(f_1,f_2)=
\int_0^1 \Bigl(f_1'(q(x))\,f_2(q(x))+f_1(q(x)) f_2'(q(x)) \Bigr) q'(x) \cdot \Bigl(-\frac 12 \frac{q''(x)}
{q'(x)^2} \Bigr)\,dx+\\+
\int_0^1 f_1(q(x)) f_2(q(x))\cdot \frac14 \frac{q''(x)^2}{q'(x)^3}\, dx
.
\label{eq:int-int}
\end{multline*}
Первый интеграл мы берем по частям и получаем,
\begin{equation}
J(f_1,f_2)=
\frac{f_1(1) f_2(1) q''(1)}{2 q'(1)^2}+ 
\frac 12\int_0^1 f_1(q(x))\,f_2(q(x)) Sq(x) \,q'(x)^{-1}\,dx
\label{eq:with-swartz} 
,\end{equation}
где $Sq(x)$ -- производная Шварца,
$$
Sq(x)= \frac{q'''(x)}{q'(x)}-\frac{3 q''(x)^2}{q'(x)^2}.
$$ 
Далее мы переходим к переменной $y=q(x)$. Для этого используем свойство коцикла
для производной Шварца,
$$
S(\tau\circ\sigma)=(S\tau\circ \sigma)\cdot (\sigma')^2+ S\sigma.
$$
В этом равенстве мы положим $\sigma=q$, $\tau=r$. Т.к. $r(q(x))=x$,  
в левой части мы получаем ноль. Поэтому (\ref{eq:with-swartz})
принимает вид
$$
J(f_1,f_2)=-\frac{r(1)''}{2r'(1)'^2}f_1(1)f_2(1) -\frac12 \int_0^1 f_1(y)\,f_2(y)\, Sr(y)\,dy.
$$

Мы замечаем, что первое слагаемое -- это квадратичная форма форма ранга 1.
Далее, квадратичная форма
$$
H(f_1,f_2):=\int_0^1  f_1(y)\,f_2(y)\,dy
$$
является ядерной в $W$, ее собственные числа в тригонометрическом базисе
образуют последовательность $1/n^2$. Оператор умножения на гладкую функцию $Sr(y)$
  ограничен
в $W$, а поэтому и форма   $J(f_1,f_2)$ является ядерной.
\hfill $\square$

\sm

{\sc Задача.} Рассмотрим пространство $Y$ функций на $[0,1]$, удовлетворяющих условию
$f(0)=f'(0)=0$ со скалярным произведением 
$$
\la f_1, f_2\ra=\int_0^1 f_1''(x)\,f_2''(x)\,dx.
$$
Рассмотрим его сигаловское расширение%
\footnote{Иными словами, мы рассматриваем движения,
скорость которых является броуновской.}
 $\wt Y$. Покажите, что 
преобразования
$$
T_2(q) f(x)=f(q(x))\,q'(x)^{-5/2}
$$ 
лежат в группе $\GLO(Y)$.
\hfill $\lozenge$

\sm

{\sc Задача.} Рассмотрим $n$-мерное броуновское движение. Это
сигаловское расширения гильбертова пространства $H$ функций
$[0,1]\to \R^n$, удовлетворяющих условию  
(\ref{eq:broun1}), со скалярным произведением
(\ref{eq:broun2}). Пусть $g$ - гладкая функция из $[0,1]$ в ортогональную
группу $\O(n)$. Покажите, что оператор
$$
A(g) f(x)= g(x) f(x)
$$
лежит в группе $\GLO(H)$. \hfill $\lozenge$

\sm

См. также выше пп.\ref{ss:diff1}, \ref{ss:hier}.

\sm
  
Вернемся к общим рассмотрениям гауссовых мер.

\sm

{\bf \punct Многочлены Эрмита.} Рассмотрим пространство $L^2$ на прямой
по гауссовой мере $e^{-x^2/2}dx$. Напомним, что многочлены Эрмита
$H_k(x)$ -- это многочлены, полученные ортогонализацией Грамма--Шмидта 
из последовательности $1$, $x$, $x^2$, \dots. Мы не будем здесь напоминать большего,
лишь отметим формулу производящей функции для многочленов  $H_k(x)$,
$$
\sum_{k=0}\frac{z^k}{k!} H_k(x)=\exp\left(\sqrt 2 xz-z^2\right)
.$$
Отметим, что система $H_k$ ортогональна, но не ортонормирована, 
$$
\|H_k\|^2=k!
$$
 
 Пусть $\bfk:=(k_1, k_2, \dots)$ -- последовательность неотрицательных
целых чисел, причем $k_j=0$ с некоторого номера. Положим
$$
H_{\bf k}(x):=\prod_{j=1}^\infty H_{k_j}(x_j) 
$$
(фактически это произведение конечно, т.к. $H_0(x)=1$).

\begin{proposition}
Функции $H_{\bf k}(x)$ образуют ортогональный базис в пространстве
$L^2(\R^\infty)$,
причем
$$
\|H_{\bf k}\|^2=\bfk!:=\prod k_j!
$$
\end{proposition} 

Ортогональность и формула для нормы здесь очевидны. Многочлены зависящие
лишь от первых $n$ переменных, образуют базис в $L^2(\R^n,e^{-x^2/2})$, а объединение
таких подпространств плотно в $L^2(\R^\infty)$. Это влечет полноту системы
$H_{\bf k}$. \hfill $\square$.

\sm

Пусть $z=(z_1,z_2,\dots)$ -- элемент комплексного $\ell_2$. Определим
функцию
\begin{equation}
\Psi_z(x)=e^{\sqrt 2 x z^t- zz^t}=\exp\left(\sqrt 2\sum x_j z_j-\sum z_j^2\right)
\label{eq:psi-fock}
.
\end{equation}
 Эта функция очевидным образом разлагается по многочленам Эрмита
 \begin{equation}
 \Psi_z(x)=\sum_{\bfk}\frac{z^\bfk}{\bfk!} H_{\bf k}
 \label{eq:psi-fock-1}
 .
 \end{equation}
 Это ряд по ортогональной системе, и он сходится в среднем квадратичном.
Легко видеть также, что
\begin{equation}
\la  \Psi_z,  \Psi_u\ra_{L^2(\R^\infty)}=e^{z\ov u^t}.
\end{equation}

{\bf\punct Голоморфные ряды.} Рассмотрим комплексную плоскость $\C$  с гауссовой
мерой $\frac 1\pi e^{-|\xi|^2}d\xi\,d\ov \xi$. Рассмотрим счетное произведение
этого пространства с мерой на себя, обозначим его через $\C^\infty$.

Понятно, что  таким способом мы не получаем нового объекта, $\C\simeq\R\times \R$,
а $\C^\infty$ -- переобозначенное пространство $\R^\infty$ (или $\C^\infty$ ).

Пусть $\bfk=(k_1,k_2,\dots)$ -- такое же, как и  в предыдущем пункте.
Рассмотрим функции 
$$
\xi^\bfk:=\prod \xi_j^{k_j}.
$$ 

\begin{proposition}
Функции $z^\bfk$ образуют {\rm(}неполную{\rm)}
ортогональную систему  в $L^2(\R^\infty)$, причем
$$
\|\xi^\bfk\|^2=\prod k_j!
$$
\end{proposition}

{\sc Доказательство.}
$$
\la \xi^\bfk, \xi^\bfl\ra= \prod_j \frac1\pi\int_\C \xi_j^{k_j}\ov \xi_j^{l_j}\, e^{-|\xi_j|^2}
\,d\xi_j\,d\ov\xi_j
$$
(фактически произведение конечно). Достаточно вычислить интеграл
$$
\frac1\pi\int_\C \xi^{k}\ov \xi^{l}\, e^{-|\xi|^2}
\,d\xi\,d\ov\xi
=
\frac1\pi \int_0^\infty \int_0^{2\pi} r^{k+l+1} e^{i(k-l)\phi}e^{-r^2}\,d\phi\, dr
,$$
и дальнейшее вычисление очевидно. \hfill $\square$

\sm

Обозначим через $\cF_\infty$ подпространство в $L^2(\C^\infty)$, порожденное
ортогональной системой $\xi^\bfk$, т.е. пространство всех рядов
\begin{equation}
F(\xi)=\sum c_{\bfk} \xi^\bfk, \qquad\text{где $\sum \bfk!\, |c_{\bfk}|^2<\infty$}.
\label{eq:ryad}
\end{equation}

{\sc Задача} (см. \cite[\S 4.1]{Ner2}) Рассмотрим такой же объект $\cF_n$, только вместо $\C^\infty$ рассмотрим
$\C^n$. Покажите, что $\cF_n$ -- это пространство всех голоморфных функций
на $\C^n$, удовлетворяющих условию
$$
\qquad\qquad\qquad\qquad\qquad
\int_{\C^n}|f(\xi)|^2 e^{-|\xi|^2}\,d\xi\,d\ov\xi<\infty.
\qquad\qquad\qquad\qquad\qquad
\lozenge
$$

Рассмотрим в $\cF_\infty$ систему функций
\begin{equation}
\label{eq:Phi-fock}
\Phi_z(\xi)=e^{\xi \ov z}=\exp\left(\sum \xi_j \ov z_j\right)
=\sum_{\bfk} \frac {\ov z^\bfk}{\bfk!} \xi^\bfk.
\end{equation}

Легко видеть, что
$$
\la \Phi_z, \Phi_u\ra=e^{u\ov z}.
$$

\begin{theorem}
Пусть $u\in \ell_2$. Тогда ряд $F$, заданный формулой {\rm(\ref{eq:ryad})},
т.е.,
$$
\sum_\bfk c_{\bfk} u^\bfk
$$
 сходится абсолютно в точке
$\xi=u$. Более того,
$$
F(u)=\la F,\Phi_u\ra_{L^2(\C^\infty)} 
.$$
\end{theorem}

{\sc Доказательство.} Посчитаем скалярное произведение $\la F,\Phi_u\ra$,
используя разложение этих функций по базису $\xi^\bfk$,
$$
\la F,\Phi_u\ra=\la \sum c_{\bfk} \xi^\bfk, \sum_{\bfk} \frac {\ov u^\bfk}{\bfk!} \xi^\bfk\ra
=\sum_\bfk c_{\bfk} \frac {u^\bfk}{\bfk!} \la\xi^\bfk,\xi^\bfk, \ra
=\sum_\bfk c_{\bfk} u^\bfk
.$$
Ряд этот сходится абсолютно по построению (так как этот ряд дает скалярное произведение двух векторов). Утверждение доказано.
\hfill$\square$

\sm

{\sc Задача.} Докажите, что этот ряд сходится равномерно на любом шаре
в $\ell_2$. \hfill $\lozenge$

\sm

Итак, каждой функции пространства $\cF_\infty$ мы поставили в соответствие
голоморфную функцию на пространстве $\ell_2$. Полученное  пространство  функций
называется {\it бозонным пространством Фока}, см. \cite[\S VI.1]{Ner1}. 
Этот объект уже находится за пределами
настоящих записок (впрочем, он фактически идентичен пространству $\cF_\infty$).

\sm

{\sc Замечание.} Хочется сказать, что мы ограничили функцию на $\C^\infty$
на $\ell^2$. Но формально это неправильно, потому что $F(\xi)$ определена с точностью до почти
всюду, а $\ell_2$ имеет меру нуль в $\C^\infty$. \hfill $\lozenge$

\sm

{\bf\punct Преобразование Сигала--Баргмана.} Напомним, что мы ввели в
$L^2(\R^\infty)$ систему функций $\Psi_z$, а в $\cF_\infty$ - систему функций
$\Phi_z$, см. (\ref{eq:psi-fock}), (\ref{eq:Phi-fock}). Обе системы нумеруются
параметром $z\in\ell_2$.

\begin{theorem}
{\rm a)} Существует единственный унитарный оператор 
{\rm (преобразование Сигала--Баргмана)}
$$\cS:L^2(\R^\infty)\to \cF_\infty,$$
такой, что
$$
\cS \Psi_z=\Phi_z\qquad\text{для всех $z$}
.$$

{\rm b)} Оператор $\cS$ переводит многочлены Эрмита $H_\bfk$ в одночлены
$\xi^\bfk$.

\sm

{\rm c)} Оператор $\cS$ коммутирует с действием ортогональной группы
$\bfO(\infty)$.
\end{theorem}

Ортогональная группа действует в $\C^\infty$ умножением ортогональной матрицы на столбец
$\xi\in\C^\infty$. Корректность этой операции не нужно отдельно доказывать, т.к.
$\C^\infty\simeq\R^{2\infty}$, в этом пространстве действует группа 
$\bfO(2\infty)$, а $\bfO(\infty)$ -- ее подгруппа.

\sm

{\sc Доказательство.} Обе системы  $\Psi_z$, $\Phi_z$ тотальны. При этом
$$
\la \Psi_z, \Psi_u\ra= e^{u\ov z}= \la \Phi_z, \Phi_u\ra.
$$
Это влечет утверждение a). 

Мы также можем рассмотреть унитарный оператор $\cS'$,
определенный равенством
$$
\cS' H_\bfk=\xi^k.
$$
Тогда в силу разложений (\ref{eq:psi-fock}),
(\ref{eq:Phi-fock}), мы имеем $\cS'\Psi_z=\Phi_z$. Мы получаем тот же оператор.
\hfill $\square$

\sm

Объясним, как получить голоморфную функцию $F$ из $f\in L^2(\R^\infty)$.
В силу унитарности оператора $\cS$,
$$
F(z)=\la F,\Phi_z\ra_{\cF_\infty}=\la f, \Psi_z\ra_{L^2(\R^\infty)}
$$
т.е., преобразование задается формулой
$$
F(z)=\int_{\R^\infty} f(x)\exp\left(\sqrt 2 xz-z^2\right)\,d\mu(x).
$$

Когда-то это преобразование вводилось Ирвингом Сигалом для потребностей квантовой теории
поля, но оно оказалось важным и содержательным объектом в конечномерном (и в одномерном) 
случае тоже. О конечномерном случае, см. \cite[Ch. 4]{Ner2}.

\sm

{\bf \punct Разложение представления ортогональной группы в $L^2(\R^\infty)$.%
	\label{ss:symmetric-sum}}
Рассмотрим в $\F_\infty$ подпространство $\F_\infty^k$, порожденное всеми 
одночленами степени $n$. По построению,
\begin{equation}
\cF_\infty=\oplus_{k=0}^\infty \cF_\infty^k.
\label{eq:fock-orthogonal}
\end{equation}
Далее, очевидно, что ортогональная группа сохраняет подпространства
$\cF_\infty^k$. Эти <<пространства однородных многочленов степени $k$>>
являются симметрическими степенями $\mathsf S^k\ell_2$ пространства $\ell^2$
(мы рассматриваем тензорную степень $\ell_2^{\otimes k}$ и берем в ней 
элементы, симметричные относительно перестановок сомножителей).

Но преобразование Сигала--Баргмана отождествляет $L^2(\R^\infty)$ с $\cF_\infty$.
Поэтому аналогичное разложение на $\bfO(\infty)$-инвариантные подпространства
имеет место и для  $L^2(\R^\infty)$,
$$L^2(\R^\infty)\simeq \oplus_{k=0}^\infty W_k,$$
слагаемое $W_k$ натянуто на многочлены Эрмита $H_\bfl(x)$ 
с $\sum l_j=k$.

\sm

Дадим другое описание тех же прямых слагаемых $W_k$.
Обозначим через $P_k\subset L^2(\R^\infty)$ замкнутое подпространство, порожденное
всеми одночленами $\prod_j x_j^{m_j}$ степени $k$. Эти подпространства, по-прежнему,
$\bfO(\infty)$-инвариантны, но они  не ортогональны.

\begin{lemma}
 Подпространство $P_k$ содержит $P_{k-2}$.
\end{lemma} 
 
{\sc Доказательство.} В силу закона больших чисел, мы имеем сходимость по мере
$$
\frac 1 N(x_1^2+\dots+x_N^2)\to 1.
$$
Нам нужна $L^2$-сходимость,  она тоже есть. Второй многочлен Эрмита равен
$H_2(x)=2(x^2-1)$. Поэтому
$$
\frac 1 N(x_1^2+\dots+x_N^2)- 1= \frac 1{2N} \sum_{k=1}^N H_2(x_k).
$$
Слагаемые  суммы в правой части попарно ортогональны, поэтому квадрат ее $L^2$-нормы
равен $(2N)^{-2} N\|H_2\|^2$. Это выражение стремится к нулю.
Таким образом $P_2\supset P_0$.

Теперь возьмем одночлен $q(x)=\prod_{j=1}^m x_j^{\alpha_j}$. 
Тогда
$$
q(x) \cdot \left(\frac 1N\sum_{m+1}^{m+N} x_j^2-1\right).
$$
стремится к $q(x)$ по тем же причинам.
\hfill $\square$

\sm

Таким образом, мы получаем две фильтрации
$$
P_0\subset P_2\subset P_4\subset\dots, \qquad
P_1\subset P_3\subset P_5\subset\dots
,$$
одна -- в пространстве четных функций ($f(x)=f(-x)$), другая -- в пространстве нечетных
функций.

\begin{proposition}
Прямое слагаемое $W_k$ совпадает с ортогональным дополнением до $P_{k-2}$ в $P_k$.
\end{proposition}

{\sc Доказательство.} Для определенности, рассмотрим случай четного $k=2m$.
Любой четный одночлен (суммарной) степени $2p$ разлагается по четным многочленам Эрмита
(суммарной) степени
$\le 2m$. Это означает, что
$\oplus_{j=0}^m W_{2j}\subset P_m$.

В силу леммы, четный многочлен  степени $\le 2m$ содержится в $P_{2m}$. Поэтому
$\oplus_{j=0}^p W_{2j}\supset P_{2m}$. Таким образом,
$$  \oplus_{j=0}^m W_{2j}= P_{2m}.
$$
Это влечет желаемое утверждение. \hfill $\square$

\sm

{\bf\punct Эргодичность действий.} Пусть $M$ -- пространство с вероятностной мерой $\mu$,
пусть группа $G$ действует на $M$ сохраняющими меру преобразованиями.
Действие называется {\it эргодичным}, если в $M$ нет инвариантных измеримых подмножеств меры,
отличной от 0 и 1. 

Рассмотрим соответствующее действие группы в $L^2(M)$. Тогда эргодичность
равносильна тому, что любая $G$-инвариантная функция в $L^2(M)$ постоянна п.в.

 \begin{theorem}
 {\rm (И.Сигал)}
 Пусть группа $G$ действует ортогональными операторами в вещественном
 гильбертовом пространстве $H$.
 Соответствующее действие $G$ в сигаловском расширении $\wt H$ 
 эргодично тогда и только тогда,
  когда $G$ не имеет конечномерных инвариантных подпространств в
 $H$.
 \end{theorem}

В частности, пусть группа $G=\Z$ порождена одним оператором $A$.  Тогда эргодичность
ее действия равносильна отсутствию у оператора $A$ дискретного спектра.

\sm

Приступим к доказательству. Пусть сначала $V$ -- конечномерное инвариантное подпространство,
пусть $W$ -- ортогональное дополнение до $V$. Выберем ортонормированный базис в $H$ так,
что первые $p$ базисных элементов  содержатся в $V$, а все остальные в $W$. Тогда
цилиндр $x_1^2+\dots+x_p^2\le 1$ является инвариантным подмножеством в $\R^\infty=\wt H$.

Пусть теперь в $H$ нет конечномерных $G$-инвариантных подпространств.  Мы
сейчас  покажем, что в тензорных степенях $H^{\otimes k}$ нет $G$-инвариантных
векторов.

\begin{lemma}
Пусть $\rho_1$, $\rho_2$ -- унитарные представления группы $G$ в бесконечномерных
 пространствах
$V_1$, $V_2$. Если представление $\rho_1\otimes \rho_2$ группы $G$ имеет ненулевой $G$--инвариантный
вектор, то в $V_1$ и в $V_2$ есть конечномерные подпредставления.
\end{lemma}

{\sc Доказательство.} Отождествим $V_1\otimes V_2$ с пространством гильберт-шмидтовских
операторов $V_1'\to V_2$, см. (\ref{eq:hom-otimes}).
Инвариантному вектору соответствует сплетающий оператор,
обозначим его через $T$. Тогда $TT^*$ -- сплетающий оператор в $V_2$. 
Пусть $\lambda$ -- ненулевое собственное значение. Тогда $\ker (TT^*-\lambda)$
-- конечномерное подпредставление (см. выше обсуждение леммы Шура).
\hfill $\square$.

\begin{lemma}
Пусть $\rho_1$, $\rho_2$ -- унитарные представления группы $G$ в бесконечномерных
 пространствах
$V_1$, $V_2$. Если представление $\rho_1\otimes \rho_2$ содержит конечномерное подпредставление,
то конечномерное подпредставления есть в обоих пространствах 
$V_1$, $V_2$.
\end{lemma}

{\sc Доказательство.}
Пусть $W$ -- конечномерное инвариантное подпространство в $V_1\otimes V_2$.
Рассмотрим двойственное действие группы $G$ в  $W'$. Тогда $W'\otimes W$ 
содержит инвариантный вектор, при изоморфизме (\ref{eq:hom-otimes})
он соответствует единичному оператору $W\to W$.

Поэтому $(V_1\otimes V_2)\otimes W'$ содержит инвариантный вектор. Осталось
применить предыдущую лемму к тому же тензорному произведению
с переставленными множителями, 
$(V_1\otimes W' )\otimes V_2$ и $(V_2\otimes W' )\otimes V_1$.
\hfill $\square$

\sm

{\sc Доказательство теоремы.}
Итак, если в $H$ не было конечномерных подпредставлений, то их нет и в тензорных степенях
$H$, а, следовательно, нет и в симметрических степенях. В частности нет и инвариантных векторов.
\hfill $\square$

\section{Меры Пуассона}

\COUNTERS

{\bf\punct Определение.} Рассмотрим лебеговское пространство $M$ с непрерывной маронечной или
$\sigma$-конечной мерой $\mu$. Сейчас мы введем вероятностную меру $\nu=\nu_M$ на пространстве $\Omega(M)$
конечных или счетных подмножеств
множества $M$ (такие подмножества мы будем называть {\it конфигурациями}). Пусть $A\subset M$--измеримое подмножество с конечной мерой. Обозначим
через $\Omega_k (A)$ множество всех $\omega\in\Omega$, таких, что $\omega\cap A$
содержит ровно $k$ точек. Положим, что 
\begin{equation}
\nu\left(\Omega_k(A)\right)= \frac{\mu(A)^k e^{-\mu(A)}}{k!}.
\label{eq:poisson-1}
\end{equation}
Отметим, что $\sum_{k=0}^\infty \nu\left(\Omega_k(A)\right)=1$, 
а поэтому вероятность того, что в $A$
попадет бесконечно много точек, равна нулю. В частности, если $M$ имеет конечную меру,
то конфигурация  $\omega$ конечна с вероятностью 1. С другой стороны,
правая часть (\ref{eq:poisson-1}) стремится к нулю при фиксированном
$k$ при $\nu(A)\to\infty$. Поэтому, если мера $M$ бесконечна, то множество
$\Omega$ с вероятностью 1 бесконечно.

Пусть далее  $A_1$, \dots, $A_n\subset M$ -- попарно непересекающиеся множества конечной меры. 
Тогда мы положим, что события $\Omega_{k_j} (A_j)$ независимы. Иными словами,
$$
\nu\Bigr(\bigcap_j\Omega_{k_j} (A_j)\Bigr)=\prod_j \nu\left( \Omega_{k_j} (A_j)\right)
.$$

Это завершает определение меры $\nu$.

\sm

{\sc Задача.} Пусть $M$ разбито на два непересекающихся подмножества,
$M=M_1\cup M_2$. Конфигурации $\omega\in M$ мы поставим в соответствие 
пару конфигураций $\omega_1=\omega\cap M_1$, $\omega_2=\omega\cap M_2$.
Покажите, что это отображение устанавливает изоморфизм пространств 
$$
\qquad\qquad\qquad\qquad\qquad
\Omega(M)\simeq \Omega(M_1)\times \Omega(M_2).
\qquad\qquad\qquad\qquad\qquad\qquad
\lozenge
$$

\sm

{\bf\punct Корректность определения.}
Пусть сначала пространство $M$ имеет конечную меру. Тогда 
$$\Omega(M)=\bigcup_{k=0}^\infty \Omega_k(M)
.$$
Каждое из множеств $\Omega_k(M)$ мы можем определить следующим образом.
Рассмотрим $k$-кратное прямое произведение $M^k=M\times\dots\times M$,
и снабдим его мерой
$$
\frac {e^{-\mu(M)}}{k!} \,\mu \times\dots\times\mu
.
$$
Мы получили меру на упорядоченных наборах $(x_1,\dots,x_k)$ из $k$ элементов $M$.
Забывая порядок, мы получаем отображение $M^k\to \Omega_k(M)$ и, соответственно,
меру на $\Omega_k(M)$. Рассматривая дизъюнктное
объединение множеств $\Omega_k(M)$, мы получаем меру $\nu=\nu_M$
 на $\Omega(M)$.

\sm

{\sc Задача.} Покажите, что построенная мера обладает желаемыми свойствами.
\hfill$\lozenge$

\sm

Пусть теперь мера на $M$ бесконечна. Рассмотрим возрастающую последовательность
измеримых подмножеств конечной меры $M_1\subset M_2\subset M_3\subset\dots$,
исчерпывающую $M$. Рассмотрим отображение
$\Omega(M_{j+1})\to \Omega(M_j)$, которое конечному подмножеству
$\omega\subset M_{j+1}$ ставит в соответствие подмножество $\omega\cap M_j$.
В силу  задачи, образ меры $\nu_{M_{j+1}}$
при этом  отображении совпадает с  $\nu_{M_{j}}$. Таким образом,
мы получаем проективную цепочку отображений пространств с мерой
$$
\dots \longleftarrow \Omega(M_{j-1}) \longleftarrow \Omega(M_{j})
 \longleftarrow \Omega(M_{j+1}) \longleftarrow\dots 
$$
и берем ее проективный предел. Если множество $A$ содержится в некотором $M_j$,
то равенство (\ref{eq:poisson-1}) мы получаем автоматически.

\sm

{\sc Задача.} Покажите, что равенство (\ref{eq:poisson-1}) выполнено для произвольных
множеств $A\subset M$. \hfill$\lozenge$

\sm

{\bf \punct Формула Кэмпбелла.}

\begin{theorem}
Пусть $h\in L^1(M)$. Тогда произведение
$$
\Phi_h(\omega)=\prod_{x_j\in\omega} (1+h(x_j))
$$
сходится абсолютно при п.в. $\omega\in \Omega(M)$ и 
$$
\int_{\Omega(M)} \Phi_h(\omega)\,d\nu_M(\omega)=\exp\Bigl\{
\int_M h(x)\,d\mu(x)
\Bigr\}.
$$
\end{theorem}

{\sc Доказательство.}
1) Пусть сначала мера пространства $M$ конечна.
Тогда
\begin{multline*}
\int_{\omega\in\Omega_k(M)}  \prod_{x_j\in\omega} (1+h(x_j))\,d\nu(\omega)
=\frac{e^{-\mu(M)}}{k!}\int_{M^k} \prod_{j=1}^k(1+h(x_k))\,\prod d\mu(x_k)
=\\=
 \frac{e^{-\mu(M)}}{k!}\biggl(\int_M (1+h(x))\,d\mu(x)\biggr)^k
 =  \frac{e^{-\mu(M)}}{k!}\biggl(\mu(M)+\int_M h(x)\,d\mu(x)\biggr)^k
\end{multline*}
Далее суммируем по $k$,
$$
\sum_k \int_{\Omega_k(M)}=e^{-\mu(M)} \exp\biggl(\mu(M)+\int_M h(x)\,d\mu(x) \biggr).
$$

2) Пусть теперь $h(x)\ge 0$. Возьмем возрастающую цепочку измеримых множеств 
$M_1\subset M_2\subset \dots$ конечной меры,
исчерпывающую $M$. Тогда соответствующая последовательность функций
$$
\Phi_p(\omega)=\prod_{x_j\in\omega\cap M_p}
(1+h(x_j))
$$
монотонно возрастает и имеет равномерно ограниченные интегралы 
$$
\int_{\Omega(M)}\Phi_p(\omega)\,d\nu_M(\omega)=\exp\Bigl\{\int_{M_p} h(x)\,d\mu(x)\Bigr).
$$
Применяя теорему о монотонной сходимости, мы получаем желаемое утверждение.

\sm

3) Если $-1\le h(x)<0$, то, рассуждая точно так же, мы можем применить теорему
о монотонной сходимости к убывающей последовательности $\Phi_p(\omega)$.

\sm

4) В общем случае мы разбиваем $M$ на три множества, $M_{--}$, $M_-$, $M_+$,
где 
$$h(x)<-1,\quad -1\le h(x)<0, \quad h(x)\ge 0.$$
Произведения
$$
\Phi^{-}(\omega)=\prod_{x_j\in M_{-}}(1+h(x_j)),\qquad
\Phi^{+}(\omega)=\prod_{x_j\in M_{+}}(1+h(x_j))
$$ 
абсолютно сходится, а аналогичное произведение $\Phi^{--}$
по $M_{--}$ конечно. Поэтому произведение $\Phi_h$ сходится абсолютно почти всюду,
причем в силу независимости
\begin{multline*}
\int\limits_{\omega\in\Omega(M)}\Phi(\omega)\,d\nu(\omega)=\\=
\int\limits_{\omega\in\Omega(M_{--})}\Phi^{--}(\omega)\,d\nu(\omega)\,\,\cdot
\int\limits_{\omega\in\Omega(M_{-})}\Phi^{-}(\omega)\,d\nu(\omega)\,\,\cdot
\int\limits_{\omega\in\Omega(M_{+})}\Phi^{+}(\omega)\,d\nu(\omega),
\end{multline*}
а все три сомножителя уже вычислены.
\hfill $\square$

\sm

{\sc Задача.} Выведите формулу Кэмпбелла следующим способом.

\sm 

a) Пусть $M$ представлено в виде дизъюнктного объединения
$M=M_1\coprod M_2$. Пусть $h_1$, $h_2$ -- ограничения функции $h$
на $M_1$, $M_2$. Тогда
$$
\int_M \Phi_{h}\,d\nu_M= \int_{M_1} \Phi_{h_1}\,d\nu_{M_1}
\cdot \int_{M_2} \Phi_{h_2}\,d\nu_{M_2}.
$$

b) Если мера $M$ конечна, а $h$ постоянна, то утверждение верно.

\sm

c) Пусть $M$  произвольно, а $h$ принимает лишь конечное число
значений (и равна нулю вне некоторого подмножества конечной меры).
Тогда формула верна.

\sm

d) Выведите общую формулу с помощью предельного перехода.
\hfill $\square$

\sm

{\sc Задача.}
Вычислите 
\begin{align*}
&\int\limits_{\Omega(M)} \sum_{x_j\in\omega} h(x_j)\,d\nu(\omega),\qquad
\int\limits_{\Omega(M)}
 \sum_{ \begin{matrix}x_i,x_j,x_k\in\omega\\ \text{попарно различны}\end{matrix}}
  h(x_i) h(x_j)h(x_k)\,d\nu(\omega)
  \\
&\int\limits_{\Omega(M)} \Bigl(\sum_{x_j\in\omega} h^2(x_j)\Bigr)
\cdot
\Bigl(\sum_{x_j\in\omega} h^3(x_j)\Bigr)
\,d\nu(\omega)  .
\end{align*}

\sm

{\bf\punct Действие группы преобразований.%
\label{ss:poisson-quasiinvariant}} Пусть мера $M$ бесконечна  Рассмотрим группу $\Gms(M)$,
состоящую из преобразований $M$, переводящих меру $\mu$ в квазиинвариантную
меру. Возьмем в ней подгруппу $\bbG$, состоящую из преобразований $g$, у которых
производная Радона--Никодима удовлетворяет условию
$$(g'(x)-1)\in L^1(M).$$

\sm

{\sc Задача.} Убедитесь, что $\bbG$ -- группа. Покажите, что отображение
$$
\sigma(g)=\int_M (g'(x)-1)\,d\mu(x)=\lim_{j\to\infty}(\mu(gM_j)-\mu(M_j))
,$$
где $M_j$ -- исчерпывающая последовательность подмножеств, является гомоморфизмом из $\bbG$
в аддитивную группу $\R$. \hfill $\lozenge$

\sm

\begin{theorem}
Пусть $g\in\bbG$. Тогда отображение $S_g:\omega\mapsto g\omega$ из
$\Omega(M)$ в себя оставляет меру $\nu$ квазиинвариантной, причем
производная Радона--Никодима равна
$$
S_g(\omega)'=
e^{-\sigma(g)}\prod_{x_j\in\omega} g'(x_j)
.$$
\end{theorem}

{\sc Доказательство.} Достаточно проверить, что для любых попарно непересекающихся
подмножеств конечной меры $A_1$, \dots, $A_m$ и любых $k_1$, \dots, $k_m$
 выполнено
$$
\int_{\bigcap_j \Omega_{k_j}(A_j)^{k_j}} S_g(\omega)' \,d\nu(\omega)
=\prod_j \frac{\mu(gA_j)e^{-\mu(gA_j)}}{k_j!}=\nu\Bigl(\bigcap_j \Omega_{k_j}(gA_j)\Bigr)
.$$
Интегрирование фактически ведется по
$$
\Omega_{k_1}(A_1)\times\dots \times \Omega_{k_m}(A_m)\times \Omega\left(M\setminus 
\cup A_j\right),
$$
а подынтегральное выражение распадается в произведение множителей
$$
\prod_{i:x_i\in \omega\cap A_j} g'(x_i)\quad\text{для каждого $j$,}\quad
\prod_{i:x_i\in \omega\cap A\setminus \cup M_j} g'(x_i)\quad
\text{и $e^{-\sigma(g)}$}
$$
Поэтому интеграл превращается в произведение
\begin{multline*}
\prod_j \frac {e^{-\mu(A_j)}}{k_j!} \biggl(\int\limits_{A_j} g'(x)\,d\mu(x)\biggr)^{k_j}
\times \int\limits_{\Omega(M\setminus \cup A_j)} \prod_{i:x_i\in (\omega\cap M\setminus \cup A_j)} g'(x_i)
\,d\nu_{M\setminus \cup A_j}
\times \\ \times e^{-\sigma(g)}
=
\prod_j
\frac{e^{-\mu(A_j)} \mu(gA_j)^{k_j} }{k_j!} \times \\\times
\exp\Bigl(\int_{M\setminus\cup A_j} (g'(x)-1)\,d\mu(x)\Bigr)\times
\exp\Bigl(-\int_M (g'(x)-1)\,d\mu(x)\Bigr),
\end{multline*}
и это, как легко убедиться, совпадает с желаемым результатом.

\sm

{\bf\punct Структура пространства Фока.} 
Скалярные произведения функций вида $\Phi_h$ легко вычисляются, потому что
$$
\Phi_h \Phi_f=\Phi_{h+f+hf}.
$$
Перенормируем функции $\Phi_h$.
Пусть $h\in L_1(M)\cap L^2(M)$. Определим функцию
$$
\wt\Phi(\omega)=\exp\left(-\int_M h(x)\,d\mu(x)\right)
\cdot 
\Phi_h(\omega).
$$
Легко видеть, что
\begin{equation}
\left\la \wt \Phi_h,\wt\Phi_f\right\ra=\exp \Bigl(\int_M h(x)\ov{ f(x)}\,d\mu(x) \Bigr).
\label{eq:boson-poisson}
\end{equation}

{\sc Задача.} Покажите, что отображение $h\mapsto \wt\Phi_h$
продолжается до непрерывного вложения $L^2(M)\to L^2(\Omega(M))$.
\hfill $\lozenge$

\sm

Теперь, сравнивая эту формулу с формулами (\ref{eq:psi-fock}) и (\ref{eq:Phi-fock}),
мы видим, что наше пространство $L^2(\Omega(M))$ отождествляется с
$L^2$ на сигаловском расширении пространства $L^2(M)$ и с пространством Фока. 
Получаются интересный бесконечномерный интегральный операторы, 
который не очень изучен.

\sm

{\sc Задача.} Рассмотрим действие группы $\bbG$
на вещественном пространстве $L^2(M)$ аффинными изометрическими преобразованиями
$$
T(g) f(x)=f(g(x))g'(x)^{1/2}+ g'(x)^{1/2}-1
$$
(убедитесь, что это в самом деле действие).
Поэтому наша группа действует и в $L^2$ на сигалоском расширении $\wt{L^2(M)}$.
 Покажите, что наше соответствие сплетает данное представление
 с представлением $\bbG$ в $L^2(M)$.

\section[Виртуальные перестановки]{Виртуальные перестановки}

\COUNTERS

{\bf \punct Проекция $S_n\to S_{n-1}$.}
Рассмотрим группу $S_n$ всех перестановок из $n$ элементов.
Оказывается, что существует естественное отображение
$$
\Upsilon^n_{n-1}: S_n\to S_{n-1}
.$$

{\sc Пример.} Пусть $n=8$.
Рассмотрим, для примера подстановки
$$
\begin{matrix}
1&2&3&4&5&6&7&8\\
3&7&4&1&6&3&5&8
\end{matrix}
\qquad\qquad
\begin{matrix}
1&2&3&4&5&6&7&8\\
7&3&6&1&8&2&5&4
\end{matrix}
$$
Вычеркнем из них цифру 8. Получим
$$
\begin{matrix}
1&2&3&4&5&6&7&\phantom{8}\\
3&7&4&1&6&3&5&\phantom{8}
\end{matrix}
\qquad\qquad
\begin{matrix}
1&2&3&4&5&6&7&\phantom{8}\\
7&3&6&1&\phantom{8}&2&5&4
\end{matrix}
$$
В первом случае мы получили элемент группы $S_7$. Во втором случае у 5 нет образа, а у  
4 нет прообраза. Поэтому мы положим, что $5$  переходит в 4,
$$
\phantom{\begin{matrix}
1&2&3&4&5&6&7&\phantom{8}\\
3&7&4&1&6&3&5&\phantom{8}
\end{matrix}}
\qquad\qquad
\begin{matrix}
1&2&3&4&\mathbf{5}&6&7&\phantom{8}\\
7&3&6&1&\mathbf{4}&2&5&\phantom{8}
\end{matrix}
$$

Общее правило ясно из этого примера.
Пусть $\sigma\in S_n$. Пусть
$\sigma i=j$:

a) если $i$, $j<n$, то мы полагаем
$\Upsilon i=j$.

b) если $\sigma i=n$ и $i<n$, то мы полагаем 
$\Upsilon  i= \sigma n$.

\sm

{\sc Задача.}  Покажите, что каждый элемент из $S_{n-1}$ имеет
ровно $n$ прообразов.
\hfill $\lozenge$

\sm

\begin{corollary}
 Образ равномерной меры на $S_n$ при отображении
 $\Upsilon$ совпадает с равномерной мерой на $S_{n-1}$.
\end{corollary}

Конечно, наше отображение не является гомоморфизмом. Однако оно удовлетворяет следующему
весьма жесткому условию эквивариантности.

\sm

{\sc Задача.} Пусть $\sigma\in S_n$, $\tau_1$, $\tau_2\in S_{n-1}$. Тогда
$$\qquad\qquad\qquad\qquad\qquad
\Upsilon(\tau_1 \sigma \tau_2^{-1})= \tau_1 \Upsilon(\sigma) \tau_2^{-1}.
\qquad\qquad\qquad\qquad\qquad
\lozenge
$$

{\bf\punct Другое определение отображения $\Upsilon$. Проективные системы мер.}
Наше отображение можно описать следующим образом. Мы раскладываем подстановку $\sigma$ в произведение 
независимых циклов и вычеркиваем $n$ из записи этого разложения.

\sm

{\sc Пример.} В разобранном выше примере мы имеем разложения
$$
(134)(2756)(8)\qquad\qquad (17584)(236).
$$
Мы вычеркиваем 8 и получаем соответственно
$$
(134)(2756)\qquad\qquad (1754)(236)
$$
Отметим, что в первом случае одноэлементный цикл исчез, во втором случае число циклов осталось тем же. 
\hfill $\lozenge$

\sm

Для подстановки $\sigma$ обозначим через $\{\sigma\}$ число ее циклов.

Теперь фиксируем $t>0$ и введем меру $\mu_t^n$ на $S_n$ положив, что
вес перестановки $\sigma$ равен 
$$
\frac{t^{\{\sigma\}}}
{t(t+1)\dots(t+n-1)}
.$$

\begin{theorem}
Мера $\mu_t^n$ является вероятностной. Образ меры $\mu_t^{n+1}$ при отображении 
$\Upsilon$ совпадает с $\mu_t^{n}$.
\end{theorem}

{\sc Доказательство.} Рассмотрим прообраз подстановки $\sigma\in S_n$.
Точка $\tau$ прообраза определяется значением $\tau(n+1)$. Если $\tau(n+1)\le n$,
то $\{\tau\}=\{\sigma\}$ и таких
точек $n$. Если $\tau(n+1)=n+1$, то $\{\tau\}=\{\sigma\}+1$. Таким образом, мера прообраза 
равна
$$
\frac{n\cdot t^{\{\sigma\}}}
{t(t+1)\dots(t+n)}
+
\frac{t\cdot t^{\{\sigma\}}}
{t(t+1)\dots(t+n)}
=\frac{ t^{\{\sigma\}}}
{t(t+1)\dots(t+n-1)}
,$$
т.е. равна мере точки  $\sigma$. \hfill $\square$

\sm

{\bf\punct Виртуальные перестановки.}
Рассмотрим цепочку отображений
$$
\dots\stackrel{\Upsilon}{\longleftarrow} S_n  \stackrel{\Upsilon}{\longleftarrow} S_{n+1}
\stackrel{\Upsilon}{\longleftarrow} \dots
$$
Для любого $n$ образ меры $\mu_t^{n+1}$ совпадает с $\mu_t^n$.
Обозначим обратный предел  этой цепочки  через
$\frS$. По теореме Колмогорова об обратных пределах, на $\frS$
определена вероятностная мера $\mu^t$, такая, что
для любого $n$ образ $\mu_t$ при естественной проекции
$\frS\to S_n$ совпадает с $\mu_t^n$.

Элементы пространства $\frS$ мы будем называть {\it виртуальными перестановками}.
Это не перестановки какого-либо бесконечного множества,
однако точка этого пространства является последовательностью
конечных перестановок
$$
\sigma_m\in S_m, \quad\text{причем $\Upsilon \sigma_n=\sigma_{n-1}$.}
$$
(ниже будет построена более прозрачная реализация пространства $\frS$).
Умножения на $\frS$ тоже нет.

Обозначим через $S_\infty$ группу всех перестановок натурального
ряда с конечным носителем.
Мы утверждаем, что 
эта группа действует на $\frS$ левыми и правыми умножениями.

Действительно, пусть $\tau\in S_\infty$. Пусть на самом деле
$\tau$ содержится в подгруппе $S_k$.
Тогда для любого $m\ge k$ определен левый сдвиг
$\sigma\mapsto \tau \sigma$ на группе $S_m$. При этом
$\Upsilon(\tau \sigma)=\tau \Upsilon (\sigma)$, а поэтому
определено и преобразование предельного объекта.

Аналогично определяются правые сдвиги.

Повторим определение немного другими словами.
Возьмем точку из $\frS$, то есть последовательность
$\sigma_m\in S_m$, такую, что $\Upsilon \sigma_m=\sigma_{m-1}$.
Далее рассмотрим последовательность $\tau \sigma_m\in S_m$.
 Мы имеем $\Upsilon (\tau \sigma_m)=\tau \sigma_{m-1}$
 при $m\ge k+1$.
 Таким образом, наша последовательность обладает желаемым
 свойством согласования, но начинается она с $m=k$. На меньшие номера
 она продолжается автоматически из условия согласования.

\sm

{\bf\punct Отступление. Определение отображения $\Upsilon$ через произведения транспозиций.}

\begin{proposition}
Любая перестановка $\sigma\in S_n$ единственным образом представима
в виде произведения транспозиций
\begin{equation}
\sigma=(i_1 j_1)\dots (i_k j_k)
,
\label{eq:in}
\end{equation}
где
$$
j_1<j_2<\dots<j_k \qquad\text{и $i_m<j_m$ для всех $m$.}
$$
\end{proposition}

{\sc Доказательство.} Пусть $\sigma\in S_n$. Если $\sigma^{-1} n= \alpha<n$,
то мы пишем $\sigma= \bigl[\sigma \cdot(\alpha n)\bigr]\cdot (\alpha n)$.
Тогда $\sigma\cdot (\alpha n)\in S_{n-1}$, и мы можем считать ее разложенной по предположению
индукции. Если же $\sigma^{-1}n=n$, то сама $\sigma$ является элементом $S_{n-1}$.
\hfill $\square$

\sm

Теперь мы можем определить отображение $\Upsilon: S_n\to S_{n-1}$ как 
отображение, <<стирающее>> множитель $(\alpha n)$ в разложении (\ref{eq:in}),
если таковой множитель есть (в противном случае мы не делаем с разложением ничего).

\sm

{\bf\punct Описание обратного предела.}
Теперь мы опишем модель множества $\frS$ (чуть ниже она будет усовершенствована). 

Сначала введем определение циклически упорядоченного множества. Неформально,
это множество, где точки идут по кругу. Теперь формально.
Рассмотрим стандартную окружность. 
Для любой упорядоченной тройки попарно различных точек окружности можно сказать,
расположены они по часовой стрелке или против часовой стрелки.
Пусть  $W$ --  абстрактное счетное множество, пусть 
множество упорядоченных троек попарно различных точек $W$
разбито на два подмножества, <<расположенных по часовой стрелке>>
и <<расположенных против часовой стрелки>>. Мы скажем, что это
{\it циклическое упорядочение}, если $W$ может быть вложено в окружность
так, что тройки, расположенные <<по часовой стрелки>>, перейдут в тройки,
расположенные <<против часовой стрелки>>, и наоборот. В принципе, можно придумать
аксиоматическое определение, но нам оно будет не нужно.

\sm

{\sc Задача.} Группа $S_n$ находится в естественном взаимно однозначном соответствии
с множеством всех разбиений множества $\{1,\dots,n\}$ на подмножества, на которых
введен циклический порядок.  \hfill $\lozenge$

\sm

Мы определим множество $\Sigma$ как множество всех счетных или конечных
разбиений $\N$ на подмножества $N_i$, при этом на каждом подмножестве
(назовем их {\it квазициклами}) введено
циклическое упорядочение.

Для каждого $n$ определим отображение $\Sigma \to S_n$.
А именно берем упомянутое выше разбиение и забываем все 
элементы с номерами $>n$. Получается разбиение множества
$\{1,2,\dots,n\}$ на циклически упорядоченные множества, а это уже
отождествляется с элементом $S_n$.

\sm

{\sc Задача.} a) Убедитесь, что множество $\Sigma$ в самом деле является проективным пределом
множеств $S_n$.

\sm

b) Опишите левое и правое действие $S_\infty$ в этих терминах.
 \hfill $\lozenge$

\sm

Проективный предел является  не просто множеством, а пространством с мерой.
Полученный объект можно понимать следующим образом. Мы <<разыгрываем>>
последовательность $\sigma_m\in S_m$ случайных перестановок следующим образом.
Пусть дана подстановка $\sigma_m$, представленная в виде произведения независимых
циклов. Между элементами этих циклов имеется $m$ промежутков.
 Равновероятно, с вероятностями $\frac{1}{t+m}$ мы вставляем $m+1$ в один из промежутков в циклах. С вероятностью $\frac t{t+m}$ мы организуем новый цикл, состоящий
 из единственного элемента $m+1$.

\sm

{\sc Задача.} Убедитесь, что число квазициклов с вероятностью 1
счетно. 
Убедитесь, что с вероятностью 1 каждый квазицикл  как циклически
упорядоченное множество устроен как $\Q/\Z$ (это равносильно тому, что для любых
$x$, $z\in N_i$ существует $y$, такой, что $x$, $y$, $z$ идут по часовой стрелке). 
\hfill $\lozenge$

\sm

Верны и более сильные высказывания. Оказывается
что у каждого квазицикла почти наверное есть длина, 
которая получается как предел нормированных длин циклов $S_n$, 
и более того, между двумя элементами квазицикла 
появляется естественное расстояние.
Мы сейчас уточним нашу конструкцию без доказательства.

\sm

{\bf \punct  <<Многолюдный ресторан>>.} Фиксируем $t>0$.
Сейчас мы опишем некоторый случайный процесс. 

Директор ресторана организует
установку счетного количества ориентированных  круглых столов с общей длиной периметров 1.
Длины столов он выбирает случайным образом по следующему закону
(<<распределение Пуассна--Дирихле>>).
Берется полупрямая $x>0$ с мерой $te^{-x}/x$. Разыгрывается пуассоновая
конфигурация (см. \S 9) на полупрямой,  обозначим ее элементы через 
$x_1$, $x_2$, \dots. В
качестве длин столов берутся числа $x_j/\sum_i x_i$, сумма этих чисел равна 1.

Далее в ресторан по очереди приходят посетители с номерами на спине и случайным 
образом садятся за столы (т.е. мы равномерно кидаем точки на окружности).
В итоге мы получаем счетное множество пронумерованных людей, сидящих за столами.
С вероятностью 1 эти точки окружности будут различны, а поэтому мы 
 получили разбиение $\N$ на счетное число квазициклов.
 
 На множестве всех таких объектов (набор столов плюс рассадка посетителей)
 по построению есть вероятностная мера. Назовем это пространство 
 {\it многолюдным рестораном}.

\begin{theorem}
Полученное вероятностное пространства совпадает с $\frS$. 
\end{theorem}

См. \cite{Tsi}.

\sm

{\bf\punct Действие группы $S_\infty\times S_\infty$.} Оно было определено выше.
Опишем его на языке предыдущего пункта. Представим себе, что каждый посетитель ресторана
состоит из левой и правой половины, которые можно отделять друг от друга. Представим себе,
что столы являются окружностями (а не кругами, что несколько противоречит аналогии с трапезой,
но по математическому смыслу столы -- именно окружности). 
Пусть $\tau \in S_\infty$, а $\frs$ - точка ресторана.
Определим точку $\tau \frs$. Пусть $\tau$ фактически содержится в $S_n$.
Мы находим всех посетителей ресторана с номерами $\le n$ и разделяем каждого
из них на левую и правую половину и в соответствующем месте разрезаем стол.
Теперь мы имеем набор круглых столов, а также набор из $n$ сегментов, каждый из которых
начинается с полупосетителя и кончается полупосетителем, причем сегменты ориентированы.
 Теперь  для каждого $j\le n$ мы приклеиваем правый конец отрезка с полупосетителем
 $j$ с левым концом отрезка с полупосетителем $\tau j$. Присваиваем
 точке склейки номер $j$. В итоге мы опять получаем набор
 круглых столов.
 
 Эта операция определяет действие группы $S_\infty$ на $\frS$ левыми умножениями.
 
 \sm
 
 {\sc Задача.} a) Покажите, что операция совпадает с операцией, введенной выше.
 
 \sm
 
 b) Опишите в тех же терминах правое действие.
 \hfill $\lozenge$ 

\sm

Теперь заметим, что хотя число столов в $\frs $ и $\tau\frs$
 счетно, оно в результате переклейки претерпевает изменение.
Пусть $\frs\in\frS$.
Рассмотрим все столы, на которых присутствует хотя
бы один посетитель с номером $\le n$.  
Обозначим  число таких столов через $a_n(\frs)$. 
Определим число $c(\tau,\frs)\in \Z$ как
$$
c(\tau,\frs)=a_n(\tau\frs) -a_n(\frs),\qquad\text{если $\tau\in S_n$}
.
$$

{\sc Задача.} Покажите, что производная Радона-Никодима преобразования
$\frs\mapsto \tau\frs$ равна $t^{c(\tau,\frs)}$.  \hfill $\lozenge$

\sm

Таким образом, мы получаем действие группы $S_\infty\times S_\infty$
на $\frS$ преобразованиями, оставляющими меру квазиинвариантной.

\sm

{\bf \punct Бисимметрическая группа.%
\label{ss:bisymmetric-1}} Обозначим через $\ov S_\infty$ 
полную группу перестановок натурального ряда, снабженную топологией как в 
п. \ref{ss:S-infty}. Группа $\ov S_\infty$  действует на многолюдном ресторане перестановками
номеров посетителей. Эти преобразования сохраняют меру. 

Поэтому на нашем пространстве действует бисимметрическая группа, введенная в п.\ref{ss:bisymmetric}.

В \cite{KOV} было построено явное спектральное разложение представления бисимметрической
группы в пространстве $L^2(\frS)$.

\section[Виртуальные унитарные матрицы]{Виртуальные унитарные матрицы}

\COUNTERS

{\bf\punct Отображение $\Upsilon$ для унитарных групп.}
Рассмотрим теперь группу $\U(n+m)$. Мы будем записывать ее элементы как
блочные матрицы размера $n+m$.
Рассмотрим следующее отображение из $\U(n+m)$ в пространство
матриц порядка $n$, опредленное почти всюду:
\begin{equation}
 \Upsilon^{n+m}_n \begin{pmatrix}
           \alpha&\beta\\ \gamma&\delta
          \end{pmatrix}
=\alpha-\beta(1+\delta)^{-1}\gamma.
\end{equation}

\begin{theorem}
\label{th:livshits}
 {\rm a)} Если $g \in \U(n+m)$, то $\Upsilon^{n+m}_n g\in \U(n)$.
 
 \sm
 
 {\rm b)} Если $h_1$, $h_2\in \U(n)$, то
 $$
 \Upsilon^{n+m}_n(h_1 g h_2)=h_1\Upsilon^{n+m}_n(g) h_2
 .$$
 
 {\rm c)} Образ меры Хаара на $\U(n+m)$ при отображении $\Upsilon^{n+m}_n$ 
 совпадает с мерой Хаара на $\U(n)$.
 
 \sm
 
 {\rm d)} $\Upsilon^q_r \Upsilon^p_q=\Upsilon^p_r$ 
\end{theorem}

{\sc Доказательство.} a) Мы записываем уравнение
$$
\begin{pmatrix}
 q\\-x
\end{pmatrix}
=
\begin{pmatrix}
 \alpha&\beta\\\gamma&\delta
\end{pmatrix}
\begin{pmatrix}
 p\\x
\end{pmatrix},
$$
где $p$, $q\in\C^n$, а $x\in \C^m$. Перепишем его в виде
\begin{align*}
 q=\alpha p+\beta x;\\
 -x=\gamma p+\delta x.
\end{align*}
Исключая из второго уравнения $x$, мы получаем $x=-(1+\delta)^{-1} \gamma p$.
Подставляя в первое уравнение, получаем 
$$
q=\bigl(\alpha-\beta(1+\delta)^{-1} \gamma) p
$$
Заметим, что для данного $p$ мы можем однозначно вычислить
$x$, а потом $q$. 
Теперь вспоминаем, что $\begin{pmatrix}
 \alpha&\beta\\\gamma&\delta
\end{pmatrix}$
унитарна, а потому
$$
\|q\|^2+\|-x\|^2=\|p\|^2+\|x\|^2
$$
или $\|q\|^2=\|p\|^2$, что и означает унитарность матрицы 
$\alpha-\beta(1+\delta)^{-1} \alpha$.

Утверждение b) очевидно. В силу этого
утверждения, образ меры Хаара является инвариантной мерой,
а потому он совпадает с мерой Хаара. Это доказывает утверждение 
c). Утверждение $d$ вытекает из конструкции следующего пункта.
\hfill $\square$

\sm

{\bf\punct Отображение $\Upsilon$ и преобразование Кэли.}
Рассмотрим преобразование Кэли
$$
r= \cC(g)=\frac{1-g}{1+g}=1+2(1+g)^{-1}.
$$
Как мы видели в \S4, оно переводит унитарные матрицы в антиэрмитовы.
Обозначим через $[r]_n$ левый верхний $n\times n$ уголок матрицы $r$. 

\begin{theorem}
Для $g\in\U(n+m)$ выполнено
 $$
 [\cC(g)]_n=\cC(\Upsilon^{n+m}_n(g)).
 $$
\end{theorem}

{\sc Доказательство.} Напомним трюк Фробениуса для обращения блочных матриц.
Мы разлагаем матрицу в произведение блочно верхнетреугольной, блочно диагональной, блочно нижнетреугольной
матриц:
$$
\begin{pmatrix}
 a&b\\c&d
\end{pmatrix}=
\begin{pmatrix}1&b d^{-1}\\ 0&1 \end{pmatrix}
\begin{pmatrix}
 a-b d^{-1} c&0\\0&d
\end{pmatrix}
\begin{pmatrix}
 1&0\\ d^{-1}c&1
\end{pmatrix}
$$
(вопрос: как это сделать, не зная заранее ответа?). 
Каждый сомножитель легко обращается, и мы получаем
\begin{multline*}
\begin{pmatrix}
 a&b\\ c&c
\end{pmatrix}^{-1}=
\begin{pmatrix}
 1&0\\ -d^{-1}c&1
\end{pmatrix}
\begin{pmatrix}
 (a-b d^{-1} c)^{-1}&0\\0&d^{-1}
\end{pmatrix}
\begin{pmatrix}1&-b d^{-1}\\ 0&1 \end{pmatrix}=
\\=
\begin{pmatrix}
 (a-b d^{-1} c)^{-1}&\dots\\\dots&\dots
\end{pmatrix}.
\end{multline*}
Вычислим преобразование Кэли
для $g=\begin{pmatrix}
 a&b\\ c&c
\end{pmatrix}$,
$$
-\begin{pmatrix}
  1&0\\0&1
 \end{pmatrix}
+2 \left[\begin{pmatrix}
  1&0\\0&1
 \end{pmatrix}+\begin{pmatrix}
 \alpha&\beta\\ \gamma&\delta
\end{pmatrix}  \right]^{-1}=
\begin{pmatrix}
 -1+2(1+\alpha-\beta (1+\delta)^{-1} \gamma)^{-1}&\dots\\\dots&\dots
\end{pmatrix},
$$
Мы получили желаемую формулу.\hfill $\square$

\sm

{\bf\punct Третье определение отображения $\Upsilon$.}
Рассмотрим евклидово пространство $\C^n$ со стандартным базисом $e_n$,
через $\C^k$ мы обозначим подпространство, натянутое на первые $k$
базисных векторов.

Мы скажем, что {\it отражение}%
\footnote{Этот странноватый термин общепринят. Заметим, что <<отражения>>
в $\C^1$ -- это в точности повороты.}-- это унитарный оператор, у которого ровно одно собственное число
отлично от 1. Любое отражение имеет вид
$$
S[v,e^{i\phi}] x= x+(e^{i\phi}-1)\la x , v\ra,\qquad \text{где $\la v,v\ra=1$, $0<\phi<2\pi$.}
$$

{\sc Задача.} Пусть $x$, $y$ -- различные единичные векторы. Тогда существует единственное отражение,
переводящее $x$ в $y$. \hfill $\lozenge$

\sm

Таким образом, отражение однозначно определяется образом вектора $e_n$, 
а множество всех отражений может быть отождествлено со сферой 
$S^{2n-1}$, из которой выброшена точка $e_n$.
 Обозначим через $L^k$ множество
всех отражений, у которых $v\in\C^k$.

\begin{theorem}
 Почти любой элемент $g$ группы $\U(n)$ единственным образом представим в виде 
 произведения
 $$
 g=S_1 S_2\dots S_n, \qquad \text{где $S_j\in L^j$}.
 $$
\end{theorem}

{\sc Доказательство.}
Пусть $ge_n=y$. В общем положении $e_n\ne y$. Берем единственное отражение $S_n$, переводящее
$e_n$ в $y$ и полагаем $g':= g S_n^{-1}$. Теперь 
$g'e_n=e_n$, а потому $g'\in \U(n-1)$. И т.д.
\hfill $\square$

\sm

Заметим, что тем самым мы построили определенную почти всюду (но не всюду!) биекцию
$$
\U(n)\simeq S^1\times S^3\times \dots \times S^{2n-1}.
$$
Очевидно, что мере Хаара на $\U(n)$ соответствует произведение равномерных мер на сферах.

\sm 

\begin{theorem}
 Стирание последнего множителя в произведении
 $$ g=S_1 S_2\dots S_n$$
 соответствует отображению 
 $$
 g\mapsto \begin{pmatrix} \Upsilon^n_{n-1} (g)&0\\0&-1\end{pmatrix}
 .$$
\end{theorem}

{\sc Замечание.} Выглядящий неэстетично минус убирается, если чуть-чуть пошевелить формулу
для $\Upsilon$, а именно положить $\Upsilon (g)=a+b(1-d)^{-1}$. 
В выбранном нами варианте отображение не имеет особенностей в окрестности единицы
(а также переводит $\SO(n+m)$ в $\SO(n)$).
\hfill $\lozenge$

\sm

{\sc Доказательство.}
Достаточно проверить, что 
$$\rk \left[g^{-1} \begin{pmatrix} \Upsilon^n_{n-1} (g)&0\\0&-1\end{pmatrix}-1\right]=1.$$
Равносильно, мы должны проверить, что ранг следующей матрицы равен 1:
\begin{multline*}
\begin{pmatrix}
 \alpha&\beta\\ \gamma&\delta
\end{pmatrix}
-
\begin{pmatrix} \alpha-\beta(1+\delta)^{-1}\gamma &0\\0&-1\end{pmatrix}=
\begin{pmatrix}
 \beta(1+\delta)^{-1} \gamma &\beta\\ \gamma& 1+\delta
\end{pmatrix}
=\\=
\begin{pmatrix}
\beta\\1+\delta 
\end{pmatrix}
\begin{pmatrix}
 (1+\delta)^{-1}&0\\0&(1+\delta)^{-1}
\end{pmatrix}
\begin{pmatrix}
 \beta&(1+\delta)
\end{pmatrix}
\end{multline*}

Мы получили произведение матрицы-столбца, скалярной матрицы и матрицы-строки,
поэтому ранг равен 1. \hfill $\square$

\sm

{\bf \punct Обратный предел.%
\label{ss:UUU}}
Теперь мы имеем цепочку отображений 
\begin{equation}
\dots \longleftarrow \U(n-1)\longleftarrow \U(n) \longleftarrow \U(n+1) \longleftarrow \dots
,
\label{eq:UUU}
\end{equation}
согласованную с мерами Хаара, и мы можем взять обратный предел $\frU$ этой цепочки,
на нем канонически определен мера -- обратный предел мер Хаара. 

\sm

a) Как и в случае симметрической
группы, здесь определено левое и  правое действия группы $\U(\infty)$
финитных унитарных матриц, а вместе --  это действие группы 
$\U(\infty)\times \U(\infty)$. Эти преобразования сохраняют меру.

\sm

b) Далее, можно показать, что
действие диагонали $\diag \U(\infty)$ продолжается по слабой непрерывности
до действия полной унитарной группы.

\sm

c) Таким образом, получается действие группы, состоящей из пар унитарных матриц 
$(g_1,g_2)$, таких, что $g_1g_2^{-1}$ финитна. 
Это условие можно ослабить до: $g_1 g_2^{-1}-1$ -- ядерный оператор.

\sm

d) Действие $\U(\infty)\times \U(\infty)$ не является индивидуальным, см. 
п.\ref{ss:individ}.

\sm

e) В случае симметрической групп мы имели явное и довольно прозрачное описание обратного
предела. В нашем случае таких описаний не известно. Конечно, наше пространство отождествляется
со счетным произведением сфер $\prod_{j=1}^\infty S^{2j-1}$, но
не видно, как писать действие группы в этой модели. Что есть точка предела
(\ref{eq:UUU}), не ясно, зато действие группы на каждом этаже понятно.
Мы можем строить обратный предел, стартуя с цепочки пространств анти-эрмитовых матриц,
$$
\dots\longleftarrow \mathrm{AHerm}_n \longleftarrow \mathrm{AHerm}_{n+1}\longleftarrow \dots,
$$
на  $n$-ном шагу мы берем меру $(1+TT^*)^{-n}\, dT$ (см. \S4).
Предельная мера, конечно, реализуется на пространстве бесконечных антиэрмитовых матриц.
Формула для действия дубля, как мы видели в Предложении \ref{pr:8.4},
оказывается несложной.
Однако непонятно, соответствуют ли эти матрицы каким-либо операторам.

\sm

{\bf\punct Меры Хуа.} Теперь возьмем функцию  на $\U(n)$, задаваемую формулой
$$
\phi^n_{\lambda,\mu}(g)=\det(1+g)^\lambda \det(1+g^{-1})^\mu:=\lim_{t\to 1-0} \det \left[(1+tg)^\lambda \right]
\cdot
\det \left[(1+tg^{-1})^\lambda \right]
, $$
где $\lambda$, $\mu$ -- коплексные параметры, $\Re(\lambda+\mu)>-1$.
В формуле комплексное число возводится в комплексную степень,
поэтому мы должны были уточнить, что это значит. 
Заметим, что степенное выражение $(1+tg)^\lambda$ разлагается в ряд Тейлора
обычным образом, а дальше мы уже можем брать определитель и переходить к пределу.

Рассмотрим
меру (заряд) $\phi_{\lambda,\mu}(g)\,d\sigma_n(g)$, где $d\sigma_n(g)$ -- 
вероятностная мера Хаара. Если $\mu=\ov\lambda$, то эта мера 
является положительной.

\begin{theorem}
{\rm a)} Образ заряда $\phi^n_{\lambda,\mu}(g)\,d\sigma_n(g)$ при отбражении $\Upsilon^n_{n-1}$
равен 
$$
\frac{\Gamma(n)\Gamma(\lambda+\mu+n)}{\Gamma(\lambda+n)\Gamma(\mu+n)}\cdot \phi^{n-1}_{\lambda,\mu} d\sigma_{n-1} .
$$
\end{theorem}

Отсюда, в частности, следует, что
$$
\int_{\U(n)} \det(1+g)^\lambda \det(1+g^{-1})^\mu \,d\sigma_n(g)=
\prod_{k=1}^n \frac{\Gamma(k)\Gamma(\lambda+\mu+k)}{\Gamma(\lambda+k)\Gamma(\mu+k)}.
$$

\begin{lemma}
Обозначим через $D$ круг $|z|\le 1$ на $\C$.
Рассмотрим отображение $\U(n)\to \U(n-1)\times D$, заданное формулой
$$
g\mapsto (\Upsilon^n_{n-1} g, g_{nn}).
$$
Образ меры Хаара при этом отображении равен
\begin{equation}
d\sigma_{n-1}(g) \cdot \pi^{-1}(n-1) (1-|z|^2)^{n-2}\, dz\,d\ov z.
\label{eq:ugolok}
\end{equation}
\end{lemma}

{\sc Доказательство леммы.} Во-первых, наше отображение коммутирует
с левыми и правыми
 умножениями на элементы группы $\U(n-1)$, поэтому образ меры инвариантен относительно 
таких умножений. Следовательно, он имеет вид 
$$
d\sigma_{n-1}(g)\cdot h(z) \, dz\,d\ov z.
$$
Нам нужно вычислить функцию $h(z)$. Она равна плотности распределения матричного элемента 
$g_{nn}$. Чтобы найти ее, мы берем отображение $\U(n)\to S^{2n-1}$, заданное
формулой 
$$g\mapsto g e_{nn}=(g_{1n},\dots, g_{nn}).$$
Образ меры Хаара, очевидно, будет равномерной мерой на сфере.
Таким образом, остается вычислить распределение последней координаты у векторов 
лежащих на сфере $S^{2n-1}\subset \C^n$. Мы оставляем это в качестве  упражнения.
\hfill $\square$

\sm

{\sc Доказательство теоремы.}
Разложим отображение $\Upsilon$ в произведение отображений
$\U(n)\to \U(n-1)\times D \to \U(n)$.

Сначала посмотрим, что прозойдет с нашим зарядом при первом отображении.
Мы применяем к $\det(1+g)$ формулу (\ref{eq:block-det-2}) для определителя блочной матрицы,
\begin{multline*}
\det\left[\begin{pmatrix}
          1&0\\0&1 
          \end{pmatrix}
 + 
 \begin{pmatrix}
 \alpha&\beta\\ \gamma&\delta
\end{pmatrix}
\right]=\\=
(1+\delta)\det(1+\alpha-\beta(1+\delta)^{-1}\gamma)=
(1+\delta)\det\bigl(1+\Upsilon^n_{n-1}(g)\bigr).
\end{multline*}
Поэтому образ меры $\det(1+g)^\lambda\det(1+g^{-1})^\mu\,d\sigma_n(g)$
при отображении 
(\ref{eq:ugolok}) равен
\begin{multline*}
\det\bigl(1+h)^\lambda\det\bigl(1+\ov h\bigr)^\mu d\sigma_{n-1}(h)
\times
\\\times
 \pi^{-1} (n-1)
(1+z)^\lambda (1+\ov z)^{\mu} (1-|z|^2)^{n-2} \,dz\,d\ov z
.
\end{multline*}
Мы видим, что полученная мера на $\U(n-1)\times D$ распадается в
прямое  произведение мер. Теперь нам нужно спроектировать
эту меру на первый сомножитель, что сводится к вычислению интеграла 

$$
 \pi^{-1} (n-1)
\int_{|z|<1} (1+z)^\lambda (1+\ov z)^{\mu} (1-|z|^2)^{n-2} \,dz\,d\ov z,
$$
что является хорошим упражнением по матану.
\hfill $\square$

\sm

Теперь мы снова можем рассмотреть цепочку (\ref{eq:UUU})
унитарных групп и ее обратный предел. Там снова действует группа
$\U(\infty)\times\U(\infty)$ и ее пополнение.

\sm

В настоящее время пространство $L^2$ на этой группе явно разложено,
хотя и сам объект и это разложение, видимо, остаются недопонятыми.
Отметим также, что аналогичная конструкция с обратными пределами есть для всех 10 серий
компактных симметрических пространств (но разложения $L^2$ не известно).

См. также
\cite[\S\S2.10, 3.6]{Ner2}

\section{Замыкания действий}

\COUNTERS

{\bf\punct Полиморфизмы.%
	\label{ss:polymorphisms}}
Напомним, что
 {\it лебеговское пространство с  мерой},
это пространство, изоморфное объединению промежутка
$[a,b]$ и конечного или счетного множества точек,
имеющих ненулевую меру (см. \cite[\S 9.4]{Bog}). Пространство может состоять
лишь из промежутка или лишь из набора точек. Стоит иметь в виду, что почти все
пространства с мерой,  встречающиеся в анализе, являются лебеговскими.

Пусть $(M,\mu)$, $(N,\nu)$ -- лебеговские пространства с вероятностной
мерой (удобно считать, что на них фиксированы борелевские структуры),
и тогда на них). {\it Полиморфимзм} $(M,\mu)\looparrowright(N,\nu)$ -- это мера $\sigma$ на произведении $M\times N$,
такая, что 

\sm

1) образ $\sigma$ при проекции $M\times N\to M$
совпадает с $\mu$;

\sm

2) образ $\sigma$ при проекции
$M\times N\to N$ совпадает с $\nu$.

\sm

{\sc Примеры.}  a) Пусть $\psi:M\to N$ -- отображение, переводящее меру $\mu$ в меру $\nu$.
Рассмотрим отображение $\wt \psi: M\to M\times N$, заданное формулой
$m\mapsto (m,\psi(m))$. Образ меры $\mu$ при этом отображении
является полиморфизмом (эта мера сосредоточена на графике отображения $\psi$).
 
\sm

b) Аналогично, отображение $N\to M$, переводящее $\nu$ в $\mu$, тоже является полиморфизмом.

\sm

c) Произведение мер $\mu\times \nu$ -- полиморфизм.

\sm

d) Пусть $s(m,n)$ -- неотрицательная функция на $M\times N$, удовлетворяющая условиям
$$
\int_M s(m,n)\, d\mu(m)=1,\qquad \int_M s(m,n)\, d\nu(n)=1.
$$
Тогда мера $s(m,n)\,d\sigma(m,n)$
является полиморфизмом. Мы будем называть такие полиморфизмы {\it абсолютно непрерывными.}

\sm

e) Пусть $\sigma\in \Pol(M,N)$  -- полиморфизм. Та же мера
может рассматриваться как полиморфизм $N\looparrowright M$.
Мы скажем, что это {\it сопряженный полиморфизм},
и обозначим его через $\sigma^*$.
\hfill $\lozenge$

\sm 

{\bf\punct Марковские операторы.}
{\it Марковский оператор} $T: L^2(M,\mu)\to L^2(N,\nu)$ -- это ограниченный
оператор, удовлетворяющий следующим свойствам:

\sm

1) $T$ переводит неотрицательные функции в неотрицательные;

\sm

2)  $T$ переводит единичную функцию в единичную, $T^*$ тоже переводит единичную функцию
в единичную.

\sm 

\begin{theorem}
Полиморфизмы и марковские операторы находятся в естественном взаимно однозначном соответствии.
\end{theorem}

{\sc Доказательство.}
{\it Построение марковского оператора по полиморфизму.}
Обозначим через $I_A$ индикаторную функцию множества $A$
(то есть функцию, равную единице на множестве $A$ и нулю вне него).
Рассмотрим полиморфизм $\sigma: M\looparrowright N$.
Рассмотрим следующую полуторалинейную форму на
$L^2(M,\mu)\times L^2(N,\nu)$:
$$
B(f,g)=\int_{M\times N} f(m) g(n)\,d\sigma(m,n)
.
$$
Согласно неравенству Коши--Буняковского,
\begin{multline}
|B(f,g)|\le \Bigl(\int_M |f(m)|^2 \,d\sigma(m,n)\Bigr)^{1/2}
\Bigl(\int_N |g(n)|^2 \,d\sigma(m,n)\Bigr)^{1/2}
=\\=
\Bigl(\int_M |f(m)|^2 \,d\mu(m)\Bigr)^{1/2}
\Bigl(\int_N |g(n)|^2 \,d\nu(n)\Bigr)^{1/2}=
\|f\|_{L^2(M)}\cdot\|g\|_{L^2(N)}
\label{eq:KBun}
\end{multline}
(мы воспользовались определением полиморфизма).
Поэтому существует оператор $T:L^2(M)\to L^2(N)$
с нормой $\le 1$, такой, что
$$
B(f,g)=\la Tf,g\ra_{L^2(N)}
$$
Осталось проверить, что $T$ -- марковский оператор.

Для неотрицательных $f$, $g$ мы имеем
$B(f,g)\ge 0$. Если бы $Tf$ была бы отрицательна
на множестве $C$ ненулевой меры, то мы имели бы
$B(f,I_C)<0$, и мы приходим к противоречию. 

Далее, если 
$$B(1,g)=\int_{M\times N} g(n)\,d\sigma(m,n)=\int_N g(n)\,d\nu(n)=
\la 1,g\ra_{L^2(N)},$$
а поэтому $T$ переводит 1 в 1.

\sm

{\it Построение полиморфизма по марковскому оператору.}
Пусть  $T$ -- марковский оператор. Определим меру $\sigma$ на
$M\times N$ по формуле 
$$
\sigma(A\times B)=\la T I_A, I_B\ra_{L^2(N)},\qquad\text{где $A\subset$, $B\subset N$}
.$$ 
Мы получили меру на полукольце, и нужно формальное доказательство счетной аддитивности.
Здесь можно избежать унылых рассуждений следующим образом. Для конечных
пространств с мерой счетная аддитивность не вызывает сомнений.
Далее мы реализуем $M$ и $N$ как проективные пределы конечных пространств,
$$
M=\lim\limits_{\longleftarrow}
\cM_j,\qquad N=\lim\limits_{\longleftarrow} \cN_j
.
$$
(например, берем отрезок $[0,1]$, делим его на$2^j$ равных частей и отображаем на множество
их $2^j$ точек). Оператор $T$ определяет полиморфизмы
$\sigma_{ij}:\cM_i\to \cN_j$,
а мы далее строим $\sigma$ как обратный предел
$$\qquad\qquad\qquad\qquad\qquad\qquad\qquad
\sigma:=\lim\limits_{\longleftarrow} \sigma_{ij}.
\qquad\qquad\qquad\qquad\qquad\qquad\qquad\square
$$  

Обозначим марковский оператор, соответствующий полиморфизму $\sigma$,
через $T(\sigma)$.  

\sm

Отметим, что для абсолютно непрерывного полиморфизма $\sigma\in\Pol(M,N)$, задаваемого функцией
$s$, мы получаем обычный интегральный оператор
$$
T(\sigma) f(n)=\int_M s(m,n) f(m) \,d\mu(m).
$$ 

\begin{proposition}
	Норма любого марковского оператора равна 1.
\end{proposition}

{\sc Доказательство.}
Утверждение вытекает из неравенства (\ref{eq:KBun}).
\hfill $\square$

\sm 

{\bf\punct  Сходимость полиморфизмов.}
Пусть $\sigma_j$, $\sigma\in\Pol(M,N)$ -- полиморфизмы.
Мы скажем, что $\sigma_j$ сходится к $\sigma$, если для 
любых измеримых подмножеств $A\subset M$, $B\subset N$ мы имеем
сходимость 
\begin{equation}
\sigma_j(A\times B)\to  \sigma(A\times B).
\label{eq:sigma-j-sigma}
\end{equation}

\begin{proposition}
Последовательность	$\sigma_j$ сходится к $\sigma$
тогда и только тогда, когда $T(\sigma_j)$ сходится
к $T(\sigma)$ слабо. 
\end{proposition}

{\sc Доказательство.} 
Сходимость (\ref{eq:sigma-j-sigma}) равносильна  сходимости
$$
\la T_j(\sigma_j) \,I_A,\,I_B\ra\to \la T_j(\sigma) \,I_A,\,I_B\ra
,$$
а это равносильно слабой сходимости $T(\sigma_j)\to T(\sigma)$.
\hfill$\square$

\begin{corollary}
	Пространство всех полиморфизмов $M\to N$ компактно.
\end{corollary}

{\sc Доказательство.} Достаточно проверить, что множество марковских
операторов слабо компактно. Пусть $T_j$ -- марковские операторы.
Т.к. $\|T_j\|=1$, из последовательности $T_j$ можно выбрать сходящуюся
подпоследовательность. Достаточно ясно, что ее предел -- марковский оператор.
\hfil $\square$

\sm

{\sc Задача.} a) Покажите, что абсолютно непрерывные полиморфизмы 
плотны в множестве всех полиморфизмов $\Pol(M,N)$.

\sm

b) Покажите, что в случае пространства $M$ c непрерывной мерой группа $\Ams(M)$
плотна в полугруппе $\Pol(M,M)$.
\hfill $\lozenge$

\sm

{\bf\punct Определение полиморфизмов в терминах условных мер.}
Пусть $\sigma:M\looparrowright N$ -- полиморфизм. Тогда  по теореме Рохлина 
на множествах $M\times n\subset M\times N$ определены условные меры
$\sigma^{\{n\}}$, они определяются из свойства: для любого
измеримого%
\footnote{Измеримость можно понимать в любом из двух смыслов -- измеримость
	относительно $\sigma$ или измеримость по Борелю.}
подмножества $V\subset W$ выполнено
$$
\sigma(V)=\int_N\sigma^{\{n\}}\bigl(V\cap (M\times a)\bigr)\, d\nu(n)
$$
или для любой измеримой функции $h$
$$
\int_{M\times N} h(m,n) \,d\sigma(m,n)=
\int_N \int_M h(m,n)\,d\sigma^{\{n\}}(m)\, d\nu(n).
$$
Мы берем $h(m,n)=f(m)g(n)$ и получаем
$$
T(\sigma) f(n)=\int_M f(m)\,d\sigma^{\{n\}}(m)
$$

В связи с этим мы можем считать полиморфизм размазывающим отображением,
он переводит точки $n\in N$ в меры $\sigma^{\{n\}}$ на $M$.

\sm

{\bf\punct Умножение полиморфизмов.}
Ясно, что произведение марковских операторов -- марковский
оператор. Поэтому определено и произведение полиморфизмов
$$
\sigma:(M,\mu)\looparrowright (N,\nu),\qquad \tau: (N,\nu)\looparrowright (K,\kappa).
$$ 
Формулу 
$$T(\sigma\tau)=T(\tau)T(\sigma)$$
можно считать определением произведения. Так как произведение операторов
с нормой $\le 1$ раздельно слабо непрерывно, то и прозведение полиморфизмов раздельно непрерывно.

\sm

{\sc Задача.} a) Пусть $\phi:M\to N$, $\psi:N\to K$
-- отображения, сохраняющие меру (как в первом примере п.\ref{ss:polymorphisms}).  Проверьте, что произведение
полиморфизмов соответствует произведению отображений, $\wt{\psi\phi}=\wt\psi\wt\phi$.
 
 \sm
 
 b) Пусть полиморфизмы $\sigma\in \Pol(M,N)$, $\tau\in\Pol(N,K)$
 aбсолютно непрерывны и задаются функциями $s(m,n)$, $t(n,k)$. Тогда
 произведение полиморфизмов абсолютно непрерывно и  соответствует
 функции 
 $$\qquad\qquad\qquad
 u(m,k)=\int_N s(m,n)\,t(m,k)\,d\nu(n).
 \qquad\qquad\qquad
 \lozenge
 $$ 
 
 Эта задача дает способ определения произведения полиморфизмов.
 В случае абсолютно непрерывных полиморфизмов произведение
 задается понятной формулой, а дальше оно продолжается по непрерыности.

Есть более прямой способ определять произведение полиморфизмов.
Пусть $\sigma\in\Pol(M,N)$, $\tau \in \Pol(N,K)$.
Тогда система условных мер, соответствующая полиморфизму
$\tau \sigma$, задается формулой
$$
(\tau\sigma)^{(k)}(m)=
\int_N \sigma^{(n)}(m)\,d\tau^{(k)}(n).
$$
Формула эта вполне наглядна и отвечает идее перемножения размазывающих
отображений.

\sm

{\bf\punct Замыкания действий.}
Пусть группа $G$ действует на пространстве $(M,\mu)$
c конечной мерой
преобразованиями, сохраняющими меру. Рассмотрим образ группы
$G$ в $\Ams(M)$ и замкнем его в полугруппе
$\Pol(M,M)$ всех полиморфизмов. Замыкание должно быть
компактно (в силу компактности полугруппы полиморфизмов)
и должно быть полугруппой (в силу раздельной непрерывности
произведения в $\Pol(M,M)$).

Соответственно возникает вопрос о явном описании таких замыканий
для конкретных действий. Он много обсуждался для групп
$\Z$ и $\R$, для групп Ли он не очень интересен, а для
бесконечномерных групп здесь известно не много.
Мы приведем простой пример (Нельсон, \cite{Nel}) на эту тему.

\sm 

{\bf \punct Пример: ортогональная группа и
	 гауссовы меры.}
Рассмотрим действие ортогональной группы $\bfO(\infty)$
на пространстве
$\R^\infty$ с гауссовой мерой, описанное в \S8.

\begin{theorem}
{\rm a)}	Замыкание группы $\bfO(\infty)$ в полугруппе полиморфизмов
	изоморфно полугруппе $B$ всех операторов в вещественном
	$\ell_2$  c нормой $\le 1$.
	
	\sm
	
	{\rm b)}  Точкой замыкания, соответствующей матрице
	$T\in B$, является гауссова мера на $\R^\infty\times \R^\infty$
	с характеристической функцией%
	\footnote{Имеется в виду преобразование Фурье меры, подробней см. в доказательстве.}
	\begin{equation}
	\exp\Bigl\{-\frac12\begin{pmatrix}\xi &\eta\end{pmatrix}
	\begin{pmatrix}
	1&T\\T^t&1
	\end{pmatrix}
	\begin{pmatrix}\xi^t \\\eta^t\end{pmatrix}
	\Bigr\}.
	\label{eq:TU1}
	\end{equation} 
\end{theorem}

{\sc Доказательство.}   Вместо замыкания в полиморфизмах
мы можем рассматривать слабое замыкание группы $\bfO(\infty)$ в пространстве
операторов в $L^2(\R^\infty)$. Как мы видели в п.\ref{ss:symmetric-sum},
представление в $L^2(\R^\infty)$ разлагается в прямую
сумму симметрических степеней пространства $\ell_2$.

Замыкание группы $\bfO(\infty)$ в пространстве операторов в 
$\ell_2$ является полугруппой $B$. Если последовательность
ортогональных операторов $U_j$ слабо сходится к 
$T$, то последовательность симметрических степеней 
$\mathsf S^k U_j$  слабо сходится к $\mathsf S^k T$. 
Поэтому искомое слабое замыкание группы $\bfO(\infty)$
совпадает с $B$.

Обозначим через
$\R^\infty_{fin}$
пространство векторов $(\xi_1,\xi_2,\dots)$,
где $\xi_j\in \R$ и $\xi_j=0$ для достаточно больших $j$.
Пусть  $\sigma$ -- произвольная вероятностная мера  
на $\R^\infty \times \R^\infty$. Мы определим ее характеристическую
функцию $\xi[\mu;\xi,\eta]$ как функцию на 
$\R^\infty_{fin}\times \R^\infty_{fin}$, заданную формулой
$$
\xi[\sigma;\xi,\eta]=\int_{\R^\infty \times \R^\infty}
e^{i(\xi x^t+ \eta y^t)}\,d\sigma(x,y)=
\la T[\sigma]e^{i\xi x^t}, e^{i\eta y^t}\ra_{L^2(\R^\infty)}.
$$
В частности, характеристические функции определены для
полиморфизмов, причем сходимость полиморфизмов влечет
поточечную сходимость характеристических функций.

Рассмотрим полиморфизм $\sigma_U$,
задаваемый ортогональным оператором
$U$. Он сосредоточен на графике отображения 
$U:\R^\infty \to \R^\infty$, а его характеристическая функция
равна
\begin{multline}
\xi[\sigma_U;\xi,\eta]=\int_{\R^\infty \times \R^\infty}
e^{i(\xi x+ \eta y)}\,d\sigma_U(x,y)
=\int_{\R^\infty} e^{i(\xi x^t+ \eta U x^t)}
\,d\mu(x)=\\=\exp\Bigl\{-\frac12(\xi+ \eta U)(\xi+ \eta U)^t\Bigr\}=
\exp\Bigl\{-\frac12(\xi\xi^t+\eta\eta^t- \eta U \xi^t-\xi U^t\eta\Bigr\}=\\=
\exp\Bigl\{-\frac 12\begin{pmatrix}
\xi&\eta
\end{pmatrix}
\begin{pmatrix}
1&U\\U^t&1
\end{pmatrix}
\begin{pmatrix}
\xi^t\\\eta^t
\end{pmatrix}
\Bigr\}.
\label{eq:TU2}
\end{multline}
Если последовательность $U_j$ cходится
слабо к $T$, то последовательность 
(\ref{eq:TU2}) сходится к (\ref{eq:TU1}).
\hfill $\square$

\sm

{\bf \punct Задача о замыкании действий с квазиинвариантной мерой.}
Действий бесконечномерных групп, сохраняющих вероятностную меру, известно
сравнительно немного. В этом смысле более интересна задача о действиях
с квазиинвариантной мерой. 

Обозначим через $\R^\times$ мультипликативную полугруппу положительных 
вещественных чисел по умножению. Обозначим через $t$ координату на 
$\R^\times$. Пусть $(M,\mu)$, $(N,\nu)$
-- лебеговские пространства с вероятностными мерами.
Мы скажем, что {\it $\R^\times$-полиморфизм} $(M,\mu)\looparrowright (N,\nu)$
-- это мера $\sigma$ на $M\times N\times \R^\times$, такая, что

\sm

1. Проекция $\sigma$ на $M$ совпадает с $\mu$;

\sm

2. Проекция $t\cdot \sigma$ на $M$ совпадает с $\nu$.

\sm

{\sc  Пример.} Пусть $M=N$ -- отрезок $[0,1]$ с мерой Лебега $\mu$, пусть $q$ -- отображение, 
переводящее меру в эквивалентную меру. Пусть $q'(m)$ -- производная Радона--Никодима.
Рассмотрим отображение $M\to M\times M\times\R^\times$,  заданное формулой
$$
m\mapsto (m, \,q(m), \,q'(m)).
$$ 
Тогда образ $\mu$ при этом отображении является $\R^\times$ -- полиморфизмом.
\hfill $\lozenge$

\sm

{\sc Пример.} Рассмотрим функцию 
$$(m,n)\mapsto u_{m,n}$$
 на $M\times N$ со значениями
в полугруппе конечных мер на $\R^\times$.
По этим данным мы можем определить меру
на $M\times N\times \R^\times$,
$$
\sigma(A\times B\times C)=
\int_A\int_B\int_C
du_{m,n}(t)\,d\nu(n)\,d\mu(m).
$$
 Допустим, что выполнены следующие тождества:
для любых измеримых $A\subset M$, $B\subset N$
\begin{align*}
\mu(A)=\int_A\int_N\int_{\R^\times} du_{m,n}(t)\,d\nu(n)\,d\mu(m);
\\
\nu(B)=\int_B\int_M  \int_{\R^\times} t\cdot du_{m,n}(t)\, d\mu(m)\,d\nu(n).
\end{align*}
Тогда $\sigma$ является $\R^\times$-полиморфизмом.
Такие полиморфизмы мы назовем {\it абсолютно непрерывными}.
\hfill $\lozenge$

\sm

Неформально, $\R^\times$-полиморфизм -- это <<отображение>>,
которое <<размазывает>> точки по пространству с мерой, причем размазывается
также производная Радона--Никодима. 

\sm

Теперь нам нужно определить умножение. Сначала определим сходимость. 
Рассмотрим $\R^\times$-полиморфизм $\sigma$. Пусть
$A\subset M$, $B\subset M$ -- измеримые подмножества.
Ограничим меру $\sigma$ на $A\times B\times \R^\times$
и спроектируем ее на $\R^\times$. Обозначим полученную меру через
$\sigma[A\times B]$. 
  
Мы скажем, что последовательность $\R^\times$-по\-ли\-мор\-физ\-мов 
$\sigma_j$
сходится к $\sigma$, если для любых измеримых подмножеств
$A\subset M$, $B\subset M$ имеют место  слабые сходимости
мер на $\R^\times$:
$$
\sigma_j[A\times B]\to \sigma[A\times B],\qquad
t\cdot \sigma_j[A\times B]\to t\cdot \sigma[A\times B].
$$
Разумеется, абсолютно непрерывные полиморфизмы
плотны во всех полиморфизмах.

Умножение легко определяется для абсолютно непрерывных
$\R^\times$-по\-ли\-мор\-физ\-мов. Пусть полиморфизмы 
$M\looparrowright N$, $N\looparrowright K$ определяются
функциями $u_{m,n}$, $v_{n,k}$. Тогда произведение
полиморфизмов определяется функцией
$$
w(m,k)=\int_N u_{m,n}*v_{n,k}\,d\nu(n),
$$
где $*$ обозначает свертку мер на $\R^\times$.

\begin{theorem}
	Это произведение допускает раздельно непрерывное продолжение на
	произвольные полиморфизмы.
\end{theorem}

Если бесконечномерная группа $G$ действует преобразованиями, оставляющими меру
квазиинвариантной, то встает вопрос об описании замыкания
группы в полугруппе полиморфизмов. В настоящее
время этот вопрос решен лишь для групп симметрий
гауссовых и пуассоновских мер, \cite{Ner-fa}, \cite{Ner-matching}.

\section*{Добавление. Лебеговские пространства}
	
Здесь собраны некоторые факты теории меры, которые используются
в разных местах этих заметок.

\sm 

{\bf Д.1. Теорема об изоморфизме.}
Мы скажем, что пространство с мерой {\it полно по Лебегу},  если 
подмножества
множеств меры ноль измеримы (и тем самым имеют меру ноль). Понятно, что
из любого пространства с мерой автоматически изготовляется полное по Лебегу пространство.

Пусть $(M,\mu)$ -- пространство 
с вероятностной  мерой ($\mu(M)=1$), полное по Лебегу. Положим, что $\sigma$-алгебра измеримых множеств
имеет счетную базу, и что мера любой точки равна нулю.  Тогда 
$(M,\mu)$ изоморфно%
\footnote{Под <<изоморфизмом>> подразумевается биективное с точностью до
	множества меры нуль отображение, сохраняющее меру.}
отрезку $[0,1]$ с мерой Лебега (см. \cite[\S41]{Hal}, \cite[теорема 9.3.4]{Bog}),
эта теорема  в разных книгах приписывается Каратеодори, Шпильрайну, Халмошу--фон Нейману,
1938-1942.

\sm

{\sc Задача.} a) Пусть $\Z_2$ -- двоеточие, причем мера каждой точки равна $1/2$. 
Рассмотрим счетное произведение $\Z_2^\infty$ пространства $\Z_2$ самого на себя.
Каждой точке $a=(a_1, a_2, \dots)\in \Z_2^\infty$ мы поставим точку отрезка 
$[0,1]$ с двоичной записью $0,a_1a_2\dots$. Покажите, что мы получили биекцию с точностью до п.в.,
переводящую меру на $\Z_2^\infty$ в меру Лебега на $[0,1]$.

\sm

b) Так как $\Z_2^\infty \times \Z_2^\infty\to \Z_2^\infty$,
мы получаем изоморфизм отрезка $[0,1]$ и квадрата $[0,1]^2$ как пространств с мерой.
Они также изоморфны кубам $[0,1]^N$ любой размерности, в том числе $[0,1]^\infty$.

\sm 

c) Пусть на двоеточии $\Z_2$  мера  с весами точек $p$, $1-p$. Покажите, что при разных
$p$ пространства с мерой   $\Z_2^\infty$ изоморфны.\hfill 
$\lozenge$

\sm 

{\bf Д.2. Лебеговские пространства.}
Напомним, что {\it лебеговское пространство} -- это пространство с вероятностной мерой, эквивалентное
объединению промежутка и конечного или счетного набора точек (<<атомов>>), имеющих ненулевую меру.

Почти все пространства, встречающиеся в анализе и теории вероятности, являются лебеговскими
(хотя руководства по функциональному анализу предлагают б\'ольшую общность,
см. \cite{RS1}, \cite{Bog}). В частности, рассмотрим произвольное
сепарабельное метрическое пространство и вероятностную меру, определенную на его борелевской
$\sigma$-алгебре. Пополним $\sigma$-алгебру, положив, что подмножества 
множеств меры ноль имеют меру ноль. Тогда получится лебеговское пространство с мерой.

См. \cite{Roh}, \cite{Hal}, \cite[\S9.4]{Bog}, \cite{Kech}.

\sm 

{\bf Д.3. Абстрактное определение Рохлина.}
Рассмотрим пространство с мерой $(M,\mu)$. Пусть $B_j\subset M$ -- счетное семейство измеримых множеств.
Назовем его {\it базой}, если для любых двух различных точек $x$, $y\in M$ существет множество
$B_n$, содержащее одну из этих точек и не содержащее другую. Обозначим через
$B^{[0]}_j$ само множество $B_j$, а через  $B^{[1]}_j$ -- дополнение до него.
Для произвольной последовательности $\omega=(\omega_1,\omega_2,\dots)$ из нулей и единиц
рассмотрим пересечение $\cap_{j=1}^\infty B^{[\omega_j]}_j$. Скажем, что $M$
 {\it полно по отношению
к базе $B_j$}, если все такие пересечения непусты. Скажем, что $M$ полно
$(\mathrm{mod}, 0)$ относительно базы $B_j$, если
$M$ можно вложить в некоторое пространство $\wt M$, полное относительно базы
$\wt B_j$, так, что $\wt M\setminus M$ имеет меру ноль, а $B_j=\wt B_j\cap M$.

{\it Лебеговское пространство} по Рохлину -- это пространство с мерой, которое полно
  $(\mathrm{mod} 0)$
относительно некоторой счетной базы. Согласно \cite{Roh}, такие пространства 
являются лебеговскими в смысле предыдущего пункта.

\sm 

{\bf Д.4. Обратные (проективные) пределы множеств.} 
Пусть $X_1$, $X_2$, \dots --множества. Пусть $\psi_j:X_j\to X_{j-1}$
-- цепочка отображений,
$$
\dots\xleftarrow{\psi_k} X_{k} \xleftarrow{\psi_{k+1}}
X_{k+1}\xleftarrow{\psi_{k+2}} X_{k+2}
\xleftarrow{\psi_{k+3}}\dots
$$
Обратным пределом системы множеств называется множество $\cX$ и набор
отображений $\pi_j:\cX\to X_j$, такие, что
$\psi_j\circ \pi_j=\psi_{j-1}$, причем для любой пары $x$, $y$ различных
элементов $\cX$ найдется номер $j$, при котором $\psi_j(x)\ne \psi_j(y)$.

Обратный предел можно описать следующим образом.
Рассмотрим прямое произведение $X_1\times X_2\times X_3\dots$. Множество
$\cX$ реализуется как подмножество в этом произведение, состоящее из последовательностей
$x_1\in X_1$, $x_2\in X_2$, \dots, таких, что $\psi_j x_j=x_{j-1}$ для
всех $j$. Отображение $\pi_j$ ставит в соответствие такой последовательности
ее $j$-ый член $x_j$.

Обратный предел единственен в следующем смысле. Если $\cX'$ --другой обратный предел,
а $\psi'_j:\cX'\to X_j$ соответствующие проекции, то существует единственное отображение
$\theta:\cX\to \cX'$, такое, что
$\psi_j=\psi'_j \theta$.   

\sm 

{\bf Д.5. Теорема Колмогорова.}
Теперь определим {\it обратный предел лебеговских пространств}.
Пусть $(M_j,\mu_j)$ -- лебеговские  пространства с вероятностными  мерами, пусть
$\psi_j:M_j\to M_{j-1}$ -- отображения, такие, что образ каждой меры
$\mu_j$ при отображении $\psi_j$ совпадает с $\mu_{j-1}$. Пусть $\cM$ -- проективный предел
пространств $M_j$.
Тогда на $\cM$ канонически определена мера $\mu$, такая, что образ $\mu$ при каждом отображении
$\pi_j$ совпадает с $\mu_j$. Фактически это означает, что для всех $j$ и всех
измеримых $A_j\subset \cM$ прообразы $\pi_j^{-1}(A)\subset \cM$ измеримы
и $\mu (\pi_j^{-1}(A))=\mu(A)$. Далее мы берем $\sigma$-алгебру, порожденную
всеми такими множествами, пополним ее по Лебегу. Это дает $\sigma$-алгебру
пространства $\cM$, полученное пространство будет лебеговским%
\footnote{Существование меры $\mu$ является одной из версий теоремы Колмогорова о проективных 
	системах мер.
	В обычной формулировке этой теоремы  в учебниках (см., например, Ширяева или \cite[\S 7.7]{Bog}),
	берутся произведения $M_1\times \dots \times M_j$ полных сепарабельных метрических пространств
	и проекции $j$-ого произведения на $(j-1)$-oe.  Пусть $M_k$ -- проективная система
	лебеговских пространств, как выше, обозначим $\psi_{j1}:=\psi_{i+1}\circ \dots \circ \psi_j$
	-- проекции $M_j\to M_i$. Определим меру $\wt \mu_j$ на каждом множестве
	$M_1\times \dots \times M_j$ как образ меры $\mu_j$ при отображении
	$M_j\to  M_1\times \dots \times M_j$, заданном формулой
	$$x\mapsto (\psi_{j1}(x), \psi_{j2}(x), \dots, x).$$
	Далее применяем обычную теорему Колмогорова к системе мер $\wt \mu_j$.}.

\sm

{\sc Задача.} Отображение $\psi_j: M_j\to M_{j-1}$ индуцирует изометричное вложение
$\psi_j^*:L^2(M_{j-1}) \to L^2(M_j)$. Таким образом мы получаем цепочку вложений
$$
\dots\xrightarrow{\psi_{j-1}^*} L^2(M_{j}) \xrightarrow{\psi_{j}^*}
L^2(M_{j+1})\xrightarrow{\psi_{j+1}^*} L^2(M_{j+2})
\xrightarrow{\psi_{j+2}^*}\dots
$$	
Покажите, что объединение всех этих  пространств плотно в $L^2(\cM)$. 	
\hfill $\lozenge$

\sm 

{\bf Д.6. Теорема Рохлина об условных мерах.}  
Пусть $(M,\mu)$, $(N,\nu)$ -- лебеговские пространства, $\psi:M\to N$
-- отображение, такое, что образ меры $\mu$ равен мере $\nu$. 
Обозначим через $M_n$ прообраз точки $n\in N$, тем самым у нас получается
разбиение множества $M$ на непересекающиеся подмножества.
Оказывается, что канонически определена {\it система условных мер},
на множествах $M_n$. А именно, на почти каждом (в смысле меры $\nu$) подмножестве
$M_n$ определена вероятностная мера $\mu_n$, причем для любого измеримого подмножества
$A\subset M$ пересечения $A\cap M_n\subset M_n$ измеримы относительно $\mu_n$ для почти всех
$n$ и
$$
\mu(A)=\int_A \mu_n(A\cap M_n)\, d\nu(n).
$$
Равносильно, для любой интегрируемой функции на $M$ ее ограничения
на почти все множества $M_n$ являются $\mu_n$-измеримыми и
$$
\int_M f(m)\,d\mu(m)=\int_N \int_{M_n} f(m)\,d\mu_n(m)\,d\nu(n).
$$

{\sc Задача.} Пусть $h:[0,1]\to \R$ -- аналитическая 
функция. Какой будет образ  $\nu$ меры Лебега при этом отображении? 
Опишите условные меры.
\hfill $\lozenge$  

\sm

См. \cite{Roh}, \cite[\S9.4, \S 10.8]{Bog}. Стоит отметить, что эта теория имеет также борелевскую
версию, см. ссылку \cite{Kech} к \S 5.


\begin{thebibliography}{GGG1}

\bibitem[Ada]{Ada}
Адамс Дж. {\it Лекции по группам Ли.}
М.: Наука, 1979

\bibitem[Bog]{Bog}
Богачев В.И., {\it Основы теории меры}, т.2, НИЦ Регулярная и хаотическая динамика, Москва-Ижевск, 2006 

\bibitem[Bou1]{Bou1}
Бурбаки Н.  {\it Элементы математики. Общая топология. Топологические группы. 
Числа и связанные с ними группы и пространства.} — М.:
Наука, 1969

\bibitem[Bou2]{Bou2}
Бурбаки Н.
{\it Интегрирование. Векторное интегрирование. Мера Хаара. Свертка и представления.}
 Наука, 1969
 
 \bibitem[GW]{GW}
 Goodman, R.; Wallach, N. R. {\it Representations and invariants of the classical groups}.
 Cambridge University Press, Cambridge, 1998.
 
 \bibitem[Haar]{Haa}
Haar, A.
{\it Der Massbegriff in der Theorie der kontinuierlichen Gruppen.}
Ann. of Math. (2) 34 (1933), no. 1, 147-169. 


\bibitem[Hal]{Hal}
Халмош П. {\it Теория меры.} М., Издательство иностранной литературы, 1953


\bibitem[HR1]{HR1}
Хьюитт Э., Росс К. {\it Абстрактный гармонический анализ. Том 1.  Структура топологических групп.
Теория интегрирования. Представления групп.}
М.: Наука, 1975.

\bibitem[HR2]{HR2}
Хьюитт Э., Росс К. {\it Абстрактный гармонический анализ. Том 2. 
 Структура и анализ компактных групп. Анализ на локально компактных абелевых группах.}
М.: Мир, 1975



\bibitem[Hua]{Hua}
 Хуа Ло-кен
{\it Гармонический анализ функций многих комплексных переменных в классических областях.}
  М., Издательство иностранной литературы, 1959 

\bibitem[Kir]{Kir}
Кириллов А.А. {\it Элементы теории представлений.} М.: Наука, 1978.

 \bibitem[Mac]{Mac}
 Macdonald, I. G. {\it Symmetric functions and Hall polynomials.} Second edition. With contributions by A. Zelevinsky.
 Oxford University Press, New York, 1995

\bibitem[Nach]{Nach} 
Nachbin L. {\it The Haar integral}, Van Nostrand, 1964.

\bibitem[NAdo]{NAdo}
Neretin, Yu. A.
{\it A construction of finite-dimensional faithful representation of Lie algebra.}
arXiv:math/0202190

\bibitem[Ner1]{Ner1}
Неретин Ю.А. 
{\it 
Категории симметрий и бесконечномерные группы.} М.:УРСС, 1998

\bibitem[Ner2]{Ner2}
Neretin, Yu. A.
{\it Lectures on Gaussian integral operators and classical groups.}
 European Mathematical Society, Z\"urich, 2011

\bibitem[Pon]{Pon}
Понтрягин Л.С. {\it Непрерывные группы, 3-е изд. М.:Наука, 1973}

\bibitem[Phe]{Phe}
Фелпс Р.
{\it Лекции о теоремах Шоке.} М. Мир, 1968. 

\bibitem[RS1]{RS1}
Рид М., Саймон Б. {\it Методы современной математической физики. Том 1. Функциональный анализ.} 
М.: Мир, 1977.

\bibitem[Ros]{Ros}
Rossmann, W.
{\it Lie groups. 
An introduction through linear groups.} 
Oxford University Press, Oxford, 2002.

\bibitem[Rud]{Rud}
Рудин У. {\it Функциональный анализ.} М.: Мир, 1975



\bibitem[Wei]{Wei}
Вейль А. {\it Интегрирование в топологических группах и его применения.}
М.: Издательство иностранной литературы. 1950.

\bibitem[Zhe1]{Zhe1}
Желобенко Д.П., Компактные группы Ли их представления. М.: Наука, 1970

\bibitem[Zhe2]{Zhe2}
Желобенко Д.П. {\it Основные структуры и методы теории представлений.} М.: Изд-во МЦНМО, 2004 -

\end{thebibliography}

\begin{thebibliography}{cc}
	
	\bibitem{Alb}
	S. Albeverio, R. H\o egh-Krohn, D. Testard und A. Vershik. 
	{\it Factorial representations of path groups.} J. Funct. Anal.,
	(1983) 51(1):115--131. 

\bibitem{Ald}
Aldous D. {\it Exchangeability and related topics.} Lect. Notes Math. 1985. V. 1117. P. 1-198. 

\bibitem{ALV}
Д. В. Алексеевский, А. М. Виноградов, В. В. Лычагин
{\it Основные идеи и понятия дифференциальной геометрии.}
 Соврем. пробл. мат. Фундам. направления, 28 (1988),  5--289
 
 \bibitem{App}
 Applebaum D. {\it Probability on compact groups.}
 Springer, 2015

\bibitem{Arn1}
Arnold, V.
{\it Sur la g\'eom\'etrie diff\'erentielle des groupes 
de Lie de dimension infinie et ses applications
\`a l'hydrodynamique des fluides parfaits}. 
Ann. Inst. Fourier (Grenoble) 16 1966 fasc. 1, 319-361.

\bibitem{Arn-ur}
Арнольд В.И.
{\it Обыкновенные дифференциальные уравнения.} Изд. 3., Ижевск, 2000.

\bibitem{Arn-mech}
Арнольд В.И. {\it Математические методы классической механики.}
3-е изд. — М.: Наука, 1989.

\bibitem{AKh}
Арнольд В.И., Хесин Б.А. {\it Топологические методы в гидродинамике.}
 М.: МЦНМО, 2007

 \bibitem{Aud}
 Audin, M., Damian, M., {\it Morse theory and Floer homology.}
 Springer,  2014.

 \bibitem{AvSy}
 Avron, J. E.; Simon, B.
{\it Almost periodic Hill's equation and the rings of Saturn.}
Phys. Rev. Lett. 46 (1981), no. 17, 1166-1168.


\bibitem{Bab}
 Бабенко И. К.
{\it Алгебра, геометрия и топология группы подстановок формальных степенных рядов.}
УМН, 2013, 68:1(409),  3-76.

 \bibitem{Ban0}
Banyaga, A.
{\it Sur la structure du groupe des diff\'eomorphismes qui pr\'eservent une forme symplectique.} 
Comment. Math. Helv. 53 (1978), no. 2, 174-227. 


 \bibitem{Ban}
 Banyaga, A. {\it The structure of classical diffeomorphism groups.}
 Kluwer, Dordrecht, 1997
 
 \bibitem{BMS}
 Bass, H.; Milnor, J.; Serre, J.-P.
{\it Solution of the congruence subgroup problem for $\SL_n$ ($\ge 3$) and $\Sp_{2n}$ ($n\ge 2$)}.
Inst. Hautes \'Etudes Sci. Publ. Math. No. 33 (1967) 59-137. 
 
 \bibitem{Ber}
 Березин   Ф.А. {\it Метод вторичного квантования}. (Изд. 2-е) М., Наука 1986 
 
  \bibitem{Ber-super}
   Березин  Ф.А. {\it Введение в алгебру и анализ сантикоммутирующими переменными.}
   Издательство МГУ, Москва, 1983; Расширенный английский перевод:
   Berezin F.A.
{\it Introduction
  to superanalysis.}  D. Reidel Publishing
  Co., Dordrecht, 1987.

 
  \bibitem{BK}
 Becker, H., Kechris, A. S. {\it The descriptive set theory of Polish group actions.} 
 Cambridge University Press, Cambridge, 1996.
 
 \bibitem{Bohr}
 Бор Г.{\it 	
Почти периодические функции.} М.-- Л., 1934

\bibitem{Boug}
Bougerol, Ph., J. Lacroix {\it Products of Random Matrices with Applications to Schr\"odinger Operators}. Birkh\"auser, Boston 1985.

\bibitem{Bre}
Brenier, Ya.
{\it Minimal geodesics on groups of volume-preserving maps 
and generalized solutions of the Euler equations.}
Comm. Pure Appl. Math. 52 (1999), no. 4, 411-452. 

\bibitem{Breu}
 Breuillard E., {\it Random walks on Lie groups},


\bibitem{Bro}
 Бродский, М. С.
{\it Унитарные операторные узлы и их характеристические функции}, 
УМН, 33:4(202) (1978), 141-168 

\bibitem{Buch}
 Buchstaber V. M., {\it $n$-valued groups: theory and applications}, Mosc. Math. J., 6:1 (2006), 57-84 

\bibitem{Cam}
Cameron, P. J.
{\it Permutation groups.}
 Cambridge University Press, Cambridge, 1999. 
 
 \bibitem{CFP}
 Cannon, J. W.; Floyd, W. J.; Parry, W. R. (1996), {]it Introductory notes on Richard Thompson's groups.}  L'Enseignement Math\'ematique. Revue Internationale. IIe S\'erie 42 (3): 215-256.
 
 \bibitem{CCTV}
 Carmeli, C., Cassinelli, G., Toigo, A., Varadarajan, V. S. {\it Unitary Representations
 of Super Lie Groups and Applications to the Classification and Multiplet
 Structure of Super Particles.} Commun. Math. Phys. 263, 217-258 (2006)

 
 \bibitem{Sophus}
  Картье П., Бланшар А., Лазар М., Брюа Ф., Берже М., Серр Ж.-П. 
  {\it Теория алгебр Ли. Топология групп Ли (Семинар <<Софус Ли>>).}
  М.: Иностранная литература, 1962.
  
 
 \bibitem{Cor}
 C. Corduneanu, {\it Almost periodic functions}, Wiley (1968).
 
 \bibitem{Mal}
 Laurence C.,  Greenleaf F. P. 
 {\it Representations of Nilpotent Lie Groups and their Applications: Volume 1, Part 1, Basic Theory and Examples}, Cambridge University Press (2004).

\bibitem{Witt} 
 De Witt, B. {\it Supermanifolds}. Cambridge University Press, 1992.
 
 \bibitem{Die}
 Дьедонне Ж. {\it Геометрия классических групп.}
 М., Мир, 1972
 
 \bibitem{Dix}
 Диксмье Ж. {\it $C^*$-алгебры и их представления}. М.: Наука, 1974. 
 
 \bibitem{Dix-U}
 Диксмье Ж. {\it Универсальные обертывающие алгебры}. М.: Мир, 1978. 
 
 \bibitem{Dym}
 Dym, H.
{\it Linear algebra in action.}
 American Mathematical Society, Providence, RI, 2013
 
 \bibitem{DRW}
 Dooley, A. H.; Repka, J.; Wildberger, N. J.
{\it Sums of adjoint orbits.}
Linear and Multilinear Algebra 36 (1993), no. 2, 79-101. 


\bibitem{EM}
Ebin, D. G.; Marsden, J. {\it Groups of diffeomorphisms and the motion of an incompressible fluid.}
Ann. of Math. (2) 92 1970, 102-163.
 
 


\bibitem{DS1}
Данфорд Н., Шварц Дж.Т.
{\it Линейные операторы
Том II. Спектральная теория}. М.: Мир, 1966

\bibitem{DZ}
Delfour, M. C., Zolesio, J.-P.
{\it Shapes and geometries.} 
 Society for Industrial and Applied Mathematics (SIAM), Philadelphia, PA, 2001.

\bibitem{FM}
Farb, B.; Margalit, D.
{\it A primer on mapping class groups.} Princeton University Press, Princeton, NJ, 2012.

  \bibitem{FTT}
  Fig\`a-Talamanca, A., Steger, T.
  {\it Harmonic analysis on trees.}  Academic Press, New York, 1987 

\bibitem{dict}
Frappat L., Sciarrino A., Sorba P.
{\it Dictionary on Lie algebras and superalgebras.}
Academic Press Inc.,
San Diego, CA (2000).

 \bibitem{GGP} 
   Гельфанд И.М., Граев М.И., Пятецкий-Шапиро И.И. 
{\it 
 Обобщенные функции, выпуск 6. Теория представлений и автоморфные функции.} 
 М.: Гос. изд-во физико-мат. лит-ры, 1966

\bibitem{GS}
 Ghys \'E.,  Sergiescu V., {\it Sur un groupe remarquable de diff\'eomorphismes du cercle}, Comment. Math. Helv.
62 (1987), 185--239.

\bibitem{Gol}
Голузин Г.М.
{\it Геометрическая теория функций комплексного переменного.} М. Наука, 1966.


\bibitem{HV}
Haller, S.; Vizman, C. {\it Non-linear Grassmannians as coadjoint orbits}, 
the abridged version is published: Math. Ann. 329 (2004), no. 4, 771-785; Preprint  arXiv:math/0305089

\bibitem{Has}
Hasimoto, R., {\it A soliton on a vortex filament}, J. Fluid Mechanics 51 (1972),
477-485.

\bibitem{Hey}
Хейер Х.
{\it Вероятностные меры на локально-компактных группах.} M. Мир, 1981.

\bibitem{HM}
Howe, R. E.; Moore, C. C. {\it Asymptotic properties of unitary representations.}
J. Funct. Anal. 32 (1979), no. 1, 72-96.

\bibitem{ology}
Iglesias-Zemmour, P.
{\it Diffeology.}
 Amer. Math. Soc., Providence, RI, 2013.
 
 \bibitem{Imb}
  Imbert M.
{\it Sur l'isomorphisme du groupe de Richard Thompson avec le groupe de
 Ptol\'em\'ee.} In {\it Geometric Galois actions,} 2, eds.
Schneps L., Lochak P., Cambrodge University Press, 1997.
 

\bibitem{IKT}
 Inci, H.; Kappeler, T.; Topalov, P. 
 {\it On the regularity of the composition of diffeomorphisms.}
 Mem. Amer. Math. Soc. 226 (2013), no. 1062.
 
 \bibitem{Ism-su2}
 Исмагилов Р. С. 
{\it Об унитарных представлениях группы $C^\infty_0(X,G)$, $G=\SU(2)$.}
 Матем. сб., 100(142):1(5) (1976),  117-131.
  
  \bibitem{Ism}
 Ismagilov, R. S.
{\it Representations of infinite-dimensional groups.}
AMS, Providence, RI, 1996

\bibitem{Iva}
 Ivanov, N. V. {\it Mapping class groups}. in {\it Handbook of geometric topology}, 523--633, North-Holland, Amsterdam, 2002.

\bibitem{Iwa}
 Iwahori, N. {\it On the structure of a Hecke ring of a Chevalley group over a finite field.} 
 J. Fac. Sci. Univ. Tokyo Sect. I 10 1964, 215-236 (1964).
 
 \bibitem{IM}
 Iwahori, N.; Matsumoto, H. {\it On some Bruhat decomposition and the structure of the Hecke rings of p-adic Chevalley groups.}
 Inst. Hautes \`Etudes Sci. Publ. Math. No. 25 1965 5-48.
 
  \bibitem{Jew}
 Jewett, R. I.
{\it Spaces with an abstract convolution of measures}.
Advances in Math. 18 (1975), no. 1, 1-101. 

\bibitem{Kap}
Kapoudjian, C.
{\it Virasoro-type extensions for the Higman-Thompson and Neretin groups.} 
Q. J. Math. 53, No. 3, 295-317 (2002).

\bibitem{KarM}
Каргаполов М. И., Мерзляков Ю. И. 
{\it Основы теории групп}.— 3-е изд., перераб. и доп.— М.: Наука, 1982.

\bibitem{Kash}
Kashaev, R. {\it Quantization of Teichm\"uller spaces and quantum dilogarithm.}
Lett. Math. Phys. 43(2), 105-115 (1998)

\bibitem{Kech}
Kechris, A. S. 
{\it Classical descriptive set theory.}
 Springer, New York, 1995.
 
\bibitem{Kech1}
Kechris, A. S. 
{\it Global aspects of ergodic group actions}, American Mathematical Society, 2010. 

\bibitem{Khe}
Khesin, B.
{\it The group and Hamiltonian descriptions of hydrodynamical systems.}
In {\it Lectures on topological fluid mechanics}, 139-155,
Lecture Notes in Math., 1973, Springer, Berlin, 2009. 
 
 \bibitem{HW}
Хесин Б.А., Вендт Р. {\it Геометрия бесконечномерных групп.} 
М.: МЦНМО, 2014  

\bibitem{Koh}
Kohno, T.
{\it S\'erie de Poincar\'e--Koszul associ?e aux groupes de tresses pures.}
Invent. Math.  82  (1985), 

\bibitem{Koo}
Koornwinder, T. H.
{\it Jacobi functions and analysis on noncompact semisimple Lie groups.}
in {\it Special functions: group theoretical aspects and applications}, 1-85,
Math. Appl., Reidel, Dordrecht, 1984.  

 \bibitem{KMR}
Kriegl A.,  Michor P.W, Rainer A.
{\it The convenient setting for quasianalytic Denjoy--Carleman differentiable mappings.}
 J. Functional Analysis 261, 7 (2011) 1799-1834.
 
 \bibitem{Kui}
 Kuiper, N. H.
{\it The homotopy type of the unitary group of Hilbert space.}
Topology 3 1965 19-30



\bibitem{McD}
McDuff, D., Salamon, D.
{\it Introduction to symplectic topology.}
Oxford University Press, New York, 1998.
 
 \bibitem{MW}
 Marsden, J.; Weinstein, A. {\it Coadjoint orbits, vortices, and Clebsch variables for incompressible fluids.}
 Phys. D 7 (1983), no. 1-3, 305-323.
 
  \bibitem{Mich}
 Michor P.W. {\it Some Geometric Evolution Equations Arising 
 as Geodesic Equations on Groups of Diffeomorphism,} in.
 {\it Progress in Non Linear Differential Equations and Their Applications}, Vol. 69. 
 Birkhauser Verlag 2006. Pages 133-215.
 
 \bibitem{MM}
 Michor, P. W.; Mumford, D.
 {\it An overview of the Riemannian metrics on spaces of curves using the Hamiltonian approach}. 
 Appl. Comput. Harmon. Anal. 23 (2007), no. 1, 74-113.
 
 \bibitem{Ner-Gamma}
  Неретин Ю. А.
{\it  Голоморфные продолжения представлений группы диффеоморфизмов окружности.}
Матем. сб., 180:5 (1989),  635-657.

\bibitem{Ner-derevo}
Неретин  Ю. А. 
{\it   О комбинаторных аналогах группы диффеоморфизмов окружности.}
Изв. РАН. Сер. матем., 56:5 (1992),  1072-1085

\bibitem{Ner-cantor}
 Неретин Ю. А.
{\it Группа диффеоморфизмов полупрямой и случайные канторовские множества.}
Матем. сб., 187:6 (1996),  73–84
 
\bibitem{Ner-MIRAN}
 Неретин Ю. А.
{\it Дробные диффузии и квазиинвариантные действия бесконечномерных групп.}
Тр. МИАН, 217 (1997),  135--181

\bibitem{Ner-trees}
Neretin, Yu.A.
{\it Groups of hierarchomorphisms of trees and related Hilbert spaces.}
J. Funct. Anal. 200, No.2, 505-535 (2003).

 \bibitem{Ner-super}
 Neretin, Yu. A. {\it
Gauss--Berezin integral operators and spinors over supergroups $\mathrm{OSp}(2p|2q)$, and Lagrangian super-Grasmannians}, Preprint arXiv:0707.0570.

 
 \bibitem{Ner-colligations}
Neretin, Yu. A. {\it Multi-operator colligations and multivariate characteristic functions.} 
Anal. Math. Phys. 1 (2011), no. 2-3, 121-138.  

\bibitem{Ner-spheric}
Неретин Ю. А. 
{\it Сферичность и умножение двойных классов смежности для бесконечномерных классических групп.}
Функц. анализ и его прил., 45:3 (2011),  79-96.
 
  \bibitem{Ner-Mal}
 Neretin, Yu.
{\it On topologies on Malcev completions of braid groups.} 
Mosc. Math. J. 12 (2012), no. 4, 803-824



 
 \bibitem{Ner-hyper}
Neretin, Yu.   {\it  On concentration of convolutions of double cosets at infinite-dimensional limit.}
arXiv:1211.6149

\bibitem{Ner-p-adic}
 Неретин Ю. А.
{\it Меры Хуа на пространстве p-адических матриц и обратные пределы грассманианов.}
Изв. РАН. Сер. матем., 77:5 (2013),  95–108

\bibitem{Ner-biinv}
 Неретин Ю. А.
{\it Биинвариантные функции на группе преобразований, оставляющих меру квазиинвариантной.}
Матем. сб., 205:9 (2014),  145--160.

\bibitem{Ner-p}
Neretin Yu.A.
{\it The space $L^2$ on semi-infinite Grassmannian over finite field.}
Adv. Math. 250 (2014), 320-350.

\bibitem{Ner-buildings}
 Неретин Ю. А.
{\it Бесконечномерные $p$-адические группы, полугруппы двойных классов смежности и внутренние функции на ансамблях Брюа–Титса.}
Изв. РАН. Сер. матем., 79:3 (2015),  87-130.


\bibitem{Ner-umn}
Neretin Yu.A.
{\it Infinite symmetric groups and combinatorial constructions of topological field theory type.}
 arXiv:1502.03472, to appear in Успехи мат. наук

 
\bibitem{Olsh-tree}
 Ольшанский Г. И.
{\it Классификация неприводимых представлений групп автоморфизмов деревьев Брюа-Титса}	
Функц. анал. и его прилож., 1977, 11:1, 32--42

\bibitem{Olsh-uni}
 Ольшанский Г. И.
{\it Унитарные представления бесконечномерных классических групп 
$\U(p,\infty)$, $\SO_0(p,\infty)$, $\Sp(p,\infty)$
и соответствующих групп движений.} Функц. анализ и его прил., 12:3 (1978),  32--44

\bibitem{Olsh-kiado}
Olshansky, G. I.
{\it Unitary representations of the infinite symmetric group: a semigroup approach}. 
In {\it Representations of Lie groups and Lie algebras {\rm(}Budapest, 1971{\rm)}}, 181-197, 
Akad. Kiad\'o, Budapest, 1985. 

\bibitem{Olsh-chip}
 Ольшанский Г. И.
{\it Унитарные представления (G,K)-пар, связанных с бесконечной симметрической группой 
$S(\infty)$.}
Алгебра и анализ, 1:4 (1989),  178-209

\bibitem{Olsh-GB}
Olshansky, G. I.
{\it Unitary representations of infinite-dimensional pairs 
$(G,K)$ and the formalism of R. Howe.} in 
{\it Representation of Lie groups and related topics}, 269-463,
 Gordon and Breach, New York, 1990. 
 
 \bibitem{Olsh-semi}
 Olshanski, G. I. {\it On semigroups related to infinite-dimensional groups}. In
 {\it Topics in representation theory.} (A. A. Kirillov, ed.),
  Amer. Math. Soc., Providence, R.I., 1991, 67-101.

 \bibitem{Pes0}
Pestov, V.
{\it Topological groups: where to from here?}
Topology Proc. 24 (1999), Summer, 421-502 (2001)
 
\bibitem{Pes}
Pestov, V.
{\it Dynamics of infinite-dimensional groups.
The Ramsey-Dvoretzky-Milman phenomenon.}
AMS, Providence, RI, 2006. 

\bibitem{Pick1}
Pickrell, D.
{\it Separable representations for automorphism groups of infinite symmetric spaces.}
J. Funct. Anal. 90 (1990), no. 1, 1-26.

\bibitem{Pick-mac}
Pickrel D. {\it Mackey analysis of infinite-dimensional motion groups.}
Pacif. J. Math. 1991. V. 150, N 1. P. 139-166. 



 \bibitem{PR}
Платонов В.П., Рапинчук А.С. {\it Алгебраические группы и теория чисел.}
М.: Наука, 1991

\bibitem{Pol1}
Polterovich, L. {\it The geometry of the group of symplectic
diffeomorphisms}. Birkh\;user, Basel, 2001.

\bibitem{Pol2}
Polterovich, L., Rosen, D. {\it Function theory on symplectic manifolds}. 
American Mathematical Society, Providence, RI, 2014. 

\bibitem{PS}
Пресли Э., Сигал Г. 
 {\it Группы петель.} М.: Мир, 1990.
 
  \bibitem{Reu}
 Reutenauer, Ch.
{\it Free Lie algebras.}
 Oxford University Press, New York, 1993.
 
 \bibitem{RZ}
 Ribes, L.; Zalesskii, P.
 {\it Profinite groups.}  Springer-Verlag, Berlin, 2010
 
 

 
  \bibitem{RosSo}
 Ch. Rosendal, Ch,  Solecki S.,
 {\it Automatic continuity of homomorphisms and fixed points on metric compacta,}
 Israel J. Math. 162 (2007), 349-371. 
 
 \bibitem{Rosl}
 R\"osler, M.
{\it Positive convolution structure for a class of Heckman-Opdam hypergeometric functions of type BC.}
J. Funct. Anal. 258 (2010), no. 8, 2779-2800
 
  \bibitem{Ser2}
 Серр Ж.-П. {\it Когомологии Галуа.}  М.: Мир, 1968

  \bibitem{Ser1}
 Серр Ж.-П. {\it Курс арифметики.}  М.: Мир, 1972
 
 
\bibitem{Ser-trees} 
 Serre, J.-P.
{\it Trees.}
 Springer Monographs in Mathematics.
Springer-Verlag, Berlin, 2003 

\bibitem{Sha}
Shale, D.
{\it Linear symmetries of free boson fields.}
Trans. Amer. Math. Soc. 103 1962 149--167.

\bibitem{Shav}
 Шавгулидзе Е. Т.
{\it Квазиинвариантные меры на группах диффеоморфизмов.}
Тр. МИАН, 217 (1997),  189–208

\bibitem{Serg}
 Сергеев А. Н.
{\it Тензорная алгебра тождественного представления как модуль над супералгебрами Ли 
	$\frG\frl(n,m)$ и $Q(n)$.}
Матем. сб., 123(165):3 (1984),  422-430


 \bibitem{SS}
 Singer, I. M.; Sternberg, S.
{\it The infinite groups of Lie and Cartan. I. The transitive groups.}
J. Analyse Math. 15 1965 1--114.

 \bibitem{Sou}
 Souriau, J.-M.
{\it Groupes diff\'erentiels.} in
{\it Differential geometrical methods in mathematical physics}, pp. 91–128,
Lecture Notes in Math., 836, Springer, Berlin-New York, 1980. 

 \bibitem{Tsa}
 Tsankov, T.
 {\it Automatic continuity for the unitary group}.
 Proc. Amer. Math. Soc. 141 (2013), no. 10, 3673-3680

 \bibitem{Viz}
Vizman, C.
{\it Induced differential forms on manifolds of functions}.
Arch. Math. (Brno) 47 (2011), no. 3, 201-215

\bibitem{Ury}
Urysohn P. S., {\it Sur un espace metrique universel}, Bull. Sci. Math., 51 (1927), 1--38; 
рус. перев. Урысон П. С., {\it Об универсальном метрическом пространстве}, П. С. Урысон,
{\it Труды по топологии},  т. 2, ГИТТЛ, М., 1951, 747-777. 

\bibitem{VaS}
 Вавилов Н. А., Степанов  А. В.
{\it Линейные группы над общими кольцами I. Общие места}
Зап. научн. сем. ПОМИ, 394 (2011),  33-139.



\bibitem{Ver-ury}
 Вершик А. М., {\it Случайные метрические пространства и универсальность}, УМН, 59:2(356) (2004),
 65-104.
 
 \bibitem{VGG}
 А. М. Вершик, И. М. Гельфанд, М. И. Граев
{\it Представления группы $\SL(2,R)$, где $R$ -- кольцо функций.}
УМН, 28:5(173) (1973),  83--128

\bibitem{VO}
 Винберг Э. Б.,  Горбацевич В. В.,  Онищик А. Л., {\it Строение групп и алгебр Ли},
 Группы Ли и алгебры Ли - 3, Итоги науки и техн. Сер. Соврем. пробл. мат.
Фундам. направления, 41, ВИНИТИ, М., 1990, 5--253 

\bibitem{You}
Younes, L.
{\it Shapes and diffeomorphisms.} 
 Springer-Verlag, Berlin, 2010.


\end{thebibliography}

\begin{thebibliography}{c}
 

\bibitem{Nai}
Наймарк М.А. {\it Нормированные кольца}.  М.: Наука, 1968

 \bibitem{Serre-finite}
 Серр Ж.-П. {\it Линейные представления конечных групп}, М. Мир, 1970

\end{thebibliography}

\begin{thebibliography}{R}
 
\bibitem{Met}
Мета М.Л. {\it Случайные матрицы}. МЦНМО, 2012.
 
 \bibitem{Ner-triangle}
 Неретин Ю. А.
 {\it Треугольники Релея и нематричная интерполяция матричных бета-интегралов.}
Матем. сб., 194:4 (2003),  49-74
 

 
  \bibitem{Wey}
Вейль Г. {\it  Классические группы, их инварианты и представления.}  М: ИЛ, 1947.


\end{thebibliography}

\begin{thebibliography}{cc}



\bibitem{Bog-gauss}
Богачев В.И. {\it Гауссовские меры,} 1997
  Москва,  Наука.
  

\bibitem{BShi}
А. В. Булинский, А. Н. Ширяев, {\it Теория случайных процессов},  Физматлит, Москва, 2005, 402 с. 

\bibitem{Gla}
E. Glasner, B. Tsirelson, B. Weiss, {\it The automorphism group of the Gaussian measure
cannot act pointwise},
Israel J. Math. 148 (2005) 305-329.


\bibitem{Ner-cmp}
Neretin Yu.A.
{\it Some remarks on quasi-invariant actions of loop groups and the group of diffeomorphisms of the circle.}
Comm. Math. Phys. 164 (1994), no. 3, 599--626.



\bibitem{ShF}
Шилов Г.Е., Фан Дык Тинь {\it Интеграл, мера и производная на линейных пространствах.}
М.: Наука, 1967

\bibitem{Shir}
Ширяев А.Н. {\it Вероятность}. М.:Наука, 1979

\end{thebibliography}

\begin{thebibliography}{cc}
 
 \bibitem{King}
 Кингман Дж. {\it Пуассоновские процессы.} — М.: МЦНМО, 2007. 
 
 \bibitem{Ner-poi}
 Неретин  Ю. А.
 {\it О соответствии между бозонным пространством Фока и пространством $L^2$ по мере Пуассона.}
Матем. сб., 188:11 (1997),  19--50
 
 \bibitem{VVG}
 Вершик  А. М.,  Гельфанд И. М.,  Граев М. И.
{\it  Представления группы диффеоморфизмов}
УМН, 30:6(186) (1975),  3-50
 
 \end{thebibliography}

\begin{thebibliography}{c}

\bibitem{KOV}
Kerov, S.; Olshanski, G.; Vershik, A.
{\it Harmonic analysis on the infinite symmetric group.}
Invent. Math. 158 (2004), no. 3, 551-642. 

\bibitem{KTs}
 Керов С. В.,  Цилевич Н. В.
{\it Случайное дробление отрезка порождает виртуальные перестановки с распределением Ювенса.}
Зап. научн. сем. ПОМИ, 223 (1995),  162-180


\bibitem{Tsi}
Цилевич Н. В.
 {\it Распределение длин циклов бесконечных перестановок.}
Зап. научн. сем. ПОМИ, 223 (1995),  148--161.

\bibitem{Olsh-NATO}
Olshanski, G. {\it An introduction to harmonic analysis on the infinite symmetric group.}
Asymptotic combinatorics with applications to mathematical physics (St. Petersburg, 2001), 127-160, 
Lecture Notes in Math., 1815, Springer, Berlin, 2003.

\bibitem{Str}
Strahov, E.
{\it Generalized characters of the symmetric group.}
Adv. Math. 212 (2007), no. 1, 109-142. 
 
\end{thebibliography}

\begin{thebibliography}{c}

\bibitem{BO}
Borodin, A.; Olshanski, G.
{\it Harmonic analysis on the infinite-dimensional unitary group
and determinantal point processes.}
Ann. of Math. (2) 161 (2005), no. 3, 1319-1422

\bibitem{Ner-hua}
  Neretin, Yu. A. {\it Hua-type integrals over unitary groups and over projective limits of unitary groups.}
  Duke Math. J. 114 (2002), no. 2, 239-266.

  \bibitem{Olsh-unitary}
Olshanski, G. {\it The problem of harmonic analysis on the infinite-dimensional unitary group.} 
J. Funct. Anal. 205 (2003), no. 2, 464-524.  
  
\bibitem{Osi}  
 Осиненко А. А.
 {\it Гармонический анализ на бесконечномерной унитарной группе.}
Зап. научн. сем. ПОМИ, 390 (2011),  237-285   

\bibitem{Pick-grass}
Pickrell D . {\it Measures on infinite-dimensional Grassmann manifolds.}
 J. Funct. Anal. 1986. V. 70, N 2. P. 323-356. 
\end{thebibliography}

\begin{thebibliography}{cc}
	
	\bibitem{Glasner}
	Glasner, E.
	{\it Ergodic theory via joinings}. 
	 American Mathematical Society, Providence, RI, 2003. 
	
	\bibitem{Kre}
Krengel U.	{\it Ergodic Theorems}, de Gruyter, 1985

\bibitem{Nel}
Nelson, E.
{\it The free Markoff field.}
J. Funct. Anal. 12 (1973), 21-–227.
	
	\bibitem{Ner-matching}
	Yu. A. Neretin
{\it	Spreading maps (polymorphisms), symmetries of Poisson processes, and matching summation.}
		Зап. научн. сем. ПОМИ, 292 (2002),  62-91.
	
	\bibitem{Ner-fa}
	Neretin, Yu.A.
	{\it Symmetries of Gaussian measures and operator colligations.}
	J. Funct. Anal. 263 (2012), no. 3, 782--802
	
	\bibitem{Ner-zanuda}
 Неретин 	Ю. А.
{\it О границе группы преобразований, оставляющих меру квазиинвариантной.}
Матем. сб., 204:8 (2013),  83--116

\bibitem{Ryzh}
  Рыжиков В. В., {\it Перемешивание, ранг и минимальное самоприсоединение действий с инвариантной мерой}, Матем. сб., 183:3 (1992), 133--160 

\bibitem{Ver-poly}
Вершик, А. М. {\it Многозначные отображения с инвариантной мерой {\rm (}полиморфизмы{\rm)} и марковские операторы},
Зап. научн. сем. ЛОМИ, 72, 1977, 26-61
\end{thebibliography}

\begin{thebibliography}{c}
	 \bibitem{Roh} 
	 Рохлин В. А. 
	 {\it	Об основных понятиях теории меры.}
	 Матем. сб., 25(67):1 (1949),  107--150.
	
\end{thebibliography}
\end{document}